\documentclass[11pt]{article}

\usepackage{amssymb,amsmath,amscd,amsthm,xy,enumitem,mathrsfs,pstricks,stmaryrd}
\xyoption{all}

\newtheorem{thm}{Theorem}[section]
\newtheorem{prp}[thm]{Proposition}
\newtheorem{lmm}[thm]{Lemma}
\newtheorem{crl}[thm]{Corollary}

\theoremstyle{definition}

\theoremstyle{remark}
\newtheorem{rmk}[thm]{Remark}

\numberwithin{equation}{section}

\newcounter{temp}

\def\lra{\longrightarrow}
\def\Lra{\Longrightarrow}
\def\Llra{\Longleftrightarrow}

\def\xlra#1{\xrightarrow{#1}}
\def\dbsqbr#1{\llbracket{#1}\rrbracket}

\def\BE#1{\begin{equation}\label{#1}}
\def\EE{\end{equation}}
\def\sm#1{\begin{small}{#1}\end{small}}

\def\wt#1{\widetilde{#1}}
\def\ov#1{\overline{#1}}
\def\eref#1{(\ref{#1})}
\def\tn#1{\textnormal{#1}}
\def\sf#1{\textsf{#1}}

\def\smsize#1{\begin{small}#1\end{small}}

\def\Ga{\Gamma}

\def\Si{\Sigma}

\def\al{\alpha}

\def\ga{\gamma}
\def\io{\iota}
\def\ka{\kappa}
\def\la{\lambda}

\def\si{\sigma}

\def\vph{\varphi}

\def\vr{\varrho}

\def\cA{\mathcal A}

\def\C{\mathbb C}
\def\cC{\mathcal C}

\def\bE{\mathbb E}

\def\nE{\tn{E}}

\def\bF{\mathbb{F}}

\def\cI{\mathcal I}

\def\fJ{\mathfrak j}

\def\cL{\mathcal L}

\def\cM{\mathcal M}
\def\cN{\mathcal N}

\def\cQ{\mathcal Q}

\def\P{\mathbb P}
\def\R{\mathbb{R}}

\def\cT{\mathcal T}

\def\nV{\textnormal{V}}
\def\cW{\mathcal W}

\def\Z{\mathbb{Z}}

\def\ff{\mathfrak f}

\def\fc{\mathfrak c}

\def\fq{\mathfrak q}

\def\au{\tn{au}}
\def\BL{\mathbb{BL}}
\def\BLau{\BL^{\!\!\fc_1}}
\def\BLauk{\BL^{\!\!\fc_1,k}}
\def\BLaua{\BL^{\!\!\fc_1,1}}
\def\BLaub{\BL^{\!\!\fc_1,2}}
\def\BLauba{\BL^{\!\!\fc_1,2;1}}
\def\BLC{\BL^{\!\!\C}}
\def\BLF{\BL^{\!\!\bF}}
\def\BLR{\BL^{\!\!\R}}
\def\CR{\tn{CR}}
\def\cha{\tn{char}}

\def\nd{\tn{d}}

\def\Edg{\tn{Edg}}
\def\End{\tn{End}}

\def\id{\tn{id}}

\def\nod{\tn{nd}}

\def\val{\tn{val}}
\def\Ver{\tn{Ver}}

\def\eset{\emptyset}
\def\i{\infty}
\def\prt{\partial}

\def\bu{\bullet}
\def\dag{\dagger}

\begin{document}

\title{Blowdowns of the Deligne-Mumford Spaces\\
of Real Rational Curves} 
\author{Xujia Chen\thanks{Partially supported by NSF grant DMS 1901979 and Simons Foundation}
~and Aleksey Zinger\thanks{Partially supported by NSF grant DMS 1901979 and Simons Foundation}}
\date{\today}
\maketitle

\begin{abstract}
\noindent   
We describe a sequence of smooth quotients of the Deligne-Mumford moduli 
space~$\R\ov\cM_{0,\ell+1}$ of real rational curves with $\ell\!+\!1$ 
conjugate pairs of marked points that terminates 
at~\hbox{$\R\ov\cM_{0,\ell}\!\times\!\C\P^1$}.
This produces an analogue of Keel's blowup construction
of the Deligne-Mumford moduli spaces~$\ov\cM_{\ell+1}$ of  
rational curves with $\ell\!+\!1$~marked points,
but with an explicit description of the intermediate spaces and
the blowups of three different types. 
The same framework readily adapts to the real moduli spaces with real points.
In a sequel, we use this inductive construction of~$\R\ov\cM_{0,\ell+1}$ to
completely determine the rational (co)homology ring of~$\R\ov\cM_{0,\ell}$.
\end{abstract}

\tableofcontents

\section{Introduction}
\label{intro_sec}

\noindent
The study of the topology of the Deligne-Mumford moduli spaces of 
stable {\it complex} marked curves has led to many remarkable results
over the past three decades.
On the other hand, very little is known about the topology of their {\it real} counterparts,
even in genus~0.
The most advanced result in this direction is perhaps the determination
of the cohomology of the Deligne-Mumford moduli space of 
stable real rational curves with real marked points in~\cite{EHKR}.
The aim of the present paper and its sequel~\cite{RDMhomol} is
to invigorate the study of the topology of 
the Deligne-Mumford moduli spaces of stable real marked curves.\\

\noindent
A key step in the determination of the cohomology of the Deligne-Mumford moduli space~$\ov\cM_{\ell}$ 
of stable (complex) rational curves with $\ell$~marked points $z_1,\ldots,z_{\ell}$ 
in~\cite{Keel} is a sequence of {\it inductively} constructed (projective) blowups 
of~$\ov\cM_{\ell}\!\times\!\ov\cM_4$ terminating at~$\ov\cM_{\ell+1}$.
These blowups are functorial in the complex analytic and projective categories.
For illustrative purposes, we first recast the intermediate blowup spaces of~\cite{Keel}
as {\it explicit} quotients of~$\ov\cM_{\ell+1}$ terminating 
at~$\ov\cM_{\ell}\!\times\!\ov\cM_4$ and show that consecutive quotients are related
by a single blowup along a smooth subvariety of complex codimension~3; see Theorem~\ref{Cblowup_thm}.
We then adapt this approach to the Deligne-Mumford moduli spaces of stable real 
curves with conjugate pairs of marked points.
The statement of Theorem~\ref{Rblowup_thm} thoroughly describes 
distinguished submanifolds of the intermediate blowups~$X_{\vr}$
that correspond to the top ``boundary" strata of the Deligne-Mumford moduli spaces.
This theorem is used in~\cite{RDMhomol} to show that such submanifolds 
generate the cohomology of~$X_{\vr}$.
A separate argument in~\cite{RDMhomol} then shows that there is 
only one ``interesting" relation among these generators in 
the rational cohomologies of the Deligne-Mumford moduli spaces, 
the one discovered in~\cite{RealEnum}.
In contrast to the complex case, the real case involves three different types of
blowups: real, complex, and ``augmented".
The blowups of the first type are functorial in the smooth category, 
but the blowups of the last two types depend on choices of compatible collections of 
coordinate charts;
in the case of Theorem~\ref{Rblowup_thm}, they arise from the cross ratio of four points 
on the Riemann sphere.
As none of the three types of blowups descends to the homeomorphism classes of smooth manifolds,
the smooth structures on the quotients~$X_{\vr}$ play a material role in
the proof of Theorem~\ref{Rblowup_thm}, even from the topological viewpoint.

\subsection{The complex case}
\label{Cintro_subs}

\noindent
The Deligne-Mumford moduli space~$\ov\cM_{\ell}$ is 
a smooth complex projective variety of (complex) dimension~$\ell\!-\!3$.
There are forgetful morphisms
\BE{ffCdfn_e}\ff_{\ell+1}\!: \ov\cM_{\ell+1}\lra\ov\cM_{\ell} 
\qquad\hbox{and}\qquad
\ff_{123,\ell+1}\!: \ov\cM_{\ell+1}\lra\ov\cM_4\EE
dropping the marked point~$z_{\ell+1}$ and 
all marked points other than $z_1,z_2,z_3,z_{\ell+1}$, respectively.
Section~1 in~\cite{Keel} decomposes the holomorphic map 
\BE{CMlind_e}\big(\ff_{\ell+1},\ff_{123,\ell+1}\big)\!: 
\ov\cM_{\ell+1}\lra\ov\cM_{\ell}\!\times\!\ov\cM_4\EE
as a sequence of inductively defined (\sf{holomorphic}) 
blowups $\pi_{\vr}\!:X_{\vr}\!\lra\!X_{\vr-1}$.\\

\noindent
The spaces $X_{\vr}$ can be explicitly described as quotients of~$\ov\cM_{\ell+1}$.
Let 
$$[\ell]=\big\{1,2,\ldots,\ell\} \qquad\hbox{and}\qquad
\cA_{\ell}=\big\{\vr\!\subset\![\ell]\!: \big|\vr\!\cap\![3]\!\big|\!\ge\!2,~
\big|[\ell]\!-\!\vr\big|\!\ge\!2\big\}.$$
The set~$\cA_{\ell}$ is partially ordered by the inclusion~$\subsetneq$ of subsets of~$[\ell]$.
We extend this partial order to a strict order~$<$ on $\{0\}\!\sqcup\!\cA_{\ell}$
so that~0 is the smallest element and define \hbox{$\vr\!-\!1\!\in\!\{0\}\!\cup\!\cA_{\ell}$}
to be the predecessor of $\vr\!\in\!\cA_{\ell}$.
Let $\vr_{\max}\!\in\!\cA_{\ell}$ be the largest element with respect to~$<$.\\

\noindent
For $\vr\!\in\![\ell]$, we denote by 
$$D_{\ell;\vr}\subset\ov\cM_{\ell}$$
the smooth divisor parametrizing the $\ell$-marked curves~$\cC$ 
with a node $\nod_{\vr}(\cC)$ that separates~$\cC$ into two topological components, 
$\cC_{\vr}'$ and~$\cC_{\vr}''$,
so that $\cC_{\vr}'$ (resp.~$\cC_{\vr}''$) carries
the marked points indexed by~$\vr$ (resp.~$[\ell]\!-\!\vr$).
Every ``boundary" divisor of~$\ov\cM_{\ell}$ is of this form for 
a unique $\vr\!\in\!\cA_{\ell}$.
The complement of the union of the ``boundary" divisors~$D_{\ell;\vr}$ 
with $\vr\!\in\!\cA_{\ell}$ is the open subspace \hbox{$\cM_{\ell}\!\subset\!\ov\cM_{\ell}$}
parametrizing smooth curves.
The divisors~$D_{4;\vr}\!\subset\!\ov\cM_4$ are points and represent the three curves
in Figure~\ref{M04rel_fig}.
We note~that
\begin{alignat}{3}
\label{ndCprp_e0} 
\vr,\vr'\in\cA_{\ell} ~~&\hbox{and}~~ D_{\ell;\vr}\!\cap\!D_{\ell;\vr'}\neq\eset
&\quad&\Lra\quad
\vr\supset\vr', &~~&\hbox{or}~~\vr\subset\vr',
~~\hbox{or}~~\vr\supset[\ell]\!-\!\vr',\\
\label{ndCprp_e} 
\vr,\vr'\in\cA_{\ell} ~~&\hbox{and}~~
D_{\ell+1;\vr}\!\cap\!D_{\ell+1;\vr'}\neq\eset
&\quad&\Lra\quad \vr\supset\vr' &~~&\hbox{or}~~ \vr\subset\vr'\,;
\end{alignat}
this follows from the assumption that $\vr$ and $\vr'$ contain 
at least two of the elements of the set $\{1,2,3\}$ and thus overlap.\\

\begin{figure}
\begin{pspicture}(-2.8,-1.7)(10,1.2)
\psset{unit=.3cm}
\pscircle*(2,0){.2}
\psline[linewidth=.07](1,-1)(5,3)\psline[linewidth=.07](1,1)(5,-3)
\pscircle*(3,1){.2}\pscircle*(4,2){.2}\pscircle*(3,-1){.2}\pscircle*(4,-2){.2}
\rput(3,2){\sm{$1$}}\rput(4,3){\sm{$2$}}\rput(3,-2){\sm{$3$}}\rput(4,-3){\sm{$4$}}
\rput(4,-5){\sm{$D_{4;\{1,2\}}$}}
\pscircle*(17,0){.2}
\psline[linewidth=.07](16,-1)(20,3)\psline[linewidth=.07](16,1)(20,-3)
\pscircle*(18,1){.2}\pscircle*(19,2){.2}\pscircle*(18,-1){.2}\pscircle*(19,-2){.2}
\rput(18,2){\sm{$1$}}\rput(19,3){\sm{$3$}}\rput(18,-2){\sm{$2$}}\rput(19,-3){\sm{$4$}}
\rput(19,-5){\sm{$D_{4;\{1,3\}}$}}
\pscircle*(32,0){.2}
\psline[linewidth=.07](31,-1)(35,3)\psline[linewidth=.07](31,1)(35,-3)
\pscircle*(33,1){.2}\pscircle*(34,2){.2}\pscircle*(33,-1){.2}\pscircle*(34,-2){.2}
\rput(33,2){\sm{$2$}}\rput(34,3){\sm{$3$}}\rput(33,-2){\sm{$1$}}\rput(34,-3){\sm{$4$}}
\rput(34,-5){\sm{$D_{4;\{2,3\}}$}}
\end{pspicture}
\caption{The three divisors $D_{4;\vr}$ in $\ov\cM_4\!\approx\!\C\P^1$.}
\label{M04rel_fig}
\end{figure}

\noindent
For $\vr\!\in\!\cA_{\ell}$ and $\wt\cC_1,\wt\cC_2\!\in\!\ov\cM_{\ell+1}$, 
we define $\wt\cC_1\!\sim_{\vr}\!\wt\cC_2$ if either $\wt\cC_1\!=\!\wt\cC_2$ or 
$$\ff_{\ell+1}(\wt\cC_1)=\ff_{\ell+1}(\wt\cC_2)\in \ov\cM_{\ell}
\quad\hbox{and}\quad \wt\cC_1,\wt\cC_2\in D_{\ell+1;\vr}.$$
By~\eref{ndCprp_e},
for every $\vr^*\!\in\!\{0\}\!\sqcup\!\cA_{\ell}$
the union of all equivalence relations~$\sim_{\vr}$ on~$\ov\cM_{\ell+1}$ with 
$\vr\!>\!\vr^*$ is again an equivalence relation.
We denote by~$X_{\vr^*}$ the quotient of~$\ov\cM_{\ell+1}$ by the last equivalence relation
and by~$[\wt\cC]_{\vr^*}$ the corresponding equivalence class of $\wt\cC\!\in\!\ov\cM_{\ell+1}$.
Since the quotient projection
\BE{qvrdfn_e}q_{\vr^*}\!:\ov\cM_{\ell+1}\lra X_{\vr^*}, \qquad 
q_{\vr^*}(\wt\cC)=[\wt\cC]_{\vr^*},\EE
is a closed map, $X_{\vr^*}$ is a compact Hausdorff space
for every \hbox{$\vr^*\!\in\!\{0\}\!\sqcup\!\cA_{\ell}$}.\\

\noindent
By definition, $X_{\vr_{\max}}\!=\!\ov\cM_{\ell+1}$.
The holomorphic map~\eref{CMlind_e} descends to a continuous bijection
$$\Psi_0\!:X_0\lra  \ov\cM_{\ell}\!\times\!\ov\cM_4\,.$$
Since $X_0$ is compact and $\ov\cM_{\ell}\!\times\!\ov\cM_4$ is Hausdorff, 
$\Psi_0$ is a homeomorphism.
By Theorem~\ref{Cblowup_thm}\ref{Cspaces_it}, 
the quotients~$X_{\vr^*}$ are complex manifolds so that
the map~\eref{qvrdfn_e} is holomorphic and the subspaces 
\BE{Yvrdfn_e}Y_{\vr^*;\vr}^0\equiv q_{\vr^*}\big(D_{\ell+1;\vr}\big)
~~\hbox{and}~~
Y_{\vr^*;\vr}^+\equiv q_{\vr^*}\big(D_{\ell+1;\vr\cup\{\ell+1\}}\big)
\qquad\hbox{with}~~\vr\!\subset\![\ell]\EE
of~$X_{\vr^*}$ are complex submanifolds.
The former implies that~$\Psi_0$ is a biholomorphism.
Each subspace $Y_{0;\vr}^0\!\subset\!X_0$ with $\vr\!\in\!\cA_{\ell}$ is a section~of
$$\pi_1\!\circ\!\Psi_0\!:X_0\lra\ov\cM_{\ell}\!\times\!\ov\cM_4 \lra\ov\cM_{\ell}$$
over the divisor $D_{\ell;\vr}\!\subset\!\ov\cM_{\ell}$.
By Theorem~\ref{Cblowup_thm}\ref{Cblowup_it}, the continuous map
\BE{bldowndfn_e}\pi_{\vr^*}\!:X_{\vr^*}\lra X_{\vr^*-1}\EE
with $\vr^*\!\in\!\cA_{\ell}$ induced by $q_{\vr^*-1}$
is the holomorphic blowup along the codimension~2 complex submanifold 
\hbox{$Y_{\vr^*-1;\vr^*}^0\!\subset\!X_{\vr^*-1}$}
with the exceptional divisor \hbox{$Y_{\vr^*;\vr^*}^0\!\subset\!X_{\vr^*}$}.
In particular,
every ``boundary" divisor $D_{\ell;\vr}\!\subset\!\ov\cM_{\ell}$ corresponds 
to a unique blowup locus~$Y_{\vr-1;\vr}^0$ in this case.
By the definition of the equivalence relations on~$\ov\cM_{\ell+1}$,
$Y_{\vr-1;\vr}^0$ is
the proper transform of~$Y_{0;\vr}^0$ under the earlier blowups.\\

\noindent
Thus, Theorem~\ref{Cblowup_thm} decomposes~\eref{CMlind_e}
into a sequence 
$$\ov\cM_{\ell+1}\!=\!X_{\vr_{\max}}\xlra{\pi_{\vr_{\max}}} 
X_{\vr_{\max}-1}\xlra{\pi_{\vr_{\max}-1}}\ldots
\lra X_0=\ov\cM_{\ell}\!\times\!\ov\cM_4$$
of holomorphic blowups with the steps indexed by the ordered set $(\cA_{\ell},<)$.
By~\eref{ndCprp_e}, 
different choices of the order~$<$ extending~$\subsetneq$ affect
this sequence by interchanging the order of blowups along disjoint loci.
The proof of this theorem in Section~\ref{Cblowup_sec} revolves around 
the cross ratio of four points on the Riemann sphere and
provides a non-inductive alternative to the inductive blowup construction
in~\cite{Keel}.

\subsection{The real case}
\label{Rintro_subs}

\noindent
The Deligne-Mumford moduli space~$\R\ov\cM_{k,\ell}$ of stable real rational curves 
with $k$~real marked points $z_1,\ldots,z_k$
and $\ell$~conjugate pairs $(z_1^+,z_1^-),\ldots,(z_{\ell}^+,z_{\ell}^-)$
of marked points is a smooth compact manifold 
of (real) dimension $k\!+\!2\ell\!-\!3$.
It is orientable if and only if either $k\!=\!0$ or $k\!+\!2\ell\!\le\!4$;
see Proposition~5.7 and Corollary~6.2 in~\cite{Cey} or
Proposition~1.5 in \cite{Penka2}.
It contains top ``boundary" strata of real codimensions~1 and~2,
which we call~\sf{boundary hypersurfaces} and \sf{boundary divisors}, respectively.
The normal bundles to the boundary divisors are orientable and thus 
admit complex structures.
There are two types of boundary hypersurfaces, $E$~and~$H$.
The normal bundles to the $E$-hypersurfaces are orientable (if $k\!>\!0$, there are no such strata); 
the normal bundles to the $H$-hypersurfaces are not orientable if $k\!+\!2\ell\!\ge\!5$
(so that they are of positive dimensions).\\

\noindent
If $\ell\!\ge\!2$, there are forgetful morphisms
\BE{ffRdfn_e}\ff_{(\ell+1)^{\pm}}^{\R}\!:\R\ov\cM_{k,\ell+1}\lra\R\ov\cM_{k,\ell} 
\quad\hbox{and}\quad
\ff_{1^{\pm}2^+(\ell+1)^+}^{\R}\!:\R\ov\cM_{k,\ell+1}\lra\ov\cM_4\EE
dropping the marked points $z_{\ell+1}^+,z_{\ell+1}^-$ and 
all marked points other than $z_1^+,z_1^-,z_2^+,z_{\ell+1}^+$,
respectively.
In this paper, we decompose the smooth map 
\BE{RMlind_e}\big(\ff_{(\ell+1)^{\pm}}^{\R},\ff_{1^{\pm}2^+(\ell+1)^+}^{\R}\big)\!: 
\R\ov\cM_{0,\ell+1}\lra \R\ov\cM_{0,\ell}\!\times\!\ov\cM_4\EE
as a sequence of blowups $\pi_{\vr}\!:X_{\vr}\!\lra\!X_{\vr-1}$ of three different types. 
The type of each blowup,  \sf{real}, \sf{complex}, or \sf{augmented}
as defined in Section~\ref{blowup_sec}, 
is determined by whether its blowup locus~$Y_{\vr-1;\vr}^-$
is associated with an $H$-hypersurface, boundary divisor, 
or $E$-hypersurface $D_{\ell;\vr}\!\subset\!\R\ov\cM_{0,\ell}$.\\

\noindent
For $\ell\!\in\!\Z^+$, we define
\BE{cAellpmdfn_e}
\big[\ell^{\pm}\big]=\big\{1^+,1^-,\ldots,\ell^+,\ell^-\big\}
~~\hbox{and}~~
\cA_{\ell}^{\pm}=\big\{\vr\!\subset\![\ell^{\pm}]\!:
\big|\vr\!\cap\!\{1^+,1^-,2^+\}\big|\!\ge\!2,\,\big|[\ell^{\pm}]\!-\!\vr\big|\!\ge\!2\big\}.\EE
The collection~$\cA_{\ell}^{\pm}$ is partially ordered by the inclusion~$\subsetneq$ 
of subsets of~$[\ell^{\pm}]$.
We extend this partial order to a strict order~$<$ on $\{0\}\!\sqcup\!\cA_{\ell}^{\pm}$
so that~0 is the smallest element and define $\vr\!-\!1\!\in\!\{0\}\!\cup\!\cA_{\ell}^{\pm}$
to be the predecessor of $\vr\!\in\!\cA_{\ell}^{\pm}$.
Let $\vr_{\max}\!\in\!\cA_{\ell}^{\pm}$ be the largest element with respect to~$<$.\\
  
\noindent
For $\vr\!\in\![\ell^{\pm}]$, we denote by 
$$D_{\ell;\vr}\subset\R\ov\cM_{0,\ell}$$
the subspace parametrizing the $[\ell^{\pm}]$-marked curves~$\cC$ with a node $\nod_{\vr}(\cC)$ 
that separates~$\cC$ into two topological components, 
$\cC_{\vr}'$ and~$\cC_{\vr}''$, so that $\cC_{\vr}'$ (resp.~$\cC_{\vr}''$) carries
the marked points indexed by~$\vr$ (resp.~$[\ell^{\pm}]\!-\!\vr$).
This subspace is empty unless~$\vr$ is an element of the subcollection 
$\cA_{\ell}^{\R}\!\subset\!\cA_{\ell}^{\pm}$ defined in Section~\ref{Rthm_subs}.
Every top ``boundary" stratum of~$\R\ov\cM_{0,\ell}$ equals~$D_{\ell;\vr}$
some $\vr\!\in\!\cA_{\ell}^{\R}$, a unique one in the case of the hypersurfaces
and with two possible choices in the case of the boundary divisors.
The complement of the union of the ``boundary" divisors~$D_{\ell;\vr}$ 
with $\vr\!\in\!\cA_{\ell}^{\R}$ is the open subspace 
\hbox{$\R\cM_{0,\ell}\!\subset\!\R\ov\cM_{0,\ell}$}
parametrizing smooth curves.
The three top boundary strata $D_{2;\vr}\!\subset\!\R\ov\cM_{0,2}$
are divisors;
all three are points and represent the first three curves in Figure~\ref{TopolRel_fig}.
There are six boundary hypersurfaces $D_{3;\vr}\!\subset\!\R\ov\cM_{0,3}$, 
each of which is diffeomorphic to the real projective space~$\R\P^1$;
two of them are represented by the last two diagrams in Figure~\ref{TopolRel_fig}.\\

\begin{figure}
\begin{pspicture}(-.5,1.8)(10,5.4)
\psset{unit=.4cm}
\psline[linewidth=.05](2,9)(6,13)\psline[linewidth=.05](2,11)(6,7)
\pscircle*(3.5,10.5){.15}\pscircle*(5,12){.15}
\pscircle*(3.5,9.5){.15}\pscircle*(5,8){.15}
\rput(3.2,8.9){\sm{$1^-$}}\rput(4.5,7.5){\sm{$2^-$}}
\rput(3.5,11.1){\sm{$1^+$}}\rput(5,12.6){\sm{$2^+$}}
\psline[linewidth=.025]{<->}(5.5,8.5)(5.5,11.5)\rput(6,10){$\si$}
\rput(4.2,5.3){\sm{$D_{2;\{1^+,2^+\}}$}}
\psline[linewidth=.05](8,9)(12,13)\psline[linewidth=.05](8,11)(12,7)
\pscircle*(9.5,10.5){.15}\pscircle*(11,12){.15}
\pscircle*(9.5,9.5){.15}\pscircle*(11,8){.15}
\rput(9.1,8.9){\sm{$1^+$}}\rput(10.5,7.5){\sm{$2^-$}}
\rput(9.5,11.1){\sm{$1^-$}}\rput(11,12.6){\sm{$2^+$}}
\psline[linewidth=.025]{<->}(11.5,8.5)(11.5,11.5)\rput(12,10){$\si$}
\rput(10.2,5.3){\sm{$D_{2;\{1^-,2^+\}}$}}
\pscircle[linewidth=.05](17,10){2}\pscircle[linewidth=.05](21,10){2}
\psline[linewidth=.025]{<->}(17,8.5)(17,11.5)\rput(17.5,10){$\si$}
\psline[linewidth=.025]{<->}(21,8.5)(21,11.5)\rput(21.5,10){$\si$}
\pscircle*(15.59,11.41){.15}\pscircle*(15.59,8.59){.15}
\pscircle*(22.41,11.41){.15}\pscircle*(22.41,8.59){.15}
\rput(15.4,12){\sm{$1^+$}}\rput(15.4,8){\sm{$1^-$}}
\rput(23,12){\sm{$2^+$}}\rput(23,8){\sm{$2^-$}}
\rput(19,5.3){\sm{$D_{2;\{1^+,1^-\}}$}}
\psline[linewidth=.05](27,12.5)(27,7.5)
\psline[linewidth=.02](26.5,12)(30,12)\psline[linewidth=.02](26.5,8)(30,8)
\pscircle*(28,12){.15}\pscircle*(29,12){.15}
\pscircle*(27,11){.15}\pscircle*(27,9){.15}
\rput(28,12.7){\sm{$1^+$}}\rput(29.4,12.7){\sm{$2^+$}}
\rput(26.4,11.3){\sm{$3^+$}}\rput(26.4,9.3){\sm{$3^-$}}
\pscircle*(28,8){.15}\pscircle*(29,8){.15}
\rput(28,7.3){\sm{$1^-$}}\rput(29.4,7.3){\sm{$2^-$}}
\psline[linewidth=.025]{<->}(29.5,8.6)(29.5,11.4)\rput(30,10){$\si$}
\rput(28,5.3){\sm{\begin{tabular}{c}$D_{3;\{1^+,2^+\}}$\\
$D_{3;\{1^+,2^+,3^+,3^-\}}$\end{tabular}}}
\psline[linewidth=.05](33,12.5)(33,7.5)
\psline[linewidth=.02](32.5,12)(36,12)\psline[linewidth=.02](32.5,8)(36,8)
\pscircle*(34,12){.15}\pscircle*(35,12){.15}
\pscircle*(33,11){.15}\pscircle*(33,9){.15}
\rput(34,12.7){\sm{$2^+$}}\rput(35.4,12.7){\sm{$3^-$}}
\rput(32.4,11.3){\sm{$1^+$}}\rput(32.4,9.3){\sm{$1^-$}}
\pscircle*(34,8){.15}\pscircle*(35,8){.15}
\rput(34,7.3){\sm{$2^-$}}\rput(35.4,7.3){\sm{$3^+$}}
\psline[linewidth=.025]{<->}(35.5,8.6)(35.5,11.4)\rput(36,10){$\si$}
\rput(34.5,5.3){\sm{\begin{tabular}{c}$D_{3;\{1^+,1^-,2^+,3^-\}}$\\
$D_{3;\{1^+,1^-,2^-,3^+\}}$\end{tabular}}}
\end{pspicture}
\caption{The three boundary hypersurfaces~$D_{2;\vr}$ in $\R\ov\cM_{0,2}\!\approx\!\R\P^1$
and two of the six boundary divisors~$D_{3;\vr}$ in the orientable threefold~$\R\ov\cM_{0,3}$.
Each line or circle represents~$\C\P^1$. 
Each double-headed arrow labeled~$\si$ indicates
the involution on the corresponding real curve.}
\label{TopolRel_fig}
\end{figure}

\noindent
For $\vr\!\subset\![\ell^{\pm}]$, let
$$\wt{D}_{\vr}''=D_{\ell+1;\vr}\!\cup\!D_{\ell+1;\vr\cup\{(\ell+1)^-\}}
\subset\R\ov\cM_{0,\ell+1}.$$
Similarly to~\eref{ndCprp_e},
\BE{ndRprp_e} 
\vr,\vr'\in\cA_{\ell}^{\pm} ~~\hbox{and}~~ 
\wt{D}_{\vr}''\!\cap\!\wt{D}_{\vr'}''\neq\eset
\quad\Lra\quad \vr\supset\vr' ~~\hbox{or}~~ \vr\subset\vr'\,.\EE
For $\vr\!\in\!\cA_{\ell}^{\pm}$ and $\wt\cC_1,\wt\cC_2\!\in\!\R\ov\cM_{0,\ell+1}$, 
we define $\wt\cC_1\!\sim_{\vr}\!\wt\cC_2$ if either $\wt\cC_1\!=\!\wt\cC_2$ or 
$$\ff_{(\ell+1)^{\pm}}^{\R}(\wt\cC_1)=\ff_{(\ell+1)^{\pm}}^{\R}(\wt\cC_2)\in\R\ov\cM_{0,\ell}
\quad\hbox{and}\quad \wt\cC_1,\wt\cC_2\in\wt{D}_{\vr}''.$$
By~\eref{ndRprp_e},
for every $\vr^*\!\in\!\{0\}\!\sqcup\!\cA_{\ell}^{\pm}$
the union of all equivalence relations~$\sim_{\vr}$ on~$\R\ov\cM_{0,\ell+1}$ with 
$\vr\!>\!\vr^*$ is an equivalence relation.
We denote by~$X_{\vr^*}$ the quotient of~$\R\ov\cM_{0,\ell+1}$ by the last equivalence relation
and by~$[\wt\cC]_{\vr^*}$ the corresponding equivalence class of
 $\wt\cC\!\in\!\R\ov\cM_{0,\ell+1}$.
Since the quotient projection
\BE{qvrdfnR_e}q_{\vr^*}\!:\R\ov\cM_{0,\ell+1}\lra X_{\vr^*}, \qquad 
q_{\vr^*}(\wt\cC)=[\wt\cC]_{\vr^*},\EE
is a closed map, $X_{\vr^*}$ is again a compact Hausdorff space
for every \hbox{$\vr^*\!\in\!\{0\}\!\sqcup\!\cA_{\ell}^{\pm}$}.\\

\noindent
By definition, $X_{\vr_{\max}}\!=\!\R\ov\cM_{0,\ell+1}$.
The smooth map~\eref{RMlind_e} descends to a continuous bijection
$$\Psi_0\!:X_0\lra  \R\ov\cM_{0,\ell}\!\times\!\ov\cM_4.$$
Since $X_0$ is compact and $\R\ov\cM_{0,\ell}\!\times\!\ov\cM_4$ is Hausdorff, 
$\Psi_0$ is a homeomorphism and thus $X_0$ is orientable.
By Theorem~\ref{Rblowup_thm}\eref{Rspaces_it}, 
the quotients~$X_{\vr^*}$ are smooth manifolds so that
the map~\eref{qvrdfnR_e} is smooth and the subspaces 
$$Y_{\vr^*;\vr}^0\equiv q_{\vr^*}(\wt{D}_{\vr}'')\subset X_{\vr^*}
\qquad\hbox{with}~~\vr\!\in\!\cA_{\ell}^{\R},~\vr\!>\!\vr^*,$$
are smooth closed submanifolds.
By the proof of Theorem~\ref{Rblowup_thm}\eref{Rspaces_it},
$\Psi_0$ is a diffeomorphism.
Each submanifold $Y_{0;\vr}^0\!\subset\!X_0$ is a section~of
$$\pi_1\!\circ\!\Psi_0\!:X_0\lra\R\ov\cM_{0,\ell}\!\times\!\ov\cM_4 
\lra\R\ov\cM_{0,\ell}$$
over the top ``boundary" stratum $D_{\ell;\vr}\!\subset\!\R\ov\cM_{0,\ell}$.
By Theorem~\ref{Rblowup_thm}\eref{Rblowup_it}, the continuous map
\BE{Rbldowndfn_e}\pi_{\vr^*}\!:X_{\vr^*}\lra X_{\vr^*-1}\EE
with $\vr^*\!\in\!\cA_{\ell}^{\R}$ induced by $q_{\vr^*-1}$
is a real, complex,  or augmented blowup 
along the  submanifold \hbox{$Y_{\vr^*-1;\vr^*}^0\!\subset\!X_{\vr^*-1}$}
with the exceptional locus 
\hbox{$Y_{\vr^*;\vr^*}^0\!\cup\!Y_{\vr^*;\vr^*}^-\!\subset\!X_{\vr^*}$}.
In particular, every boundary hypersurface \hbox{$D_{\ell;\vr}\!\subset\R\ov\cM_{0,\ell}$} 
corresponds 
to a unique blowup locus~$Y_{\vr-1;\vr}^0$ in this case,
while every boundary divisor~$D_{\ell;\vr}$ corresponds to two blowups loci.
By the definition of the equivalence relations on~$\R\ov\cM_{0,\ell+1}$,
$Y_{\vr-1;\vr}^0$ is the proper transform of~$Y_{0;\vr}^0$ under the earlier blowups
in most, but not all, cases
(the exceptional cases correspond to the elements~$\vr$
of the subcollection \hbox{$\cA_{\ell;2}^D\!\subset\!\cA_{\ell}^{\R}$}
defined in Section~\ref{Rthm_subs}).\\

\noindent
Thus, Theorem~\ref{Rblowup_thm} decomposes~\eref{RMlind_e}
into a sequence 
$$\R\ov\cM_{0,\ell+1}\!=\!X_{\vr_{\max}}\xlra{\pi_{\vr_{\max}}} 
X_{\vr_{\max}-1}\xlra{\pi_{\vr_{\max}-1}}\ldots
\lra X_0=\R\ov\cM_{0,\ell}\!\times\!\ov\cM_4$$
of blowups of the three types defined in Section~\ref{blowup_sec}.
By~\eref{ndRprp_e}, 
different choices of the order~$<$ extending~$\subsetneq$ affect
this sequence by interchanging the order of blowups along disjoint loci.
Theorem~\ref{Rblowup_thm} is proved in Section~\ref{Rblowup_sec}
also using the cross ratio of four points on the Riemann sphere.
The proof in this case is significantly more technical than in the complex case
because there are now three different types of blowups and five types
of blowup loci involved.

\section{Main statements}
\label{thm_sec}

\noindent
Let $\ell\!\in\!\Z$. 
For $\vr\!\subset\![\ell]$, let $\vr_{\ell}^c\!=\![\ell]\!-\!\vr$.
For $\vr\!\subset\![\ell^{\pm}]$, define
$$ \vr_{\ell^{\pm}}^c=\big[\ell^{\pm}\big]\!-\!\vr \quad\hbox{and}\quad
\ov\vr=\big\{i^+\!:i^-\!\in\!\vr\big\}\!\cup\!\big\{i^-\!:i^+\!\in\!\vr\big\}.$$

\vspace{.15in}

\noindent
We say that smooth submanifolds $Y,Y'\!\subset\!X$ of a smooth manifold
\sf{intersect cleanly} if 
\hbox{$Y\!\cap\!Y'\!\subset\!X$}  is also a submanifold~and 
$$T(Y\!\cap\!Y')=TY|_{Y\cap Y'}\!\cap\!TY'|_{Y\cap Y'}\subset TX|_{Y\cap Y'}.$$
This is equivalent to the condition that the homomorphism
$$\cN_{Y'}(Y\!\cap\!Y')\!\equiv\!\frac{TY'|_{Y\cap Y'}}{T(Y\!\cap\!Y')}
\lra \frac{TX|_{Y\cap Y'}}{TY|_{Y\cap Y'}}
\!\equiv\!\cN_XY|_{Y\cap Y'}$$
from the normal bundle of~$Y\!\cap\!Y'$ in~$Y'$ to 
the normal bundle of~$Y$ in~$X$ induced by the obvious inclusions is injective.
If $Y$ contains~$Y'$ or intersects~it transversely, then
$Y$ and~$Y'$ intersect cleanly, but these are just extreme examples of clean intersection.\\
 
\noindent
Suppose $\pi\!:\wt{X}\!\lra\!X$ is a smooth map and 
$Y'\!\subset\!X$ and $\wt{Y}'\!\subset\!\wt{X}$ are submanifolds.
If $R$ is a commutative ring with unity, 
we call $Y'$ and~$\wt{Y}'$ \sf{$(\pi,R)$-related} if
\begin{enumerate}[label=$\bullet$,leftmargin=*]

\item either~$Y'$ and of~$\wt{Y}'$ are oriented or $\cha(R)\!=\!2$ and 

\item $Y'$ and~$\wt{Y}'$ are compact and $\pi_*[\wt{Y}']\!=\![Y']\!\in\!H_*(X;R)$.

\end{enumerate}
If $\pi\!:\wt{X}\!\lra\!X$ is a blowup along a closed submanifold $Y\!\subset\!X$
as defined in Sections~\ref{StanBl_subs} and~\ref{AugblowupDfn_subs},
$\bE$ is its exceptional locus, and $Y'$ intersects~$Y$ cleanly,
we call $Y'$ and~$\wt{Y}'$ \sf{$\pi$-equivalent}
if $\pi|_{\wt{Y}'}$ is a bijection onto~$Y'$, 
$\wt{Y}'$ intersects~$\bE$ cleanly, and 
$$T_{\wt{x}}\wt{Y}'\cap\ker\nd_{\wt{x}}\big\{\pi|_{\bE}\big\}=\{0\}
\quad\forall\,\wt{x}\!\in\!\bE\!\cap\!\wt{Y}'.$$
These conditions imply that $\pi|_{\wt{Y}'}$ is a diffeomorphism onto~$Y'$
and thus the submanifolds~$Y'$ and~$\wt{Y}'$ are $(\pi,R)$-related for 
any admissible ring~$R$, provided either is compact.
If $Y'$ intersects~$Y$ cleanly without a topological component contained in~$Y$,
then its proper transform~$\wt{Y}'$ under the blowup is  $(\pi,R)$-related to~$Y'$,
but usually is not $\pi$-equivalent.

\subsection{The complex case}
\label{Cthm_subs}

\noindent
With $\ell\!\in\!\Z^+$, let
$$\wt\cA_{\ell}=\cA_{\ell}\!\sqcup\!\big\{[\ell]\!-\!\{i\}\!:i\!\in\![\ell]\big\}.$$
For \hbox{$\vr^*\!\in\!\{0\}\!\sqcup\!\cA_{\ell}$}, define 
$$\cA_{\ell}(\vr^*)=\big\{\vr\!\in\!\cA_{\ell}\!:\vr\!>\!\vr^*\big\},
\qquad 
\wt\cA_{\ell}(\vr^*)=\cA_{\ell}(\vr^*)\!\sqcup\!\big\{\![\ell]\!-\!\{i\}\!:i\!\in\![\ell]\big\}
\subset\wt\cA_{\ell}.$$
With the notation as in~\eref{Yvrdfn_e},
\BE{Yvrmp_e}
Y_{\vr^*;\vr}^0=q_{\vr^*}
\big(D_{\ell+1;\vr\cup\{\ell+1\}}\!\cap\!D_{\ell+1;\vr}\big)
=  Y_{\vr^*;\vr}^+\!\cap\!Y_{\vr^*;[\ell]-\{i\}}^0
\quad\forall~\vr\!\in\!\cA_{\ell}(\vr^*),~i\!\in\![\ell]\!-\!\vr.\EE
The next theorem is proved in Section~\ref{Cblowup_sec},
with Proposition~\ref{Cblowup_prp} playing a central role.

\begin{thm}\label{Cblowup_thm}
Suppose $\ell\!\in\!\Z^+$ with $\ell\!\ge\!3$.
\begin{enumerate}[label=($\C\arabic*$),leftmargin=*]

\item\label{Cspaces_it} For every $\vr^*\!\in\!\{0\}\!\sqcup\!\cA_{\ell}$,
the space $X_{\vr^*}$ is a complex manifold so that the quotient projection~$q_{\vr^*}$
in~\eref{qvrdfn_e} is holomorphic and 
the compact subspaces $Y_{\vr^*;\vr}^+,Y_{\vr^*;\vr}^0\!\subset\!X_{\vr^*}$ 
with $\vr\!\in\!\wt\cA_{\ell}$ are complex submanifolds.

\item\label{Cblowup_it} If $\vr^*\!\in\!\cA_{\ell}$, 
the map~$\pi_{\vr^*}$ in~\eref{bldowndfn_e}
is the holomorphic blowup of~$X_{\vr^*-1}$ along~$Y_{\vr^*-1;\vr^*}^0$
with the exceptional divisor~$Y_{\vr^*;\vr^*}^0$.

\item\label{Csubspaces_it} If $\vr^*\!\in\!\cA_{\ell}$, 
the submanifolds $Y_{\vr^*;\vr}^{\bu}\!\subset\!X_{\vr^*}$ 
with \hbox{$\vr\!\in\!\wt\cA_{\ell}$}, \hbox{$\bu\!=\!+,0$}, and $(\vr,\bu)\!\neq\!(\vr^*,0)$ 
are $(\pi_{\vr^*},\Z)$-related to the corresponding submanifolds~$Y_{\vr^*-1;\vr}^{\bu}$.
The submanifolds \hbox{$Y_{\vr^*;\vr}^0\!\subset\!X_{\vr^*}$} with 
\hbox{$\vr\!\in\!\cA_{\ell}(\vr^*)$}
are $\pi_{\vr^*}$-equivalent to the corresponding submanifolds~$Y_{\vr^*-1;\vr}^0$.

\end{enumerate}
\end{thm}

\vspace{.15in}

\noindent
An implicit part of the last claim of~\ref{Csubspaces_it} is that 
the submanifold $Y_{\vr^*-1;\vr^*}^0\!\subset\!X_{\vr^*-1}$
intersects each submanifold~$Y_{\vr^*-1;\vr}^0$ with $\vr\!\in\!\cA_{\ell}(\vr^*)$
cleanly.
In fact,
the submanifolds $Y_{\vr^*;\vr}^+\!\subset\!X_{\vr^*}$ with $\vr\!\in\!\cA_{\ell}$
and $Y_{\vr^*;\vr}^0$ with $\vr\!\in\!\cA_{\ell}\!-\!\cA_{\ell}(\vr^*)$
intersect~$Y_{\vr^*;\vr'}^{\bu}$ with $\vr\!\neq\!\vr'$ transversely and thus cleanly
(if $\vr\!\not\in\!\cA_{\ell}(\vr^*)$, the submanifolds $Y_{\vr^*;\vr}^+$ and 
$Y_{\vr^*;\vr}^0$ are transverse as~well).
The submanifolds $Y_{\vr^*;[\ell]-\{i\}}^0\!\subset\!X_{\vr^*}$ with $i\!\in\![\ell]$
intersect each of the submanifolds~$Y_{\vr^*;\vr'}^{\bu}$ transversely or contain~it.
However, \hbox{$Y_{\vr^*;\vr}^0,Y_{\vr^*;\vr'}^0\!\subset\!X_{\vr^*}$} with 
$\vr,\vr'\!\in\!\cA_{\ell}(\vr^*)$ may not intersect cleanly.
For example, if $\ell\!=\!6$, 
\hbox{$Y_{0;1234}^0\!\cap\!Y_{0;1256}^0$} is a single point~$[\wt\cC]_0$,
but
$$T_{[\wt\cC]_0}Y_{0;1234}^0\cap T_{[\wt\cC]_0}Y_{0;1256}^0\neq\{0\}\,.$$

\subsection{The real case}
\label{Rthm_subs}

\noindent
For $\ell\!\in\!\Z^+$, define
\begin{gather*}
\cA_{\ell}^H=\big\{\vr\!\in\!\cA_{\ell}^{\pm}\!: \ov\vr\!=\!\vr\big\}, \quad
\cA_{\ell}^E=\big\{\vr\!\in\!\cA_{\ell}^{\pm}\!:\ov\vr\!=\!\vr_{\ell^{\pm}}^c\big\},\\
\cA_{\ell;1}^D=
\big\{\vr\!\in\!\cA_{\ell}^{\pm}\!:\ov\vr\!\subsetneq\!\vr_{\ell^{\pm}}^c\big\},
~\cA_{\ell;2}^D=\big\{\ov\vr_{\ell^{\pm}}^c\!:\vr\!\in\!\cA_{\ell;1}^D\big\}
\subset\cA_{\ell}^{\pm},
~\cA_{\ell;3}^D=
\big\{\vr\!\in\!\cA_{\ell}^{\pm}\!: \ov\vr\!\supsetneq\!\vr_{\ell^{\pm}}^c\big\}
\!-\!\cA_{\ell;2}^D,\\
\cA_{\ell}^D=\cA_{\ell;1}^D\!\cup\!\cA_{\ell;2}^D\!\cup\!\cA_{\ell;3}^D, \quad
\cA_{\ell}^{\R}=
\cA_{\ell}^H\!\cup\!\cA_{\ell}^E\!\cup\!\cA_{\ell}^D, \quad
\wt\cA_{\ell}^{\R}=\cA_{\ell}^{\R}\!\sqcup\!
\big\{\!\big[\ell^{\pm}\big]\!-\!\{i\}\!:i\!\in\!\big[\ell^{\pm}\big]\!\big\}.
\end{gather*}
If $\vr\!\in\!\cA_{\ell;3}^D$, then $\ov\vr\!\in\!\cA_{\ell;3}^D$ as well.
For \hbox{$\vr^*\!\in\!\{0\}\!\sqcup\!\cA_{\ell}^{\R}$}, let 
$$\cA_{\ell}^{\R}(\vr^*)=\big\{\vr\!\in\!\cA_{\ell}^{\R}\!:\vr\!>\!\vr^*\big\},
\qquad
\wt\cA_{\ell}^{\R}(\vr^*)=\cA_{\ell}^{\R}(\vr^*)\!\sqcup\!
\big\{\!\big[\ell^{\pm}\big]\!-\!\{i\}\!:i\!\in\!\big[\ell^{\pm}\big]\!\big\}
\subset \wt\cA_{\ell}^{\R}.$$

\vspace{.15in}

\noindent
The subspace $D_{\ell;\vr}\!\subset\!\R\ov\cM_{0,\ell}$ defined in Section~\ref{Rintro_subs}
is an $H$-hypersurface (resp.~$E$-hypersurface, boundary divisor) if and only if 
$\vr\!\in\!\cA_{\ell}^H$ (resp.~$\vr\!\in\!\cA_{\ell}^E$, $\vr\!\in\!\cA_{\ell}^D$);
see the first column of diagrams in Figure~\ref{DT_fig}.
If $\vr\!\in\!\cA_{\ell;1}^D\!\cup\!\cA_{\ell;2}^D$, 
then $D_{\ell;\vr}\!=\!D_{\ell;\ov\vr_{\ell^{\pm}}^c}$;
see the penultimate diagram in Figure~\ref{TopolRel_fig}.
If $\vr\!\in\!\cA_{\ell;3}^D$, then $D_{\ell;\vr}\!=\!D_{\ell;\ov\vr}$;
see the last diagram in Figure~\ref{TopolRel_fig}.
These are the only cases of distinct $\vr,\vr'\!\in\cA_{\ell}^{\R}$ 
with $D_{\ell;\vr}\!=\!D_{\ell;\vr'}$.\\

\noindent
We now describe top boundary strata 
$$\wt{D}_{\vr}^{\bu}\subset \big\{\ff_{(\ell+1)^{\pm}}^{\R}\big\}^{-1}(D_{\ell;\vr})
\subset\R\ov\cM_{0,\ell+1} \qquad\hbox{with}\quad \bu=+,0,-.$$
These strata are illustrated in Figure~\ref{DT_fig}.
For $\vr\!\in\!\cA_{\ell}^E\!\cup\!\cA_{\ell;1}^D$, let
$$\wt{D}_{\vr}^+=D_{\ell+1;\vr\cup\{(\ell+1)^+\}},\quad
\wt{D}_{\vr}^0=D_{\ell+1;\vr}, \quad
\wt{D}_{\vr}^-=D_{\ell+1;\vr\cup\{(\ell+1)^-\}}.$$
If $\vr\!\in\!\cA_{\ell}^H$ (resp.~$\vr\!\in\!\cA_{\ell;2}^D\!\cup\!\cA_{\ell;3}^D$), 
\hbox{$D_{\ell+1;\vr\cup\{(\ell+1)^-\}}\!=\!\eset$}
(resp.~$D_{\ell+1;\vr}\!=\!\eset$).
In these cases, we~set
$$\wt{D}_{\vr}^+=D_{\ell+1;\vr\cup\{(\ell+1)^+,(\ell+1)^-\}}\,,
\quad
\wt{D}_{\vr}^0,\wt{D}_{\vr}^-= \wt{D}_{\vr}''=
\begin{cases}
D_{\ell+1;\vr},&\hbox{if}~\vr\!\in\!\cA_{\ell}^H;\\
D_{\ell+1;\vr\cup\{(\ell+1)^-\}},&\hbox{if}~
\vr\!\in\!\cA_{\ell;2}^D\!\cup\!\cA_{\ell;3}^D.
\end{cases}$$
If $\vr\!\in\!\cA_{\ell;1}^D$, 
$\wt{D}_{\vr}^0\!=\!\wt{D}_{\ov\vr_{\ell^{\pm}}^c}^+$ and
$\wt{D}_{\vr}^-\!=\!\wt{D}_{\ov\vr_{\ell^{\pm}}^c}^0,\wt{D}_{\ov\vr_{\ell^{\pm}}^c}^-$.
If $\vr\!\in\!\cA_{\ell;3}^D$, then $\wt{D}_{\vr}^+\!=\!\wt{D}_{\ov\vr}^+$.
These are the only cases of pairs~$(\vr,\bu)$ and~$(\vr',\circ)$ 
with distinct $\vr,\vr'\!\in\!\cA_{\ell}^{\R}$ such that 
$\wt{D}_{\vr}^{\bu}\!=\!\wt{D}_{\vr'}^{\circ}$.
For $i\!\in\![\ell^{\pm}]$, let 
$$\wt{D}_{[\ell^{\pm}]-\{i\}}^+=\eset\subset\ov\cM_{0,\ell+1}\,,\quad
\wt{D}_{[\ell^{\pm}]-\{i\}}^0,\wt{D}_{[\ell^{\pm}]-\{i\}}^-
=D_{\ell+1;([\ell^{\pm}]-\{i\})\cup\{(\ell+1)^-\}}\,.$$

\vspace{.15in}

\begin{figure}
\begin{pspicture}(-1,-9.2)(10,5.6)
\psset{unit=.4cm}
\rput(0,7.5){\framebox{$\vr\!\in\!\cA_{\ell}^H$}}
\pscircle[linewidth=.05](2.5,11){1.5}\pscircle[linewidth=.05](5.5,11){1.5}
\pscircle*(1.44,12.06){.15}\pscircle*(1.44,9.94){.15}
\rput(1.1,12.7){\sm{$1^+$}}\rput(1.1,9.3){\sm{$1^-$}}
\rput(1.6,11.1){\sm{$\vr$}}\rput(6.3,11.1){\sm{$\vr_{\ell^{\pm}}^c$}}
\psline[linewidth=.025]{<->}(7.5,9.5)(7.5,12.5)\rput(8,12){$\si$}
\rput(5.5,7.5){\sm{$D_{\ell;\vr}$}}
\pscircle[linewidth=.05](12.5,11){1.5}\pscircle[linewidth=.05](15.5,11){1.5}
\pscircle*(11.44,12.06){.15}\pscircle*(11.44,9.94){.15}
\pscircle*(12.5,12.5){.15}\pscircle*(12.5,9.5){.15}
\rput(11.1,12.7){\sm{$1^+$}}\rput(11.1,9.3){\sm{$1^-$}}
\rput(11.6,11.1){\sm{$\vr$}}\rput(16.3,11.1){\sm{$\vr_{\ell^{\pm}}^c$}}
\psline[linewidth=.025]{<->}(17.5,9.5)(17.5,12.5)\rput(18,11){$\si$}
\rput(13.5,13.2){\sm{$(\ell\!+\!1)^+$}}\rput(13.5,8.8){\sm{$(\ell\!+\!1)^-$}}
\rput(14.5,7.5){\sm{$\wt{D}_\vr^+\!\subset\!\R\ov\cM_{0,\ell+1}$}}
\pscircle[linewidth=.05](22.5,11){1.5}\pscircle[linewidth=.05](25.5,11){1.5}
\pscircle*(21.44,12.06){.15}\pscircle*(21.44,9.94){.15}
\pscircle*(25.5,12.5){.15}\pscircle*(25.5,9.5){.15}
\rput(21.1,12.7){\sm{$1^+$}}\rput(21.1,9.3){\sm{$1^-$}}
\rput(21.6,11.1){\sm{$\vr$}}\rput(26.3,11.1){\sm{$\vr_{\ell^{\pm}}^c$}}
\psline[linewidth=.025]{<->}(27.5,9.5)(27.5,12.5)\rput(28,11){$\si$}
\rput(25.5,13.2){\sm{$(\ell\!+\!1)^+$}}\rput(25.5,8.8){\sm{$(\ell\!+\!1)^-$}}
\rput(24.5,7.5){\sm{$\wt{D}_\vr^0,\wt{D}_\vr^-\!\subset\!\R\ov\cM_{0,\ell+1}$}}
\rput(0,-2.5){\framebox{$\vr\!\in\!\cA_{\ell}^E$}}
\psline[linewidth=.05](2,1)(6,5)\psline[linewidth=.05](2,3)(6,-1)
\pscircle*(3.5,2.5){.15}\pscircle*(5,4){.15}
\pscircle*(3.5,1.5){.15}\pscircle*(5,0){.15}
\rput(3.5,3.1){\sm{$1^{\pm}$}}\rput(3.2,.9){\sm{$1^{\mp}$}}
\rput(4.9,4.6){\sm{$2^+$}}\rput(4.7,-.7){\sm{$2^-$}}
\rput(6.5,5.1){\sm{$\vr$}}\rput(6.7,-1){\sm{$\vr_{\ell^{\pm}}^c$}}
\psline[linewidth=.025]{<->}(5.5,.5)(5.5,3.5)\rput(6,2){$\si$}
\rput(5.5,-2.5){\sm{$D_{\ell;\vr}$}}
\psline[linewidth=.05](12,1)(16,5)\psline[linewidth=.05](12,3)(16,-1)
\pscircle*(13.5,2.5){.15}\pscircle*(15,4){.15}\pscircle*(14.25,3.25){.15}
\pscircle*(13.5,1.5){.15}\pscircle*(15,0){.15}\pscircle*(14.25,.75){.15}
\rput(13.3,3){\sm{$1^{\pm}$}}\rput(13.2,.9){\sm{$1^{\mp}$}}
\rput(15,3){\sm{$2^+$}}\rput(15,1.2){\sm{$2^-$}}
\rput(13.5,-.5){\sm{$(\ell\!+\!1)^-$}}
\rput(14,4.6){\sm{$(\ell\!+\!1)^+$}}
\psline[linewidth=.025]{<->}(16.5,.5)(16.5,3.5)\rput(17,2){$\si$}
\rput(16.5,5.1){\sm{$\vr$}}\rput(16.7,-1){\sm{$\vr_{\ell^{\pm}}^c$}}
\rput(14.5,-2.5){\sm{$\wt{D}_\vr^+\!\subset\!\R\ov\cM_{0,\ell+1}$}}
\psline[linewidth=.05](23,4.5)(23,-.5)
\psline[linewidth=.02](22.5,4)(26,4)\psline[linewidth=.02](22.5,0)(26,0)
\pscircle*(24,4){.15}\pscircle*(25.5,4){.15}\pscircle*(25.5,0){.15}
\pscircle*(23,3){.15}\pscircle*(23,1){.15}\pscircle*(24,0){.15}
\rput(24,4.7){\sm{$1^{\pm}$}}\rput(24,-.7){\sm{$1^{\mp}$}}
\rput(25.5,4.7){\sm{$2^+$}}\rput(25.5,-.7){\sm{$2^-$}}
\psline[linewidth=.025]{<->}(25.5,.6)(25.5,3.4)\rput(26,2){$\si$}
\rput(26.5,4.1){\sm{$\vr$}}\rput(26.7,0){\sm{$\vr_{\ell^{\pm}}^c$}}
\rput(21.5,3.2){\sm{$(\ell\!+\!1)^+$}}\rput(21.5,1.2){\sm{$(\ell\!+\!1)^-$}}
\rput(24.5,-2.5){\sm{$\wt{D}_\vr^0\!\subset\!\R\ov\cM_{0,\ell+1}$}}
\psline[linewidth=.05](32,1)(36,5)\psline[linewidth=.05](32,3)(36,-1)
\pscircle*(33.5,2.5){.15}\pscircle*(35,4){.15}\pscircle*(34.25,3.25){.15}
\pscircle*(33.5,1.5){.15}\pscircle*(35,0){.15}\pscircle*(34.25,.75){.15}
\rput(33.3,3){\sm{$1^{\pm}$}}\rput(33.2,.9){\sm{$1^{\mp}$}}
\rput(35,3){\sm{$2^+$}}\rput(35,1.2){\sm{$2^-$}}
\rput(33.5,-.5){\sm{$(\ell\!+\!1)^+$}}
\rput(34,4.6){\sm{$(\ell\!+\!1)^-$}}
\psline[linewidth=.025]{<->}(36.5,.5)(36.5,3.5)\rput(37,2){$\si$}
\rput(36.5,5.1){\sm{$\vr$}}\rput(36.7,-1){\sm{$\vr_{\ell^{\pm}}^c$}}
\rput(34.5,-2.5){\sm{$\wt{D}_\vr^-\!\subset\!\R\ov\cM_{0,\ell+1}$}}
\rput(0,-12.5){\framebox{$\vr\!\in\!\cA_{\ell;1}^D$}}
\psline[linewidth=.05](3,-5.5)(3,-10.5)
\psline[linewidth=.02](2.5,-6)(6,-6)\psline[linewidth=.02](2.5,-10)(6,-10)
\pscircle*(3.5,-6){.15}\pscircle*(4.5,-6){.15}
\pscircle*(3.5,-10){.15}\pscircle*(4.5,-10){.15}
\rput(3.7,-5.3){\sm{$1^{\pm}$}}\rput(3.7,-10.7){\sm{$1^{\mp}$}}
\rput(4.7,-6.6){\sm{$2^+$}}\rput(4.7,-9.3){\sm{$2^-$}}
\psline[linewidth=.025]{<->}(5.8,-9.4)(5.8,-6.6)\rput(6.3,-8){$\si$}
\rput(6.5,-5.9){\sm{$\vr$}}\rput(6.7,-10){\sm{$\ov\vr$}}
\rput(1.6,-8){\sm{$\vr_{\ell^{\pm}}^c\!-\!\ov\vr$}}
\rput(5.5,-12.5){\sm{$D_{\ell;\vr}$}}
\psline[linewidth=.05](13,-5.5)(13,-10.5)
\psline[linewidth=.02](12.5,-6)(16,-6)\psline[linewidth=.02](12.5,-10)(16,-10)
\pscircle*(13.5,-6){.15}\pscircle*(14.5,-6){.15}\pscircle*(15.5,-6){.15}
\pscircle*(13.5,-10){.15}\pscircle*(14.5,-10){.15}\pscircle*(15.5,-10){.15}
\rput(13.7,-5.3){\sm{$1^{\pm}$}}\rput(13.7,-10.7){\sm{$1^{\mp}$}}
\rput(14.7,-6.6){\sm{$2^+$}}\rput(14.7,-9.3){\sm{$2^-$}}
\rput(16.5,-5.3){\sm{$(\ell\!+\!1)^+$}}\rput(16.5,-10.7){\sm{$(\ell\!+\!1)^-$}}
\psline[linewidth=.025]{<->}(15.8,-9.4)(15.8,-6.6)\rput(16.3,-8){$\si$}
\rput(16.5,-6.2){\sm{$\vr$}}\rput(16.7,-9.7){\sm{$\ov{\vr}$}}
\rput(11.6,-8){\sm{$\vr_{\ell^{\pm}}^c\!-\!\ov{\vr}$}}
\rput(14.5,-12.5){\sm{$\wt{D}_\vr^+\!\subset\!\R\ov\cM_{0,\ell+1}$}}
\psline[linewidth=.05](23,-5.5)(23,-10.5)
\psline[linewidth=.02](22.5,-6)(26,-6)\psline[linewidth=.02](22.5,-10)(26,-10)
\pscircle*(23.5,-6){.15}\pscircle*(24.5,-6){.15}
\pscircle*(23.5,-10){.15}\pscircle*(24.5,-10){.15}
\pscircle*(23,-6.5){.15}\pscircle*(23,-9.5){.15}
\rput(23.7,-5.3){\sm{$1^{\pm}$}}\rput(23.7,-10.7){\sm{$1^{\mp}$}}
\rput(24.7,-6.6){\sm{$2^+$}}\rput(24.7,-9.3){\sm{$2^-$}}
\psline[linewidth=.025]{<->}(25.8,-9.4)(25.8,-6.6)\rput(26.3,-8){$\si$}
\rput(26.5,-5.9){\sm{$\vr$}}\rput(26.7,-10){\sm{$\ov{\vr}$}}
\rput(21.6,-8){\sm{$\vr_{\ell^{\pm}}^c\!-\!\ov{\vr}$}}
\rput(21.3,-6.8){\sm{$(\ell\!+\!1)^+$}}\rput(21.4,-9.2){\sm{$(\ell\!+\!1)^-$}}
\rput(24.5,-12.5){\sm{$\wt{D}_\vr^0\!\subset\!\R\ov\cM_{0,\ell+1}$}}
\psline[linewidth=.05](33,-5.5)(33,-10.5)
\psline[linewidth=.02](32.5,-6)(36,-6)\psline[linewidth=.02](32.5,-10)(36,-10)
\pscircle*(33.5,-6){.15}\pscircle*(34.5,-6){.15}\pscircle*(35.5,-6){.15}
\pscircle*(33.5,-10){.15}\pscircle*(34.5,-10){.15}\pscircle*(35.5,-10){.15}
\rput(33.7,-5.3){\sm{$1^{\pm}$}}\rput(33.7,-10.7){\sm{$1^{\mp}$}}
\rput(34.7,-6.6){\sm{$2^+$}}\rput(34.7,-9.3){\sm{$2^-$}}
\rput(36.5,-5.3){\sm{$(\ell\!+\!1)^-$}}\rput(36.5,-10.7){\sm{$(\ell\!+\!1)^+$}}
\psline[linewidth=.025]{<->}(35.8,-9.4)(35.8,-6.6)\rput(36.3,-8){$\si$}
\rput(36.5,-6.2){\sm{$\vr$}}\rput(36.7,-9.7){\sm{$\ov{\vr}$}}
\rput(31.6,-8){\sm{$\vr_{\ell^{\pm}}^c\!-\!\ov{\vr}$}}
\rput(34.5,-12.5){\sm{$\wt{D}_\vr^-\!\subset\!\R\ov\cM_{0,\ell+1}$}}
\rput(1,-22.5){\framebox{$\vr\!\in\!\cA_{\ell;2}^D\!\cup\!\cA_{\ell;3}^D$}}
\psline[linewidth=.05](3,-15.5)(3,-20.5)
\psline[linewidth=.02](2.5,-16)(6,-16)\psline[linewidth=.02](2.5,-20)(6,-20)
\psline[linewidth=.025]{<->}(5.8,-19.4)(5.8,-16.6)\rput(6.3,-18){$\si$}
\rput(6.7,-15.9){\sm{$\ov{\vr}_{\ell^{\pm}}^c$}}\rput(6.7,-20){\sm{$\vr_{\ell^{\pm}}^c$}}
\rput(1.6,-18){\sm{$\ov{\vr}\!-\!\vr_{\ell^{\pm}}^c$}}
\rput(5.8,-22.5){\sm{$D_{\ell;\vr}$}}
\psline[linewidth=.05](13,-15.5)(13,-20.5)
\psline[linewidth=.02](12.5,-16)(16,-16)\psline[linewidth=.02](12.5,-20)(16,-20)
\pscircle*(13,-16.5){.15}\pscircle*(13,-19.5){.15}
\psline[linewidth=.025]{<->}(15.8,-19.4)(15.8,-16.6)\rput(16.3,-18){$\si$}
\rput(16.7,-16.5){\sm{$\ov{\vr}_{\ell^{\pm}}^c$}}\rput(16.7,-19.5){\sm{$\vr_{\ell^{\pm}}^c$}}
\rput(11.6,-18){\sm{$\ov{\vr}\!-\!\vr_{\ell^{\pm}}^c$}}
\rput(11.3,-16.8){\sm{$(\ell\!+\!1)^+$}}\rput(11.4,-19.2){\sm{$(\ell\!+\!1)^-$}}
\rput(14.5,-22.5){\sm{$\wt{D}_\vr^+\!\subset\!\R\ov\cM_{0,\ell+1}$}}
\psline[linewidth=.05](33,-15.5)(33,-20.5)
\psline[linewidth=.02](32.5,-16)(36,-16)\psline[linewidth=.02](32.5,-20)(36,-20)
\pscircle*(34.5,-16){.15}\pscircle*(34.5,-20){.15}
\rput(35.5,-15.3){\sm{$(\ell\!+\!1)^-$}}\rput(35.5,-20.7){\sm{$(\ell\!+\!1)^+$}}
\psline[linewidth=.025]{<->}(35.8,-19.4)(35.8,-16.6)\rput(36.3,-18){$\si$}
\rput(36.7,-16.5){\sm{$\ov\vr_{\ell^{\pm}}^c$}}\rput(36.7,-19.5){\sm{$\vr_{\ell^{\pm}}^c$}}
\rput(31.6,-18){\sm{$\ov\vr\!-\!\vr_{\ell^{\pm}}^c$}}
\rput(34.5,-22.5){\sm{$\wt{D}_\vr^0,\wt{D}_\vr^-\!\subset\!\R\ov\cM_{0,\ell+1}$}}
\end{pspicture}
\caption{Generic representatives of the top boundary strata $D_{\ell;\vr}$ in $\R\ov\cM_{0,\ell}$ 
and $\wt{D}_{\vr}^{\bu}$ in $\R\ov\cM_{0,\ell+1}$.
Each line or circle represents~$\C\P^1$.
Each double-headed arrow labeled~$\si$ indicates
the involution on the corresponding real curve.}
\label{DT_fig}
\end{figure}

\noindent
For $\vr^*\!\in\!\{0\}\!\sqcup\!\cA_{\ell}^{\R}$, $\vr\!\in\!\wt\cA_{\ell}^{\R}$, 
and $\bu\!=\!+,-,0$, let
$$Y_{\vr^*;\vr}^{\bu}=q_{\vr^*}(\wt{D}_{\vr}^{\bu})\subset X_{\vr^*}.$$
If $\vr\!\in\!\cA_{\ell}(\vr^*)$, then $Y_{\vr^*;\vr}^0\!=\!Y_{\vr^*;\vr}^-$.
Furthermore, 
\begin{gather}\label{RYvrmp_e}
Y_{\vr^*;\vr}^0
=q_{\vr^*}\big(\wt{D}_{\vr}^+\!\cap\!\wt{D}_{\vr}^0\big)
=  Y_{\vr^*;\vr}^+\!\cap\!Y_{\vr^*;[\ell^{\pm}]-\{i\}}^0
 ~~\forall~\vr\!\in\!\cA_{\ell}^{\R}(\vr^*),~i\!\in\!\vr_{\ell^{\pm}}^c,\\
\notag
Y_{\vr^*;\vr}^0=Y_{\vr^*;\ov\vr_{\ell}^c}^+~~\forall~\vr\!\in\!\cA_{\ell;1}^D,
\quad
Y_{\vr^*;\vr}^-=q_{\vr^*}\big(\wt{D}_{\vr}^0\!\cap\!\wt{D}_{\vr}^-\big)
\subset Y_{\vr^*;\vr}^0~~\forall~\vr\!\in\!\cA_{\ell;1}^D,~
\ov\vr_{\ell^{\pm}}^c\!\in\!\cA_{\ell}(\vr^*),\\
\notag
Y_{\vr^*;\vr}^0\!\cap\!Y_{\vr^*;\ov\vr_{\ell}^c}^+=\eset
~~\forall~\vr\!\in\!\cA_{\ell;2}^D,~\ov\vr_{\ell^{\pm}}^c\!\le\!\vr^*.
\end{gather}
The next theorem is proved in Section~\ref{Rblowup_sec},
with Proposition~\ref{Rblowup_prp} playing a central role.

\begin{thm}\label{Rblowup_thm}
Suppose $\ell\!\in\!\Z^+$ with $\ell\!\ge\!2$.
\begin{enumerate}[ref=$\R\arabic*$,label=($\R\arabic*$),leftmargin=*]

\item\label{Rspaces_it} For every $\vr^*\!\in\!\{0\}\!\sqcup\!\cA_{\ell}^{\R}$,
the space $X_{\vr^*}$ is a smooth manifold so that the quotient projection~$q_{\vr^*}$
in~\eref{qvrdfnR_e} is smooth
and the compact subspaces~$Y_{\vr^*;\vr}^{\bu}\!\subset\!X_{\vr^*}$ 
with $\vr\!\in\!\wt\cA_{\ell}^{\R}$ and \hbox{$\bu\!=\!+,0,-$} are submanifolds.

\item\label{Rblowup_it} If $\vr^*\!\in\!\cA_{\ell}^H$
(resp.~$\vr^*\!\in\!\cA_{\ell}^D$, $\vr^*\!\in\!\cA_{\ell}^E$), 
the map~$\pi_{\vr^*}$ in~\eref{Rbldowndfn_e}  
is the real blowup (resp.~a complex blowup, a 1-augmented blowup)
of~$X_{\vr^*-1}$ along~$Y_{\vr^*-1;\vr^*}^0$
with the exceptional locus \hbox{$Y_{\vr^*;\vr^*}^0\!\cup\!Y_{\vr^*;\vr^*}^-$}.
If $\vr^*\!\in\!\cA_{\ell;1}^D$, 
$Y_{\vr^*;\ov{\vr^*}^c_{\ell^{\pm}}}^0$ is a section of 
$\pi_{\vr^*}\!:Y_{\vr^*;\vr^*}^0\!\lra\!Y_{\vr^*-1;\vr^*}^0$.

\item\label{Rsubspaces_it} 
If $\vr^*\!\in\!\cA_{\ell}^{\R}$, 
$\vr\!\in\!\wt\cA_{\ell}^{\R}\!-\!\cA_{\ell}^H$ (resp.~$\vr\!\in\!\cA_{\ell}^H$),
and $\bu\!=\!+,0,-$ with 
$$\!\!
(\vr,\bu)\neq(\vr^*,0),(\vr^*,-), ~~
(\vr,\bu)\neq(\ov{\vr^*}_{\ell^{\pm}}^c,+)~\hbox{if}~\vr^*\!\in\!\cA_{\ell;1}^D, 
~~\hbox{and}~~
(\vr,\bu)\neq(\ov{\vr^*}_{\ell^{\pm}}^c,-)~\hbox{if}~\vr^*\!\in\!\cA_{\ell;2}^D,$$
then the submanifold $Y_{\vr^*;\vr}^{\bu}\!\subset\!X_{\vr^*}$ is
$(\pi_{\vr^*},\Z)$-related to~$Y_{\vr^*-1;\vr}^{\bu}$
(resp.~$(\pi_{\vr^*},R)$-related to~$Y_{\vr^*-1;\vr}^{\bu}$ for any admissible ring~$R$
with $\cha(R)\!=\!2$).
The submanifolds \hbox{$Y_{\vr^*;\vr}^0\!\subset\!X_{\vr^*}$} 
with \hbox{$\vr\!\in\!\cA_{\ell}^{\R}(\vr^*)$}
are $\pi_{\vr^*}$-equivalent to the corresponding submanifolds~$Y_{\vr^*-1;\vr}^0$.

\end{enumerate}
\end{thm}

\vspace{.15in}

\noindent
In the case of the last statement of Theorem~\ref{Rblowup_thm}\eref{Rsubspaces_it}
with $\vr^*\!\in\!\cA_{\ell;1}^D$ and $\vr\!=\!\ov{\vr^*}_{\ell^{\pm}}^c$,
$Y_{\vr^*-1;\vr}^0\!=\!Y_{\vr^*-1;\vr^*}^0$
and $Y_{\vr^*;\vr}^0\!\subset\!Y_{\vr^*;\vr^*}^0$ is a section of
$$\pi_{\vr^*}\!:Y_{\vr^*;\vr^*}^0\lra Y_{\vr^*-1;\vr^*}^0.$$
In the remaining cases of the last statement of Theorem~\ref{Rblowup_thm}\eref{Rsubspaces_it}
and in all cases of the last statement of Theorem~\ref{Cblowup_thm}\ref{Csubspaces_it},
$\cN_{Y_{\vr^*-1;\vr}^0\cap Y_{\vr^*-1;\vr^*}^0}\!Y_{\vr^*-1;\vr}^0$
is of rank~1 (real or complex, depending on~$\vr^*$) 
and $Y_{\vr^*;\vr}^0$ is the proper transform of~$Y_{\vr^*-1;\vr}^0$
under the blowup~$\pi_{\vr^*}$.\\

\noindent
The approach of this paper readily applies to decompose the analogues
$$\big(\ff_{(\ell+1)^{\pm}}^{\R},\ff_{1^{\pm}2^+(\ell+1)^+}^{\R}\big)\!: 
\R\ov\cM_{k,\ell+1}\lra\R\ov\cM_{k,\ell}\!\times\!\ov\cM_4$$
of the morphism~\eref{RMlind_e} for $k\!>\!0$, provided $\ell\!\ge\!2$.
If $\ell\!\ge\!1$ and $k\!\ge\!1$, 
we can also replace $\ff_{1^{\pm}2^+(\ell+1)^+}^{\R}$ above
with the forgetful morphism
$$\ff_{1^{\pm}1(\ell+1)^+}^{\R}\!:\R\ov\cM_{k,\ell+1}\lra\ov\cM_4$$
dropping all marked points other than $z_1^+,z_1^-$, the first real marked point~$z_1$,
and~$z_{\ell+1}^+$.
In either case, only the real and complex blowups appear in
the resulting decomposition. 
If $k\!\ge\!3$, we can similarly decompose the smooth~map
$$\big(\ff_{k+1}^{\R},\ff_{123,k+1}^{\R}\big)\!: 
\R\ov\cM_{k+1,\ell}\lra X_0\!\equiv\!\R\ov\cM_{k,\ell}\!\times\!\R\ov\cM_{4,0},$$
where
$$\ff_{k+1}^{\R}\!:\R\ov\cM_{k+1,\ell}\lra\R\ov\cM_{k,\ell}  \qquad\hbox{and}\qquad
\ff_{123,k+1}^{\R}\!:\R\ov\cM_{k+1,\ell}\lra\R\ov\cM_{4,0}$$
are the morphisms dropping the last real marked point and 
all marked points other than the real marked points $z_1,z_2,z_3,z_{k+1}$, respectively.
In this case, only real blowups occur.

\section{Blowups: a local perspective}
\label{blowup_sec}

\noindent
We describe the real, complex, and augmented blowups
of a smooth manifold~$X$ along a closed codimension~$\fc$ submanifold~$Y$ 
appearing in Theorems~\ref{Cblowup_thm} and~\ref{Rblowup_thm},
i.e.~smooth proper surjections
\BE{blowdownmaps_e}\pi\!:\BLR_YX\lra X, \qquad \pi\!:\BLC_YX\lra X, 
\qquad\hbox{and}\qquad \pi\!:\BLau_YX\lra X\EE
with certain properties,
via charts for~$X$ that cover~$Y$ in Sections~\ref{StanBl_subs} and~\ref{AugblowupDfn_subs};
the augmented blowup involves a choice of $\fc_1\!\in\![\fc\!-\!1]$.
This adapts the local perspective of~\cite[p603]{GH} from the complex blowup 
in the holomorphic category to the three blowups in the smooth category;
a~global perspective on these blowups, which makes the changes in the topology
more evident, appears in~\cite{blowups}.
In the latter perspective, an augmented blowup is actually two real blowups
followed by a real blowdown in a different direction.
The complex and augmented blowups depend on the choice of a compatible collection
of local charts for~$X$ that cover~$Y$.
Such a collection determines a complex structure on 
the normal bundle~$\cN_XY$ of~$Y$ in~$X$ in the complex case 
and a distinguished subbundle $\cN_X^{\fc_1}\!Y\!\subset\!\cN_XY$ of corank~$\fc_1$
in the $\fc_1$-augmented case.
In the settings of Theorems~\ref{Cblowup_thm} and~\ref{Rblowup_thm}, 
the local charts are provided by the cross ratios of quadruples of
marked points.
For $k\!=\!0,1$, let \hbox{$\tau_Y^k\!\lra\!Y$} 
denote the trivial rank~$k$ real line bundle.\\

\noindent
In all three cases of~\eref{blowdownmaps_e}, the restrictions
\begin{gather}
\label{standdiff_e}
\pi\!:\BLR_YX\!-\!\pi^{-1}(Y)\lra X\!-\!Y, \qquad 
\pi\!:\BLC_YX\!-\!\pi^{-1}(Y)\lra X\!-\!Y,\\
\label{augdiff_e}
\hbox{and}\qquad 
\pi\!:\BLau_YX\!-\!\pi^{-1}(Y)\lra X\!-\!Y
\end{gather}
are diffeomorphisms.
The restrictions of~$\pi$ to the \sf{exceptional loci}
$$\pi\!:\bE_Y^{\R}X\!\equiv\!\pi^{-1}(Y)\lra Y \qquad\hbox{and}\qquad
\pi\!:\bE_Y^{\C}X\!\equiv\!\pi^{-1}(Y)\lra Y$$
in the first two cases
are isomorphic to the real and complex projectivizations, 
$\R\P(\cN_XY)$ and $\C\P(\cN_XY)$, respectively, of~$\cN_XY$.
The normal bundles to the exceptional loci in these two cases are
isomorphic to the real and complex tautological line bundles:
\BE{cNEisom_e}\cN_{\BLR_YX}\bE_Y^{\R}X\approx\ga_{\cN_XY}^{\R}\lra\R\P(\cN_XY)
\quad\hbox{and}\quad
\cN_{\BLC_YX}\bE_Y^{\C}X\approx\ga_{\cN_XY}^{\C}\lra\C\P(\cN_XY)\,.\EE

\vspace{.15in}

\noindent
In the last case in~\eref{blowdownmaps_e}, the exceptional locus
$$\bE_Y^{\fc_1}X\equiv\pi^{-1}(Y)\subset\wt{X}$$ 
is the union of two closed submanifolds, $\bE^0_YX$ and~$\bE^-_YX$, so that 
$$\big(\bE^0_YX,\bE^0_YX\!\cap\!\bE^-_YX\big)\approx
\big(\R\P\big(\cN_YX/\cN_Y^{\fc_1}X\!\oplus\!\tau_Y^1\big),
\R\P\big(\cN_YX/\cN_Y^{\fc_1}X\!\oplus\!\tau_Y^0\big)\big)$$
as fiber bundle pairs over~$Y$, while 
$$\pi\!:\bE^-_YX\lra Y$$
is a fiber bundle with fibers that themselves form a fiber bundle over~$\R\P^{\fc_1-1}$
with fibers $S^{\fc-\fc_1}$ and two sections, one of which is~$\bE^0_YX\!\cap\!\bE^-_YX$.\\

\noindent
The cases relevant for Theorems~\ref{Cblowup_thm} and~\ref{Rblowup_thm}
are the real blowup with $\fc\!=\!3$, complex blowup with $\fc/2\!=\!2$,
and augmented blowup with $(\fc,\fc_1)\!=\!(3,1)$.
In the last case, 
$\bE^-_YX$ is the $S^2$-sphere bundle over~$Y$ obtained by collapsing 
each fiber of the $\R\P^1$-fiber subbundle 
\hbox{$\R\P(\cN_Y^{\fc_1}X)\!\subset\!\R\P(\cN_YX)$} to a point.
The section of \hbox{$\bE_Y^{\fc_1}X\!\lra\!Y$}
determined by this collapse is~$\bE^0_YX\!\cap\!\bE^-_YX$.\\

\noindent
The constructions of the blowups out of local charts are somewhat involved
in the real and complex cases and quite technical in the augmented case.
However, verifying that a given smooth map 
\BE{pipr_e}\pi'\!:\wt{X}\lra X\EE
is a blowup of one of the three kinds involves only a comparison of 
local charts on~$X$ and~$\wt{X}$ as in Lemmas~\ref{StanBl_lmm} and~\ref{AugBl_lmm}. 
On the other hand, the cross ratios do not readily restrict to coordinate charts 
on the real moduli space
\hbox{$\R\ov\cM_{0,\ell}\!\subset\!\ov\cM_{2\ell}$}.
For this reason, we give an alternative characterization of the real blowup
in the $\fc\!=\!3$ case in Lemma~\ref{Rblowup_lmm}.\\

\noindent
An \sf{involution} on a set $W$ is a map $F\!:W\!\lra\!W$ such that $F^2\!=\!\id_W$.
For such a map, let
\BE{WFdfn_e}W_F=\big\{w\!\in\!W\!:F(w)\!=\!w\big\}.\EE
If $W$ is a smooth manifold and $F$ is smooth,
$W_F\!\subset\!W$ is a smooth submanifold and 
$$T_wW_F=\big\{\dot{w}\!\in\!T_wW\!:\nd_wF(\dot{w})\!=\!\dot{w}\big\}
\qquad\forall~w\!\in\!W_F;$$
see \cite[Lemma~3]{Meyer81}.

\subsection{Standard blowups}
\label{StanBl_subs}

\noindent
Let $\fc\!\in\!\Z^+$ and $\bF\!=\!\R,\C$.
We denote~by
$$\ga^{\bF}_{\fc}\equiv\big\{\!(L,v)\!\in\!\bF\P^{\fc-1}\!\times\!\bF^{\fc}\!:
v\!\in\!L\!\subset\!\bF^{\fc}\big\}$$
the \sf{tautological line bundle} over the $\bF$-projective space $\bF\P^{\fc-1}$.
For each $i\!\in\![\fc]$, let
\begin{gather}\label{vphfcidfn_e}
\wt\vph_{\fc;i}\!\equiv\!
\big(\wt\vph_{\fc;i;1},\ldots,\wt\vph_{\fc;i;\fc}\big)\!:
\wt{U}_{\fc;i}^{\bF}\!\equiv\!\big\{\!
\big([r_1,\ldots,r_{\fc}],v\big)\!\in\!\ga^{\bF}_{\fc}\!:
r_i\!\neq\!0\big\}\lra\bF^{\fc},\\
\notag
\wt\vph_{\fc;i;j}\big([r_1,\ldots,r_{\fc}],(v_1,\ldots,v_{\fc})\!\big)
=\begin{cases}r_j/r_i,&\hbox{if}~j\!\in\![\fc]\!-\!\{i\};\\
v_i,&\hbox{if}~j\!=\!i;
\end{cases}
\end{gather}
be the $i$-th standard coordinate chart on~$\ga^{\bF}_{\fc}$.\\

\noindent
Let $X$ be a smooth manifold and $Y\!\subset\!X$
be a closed submanifold of $\bF$-codimension~$\fc$ and real dimension~$m$.
We call a coordinate chart 
$$\vph\!\equiv\!(\vph_1,\ldots,\vph_{\fc+m})\!:U\!\lra\!\bF^{\fc}\!\times\!\R^m$$
on~$X$ a \sf{chart for~$Y$ in~$X$} if 
\BE{vphYcond_e}U\!\cap\!Y=\big\{x\!\in\!U\!:
\vph_1(x),\ldots,\vph_{\fc}(x)\!=\!0\big\}.\EE
For such a chart, the subspace
$$\BLF_Y\vph\equiv \big\{\!(L,x)\!\in\!\bF\P^{\fc-1}\!\times\!U\!:
\big(\vph_1(x),\ldots,\vph_{\fc}(x)\!\big)\!\in\!L\!\subset\!\bF^{\fc}\big\}$$
of $\bF\P^{\fc-1}\!\times\!U$ is a closed submanifold.
We denote by 
$\pi_{\bF\P^{\fc-1}},\pi_U\!:\BLF_Y\vph\!\lra\!\bF\P^{\fc-1},U$
the two projections.
The latter restricts to a diffeomorphism
$$\pi_U\!:\big\{\!(L,x)\!\in\!\bF\P^{\fc-1}\!\times\!U\!:x\!\not\in\!Y\big\}\lra U\!-\!Y.$$
For each $i\!\in\![\fc]$, the smooth map
\begin{gather}\label{BLFchart_e}
\wt\vph_i\!\equiv\!(\wt\vph_{i;1},\ldots,\wt\vph_{i;\fc+m})\!:
\BLF_{Y;i}\vph\!\equiv\!\big\{\!([r_1,\ldots,r_{\fc}],x)\!\in\!\BLF_Y\vph\!:
r_i\!\neq\!0\big\}
\lra\bF^{\fc}\!\times\!\R^m,\\
\notag
\wt\vph_{i;j}=\begin{cases}
\wt\vph_{\fc;i;j}\!\circ\!\pi_{\bF\P^{\fc-1}},
&\hbox{if}~j\!\in\![\fc]\!-\!\{i\};\\
\vph_j\!\circ\!\pi_U,
&\hbox{if}~j\!\in\!\{i\}\!\cup\!\big([\fc\!+\!m]\!-\![\fc]\big);
\end{cases}
\end{gather}
is a coordinate chart on~$\BLF_Y\vph$.\\

\noindent
Let $\{\vph_{\al}\!:U_{\al}\!\lra\!\bF^{\fc}\!\times\!\R^m\}_{\al\in\cI}$ be a collection
of charts for~$Y$ in~$X$.
For $\al,\al'\!\in\!\cI$, let
$$\vph_{\al\al'}\!\equiv\!\big(\vph_{\al\al';1},\ldots,\vph_{\al\al';\fc+m}\big)
\!\equiv\!\vph_{\al}\!\circ\!\vph_{\al'}^{-1}\!:
\vph_{\al'}\big(U_{\al}\!\cap\!U_{\al'}\big)\lra \vph_{\al}\big(U_{\al}\!\cap\!U_{\al'}\big)$$
be the overlap map between the charts $\vph_{\al}$ and $\vph_{\al'}$.
By~\eref{vphYcond_e},
\BE{Rslice_e}
\vph_{\al\al'}^{~-1}\big(0^{\fc}\!\times\R^m\big)
=\vph_{\al'}\big(U_{\al}\!\cap\!U_{\al'}\big)\!\cap\!\big(0^{\fc}\!\times\!\R^m\big).\EE
We call $\{\vph_{\al}\}_{\al\in\cI}$ an \sf{$\bF$-atlas for~$Y$ in~$X$}
if the domains~$U_{\al}$ of~$\vph_{\al}$ cover~$Y$ and
 for all $\al,\al'\!\in\cI$ there exists a smooth map 
\begin{gather}
\label{halalprdfn_e}
h_{\al\al'}\!\equiv\!\big(h_{\al\al';jj'}\big)_{j,j'\in[\fc]}\!:
\vph_{\al'}\big(U_{\al}\!\cap\!U_{\al'}\big)\lra\End_{\bF}(\bF^{\fc})
\qquad\hbox{s.t.}\\
\notag
\vph_{\al\al';j}(r,s)=\sum_{j'=1}^{\fc}h_{\al\al';jj'}(r,s)r_{j'}
~~\forall~j\!\in\![\fc],~
(r,s)\!\in\!\vph_{\al'}(U_{\al}\!\cap\!U_{\al'}),~
r\!\equiv\!(r_{j'})_{j'\in[\fc]}\!\in\!\bF^{\fc},~s\!\in\!\R^m.
\end{gather}
By~\eref{Rslice_e}, 
$$h_{\al\al'}(r,s)r\neq0 \quad\forall~
(r,s)\!\in\!\vph_{\al'}\big(U_{\al}\!\cap\!U_{\al'}\big),~
r\!\in\!\bF^{\fc}\!-\!\{0\},~s\!\in\!\R^m.$$
Since $h_{\al\al';jj'}$ equals $\prt\vph_{\al\al';j}/\prt r_{j'}$
along the right-hand side in~\eref{Rslice_e}, 
$h_{\al\al'}$ is an isomorphism along this subspace.
Thus, the~map
\BE{wtvphalalprdfn_e}\wt\vph_{\al\al'}\!:
\BLF_Y\big(\vph_{\al'}|_{U_{\al}\cap U_{\al'}}\big)
\lra \BLF_Y\big(\vph_{\al}|_{U_{\al}\cap U_{\al'}}\big), \quad
\wt\vph_{\al\al'}(L,x)=\big(h_{\al\al'}(\vph_{\al'}(x)\!)L,x\big),\EE
is well-defined and smooth.
It satisfies
$$\pi_{U_{\al}}\!\circ\!\wt\vph_{\al\al'}=
\pi_{U_{\al'}}\big|_{U_{\al}\cap U_{\al'}}\,.$$
By the uniqueness of continuous extensions and the cocycle condition
for the overlap maps~$\vph_{\al\al'}$,
\BE{cocyclecond_e}
\wt\vph_{\al\al''}\!=\!\wt\vph_{\al\al'}\!\circ\!\wt\vph_{\al'\al''}\!:
\BLF_Y\big(\vph_{\al''}|_{U_{\al}\cap U_{\al'}\cap U_{\al''}}\big)
\lra \BLF_Y\big(\vph_{\al}|_{U_{\al}\cap U_{\al'}\cap U_{\al''}}\big)
~~\forall\,\al,\al',\al''\!\in\!\cI.\EE

\vspace{.15in}

\noindent
Suppose $\{\vph_{\al}\}_{\al\in\cI}$ is an $\bF$-atlas for~$Y$ in~$X$.
An  \sf{$\bF$-blowup of~$X$ along~$Y$} is the smooth manifold
\begin{gather*}
\BLF_YX\equiv\Big(\!\!(X\!-\!Y)\!\sqcup\!\bigsqcup_{\al\in\cI}\!\!\BLF_Y\vph_{\al}\!
\Big)\!\!\Big/\!\!\!\sim,\quad
\BLF_Y\vph_{\al}\ni(L,x)\sim x\in X\!-\!Y~\forall\,x\!\in\!U_{\al}\!-\!Y,\,\al\!\in\!\cI,\\
\BLF_Y\big(\vph_{\al'}|_{U_{\al}\cap U_{\al'}}\big)\ni(L,x)\sim
\wt\vph_{\al\al'}(L,x)\in
\BLF_Y\big(\vph_{\al}|_{U_{\al}\cap U_{\al'}}\big)
~~\forall~\al,\al'\!\in\!\cI.
\end{gather*}
The \sf{blowdown map}
\BE{blowmapdfn_e}\pi\!:\BLF_YX\lra X, \qquad
\pi\big([\wt{x}]\big)=\begin{cases} 
\pi_{U_{\al}}(\wt{x}),&\hbox{if}~\wt{x}\!\in\!\BLF_Y\vph_{\al},\,\al\!\in\!\cI;\\
\wt{x},&\hbox{if}~\wt{x}\!\in\!X\!-\!Y;
\end{cases}\EE
is well-defined, surjective, proper, and smooth.\\

\noindent
The \sf{exceptional locus} for the blowdown map~$\pi$ in~\eref{blowmapdfn_e},
\begin{gather*}
\bE^{\bF}_YX\equiv\pi^{-1}(Y)=
\bigg(\bigsqcup_{\al\in\cI}\!\!
\big\{\!(L,x)\!\in\!\BLF_Y\vph_{\al}\!:x\!\in\!Y\big\}\!\!\bigg)\!\!\Big/\!\!\!\sim,\\
\BLF_Y\big(\vph_{\al'}|_{U_{\al}\cap U_{\al'}}\big)\ni(L,x)\sim
\wt\vph_{\al\al'}(L,x)\in
\BLF_Y\big(\vph_{\al}|_{U_{\al}\cap U_{\al'}}\big)
~~\forall~\al,\al'\!\in\!\cI,
\end{gather*}
is a smooth submanifold of~$\BLF_YX$ of $\bF$-codimension~1.
The blowdown map~\eref{blowmapdfn_e} restricts to a diffeomorphism~\eref{standdiff_e}.
The vector bundle isomorphisms
\BE{NXYisom_e}\bF^{\fc}\!\times\!(U_{\al}\!\cap\!Y)\lra \cN_YX|_{U_{\al}\cap Y}, \quad
\big(\!(r_1,\ldots,r_{\fc}),x\big)\lra\sum_{j=1}^{\fc}r_j
 \frac{\prt}{\prt\vph_{\al;j}}\Bigg|_x\!+\!T_xY,
\quad\al\!\in\!\cI,\EE
determine a complex structure on  the normal bundle~$\cN_XY$ of~$Y$ in~$X$
in the $\bF\!=\!\C$ case.
In both cases, they induce an identification of $\bE^{\bF}_YX$ with 
the $\bF$-projectivization~$\bF\P(\cN_XY)$ of~$\cN_XY$
as fiber bundles over~$Y$ and
an identification of~$\cN_{\BLF_YX}\bE^{\bF}_YX$ with 
the tautological line bundle~$\ga_{\cN_XY}^{\bF}$ over~$\bF\P(\cN_XY)$.\\

\noindent
Suppose $\bF\!=\!\R$.
By~\eref{Rslice_e}, smooth maps~$h_{\al\al'}$ as in~\eref{halalprdfn_e}
then exist for any collection $\{\vph_{\al}\}_{\al\in\cI}$ of charts for~$Y$ in~$X$
and for all $\al,\al'\!\in\!\cI$.
Since the above construction can be carried out with the maximal collection
$\{\vph_{\al}\}_{\al\in\cI}$ of charts for~$Y$ in~$X$ 
in this case,
the diffeomorphism class of the projection~$\pi$ in~\eref{blowmapdfn_e}
is independent of the choice of such a collection.
We call~$\pi$ \sf{the real blowup of~$X$ along~$Y$}.\\

\noindent
Suppose $\bF\!=\!\C$.
By~\eref{Rslice_e}, smooth maps~$h_{\al\al'}$ as in~\eref{halalprdfn_e}
then exist if the maps
$$\vph_{\al'}\big(U_{\al}\!\cap\!U_{\al'}\big)\!\cap\!
\big(\C^{\fc}\!\times\!\{s\}\!\big)\lra\C, 
~~r\lra \vph_{\al\al';j}(r,s),
\quad j\!\in\![\fc],\,s\!\in\!\R^m,$$
are holomorphic.
A \sf{co-holomorphic atlas for~$Y$ in~$X$} is
a collection $\{\vph_{\al}\}_{\al\in\cI}$ of charts for~$Y$ in~$X$
that covers~$Y$ and satisfies this condition for all $\al,\al'\!\in\!\cI$.
We define two such atlases to be \sf{equivalent} if their union is still  
a co-holomorphic atlas for~$Y$ in~$X$.
A \sf{co-holomorphic structure for~$Y$ in~$X$} is an equivalence class of 
co-holomorphic atlases for~$Y$ in~$X$.
Since the above construction can be carried out with the maximal collection of charts
in a co-holomorphic structure for~$Y$ in~$X$,
the diffeomorphism class of the projection~$\pi$ in~\eref{blowmapdfn_e}
is determined by such a structure, but may depend on~it.
We call~$\pi$ \sf{the complex blowup of~$X$ along~$Y$ determined} 
by the given co-holomorphic structure for~$Y$ in~$X$ or simply  
\sf{a complex blowup of~$X$ along~$Y$}.
If~$Y$ is a complex submanifold of a complex manifold~$X$,
an atlas of holomorphic charts on~$X$ satisfying~\eref{vphYcond_e} 
determines a co-holomorphic structure for~$Y$ in~$X$
and the complex blowup of~$X$ along~$Y$ in the sense of complex geometry,
which we call \sf{the holomorphic blowup of~$X$ along~$Y$}.\\

\noindent
The $\bF\!=\!\R$ case of the next lemma is equivalent to a special case 
of \cite[Lemma~2.1]{AkKi85} and \cite[Theorem~4.1]{ArKa}.
However, the direct analogues of \cite[Lemma~2.1]{AkKi85} and \cite[Theorem~4.1]{ArKa}
in the \hbox{$\bF\!=\!\C$} case do not hold, even if the differential between the normal
bundles of the submanifolds in their statements is assumed to be $\C$-linear,
because the complex blowup depends on the charts in an atlas
overlapping ``holomorphically in the normal direction" 
to the submanifold. 

\begin{lmm}\label{StanBl_lmm}
Suppose $X$ is a smooth manifold, 
$Y\!\subset\!X$ is a closed submanifold of $\bF$-codimension~$\fc$,
$\pi'$ is a smooth proper map as in~\eref{pipr_e} so~that 
\BE{StanBl_e}\pi'\!:\wt{X}\!-\!\pi'^{-1}(Y)\lra X\!-\!Y\EE
is a diffeomorphism, and
$\{\vph_{\al}\!\equiv\!(\vph_{\al;1},\ldots,\vph_{\al;\fc+m})\!:
U_{\al}\!\lra\!\bF^{\fc}\!\times\!\R^m\}_{\al\in\cI}$ is an $\bF$-atlas of charts for~$Y$ in~$X$.
Let \hbox{$\pi\!:\BLF_YX\!\lra\!X$} be the $\bF$-blowup of~$X$ along~$Y$
determined by this atlas.
If there exists a collection 
$$\big\{\wt\phi_{\al;i}\!\equiv\!
\big(\wt\phi_{\al;i;1},\ldots,\wt\phi_{\al;i;\fc+m}\big)
\!:\wt{U}_{\al;i}\!\lra\!\bF^{\fc}\!\times\!\R^m
\big\}_{\al\in\cI,i\in[\fc]}$$
of charts on~$\wt{X}$ so that 
$$\pi'^{-1}(U_{\al})=\bigcup_{i=1}^{\fc}\!\wt{U}_{\al;i}\,,
~~
\vph_{\al;j}\!\circ\!\pi'|_{\wt{U}_{\al;i}}=\begin{cases}
\wt\phi_{\al;i;j}\!\cdot\!\{\vph_{\al;i}\!\circ\!\pi'|_{\wt{U}_{\al;i}}\},
&\hbox{if}~i\!\in\![\fc],\,j\!\in\![\fc]\!-\!\{i\};\\
\wt\phi_{\al;i;j},&\hbox{if}~i\!\in\![\fc],\,
j\!\in\!\{i\}\!\cup\!\big([\fc\!+\!m]\!-\![\fc]);
\end{cases}$$
for every~$\al\!\in\!\cI$, then the maps $\pi$ and $\pi'$ are diffeomorphic.
\end{lmm}

\begin{proof}
By  \eref{vphYcond_e} and the second condition on the charts~$\wt\phi_{\al;i}$,
\BE{StanBl_e1}
\wt{U}_{\al;i}\!\cap\!\pi'^{-1}(Y)=\wt\phi_{\al;i;i}^{-1}(0).\EE
For $\al\!\in\!\cI$ and $i\!\in\![\fc]$, define
$$\wt\psi_{\al;i}\!\equiv\!\big(\wt\psi_{\al;i;1},\ldots,\wt\psi_{\al;i;\fc}\big)
\!:\wt{U}_{\al;i}\lra\bF^{\fc}\!-\!\{0\}, \quad
\wt\psi_{\al;i;j}(\wt{x})=\begin{cases}
\wt\phi_{\al;i;j}(\wt{x}),&\hbox{if}~j\!\in\![\fc]\!-\!\{i\};\\
1,&\hbox{if}~j\!=\!i.
\end{cases}$$
By the second condition on the charts~$\wt\phi_{\al;i}$, the map
$$\wt\Psi_{\al;i}\!:\wt{U}_{\al;i}\lra\BLF_{Y;i}\vph_{\al}\!\subset\!\BLF_Y\vph_{\al}, \quad
\wt\Psi_{\al;i}(\wt{x})=\big(\bF\wt\psi_{\al;i}(\wt{x}),\pi'(\wt{x})\!\big),$$
is well-defined and
$$\pi\big([\wt\Psi_{\al;i}(\wt{x})]\big)=\pi'(\wt{x})\in X
\quad\forall\,\wt{x}\!\in\!\wt{U}_{\al;i}.$$
Since
$$\wt\vph_{\al;i}\!\circ\!\wt\Psi_{\al;i}\!=\!\wt\phi_{\al;i}\!:
\wt{U}_{\al;i}\lra \bF^{\fc}\!\times\!\R^m,$$
where $\wt\vph_{\al;i}$ is a coordinate chart on~$\BLF_Y\vph_{\al}$
as in~\eref{BLFchart_e}, 
$\wt\Psi_{\al;i}$ is a diffeomorphism onto 
an open subset of~$\BLF_{Y;i}\vph_{\al}$ for every $i\!\in\![\fc]$.\\

\noindent
Since the restriction~\eref{standdiff_e} is injective,
\BE{StanBl_e6}\big[\wt\Psi_{\al;i}(\wt{x})\big]=\big[\wt\Psi_{\al';i'}(\wt{x})\big]\in\BLF_YX
\quad\forall~\wt{x}\!\in\!\wt{U}_{\al;i}\!\cap\!\wt{U}_{\al';i'}\!-\!\pi'^{-1}(Y),\,
\al,\al'\!\in\!\cI,\,i,i'\!\in\![\fc].\EE
Along with~\eref{StanBl_e1}, this implies that the equality in~\eref{StanBl_e6}
holds on all of~$\wt{U}_{\al;i}\!\cap\!\wt{U}_{\al';i'}$.
Thus, the~map
$$\wt\Psi\!:\wt{X}\lra\BLF_YX, \qquad
\wt\Psi(\wt{x})=\begin{cases}[\wt\Psi_{\al;i}(\wt{x})],
&\hbox{if}~\wt{x}\!\in\!\wt{U}_{\al;i},\,\al\!\in\!\cI,\,i\!\in\![\fc];\\
[\pi'(\wt{x})],&\hbox{if}~\wt{x}\!\in\!\wt{X}\!-\!\pi'^{-1}(Y);
\end{cases}$$
is a well-defined smooth submersion onto an open subset of~$\BLF_YX$.
It satisfies \hbox{$\pi'\!=\!\pi\!\circ\!\wt\Psi$}. 
Since the restriction~\eref{StanBl_e} is injective, 
\eref{StanBl_e1} implies that $\wt\Psi$ is injective as~well.
Since the restriction~\eref{StanBl_e} is surjective, 
the image of~$\wt\Psi$ contains $\BLF_YX\!-\!\bE^{\bF}_YX$. 
Since $\pi'$ is a proper map, 
$$\wt\Psi\big(\pi'^{-1}(x)\!\big)\subset \pi^{-1}(x)
\subset \bE^{\bF}_YX\subset \BLF_YX$$
is a nonempty, compact, open subspace of~$\pi^{-1}(x)$ for every $x\!\in\!Y$.
Thus, $\pi'$ is onto.
\end{proof}

\subsection{Augmented blowup}
\label{AugblowupDfn_subs}

\noindent
Let $\fc\!\in\!\Z^+$, $\fc_1\!\in\![\fc\!-\!1]$, 
$\fc_2\!=\!\fc\!-\!\fc_1$, and $\dbsqbr{\fc_1}\!=\!\{0\}\!\sqcup\![\fc_1]$.
For $r\!\equiv\!(r_j)_{j\in[\fc_2]},v\!\equiv\!(v_j)_{j\in[\fc_2]}\!\in\!\R^{\fc_2}$, 
denote~by
$$r\!\cdot\!v\equiv r_1v_1\!+\!\ldots\!+\!r_{\fc_2}v_{\fc_2}\in\R
\qquad\hbox{and}\qquad
|r|=(r\!\cdot\!r)^{1/2}$$
the usual inner-product of~$r$ and~$v$ and the corresponding norm of~$r$. 
For \hbox{$r\!\equiv\!(r_i)_{i\in[\fc_1]}\!\in\!\R^{\fc_1}$} and 
$\la\!\equiv\!(\la_j)_{j\in[\fc_2]}\!\in\!\R^{\fc_2}$, let
$$r\!\times\!\la=\big(r_i\la_j)_{i\in[\fc_1],j\in[\fc_2]}\in
\R^{[\fc_1]\times[\fc_2]}\,.$$

\vspace{.15in}

\noindent
We define
\begin{gather*}
\pi_{\ga_{\fc;\fc_1}^1}\!\!:
\ga_{\fc;\fc_1}^1\!\equiv\!\big\{\!
\big([r_1,\ldots,r_{\fc}],v\big)\!\in\!\ga^{\R}_{\fc}\!:
(r_1,\ldots,r_{\fc_1})\!\neq\!0\big\}\lra \R\P^{\fc-1}\\
\hbox{and}\qquad
\pi_{\ga_{\fc;\fc_1}^2}\!\!:
\ga_{\fc;\fc_1}^2\!\equiv\!\big(\ga_{\fc_1+1}^{\R}\big)^{\!\oplus\fc_2}
\lra\R\P^{\fc_1}
\end{gather*}
to be the obvious bundle projections.
For each $j\!\in\![\fc]$, let
\begin{alignat*}{2}
\pi_{\fc;\fc_1;j}^1\!:\ga_{\fc;\fc_1}^1&\lra\R, &\quad 
\pi_{\fc;\fc_1;j}^1\big(L,(v_{j'})_{j'\in[\fc]}\big)&=v_j,\\
\pi_{\fc;\fc_1;j}^2\!:\ga_{\fc;\fc_1}^2&\lra\R, &\quad 
\pi_{\fc;\fc_1;j}^2
\big(L,(v_{ij'})_{i\in\dbsqbr{\fc_1},j'\in[\fc_2]}\big)
&=\begin{cases}
(v_{0j'})_{j'\in[\fc_2]}\!\cdot\!(v_{jj'})_{j'\in[\fc_2]},
&\hbox{if}~j\!\in\![\fc_1];\\
v_{0(j-\fc_1)},&\hbox{if}~j\!\in\![\fc]\!-\![\fc_1].
\end{cases}
\end{alignat*}
For any $(L,(v_{ij})_{i\in\dbsqbr{\fc_1},j\in[\fc_2]})\!\in\!\ga_{\fc;\fc_1}^2$
with $L\!=\![(r_i)_{i\in\dbsqbr{\fc_1}}]$, 
there exists $\la\!\in\!\R^{\fc_2}$ so that 
$$v_{ij}=r_i\la_j \qquad\forall~i\!\in\!\dbsqbr{\fc_1},\,j\!\in\![\fc_2].$$
For any homomorphisms 
$h^{11}\!:\R^{\fc_1}\!\lra\!\R^{\fc_1}$,
$h^{21}\!:\R^{\fc_1}\!\lra\!\R^{\fc_2}$, and $h^{22}\!:\R^{\fc_2}\!\lra\!\R^{\fc_2}$,
\begin{alignat*}{2}
\{h^{11}\!\!\times\!h^{22}\}\!:\ga_{\fc;\fc_1}^2&\lra\R^{[\fc_1]\times[\fc_2]}, &~~
\{h^{11}\!\!\times\!h^{22}\}\big(L,(r_i\la_j)_{i\in\dbsqbr{\fc_1},j\in[\fc_2]}\big)
&=h^{11}\big(\!(r_i)_{i\in[\fc_1]}\big)\!\times\!
h^{22}\big(\!(\la_j)_{j\in[\fc_2]}\big),\\
\{h^{21}\!\cdot\!h^{22}\}\!:\ga_{\fc;\fc_1}^2&\lra \R, &~~
\{h^{21}\!\cdot\!h^{22}\}\big(L,(r_i\la_j)_{i\in\dbsqbr{\fc_1},j\in[\fc_2]}\big)
&=h^{21}\big(\!(r_i)_{i\in[\fc_1]}\big)\!\cdot\!
h^{22}\big(\!(\la_j)_{j\in[\fc_2]}\big),
\end{alignat*}
are thus well-defined maps which are linear on each fiber of~$\pi_{\ga_{\fc;\fc_1}^2}$.\\

\noindent
The map given~by
\begin{gather*}
\phi_{\fc;\fc_1}^{12}\!:
\ga_{\fc;\fc_1}^{2;1}\!\equiv\!\big\{\!
\big(L,(v_{ij})_{i\in\dbsqbr{\fc_1},j\in[\fc_2]}\big)
\!\in\!\ga_{\fc;\fc_1}^2\!:
(v_{ij})_{i\in[\fc_1],j\in[\fc_2]}\!\neq\!0\!\in\!\R^{\fc_1\fc_2}\big\}
\lra \R\P^{\fc_2-1},\\
\phi_{\fc;\fc_1}^{12}
\big(L,(r_i\la_j)_{i\in\dbsqbr{\fc_1},j\in[\fc_2]}\big)
=\big[(\la_j)_{j\in[\fc_2]}\big],
\end{gather*} 
is well-defined and smooth.
It satisfies
$$(v_{ij})_{j\in[\fc_2]}\in\phi_{\fc;\fc_1}^{12}
\big(L,(v_{i'j})_{i'\in\dbsqbr{\fc_1},j\in [\fc_2]}\big)
\subset\R^{\fc_2} \quad\forall~i\!\in\!\dbsqbr{\fc_1}.$$
Thus, the map
\begin{gather*}
\wt\phi_{\fc;\fc_1;0}^{12}\!:
\ga_{\fc;\fc_1}^{2;1}\lra \R\P^{\fc-1},\\
\wt\phi_{\fc;\fc_1;0}^{12}
\big(L,(v_{ij})_{i\in\dbsqbr{\fc_1},j\in[\fc_2]}\big)
=\big[\!\big(w\!\cdot\!(v_{ij})_{j\in[\fc_2]}\big)_{i\in[\fc_1]},w\big]~
\forall\,w\!\in\!
\phi_{\fc;\fc_1}^{12}
\big(L,(v_{ij})_{i\in\dbsqbr{\fc_1},j\in [\fc_2]}\big)\!-\!\{0\},
\end{gather*} 
is well-defined and smooth.
The map
$$\wt\phi_{\fc;\fc_1}^{12}\!\equiv\!
\big(\wt\phi_{\fc;\fc_1;0}^{12},(\pi_{\fc;\fc_1;j}^2)_{j\in[\fc]}\big)
\!:\ga_{\fc;\fc_1}^{2;1}\lra 
\ga_{\fc;\fc_1}^{1;2}
\!\equiv\!\big\{\!
\big([(r_j)_{j\in[\fc]}],v\big)\!\in\!\ga^1_{\fc;\fc_1}\!:
(r_j)_{j\in[\fc]-[\fc_1]}\!\neq\!0\big\}$$
is a well-defined diffeomorphism onto an open subspace of~$\ga_{\fc;\fc_1}^1$.
It satisfies
\BE{pifcfc1dfn_e}\pi_{\fc;\fc_1;j}^1\!\circ\!\wt\phi_{\fc;\fc_1}^{12}
=\pi_{\fc;\fc_1;j}^2\big|_{\ga_{\fc;\fc_1}^{2;1}}
\qquad\forall\,j\!\in\![\fc].\EE

\vspace{.15in}

\noindent
The space
$$\ga_{\fc;\fc_1}^{\au}=\big(\ga_{\fc;\fc_1}^1\!\sqcup\!\ga_{\fc;\fc_1}^2
\big)\!\big/\!\sim,\quad 
\ga_{\fc;\fc_1}^1\ni \wt\phi_{\fc;\fc_1}^{12}(\wt{x})
\sim \wt{x}\in \ga_{\fc;\fc_1}^{2;1}\subset\ga_{\fc;\fc_1}^2,$$
is a smooth manifold.
By~\eref{pifcfc1dfn_e}, the map
$$\pi\!:\ga_{\fc;\fc_1}^{\au}\lra\R^{\fc}, \qquad
\pi\big([\wt{r}]\big)=\begin{cases} 
\big(\pi_{\fc;\fc_1;j}^1(\wt{r})\!\big)_{j\in[\fc]},
&\hbox{if}~\wt{r}\!\in\!\ga_{\fc;\fc_1}^1;\\
\big(\pi_{\fc;\fc_1;j}^2(\wt{r})\!\big)_{j\in[\fc]},
&\hbox{if}~\wt{r}\!\in\!\ga_{\fc;\fc_1}^2;
\end{cases}$$
is well-defined and smooth.
With $\wt{U}_{\fc;i}^{\R}\!\subset\!\ga_{\fc;\fc_1}^1$
and $\wt\vph_{\fc;i}$ as in~\eref{vphfcidfn_e} for each $i\!\in\![\fc_1]$, 
the~map
$$\wt\vph_{\fc;\fc_1;i}^1\!\equiv\!
(\wt\vph_{\fc;\fc_1;i;j}^1)_{j\in[\fc]}\!:
\wt{U}_{\fc;\fc_1;i}^1\!\equiv\!\big\{[\wt{x}]\!\in\!\ga_{\fc;\fc_1}^{\au}\!:
\wt{x}\!\in\!\wt{U}_{\fc;i}^{\R}\big\}\lra\R^{\fc},\quad
\wt\vph_{\fc;\fc_1;i}^1\big([\wt{x}]\big)=\wt\vph_{\fc;i}\big(\wt{x}\big),$$
is a coordinate chart on~$\ga_{\fc;\fc_1}^{\au}$.
For $i\!\in\!\dbsqbr{\fc_1}$, the~map
\begin{gather*}
\wt\vph_{\fc;\fc_1;i}^2\!\equiv\!
(\wt\vph_{\fc;\fc_1;i;j}^2)_{j\in[\fc]}\!:
\wt{U}_{\fc;\fc_1;i}^2\!\equiv\!\big\{[\wt{x}]\!\in\!\ga_{\fc;\fc_1}^{\au}\!:
\wt{x}\!=\!\big([(r_{i'})_{i'\in\dbsqbr{\fc_1}}],v\big)\!\in\!
\ga_{\fc;\fc_1}^2,~r_i\!\neq\!0\big\}\lra\R^{\fc},\\
\wt\vph_{\fc;\fc_1;i;j}^2\big(\big[[(r_{i'})_{i'\in\dbsqbr{\fc_1}}],
(v_{i'j'})_{i'\in\dbsqbr{\fc_1},j'\in[\fc_2]}\big]\big)
=\begin{cases}r_{j-1}/r_i,&\hbox{if}~j\!\in\![i];\\
r_j/r_i,&\hbox{if}~j\!\in\![\fc_1]\!-\![i];\\
v_{i(j-\fc_1)},&\hbox{if}~j\!\in\![\fc]\!-\![\fc_1];
\end{cases}
\end{gather*}
is also a coordinate chart on~$\ga_{\fc;\fc_1}^{\au}$.
These charts cover~$\ga_{\fc;\fc_1}^{\au}$.\\

\noindent
Let $X$ be a smooth manifold and $Y\!\subset\!X$
be a closed submanifold of real codimension~$\fc$ and dimension~$m$.
For a coordinate chart 
$$\vph\!\equiv\!(\vph_1,\ldots,\vph_{\fc+m})\!:U\lra\R^{\fc}\!\times\!\R^m$$
for~$Y$ in~$X$, the subspaces
\begin{equation*}\begin{split}
\BLauk_Y\vph&\equiv\big\{\!(\wt{r},x)\!\in\!\ga^k_{\fc;\fc_1}\!\times\!U\!:
\pi_{\fc;\fc_1;j}^k(\wt{r})\!=\!\vph_j(x)\!\in\!\R~\forall\,j\!\in\![\fc]\big\}
~~\hbox{with}~~k=1,2 \quad\hbox{and}\\
\BLau_Y\vph&\equiv \big\{\!(\wt{r},x)\!\in\!\ga^{\au}_{\fc;\fc_1}\!\times\!U\!:
\pi(\wt{r})\!=\!\big(\vph_1(x),\ldots,\vph_{\fc}(x)\!\big)
\!\in\!\R^{\fc}\big\} 
\end{split}\end{equation*}
are closed submanifolds of~$\ga^k_{\fc;\fc_1}\!\times\!U$
 and $\ga^{\au}_{\fc;\fc_1}\!\times\!U$, respectively.
Let
$$\BLauba_Y\vph\equiv\big\{\!(\wt{r},x)\!\in\!\BLaub_Y\vph\!:
\wt{r}\!\in\!\ga^{2;1}_{\fc;\fc_1}\big\}.$$
We denote~by
$$\pi^k_U\!:\BLauk_Y\vph\lra U \quad\hbox{and}\quad
\pi_U\!:\BLau_Y\vph\lra U$$
the projection maps.
The restriction
\BE{AugpiUrestr_e}\pi_U\!:\BLau_Y\vph\!-\!\pi_U^{-1}(Y) \lra U\!-\!Y\EE
is a diffeomorphism.
For $k\!=\!1,2$, $i\!\in\![\fc_1]$ if $k\!=\!1$, and $i\!\in\!\dbsqbr{\fc_1}$ if $k\!=\!2$,
the smooth map
\begin{gather}\label{augBLFchart_e}
\wt\vph_i^k\!\equiv\!\big(\wt\vph_{i;1}^k,\ldots,\wt\vph_{i;\fc+m}^k\big)\!:
\BLauk_{Y;i}\vph\!\equiv\!\big\{\!(\wt{r},x)\!\in\!\BLau_Y\vph\!:
\wt{r}\!\in\!\wt{U}_{\fc;\fc_1;i}^k\big\}
\lra\R^{\fc}\!\times\!\R^m,\\
\notag
\wt\vph_{i;j}^k(\wt{r},x)=\begin{cases}
\wt\vph_{\fc;\fc_1;i;j}^k(\wt{r}),&\hbox{if}~j\!\in\![\fc];\\
\vph_j(x),&\hbox{if}~j\!\in\![\fc\!+\!m]\!-\![\fc];
\end{cases}
\end{gather}
is a coordinate chart on~$\BLau_Y\vph$.\\

\noindent
The subspaces  
\begin{equation*}\begin{split}
\bE_Y^0\vph&=\big\{\![\wt{r},x]\!\in\!\BLau_Y\vph\!:
\wt{r}\!\in\!\R\P^{\fc_1}\!\subset\!\ga^2_{\fc;\fc_1}\big\} \qquad\hbox{and}\\
\bE_Y^-\vph&=\big\{\![\wt{r},x]\!\in\!\BLau_Y\vph\!:
\wt{r}\!\in\!\R\P^{\fc-1}\!\cap\!\ga^1_{\fc;\fc_1}\big\}
\!\cup\!
\big\{\!\big[\big([(r_i)_{i\in\dbsqbr{\fc_1}}],0\big),x\big]\!\in\!\bE_Y^0\!:r_0\!=\!0\big\}
\end{split}\end{equation*}
are closed submanifolds of~$\BLau_Y\vph$ identified along
an $\R\P^{\fc_1-1}$-fiber bundle over~$U\!\cap\!Y$ so that 
$$\pi_U^{-1}(Y)=\bE_Y^0\vph\!\cup\!\bE_Y^-\vph\,.$$
The first submanifold is an $\R\P^{\fc_1}$-fiber bundle.
The second one is a fiber bundle with fibers that themselves form a fiber bundle
over~$\R\P^{\fc_1-1}$ with fibers~$S^{\fc_2}$ and two sections.\\

\noindent
We call an $\R$-atlas 
$\{\vph_{\al}\!:U_{\al}\!\lra\!\R^{\fc}\!\times\!\R^m
\big\}_{\al\in\cI}$ 
for~$Y$ in~$X$ \sf{$\fc_1$-augmented} if for all $\al,\al'\!\in\!\cI$ there exists
$h_{\al\al'}$ as in~\eref{halalprdfn_e} with $h_{\al\al';jj'}\!=\!0$
for $j\!\le\!\fc_1$ and $j'\!>\!\fc_1$ and a smooth function
\BE{augatlcond_e}\begin{split}
&\hspace{1in}f_{\al\al'}\!:\vph_{\al'}\big(U_{\al}\!\cap\!U_{\al'}\big)\lra\R^+
\qquad\hbox{s.t.}\\
&\big|\big(h_{\al\al';jj'}(r)\!\big)_{j,j'\in[\fc]-[\fc_1]}v\big|^2
=f_{\al\al'}(r)|v|^2
\quad\forall~r\!\in\!\vph_{\al'}\big(U_{\al}\!\cap\!U_{\al'}\big),\,
v\!\in\!\R^{\fc_2}.
\end{split}\EE
For such an atlas, $\al,\al'\!\in\!\cI$, and $x\!\in\!U_{\al}\!\cap\!U_{\al'}$, let
\begin{gather*}
h_{\al\al'}^{11}(x)=
\big(h_{\al\al';jj'}(\vph_{\al'}(x)\!)\!\big)_{j,j'\in[\fc_1]}, \quad
h_{\al\al'}^{21}(x)=
\big(h_{\al\al';jj'}(\vph_{\al'}(x)\!)\!\big)_{j\in[\fc]-[\fc_1],j'\in[\fc_1]},\\
h_{\al\al'}^{22}(x)=
\big(h_{\al\al';jj'}(\vph_{\al'}(x)\!)\!\big)_{j,j'\in[\fc]-[\fc_1]}.
\end{gather*}
The endomorphism $h_{\al\al'}^{11}(x)$ of $\R^{\fc_1}$ is invertible for all~$x$
in a neighborhood~$U_{\al\al'}$ of~$U_{\al}\!\cap\!U_{\al'}\!\cap\!Y$ 
in~$U_{\al}\!\cap\!U_{\al'}$.
Let
$$\wt{U}_{\al\al'}^1=\big\{\pi_{U_{\al'}}^1\big\}^{-1}\big(U_{\al\al'}\big)
\subset 
\BLaua_Y\big(\vph_{\al'}|_{U_{\al}\cap U_{\al'}}\big).$$

\vspace{.15in}

\noindent
The map
\begin{gather*}
\psi_{\al\al';0}^2\!\equiv\!\big(\psi_{\al\al';0j}^2\big)_{j\in[\fc_2]}:
\big(U_{\al}\!\cap\!U_{\al'}\big)
\!\times\!\R^{\dbsqbr{\fc_1}\times[\fc_2]}\lra\R^{\fc_2},\\
\psi_{\al\al';0}^2
\big(x,(v_{ij'})_{i\in\dbsqbr{\fc_1},j'\in[\fc_2]}\big)
=h_{\al\al'}^{22}(x)(v_{0j'})_{j'\in[\fc_2]}
\!+\!
\sum_{j'=1}^{\fc_2}v_{0j'}h_{\al\al'}^{21}(x)(v_{ij'})_{i\in[\fc_1]}\,,
\end{gather*}
is smooth.
The smooth function 
\begin{gather*}
F_{\al\al'}\!:\BLaub_Y\big(\vph_{\al'}|_{U_{\al}\cap U_{\al'}}\big)\lra \R,\\
\begin{split}
F_{\al\al'}\big(\!\big(L,(v_{ij'})_{i\in\dbsqbr{\fc_1},j'\in[\fc_2]}\big),x\big)
=f_{\al\al'}\big(\vph_{\al'}(x)\!\big)+
&\sum_{j=1}^{\fc_2}\big|h_{\al\al'}^{21}(x)(v_{ij})_{i\in[\fc_1]}\big|^2\\
+&2\big\{h_{\al\al'}^{21}(x)\!\cdot\!h_{\al\al'}^{22}(x)\!\big\}
\!\big(L,(v_{ij'})_{i\in\dbsqbr{\fc_1},j'\in[\fc_2]}\big),
\end{split}\end{gather*}
takes values in $\R^+$ along 
$\BLaub_Y(\vph_{\al'}|_{U_{\al}\cap U_{\al'}})
\!-\!\BLauba_Y(\vph_{\al'}|_{U_{\al}\cap U_{\al'}})$
and satisfies
\BE{Fprp_e}\begin{split}
\big|\psi_{\al\al';0}^2
\big(x,(v_{ij})_{i\in\dbsqbr{\fc_1},j\in[\fc_2]}\big)\big|^2
&=F_{\al\al'}\big(\!\big(L,(v_{ij})_{i\in\dbsqbr{\fc_1},j\in[\fc_2]}\big),x\big)
\big|(v_{0j})_{j\in[\fc_2]}\big|^2\\
&\quad\forall~
\big(\!\big(L,(v_{ij})_{i\in\dbsqbr{\fc_1},j\in[\fc_2]}\big),x\big)
\in\BLaub_Y\big(\vph_{\al'}|_{U_{\al}\cap U_{\al'}}\big)
\end{split}\EE
by~\eref{augatlcond_e}.\\

\noindent
The smooth~map
\begin{gather*}
\phi_{\al\al'}^2\!:
\BLaub_Y\big(\vph_{\al'}|_{U_{\al}\cap U_{\al'}}\big)\!\times\!\R^{\dbsqbr{\fc_1}}
\lra\R\!\times\!\R^{\fc_1},\\
\begin{split}
&\phi_{\al\al'}^2
\big(\!\big(L,(v_{ij'})_{i\in\dbsqbr{\fc_1},j'\in[\fc_2]}\big),x,
(r_i)_{i\in\dbsqbr{\fc_1}}\big)
=\Big(\!F_{\al\al'}\big(\!\big(L,(v_{ij'})_{i\in\dbsqbr{\fc_1},j'\in[\fc_2]}\big),x\big)r_0,
h_{\al\al'}^{11}(x)(r_i)_{i\in[\fc_1]}\Big),
\end{split}\end{gather*}
is linear in the $\R^{\dbsqbr{\fc_1}}$-factor everywhere and injective over 
a neighborhood~$\wt{U}_{\al\al'}^2$ of 
$$\BLaub_Y\big(\vph_{\al'}|_{U_{\al\al'}}\big)
\!-\!\BLauba_Y\big(\vph_{\al'}|_{U_{\al\al'}}\big)
\subset \BLaub_Y\big(\vph_{\al'}|_{U_{\al\al'}}\big).$$
Since $F_{\al\al'}$ does not vanish over~$\wt{U}_{\al\al'}^2$,
the~map 
\begin{gather*}
\psi_{\al\al'}^2\!\equiv\!\big(\psi_{\al\al';ij}^2\big)_{i\in[\fc_1],j\in[\fc_2]}:
\wt{U}_{\al\al'}^2\lra\R^{[\fc_1]\times[\fc_2]},\\
\begin{split}
&\psi_{\al\al'}^2
\big(\!\big(L,(v_{i'j'})_{i'\in\dbsqbr{\fc_1},j'\in[\fc_2]}\big),x\big)
=\frac{1}{F_{\al\al'}(\!(L,(v_{ij'})_{i\in\dbsqbr{\fc_1},j'\in[\fc_2]}),x)}\\
&\hspace{.1in}
\times\!\!
\bigg(\!\!\big\{\!h^{11}(x)\!\times\!h^{22}(x)\!\big\}\!
\big(L,(v_{i'j'})_{i'\in\dbsqbr{\fc_1},j'\in[\fc_2]}\big)
+\!\sum_{j'=1}^{\fc_2}\{h^{11}(x)\!\}\big(\!(v_{i'j'})_{i'\in[\fc_1]}\big)
\!\times\!\{h^{21}(x)\!\}\big(\!(v_{i'j'})_{i'\in[\fc_1]}\big)\!\!\bigg),
\end{split}
\end{gather*}
is well-defined.
We note~that
\begin{equation*}\begin{split}
&\Big(\psi_{\al\al';0j}^2\big(x,(v_{i'j'})_{i'\in\dbsqbr{\fc_1},j'\in[\fc_2]}\big),
\big(\psi_{\al\al';ij}^2
\big(\!\big(L,(v_{i'j'})_{i'\in\dbsqbr{\fc_1},j'\in[\fc_2]}\big),x\big)\!\big)_{i\in[\fc_1]}\Big)\\
&\hspace{.8in}
\in \phi_{\al\al'}^2\big(\!\big(L,(v_{i'j'})_{i'\in\dbsqbr{\fc_1},j'\in[\fc_2]}\big),x,L\big)
\quad\forall~j\!\in\![\fc_2],\,
\big(\!\big(L,(v_{i'j'})_{i'\in\dbsqbr{\fc_1},j'\in[\fc_2]}\big),x\big)
\!\in\!\wt{U}_{\al\al'}^2.
\end{split}\end{equation*}
Let $\wt{U}_{\al\al'}^{2;1}\!=\!\wt{U}_{\al\al'}^2\!\cap\!
\BLauba_Y(\vph_{\al'}|_{U_{\al\al'}})$.\\

\noindent
Similarly to~\eref{wtvphalalprdfn_e}, the~map
\begin{gather*}
\wt\vph_{\al\al'}^1\!:\wt{U}_{\al\al'}^1
\lra \BLaua_Y\big(\vph_{\al}|_{U_{\al\al'}}\big), \\
\wt\vph_{\al\al'}^1\big(\!(L,(v_j)_{j\in[\fc]}),x\big)
=\big( \big(h_{\al\al'}(\vph_{\al'}(x)\!)L,
h_{\al\al'}(\vph_{\al'}(x)\!)(v_j)_{j\in[\fc]}\big),x\big),
\end{gather*}
is well-defined and smooth.
By~\eref{halalprdfn_e} and~\eref{Fprp_e}, the map 
\begin{gather*}
\wt\vph_{\al\al'}^2\!:\wt{U}_{\al\al'}^2
\lra \BLaub_Y\big(\vph_{\al}|_{U_{\al\al'}}\big),\\
\begin{split}
&\wt\vph_{\al\al'}^2\big(\!\big(L,(v_{i'j'})_{i'\in\dbsqbr{\fc_1},j'\in[\fc_2]}\big),x\big)
=\bigg(\!\!\!\Big(
\phi_{\al\al'}^2\big(\!\big(L,(v_{i'j'})_{i'\in\dbsqbr{\fc_1},j'\in[\fc_2]}\big),x,L\big),\\
&\hspace{.7in}
\big(\psi_{\al\al';0j}^2(x,(v_{i'j'})_{i'\in\dbsqbr{\fc_1},j'\in[\fc_2]}),
\big(\psi_{\al\al';ij}^2
\big(\!\big(L,(v_{i'j'})_{i'\in\dbsqbr{\fc_1},j'\in[\fc_2]}\big),x\big)\!\big)_{i\in[\fc_1]}
\big)_{j\in[\fc_2]}\Big),x\!\!\bigg),
\end{split}\end{gather*}
is also well-defined and smooth.
The two maps are diffeomorphisms onto their images and satisfy
\begin{gather*}
\pi_{U_{\al}}^k\!\circ\!\wt\vph_{\al\al'}^k=
\pi_{U_{\al'}}^k\big|_{\wt{U}_{\al\al'}^k}\,,\quad
\big\{\wt\phi_{\fc;\fc_1}^{12}\!\times\!\id_{U_{\al}}\big\}
\!\circ\!\wt\vph_{\al\al'}^2\big|_{\wt{U}_{\al\al'}^{2;1}}
=\wt\vph_{\al\al'}^1\!\circ\!
\big\{\wt\phi_{\fc;\fc_1}^{12}\!\times\!\id_{U_{\al'}}\big\}
\big|_{\wt{U}_{\al\al'}^{2;1}}.
\end{gather*}
The second equality above on $\wt{U}_{\al\al'}^{2;1}\!-\!\{\pi_{U_{\al'}}^2\}^{-1}(Y)$
follows immediately from the first and the injectivity of the restriction~\eref{AugpiUrestr_e}.
By the continuity of the two sides, this equality must then hold on all 
of~$\wt{U}_{\al\al'}^{2;1}$.
This can also be verified through a direct computation.\\

\noindent
Thus, the map
\begin{gather*}
\wt\vph_{\al\al'}\!:
\BLau_Y\big(\vph_{\al'}|_{U_{\al}\cap U_{\al'}}\big)
\lra \BLau_Y\big(\vph_{\al}|_{U_{\al}\cap U_{\al'}}\big),\\
\wt\vph_{\al\al'}\big([\wt{x}']\big)=
\begin{cases}\big[\wt\vph_{\al\al'}^k(\wt{x}')\big],\!\!
&\hbox{if}~\wt{x}'\!\in\!\wt{U}_{\al\al'}^k,\,k\!=\!1,2;\\
[\wt{x}],&
\hbox{if}~\pi_{U_{\al'}}([\wt{x}'])\!=\!\pi_{U_{\al}}([\wt{x}])\!\not\in\!Y;
\end{cases}
\end{gather*}
is a well-defined diffeomorphism that satisfies
$$\pi_{U_{\al}}\!\circ\!\wt\vph_{\al\al'}=
\pi_{U_{\al'}}\big|_{\BLau_Y\big(\vph_{\al'}|_{U_{\al}\cap U_{\al'}}\big)}.$$
By the uniqueness of continuous extensions and the cocycle condition
for the overlap maps~$\vph_{\al\al'}$,
the maps~$\wt\vph_{\al\al'}$ satisfy the cocycle condition~\eref{cocyclecond_e}
with~$\bF\!=\!\fc_1$.\\

\noindent
For a $\fc_1$-augmented atlas $\{\vph_{\al}\}_{\al\in\cI}$ for~$Y$ in~$X$,
we define the corresponding $\fc_1$-augmented blowup
\BE{blowmapaugdfn_e}\pi\!:\BLau_YX\!\equiv\!\Big(\!\!
(X\!-\!Y)\!\sqcup\!\bigsqcup_{\al\in\cI}\!\!\BLau_Y\vph_{\al}\!\Big)\!\!\Big/\!\!\!\sim
\lra X\EE
as in~\eref{blowmapdfn_e} and just above with~$\bF$ replaced by~$\fc_1$.
This map is again surjective, proper, and smooth.
The subspaces
$$\bE^0_YX\equiv
\bigg(\bigsqcup_{\al\in\cI}\!\!\bE^0_Y\vph_{\al}\!\!\bigg)\!\!\Big/\!\!\!\sim
\,\subset\BLau_YX 
\quad\hbox{and}\quad
\bE^-_YX\equiv
\bigg(\bigsqcup_{\al\in\cI}\!\!\bE^-_Y\vph_{\al}\!\!\bigg)\!\!\Big/\!\!\!\sim
\,\subset\BLau_YX$$
are smooth closed submanifolds of codimensions~$\fc_2$ and~1, respectively,
intersecting transversely.
The \sf{exceptional locus} for the \sf{blowdown map}~$\pi$ in~\eref{blowmapaugdfn_e} is
$$\bE^{\fc_1}_YX\equiv
\pi^{-1}(Y)=\bE^0_YX\!\cup\!\bE^-_YX\,.$$
The blowdown map~\eref{blowmapaugdfn_e} restricts to a diffeomorphism~\eref{augdiff_e}.\\

\noindent
The restrictions of the vector bundle isomorphisms~\eref{NXYisom_e} with $\bF\!=\!\R$
to the subbundles 
\BE{NXYisomsub_e}0^{\fc_1}\!\!\times\!\R^{\fc_2}\!\times\!(U_{\al}\!\cap\!Y)
\subset \R^{\fc_1}\!\!\times\!\R^{\fc_2}\!\times\!(U_{\al}\!\cap\!Y)\EE
determine a corank~$\fc_1$-subbundle $\cN_Y^{\fc_1}X\!\subset\!\cN_YX$.
Let $\wt{\cN}_Y^{\fc_1}X\!\subset\!TX|_Y$ be its preimage under the quotient~map.
The restrictions of these isomorphisms to the complements of the subbundles
in~\eref{NXYisomsub_e} induce an identification of
$\bE^-_YX\!-\!\bE^0_YX$ with $\R\P(\cN_YX)\!-\!\R\P(\cN_Y^{\fc_1}X)$.
The smooth functions in~\eref{augatlcond_e} determine a (trivializable) line bundle
\begin{gather*}
\cL_Y^{\fc_1}X\equiv
\bigg(\bigsqcup_{\al\in\cI}\!\!\{\al\}\!\times\!(U_{\al}\!\cap\!Y)\!\times\!\R\!\!\bigg)
\!\!\Big/\!\!\!\sim\lra Y, \\
\big(\al,x,f_{\al\al'}(x)r_0\big)\sim(\al',x,r_0)
\quad \forall\,x\!\in\!U_{\al}\!\cap\!U_{\al'}\!\cap\!Y,\,
\al,\al'\!\in\!\cI.
\end{gather*}
The vector bundle isomorphisms
\begin{gather*}
\R^{\dbsqbr{\fc_1}}\!\times\!(U_{\al}\!\cap\!Y)\lra 
\big(\cN_YX/\cN_Y^{\fc_1}X\!\oplus\!\cL_Y^{\fc_1}X\big)\!\big|_{U_{\al}\cap Y},\\ 
\big(\!(r_i)_{i\in\dbsqbr{\fc_1}},x\big)\lra
\bigg(\sum_{i=1}^{\fc_1}r_i\frac{\prt}{\prt\vph_{\al;i}}\Bigg|_x
\!+\!\wt{\cN}_Y^{\fc_1}X|_x,r_0\!\!\bigg),
\quad\al\!\in\!\cI,
\end{gather*}
induce identifications of $\bE^0_YX$ and $\bE^0_YX\!\cap\!\bE^-_YX$ with 
$$\R\P\big(\cN_YX/\cN_Y^{\fc_1}X\!\oplus\!\cL_Y^{\fc_1}X\big)
\quad\hbox{and}\quad 
\R\P\big(\cN_YX/\cN_Y^{\fc_1}X\!\oplus\!\tau_Y^0\big),$$ 
respectively, as fiber bundles over~$Y$.
If $\fc_1\!=\!1$, the isomorphisms~\eref{NXYisom_e} induce an identification 
of~$\bE^-_YX$ with the $S^{\fc_2}$-sphere bundle over~$Y$ obtained by collapsing 
each fiber of the $\R\P^{\fc_2-1}$-fiber bundle 
\hbox{$\R\P(\cN_Y^{\fc_1}X)\!\subset\!\R\P(\cN_YX)$} to a point.
The section of \hbox{$\pi^{-1}(Y)\!\lra\!Y$}
determined by this collapse is~$\bE^0_YX\!\cap\!\bE^-_YX$.\\

\noindent
We define two $\fc_1$-augmented atlases for~$Y$ in~$X$
to be \sf{$\fc_1$-augmented equivalent} if their union is still  
a $\fc_1$-augmented atlas for~$Y$ in~$X$.
A \sf{co-augmented structure for~$Y$ in~$X$} is a $\fc_1$-augmented equivalence class of 
augmented atlases for~$Y$ in~$X$.
Since the above construction can be carried out with the maximal collection of charts
in a co-augmented structure for~$Y$ in~$X$,
the diffeomorphism class of the projection~$\pi$ in~\eref{blowmapaugdfn_e}
is determined by such a structure, but may depend on~it.

\begin{lmm}\label{AugBl_lmm}
Suppose $X$, $Y$, $\fc$, and $\pi'$ are as in Lemma~\ref{StanBl_lmm} with $\bF\!=\!\R$,
$\fc_1\!\in\![\fc\!-\!1]$, and
$\{\vph_{\al}\!\equiv\!(\vph_{\al;1},\ldots,\vph_{\al;\fc+m})\!:
U_{\al}\!\lra\!\R^{\fc+m}\}_{\al\in\cI}$ 
is a \hbox{$\fc_1$-augmented} atlas of charts for~$Y$ in~$X$.
Let \hbox{$\pi\!:\BLau_YX\!\lra\!X$} be the $\fc_1$-augmented blowup of~$X$ along~$Y$
determined by this atlas.
If there exists a collection 
$$\big\{\big(\wt\phi_{\al;i}^k\!\equiv\!
(\wt\phi_{\al;i;1}^k,\ldots,\wt\phi_{\al;i;\fc+m}^k)
\!:\wt{U}_{\al;i}^k\!\lra\!\R^{\fc+m}\big)\!:
k\!=\!1,2,\,i\!\in\![\fc_1]\,\hbox{if}\,k\!=\!1,\,
i\!\in\!\dbsqbr{\fc_1}\,\hbox{if}\,k\!=\!2\big\}_{\al\in\cI}$$
of charts on~$\wt{X}$ so that 
\begin{gather*}
\pi'^{-1}(U_{\al})=\bigcup_{i=1}^{\fc_1}\!\wt{U}_{\al;i}^1\cup
\bigcup_{i=0}^{\fc_1}\!\wt{U}_{\al;i}^2, \quad
\wt\phi_{\al;i;j}^k=\vph_{\al;j}\!\circ\!\pi'|_{\wt{U}_{\al;i}^k}
~~\hbox{if}~j\!>\!\fc~\hbox{or}~(k,j)\!=\!(1,i),\\
\vph_{\al;j}\!\circ\!\pi'|_{\wt{U}_{\al;i}^1}=
\wt\phi_{\al;i;j}^1\!\cdot\!\{\vph_{\al;i}\!\circ\!\pi'|_{\wt{U}_{\al;i}^1}\}
~\hbox{if}\,j\!\in\![\fc]\!-\!\{i\},\\
\vph_{\al;j}\!\circ\!\pi'|_{\wt{U}_{\al;i}^2}=
\wt\phi_{\al;i;j}^2\!\cdot\!
\begin{cases}1,&\hbox{if}~i\!=\!0,\,j\!\in\![\fc]\!-\![\fc_1];\\
\wt\phi_{\al;i;1}^2,&\hbox{if}~i\!\in\![\fc_1],\,j\!\in\![\fc]\!-\![\fc_1];
\end{cases}\\
\vph_{\al;j}\!\circ\!\pi'|_{\wt{U}_{\al;i}^2}=
\bigg(\!\sum_{j'=\fc_1+1}^{\fc}\!\!\!\!(\wt\phi_{\al;i;j'}^2)^2\!\!\bigg)
\!\cdot\!
\begin{cases}\wt\phi_{\al;i;j}^2,&\hbox{if}~i\!=\!0,\,j\!\le\!\fc_1;\\
\wt\phi_{\al;i;1}^2\wt\phi_{\al;i;j}^2,&\hbox{if}~1\!\le\!i\!<\!j\!\le\!\fc_1;\\
\wt\phi_{\al;i;1}^2,&\hbox{if}~i\!=\!j;\\
\wt\phi_{\al;i;1}^2\wt\phi_{\al;i;j+1}^2,&\hbox{if}~i\!>\!j;
\end{cases}
\end{gather*} 
for every~$\al\!\in\!\cI$, then the maps $\pi$ and $\pi'$ are diffeomorphic.
\end{lmm}

\begin{proof}
Similarly to~\eref{StanBl_e1},
\BE{AugBl_e1}\begin{split}
\wt{U}_{\al;i}^1\!\cap\!\pi'^{-1}(Y)&=\big\{\wt\phi_{\al;i;i}^1\big\}^{-1}(0),\\
\wt{U}_{\al;i}^2\!\cap\!\pi'^{-1}(Y)&=
\bigcap_{j=\fc_1+1}^{\fc}\!\!\!\!\!\big\{\wt\phi_{\al;i;j}^2\big\}^{-1}(0)
\!\cup\!\begin{cases}\eset,&\hbox{if}~i\!=\!0;\\
\{\wt\phi_{\al;i;1}^2\}^{-1}(0),&\hbox{if}~i\!\in\![\fc_1].
\end{cases}
\end{split}\EE
For $\al\!\in\!\cI$ and $i\!\in\![\fc_1]$, let
$$\wt\psi_{\al;i}^1\!\equiv\!
\big(\wt\psi_{\al;i;j}^1\big)_{j\in[\fc_1]}\!:
\wt{U}_{\al;i}^1\lra\R^{\fc}\!-\!\{0\}, \quad
\wt\psi_{\al;i;j}^1(\wt{x})=\begin{cases}
\wt\phi_{\al;i;j}^1(\wt{x}),&\hbox{if}~j\!\in\![\fc]\!-\!\{i\};\\
1,&\hbox{if}~j\!=\!i.
\end{cases}$$
For $i\!\in\!\dbsqbr{\fc_1}$, let
$$\wt\psi_{\al;i}^{(2)}\!\equiv\!
\big(\wt\psi_{\al;i;i'}^2\big)_{i'\in\dbsqbr{\fc_1}}
\!:\wt{U}_{\al;i}^2\lra\R^{\dbsqbr{\fc_1}}\!-\!\{0\}, \quad
\wt\psi_{\al;i;i'}^2(\wt{x})=\begin{cases}
\wt\phi_{\al;i;i'+1}^2(\wt{x}),&\hbox{if}~i'\!<\!i;\\
1,&\hbox{if}~i'\!=\!i;\\
\wt\phi_{\al;i;i'}^2(\wt{x}),&\hbox{if}~i'\!>\!i.
\end{cases}$$
Combining the conditions on the coordinate charts~$\{\wt\phi_{\al;i}^k\}$
with the coordinate charts~\eref{augBLFchart_e}, we find that the~maps
\begin{alignat*}{3}
\wt\Psi_{\al;i}^1\!:\wt{U}_{\al;i}^1&\lra\BLaua_{Y;i}\vph_{\al},
&~~
\wt\Psi_{\al;i}^1(\wt{x})&=\big(\!\big(\R\wt\psi_{\al;i}^1(\wt{x}),
\wt\phi_{\al;i;i}^1(\wt{x})\wt\psi_{\al;i}^1(\wt{x})\!\big),
\pi'(\wt{x})\!\big), &\quad &i\!\in\![\fc_1],\\
\wt\Psi_{\al;i}^2\!:\wt{U}_{\al;i}^2&\lra\BLaub_{Y;i}\vph_{\al},
&~~
\wt\Psi_{\al;i}^2(\wt{x})&=\big(\!\big[\R\wt\psi_{\al;i}^2(\wt{x}),
\big(\wt\phi_{\al;i;\fc_1+j}^2(\wt{x})\wt\psi_{\al;i}^2(\wt{x})\!
\big)_{j\in[\fc-\fc_1]}\big],\pi'(\wt{x})\!\big),
&\quad &i\!\in\!\dbsqbr{\fc_1},
\end{alignat*}
are well-defined diffeomorphisms onto open subsets of~$\BLaua_{Y;i}\vph_{\al}$
and~$\BLaub_{Y;i}\vph_{\al}$, respectively, so~that
$$\pi\big([\wt\Psi_{\al;i}^k(\wt{x})]\big)=\pi'(\wt{x})\in X
\quad\forall\,\wt{x}\!\in\!\wt{U}_{\al;i}^k.$$

\vspace{.15in}

\noindent
Similarly to the proof of Lemma~\ref{StanBl_lmm},
the injectivity of the restriction~\eref{augdiff_e} and~\eref{AugBl_e1} imply~that
\begin{equation*}\begin{split}
\big[\wt\Psi_{\al;i}^{(k)}(\wt{x})\big]
=\big[\wt\Psi_{\al';i'}^{k'}(\wt{x})\big]\in\BLau_YX
\quad\forall~
&\wt{x}\!\in\!\wt{U}_{\al;i}^k\!\cap\!\wt{U}_{\al';i'}^{k'},\,
\al,\al'\!\in\!\cI,\\
&(k,i),(k',i')\!\in\!\big\{\!(1,j)\!:j\!\in\![\fc_1]\big\}\!\cup\!
\big\{\!(2,j)\!:j\!\in\!\dbsqbr{\fc_1}\big\}.
\end{split}\end{equation*}
Thus, the~map
$$\wt\Psi\!:\wt{X}\lra\BLau_YX, \qquad
\wt\Psi(\wt{x})=\begin{cases}[\wt\Psi_{\al;i}^1(\wt{x})],
&\hbox{if}~\wt{x}\!\in\!\wt{U}_{\al;i}^1,\,\al\!\in\!\cI,\,i\!\in\![\fc_1];\\
[\wt\Psi_{\al;i}^2(\wt{x})],
&\hbox{if}~\wt{x}\!\in\!\wt{U}_{\al;i}^2,\,\al\!\in\!\cI,\,i\!\in\!\dbsqbr{\fc_1};\\
[\pi'(\wt{x})],&\hbox{if}~\wt{x}\!\in\!\wt{X}\!-\!\pi'^{-1}(Y);
\end{cases}$$
is a well-defined smooth submersion onto an open subset of~$\BLau_YX$.
It satisfies \hbox{$\pi'\!=\!\pi\!\circ\!\wt\Psi$}. 
Since the restriction~\eref{StanBl_e} is injective, 
\eref{AugBl_e1} implies that $\wt\Psi$ is injective as~well.
Since the restriction~\eref{StanBl_e} is surjective, 
the image of~$\wt\Psi$ contains \hbox{$\BLau_YX\!-\!\pi^{-1}(Y)$}. 
Since $\pi'$ is a proper map, 
$$\wt\Psi\big(\pi'^{-1}(x)\!\big)\subset \pi^{-1}(x)
\subset \pi^{-1}(Y)\subset \BLau_YX$$
is a nonempty, compact, open subspace of~$\pi^{-1}(x)$ for every $x\!\in\!Y$.
Thus, $\pi'$ is onto.
\end{proof}

\subsection{Real blowups and involutions}
\label{Rblowup_subs}

\noindent
Let $X$ be a smooth manifold and $Y\!\subset\!X$
be a closed submanifold of real codimension~3 and real dimension~$m$.
We now give a different characterization of the real blowup~$\BLR_YX$
of~$X$ along~$Y$, viewing the normal direction~$\R^3$ of~$Y$ in~$X$
as $\C\!\times\!\R$.
In particular, an element of~$\R\P^2$ is denoted below by~$[z,r]$
with \hbox{$(z,r)\!\in\!\C\!\times\!\R$} nonzero.\\

\noindent
Suppose 
$$\vph\!\equiv\!(\vph_1,\ldots,\vph_{2+m})\!:U\!\lra\!
\C\!\times\!\R\!\times\!\R^m$$
is a chart for~$Y$ in~$X$.
Let
\BE{BLRchartR_e}
\wt\vph_{\R}\!\equiv\!
\big(\wt\vph_{\R;j}\big)_{j\in[2+m]}\!\equiv\!
\big(\wt\vph_{3;j}\big)_{j\in[2+m]}\!:
\BLR_{Y;\R}\vph
\!\equiv\!\big\{\!([z,r],x)\!\in\!\BLR_Y\vph\!:r\!\neq\!0\big\}
\lra\C\!\times\!\R\!\times\!\R^m\EE
be a coordinate chart on~$\BLR_Y\vph$ as in~\eref{BLFchart_e}.
The smooth~map
\begin{gather}\label{BLRchartC_e}
\wt\vph_{\C}\!\equiv\!\big(\wt\vph_{\C;j}\big)_{j\in\dbsqbr{2+m}}\!:
\BLR_{Y;\C}\vph\!\equiv\!\big\{\!([z,r],x)\!\in\!\BLR_Y\vph\!:
z\!\neq\!0\big\}
\lra S^1\!\times\!\C^2\!\times\!\R^m,\\
\notag
\wt\vph_{\C;j}\big([z,r],x\big)
=\begin{cases}
|z|^2/z^2,&\hbox{if}~j\!=\!0;\\
r/z,&\hbox{if}~j\!=\!2;\\
\vph_j(x),&\hbox{if}~j\!\in\!\{1\}\!\cup\!\big([2\!+\!m]\!-\![2]\big);
\end{cases}
\end{gather}
is an embedding.
The open subsets $\BLR_{Y;\R}\vph,\BLR_{Y;\C}\vph\!\subset\!\BLR_Y\vph$
cover~$\BLR_Y\vph$.\\

\noindent
If $W\!\subset\!\C\!\times\!\R^m$ and \hbox{$g\!:W\!\lra\!\C^*$}, let
$$F_g\!:W\lra\C\!\times\!\R^m,  ~~
F_g(z,r)=\big(\ov{zg(z,r)},r\big),\quad
W_g=\big\{\!(z,r)\!\in\!W\!:F_g(z,r)\!=\!(z,r)\!\big\}.$$
If $W$ is an open subspace, $g$ is a smooth map, and $F_g$ is an involution, then
$g(0,r)\!\in\!S^1$ for all $(0,r)\!\in\!W$ and
\hbox{$W_g\!\subset\!W$} is a submanifold so~that
\begin{gather*}
W\!\cap\!\big(0\!\times\!\R^m\big)\subset W_g,\\
T_{(0,r)}W_g=
\big\{\dot{z}\!\in\!\C\!:
\sqrt{g(0,r)}\dot{z}\!\in\!\R\big\}\!\oplus\!\R^m\subset 
T_{(0,r)}\big(\C\!\times\!\R^m\big)
~~\forall~(0,r)\!\in\!W\!\subset\!\C\!\times\!\R^m\,;
\end{gather*}
see \cite[Lemma~3]{Meyer81}.
We call a~pair
\BE{quasichartdfn_e}\vph\!\equiv\!(\vph_j)_{j\in[2+\fc]}\!: U\lra\C^2\!\times\!\R^m
\quad\hbox{and}\quad g\!:W\lra\C^*\EE
of smooth maps from open subspaces of~$U\!\subset\!X$ 
and $W\!\subset\!\C\!\times\!\R^m$ a \sf{type~1 pre-chart on~$X$}
if $F_g$ is an involution on~$W$ and 
$$\vph\!:U\lra\C\!\times\!W_g\subset \C^2\!\times\!\R^m$$
is a diffeomorphism onto an open subset of~$\C\!\times\!W_g$.
We call the pair~\eref{quasichartdfn_e} 
a \sf{type~1 pre-chart for~$Y$ in~$X$}
if in addition the condition~\eref{vphYcond_e} with $\fc\!=\!2$ is satisfied.
If~$g$ is the constant function~1, then $\vph$ is a chart in the usual sense.
We call a collection 
\BE{quasiatldfn_e}\big\{\!(\vph_{\al}\!:U_{\al}\!\lra\!\C^2\!\times\!\R^m,
g_{\al}\!:W_{\al}\!\lra\!\C^*)\!\big\}_{\al\in\cI}\EE
of type~1 pre-charts for~$Y$ in~$X$ a \sf{quasi-atlas for $Y$ in~$X$}
if the open subsets $U_{\al}\!\subset\!X$ cover~$Y$.\\

\noindent
If $W\!\subset\!\C^*\!\times\!\C^2\!\times\!\R^m$ 
and \hbox{$g\!\equiv\!(g_0,g_1,g_2)\!:W\!\lra\!(\C^*)^3$}, 
we similarly define
\begin{gather*}
F_g^{\circ}\!:W\lra\C^3\!\times\!\R^m, \quad
F_g^{\circ}\big(z\!\equiv\!(z_j)_{j\in[3]},r\big)
=\big(\ov{z_0^{-1}g_0(z_0,r)},
\ov{z_1z_0g_1(z,r)},\ov{z_2z_0^{-1}g_2(z,r)},r\big),\\
W_g^{\circ}=\big\{\!(z,r)\!\in\!W\!:F_g^{\circ}(z,r)\!=\!(z,r)\!\big\}.
\end{gather*}
If $W$ is an open subspace, $g$ is a smooth map, 
$F_g^{\circ}$ is an involution, and $g_0(z_0,0,r)\!=\!1$ 
for all $(z_0,0,r)\!\in\!W$, then
$g_j(z_0,0,r)\!\in\!S^1$ for all \hbox{$(z_0,0,r)\!\in\!W$} 
with $z_0\!\in\!S^1$ and
\hbox{$W_g^{\circ}\!\subset\!W$} is a submanifold so~that
\begin{gather}\label{pre12chart_e3}
W\!\cap\!\big(S^1\!\times\!0^2\!\times\!\R^m\big)
= W_g^{\circ}\!\cap\!\big(\C^*\!\times\!0^2\!\times\!\R^m\big),\\
\notag
T_zW_g^{\circ}=
T_{z_0}S^1\!\oplus\!
\big\{\!(\dot{z}_1,\dot{z}_2)\!\in\!\C^2\!:
\sqrt{z_0g_1(z)}\dot{z}_1,\sqrt{z_0^{-1}g_2(z)}\dot{z}_2\!\in\!\R\big\}\!\oplus\!\R^m
\subset T_zW~~\forall~z\!\equiv\!(z_0,0,r)\!\in\!W_g^{\circ};
\end{gather}
see \cite[Lemma~3]{Meyer81}.
We call a~pair
\BE{quasichartdfn3_e}\vph\!\equiv\!(\vph_j)_{j\in\dbsqbr{2+\fc}}\!: 
U\lra\C^*\!\times\!\C^2\!\times\!\R^m
\quad\hbox{and}\quad g\!\equiv\!(g_j)_{j\in\dbsqbr{2}}:W\lra(\C^*)^3 \EE
of smooth maps from open subspaces of~$U\!\subset\!X$ 
and $W\!\subset\!\C^*\!\times\!\C^2\!\times\!\R^m$ a \sf{type~$(1,2)$ pre-chart on~$X$}
if  $F_g^{\circ}$ is an involution on~$W$,
$g_0(z_0,0,r)\!=\!1$ for all $(z_0,0,r)\!\in\!W$, 
the functions
\BE{g1cond_e}S^1\lra S^1, \qquad z_0\lra g_1(z_0,0,r),\EE
with $r\!\in\!\R^m$ are constant, and 
$$\vph\!:U\lra W_g^{\circ}\subset \C^*\!\times\!\C^2\!\times\!\R^m$$
is a diffeomorphism onto an open subset of~$W_g^{\circ}$.

\begin{lmm}\label{Rblowup_lmm}
Suppose $X$, $Y$, and $\pi'$ are as in Lemma~\ref{StanBl_lmm}
with \hbox{$(\bF,\fc)\!=\!(\R,3)$} and~\eref{quasiatldfn_e}
is a quasi-atlas for $Y$ in~$X$.
Let \hbox{$\pi\!:\BLR_YX\!\lra\!X$} be the real blowup of~$X$ along~$Y$.
If for every $\al\!\in\!\cI$ there exist open subsets 
\hbox{$\wt{U}_{\al;1},\wt{U}_{\al;2}\!\subset\!\pi'^{-1}(U_{\al})$},
\hbox{$\wt{W}_{\al}\!\subset\!\C^*\!\times\!\C^2\!\times\!\R^m$},
and smooth maps
\begin{gather*}
\wt\phi_{\al;1}\!\equiv\!\big(\wt\phi_{\al;1;j}\big)_{j\in\dbsqbr{2+m}}\!:
\wt{U}_{\al;1}\lra\C^*\!\times\!\C^2\!\times\!\R^m, \quad
\wt\phi_{\al;2}\!\equiv\!\big(\wt\phi_{\al;2;j}\big)_{j\in[2+m]}\!:
\wt{U}_{\al;2}\lra\C^2\!\times\!\R^m, \\
\hbox{and}\qquad
\wt{g}_{\al}\!\equiv\!(\wt{g}_{\al;j})_{j\in\dbsqbr{2}}\!:\wt{W}_{\al}\lra(\C^*)^3
\end{gather*}
such that  
$(\wt\phi_{\al;1},\wt{g}_{\al})$ is a type~$(1,2)$ pre-chart on~$\wt{X}$,
$(\wt\phi_{\al;2},g_{\al})$ is a type~1 pre-chart on~$\wt{X}$,
\begin{gather}\label{Rblowup_e}
\pi'^{-1}\big(U_{\al}\big)\!-\!\wt{U}_{\al;2}
=\big\{\wt{x}\!\in\!\wt{U}_{\al;1}\!:
\wt\phi_{\al;1;1}(\wt{x}),\wt\phi_{\al;1;2}(\wt{x})\!=\!0\big\},\\
\notag
\vph_{\al;j}\!\circ\!\pi'|_{\wt{U}_{\al;i}}=\begin{cases}
\wt\phi_{\al;i;j}\!\cdot\!\{\vph_{\al;i}\!\circ\!\pi'|_{\wt{U}_{\al;i}}\},
&\hbox{if}~j\!=\!3\!-\!i;\\
\wt\phi_{\al;i;j},&\hbox{if}~j\!\in\!\{i\}\!\cup\!\big([2\!+\!m]\!-\![2]);
\end{cases}\\
\wt{g}_{\al;1}\big(\!(z_j)_{j\in\dbsqbr{2}},r\big)
\wt{g}_{\al;2}\big(\!(z_j)_{j\in\dbsqbr{2}},r\big)
=g_{\al}\big(z_1z_2,r\big)
\quad\forall\,\big(\!(z_j)_{j\in\dbsqbr{2}},r)\!\in\!\wt{W}_{\al},
\end{gather}
then the maps $\pi$ and $\pi'$ are diffeomorphic.
\end{lmm}

\begin{proof}
Let $\al\!\in\!\cI$.
The property~\eref{StanBl_e1} still holds.
For every $x\!\in\!U_{\al}\!\cap\!Y$, there exist 
a neighborhood $U_{\al;x}\!\subset\!U_{\al}$ of~$x$
and a smooth function $g_{\al;x}\!:U_{\al;x}\!\lra\!S^1$ so~that 
$$g_{\al;x}(x')^2=g_{\al}\big(\!(\vph_{\al;j}(x')\!)_{j\in[2+m]-\{1\}}\big)
\qquad\forall~x'\!\in\!U_{\al;x}.$$
For $i\!=\!1,2$, let $\wt{U}_{\al;x;i}\!=\!\wt{U}_{\al;i}\!\cap\!\pi'^{-1}(U_{\al;x})$.\\

\noindent
If $U_{\al;x}$ is sufficiently small, the smooth maps
\begin{alignat*}{2}
\vph_{\al;x}\!\equiv\!(\vph_{\al;x;j})_{j\in[2+\fc]}\!: 
U_{\al;x}&\lra\C\!\times\!\R\!\times\!\R^m,
&~~
\vph_{\al;x;j}(x')&=
\begin{cases}
\vph_{\al;2}(x')g_{\al;x}(x'),&\hbox{if}~j\!=\!2;\\
\vph_{\al;j}(x'),&\hbox{if}~j\!\neq\!2;
\end{cases}\\
\wt\phi_{\al;x;2}\!\equiv\!(\wt\phi_{\al;x;2;j})_{j\in[2+\fc]}\!: 
\wt{U}_{\al;x;2}&\lra\C\!\times\!\R\!\times\!\R^m,
&~~
\wt\phi_{\al;x;2;j}(\wt{x}')&=
\begin{cases}
\wt\phi_{\al;2;1}(\wt{x}')/g_{\al;x}(\pi'(\wt{x}')\!),&\hbox{if}~j\!=\!1;\\
\wt\phi_{\al;2;2}(\wt{x}')g_{\al;x}(\pi'(\wt{x}')\!),&\hbox{if}~j\!=\!2;\\
\wt\phi_{\al;2;j}(\wt{x}'),&\hbox{if}~j\!>\!2;
\end{cases}
\end{alignat*}
are a coordinate chart for~$Y$ in~$X$ and a coordinate chart on~$\wt{X}$, respectively,
such that
$$\vph_{\al;x;j}\!\circ\!\pi'|_{\wt{U}_{\al;x;2}}=\begin{cases}
\wt\phi_{\al;x;2;1}\!\cdot\!\{\vph_{\al;x;2}\!\circ\!\pi'|_{\wt{U}_{\al;x;2}}\},
&\hbox{if}~j\!=\!1;\\
\wt\phi_{\al;x;2;j},&\hbox{if}~j\!\in\![2\!+\!m]\!-\!\{1\}.
\end{cases}$$
By~\eref{Rblowup_e}, \eref{pre12chart_e3} with $g_j\!=\!\wt{g}_{\al;j}$, 
and the assumption that the functions~\eref{g1cond_e} with 
$g_1\!=\!\wt{g}_{\al;1}$ are constant,
the smooth~map
\begin{gather*}
\wt\phi_{\al;x;1}\!\equiv\!\big(\wt\phi_{\al;x;1;j}\big)_{j\in\dbsqbr{2+m}}\!:
\wt{U}_{\al;x;1}\lra S^1\!\times\!\C^2\!\times\!\R^m, \\
\wt\phi_{\al;x;1;j}(\wt{x})=\begin{cases}
\wt\phi_{\al;1;0}(\wt{x})\wt{g}_{\al;1}(\wt\phi_{\al;1}(\wt{x})),&\hbox{if}~j\!=\!0;\\
\wt\phi_{\al;1;2}(\wt{x})g_{\al;x}\big(\pi'(\wt{x})\!\big),&\hbox{if}~j\!=\!2;\\
\wt\phi_{\al;1;j}(\wt{x}),&\hbox{if}~j\!\in\!\{1\}\!\cup\!([2\!+\!m]\!-\![2]\big);
\end{cases}
\end{gather*}
restricts to a diffeomorphism from a neighborhood of 
$$\pi'^{-1}\big(U_{\al;x}\big)\!-\!\wt{U}_{\al;x;2}\subset\wt{U}_{\al;x;1}$$
onto a codimension~2 submanifold of the target.
Let
$$\wt\vph_{\al;x;\R}\!:\BLR_{Y;\R}\vph_{\al;x}\lra\C\!\!\times\!\R\!\times\!\R^m
\quad\hbox{and}\quad
\wt\vph_{\al;x;\C}\!:\BLR_{Y;\C}\vph_{\al;x}\lra S^1\!\times\!\C^2\!\times\!\R^m$$
be the coordinate chart as  in~\eref{BLRchartR_e} and
 embedding as in~\eref{BLRchartC_e}, respectively, 
corresponding to~$\vph_{\al;x}$.\\

\noindent
For $i\!=\!1,2$, define
\begin{gather*}
\wt\psi_{\al;x;i}\!\equiv\!\big(\wt\psi_{\al;x;i;1},\wt\psi_{\al;x;i;2}\big)
\!:\wt{U}_{\al;x;i}\lra\C\!\times\!\R\!-\!\{0\}, \\
\begin{aligned}
\wt\psi_{\al;x;1;1}(\wt{x})&=
\big(\wt\phi_{\al;1;0}(\wt{x})\wt{g}_{\al;1}\!\big(\wt\phi_{\al;1}(\wt{x})\!\big)\!\big)^{-1/2},
&\quad
\wt\psi_{\al;x;2;1}(\wt{x})&=\wt\phi_{\al;2;1}(\wt{x})
g_{\al;x}\big(\pi'(\wt{x})\!\big)^{-1},\\
\wt\psi_{\al;x;1;2}(\wt{x})&=
\wt\psi_{\al;x;1;1}(\wt{x})\wt\phi_{\al;1;2}(\wt{x})
g_{\al;x}\big(\pi'(\wt{x})\!\big),
&\quad
\wt\psi_{\al;x;2;2}(\wt{x})&=1.
\end{aligned}
\end{gather*}
While $\wt\psi_{\al;x;1;1}(\wt{x})\!\in\!S^1$ is defined only up to the multiplication
by~$-1$, the~map
$$\wt\Psi_{\al;x;1}\!:\wt{U}_{\al;x;1}\lra\BLR_{Y;\C}\vph_{\al;x}
\!\subset\!\BLR_Y\vph_{\al;x}, \quad
\wt\Psi_{\al;x;1}(\wt{x})=\big(\R\wt\psi_{\al;x;1}(\wt{x}),\pi'(\wt{x})\!\big),$$
is well-defined.
Since
$$\wt\vph_{\al;x;\C}\!\circ\!\wt\Psi_{\al;x;1}\!=\!\wt\phi_{\al;x;1}\!:
\wt{U}_{\al;x;1}\lra  S^1\!\times\!\C^2\!\times\!\R^m,$$
$\wt\Psi_{\al;x;1}$ restricts to a diffeomorphism
from a neighborhood of $\pi'^{-1}(U_{\al;x})\!-\!\wt{U}_{\al;x;2}$
onto an open subset of~$\BLR_{Y;\C}\vph_{\al;x}$.
The~map
$$\wt\Psi_{\al;x;2}\!:\wt{U}_{\al;x;2}\lra\BLR_{Y;\R}\vph_{\al;x}
\!\subset\!\BLR_Y\vph_{\al;x}, \quad
\wt\Psi_{\al;x;2}(\wt{x})=\big(\R\wt\psi_{\al;x;2}(\wt{x}),\pi'(\wt{x})\!\big),$$
is similarly a diffeomorphism onto an open subset of~$\BLR_{Y;\R}\vph_{\al;x}$,
since
$$\wt\vph_{\al;x;\R}\!\circ\!\wt\Psi_{\al;x;2}\!=\!\wt\phi_{\al;x;2}\!:
\wt{U}_{\al;x;2}\lra  \C\!\times\!\R\!\times\!\R^m\,.$$
By the same reasoning as in the proof of Lemma~\ref{StanBl_lmm}, the~map
$$\wt\Psi\!:\wt{X}\lra\BLR_YX, \qquad
\wt\Psi(\wt{x})=\begin{cases}[\wt\Psi_{\al;x;i}(\wt{x})],
&\hbox{if}~\wt{x}\!\in\!\wt{U}_{\al;x;i},\,\al\!\in\!\cI,\,x\!\in\!U_{\al},\,i\!\in\![2];\\
[\pi'(\wt{x})],&\hbox{if}~\wt{x}\!\in\!\wt{X}\!-\!\pi'^{-1}(Y);
\end{cases}$$
is a well-defined diffeomorphism onto an open subset of~$\BLR_YX$
which satisfies \hbox{$\pi'\!=\!\pi\!\circ\!\wt\Psi$}. 
\end{proof}

\section{The complex quotient sequence}
\label{Cblowup_sec}

\noindent
We identify the moduli space $\ov\cM_4$ with $\C\P^1$ via the cross ratio
\BE{CRdfn_e}\CR\!:\ov\cM_4\lra\C\P^1, \qquad
\CR(\cC)=\frac{z_1(\cC)\!-\!z_3(\cC)}{z_1(\cC)\!-\!z_4(\cC)}\!:
\!\frac{z_2(\cC)\!-\!z_3(\cC)}{z_2(\cC)\!-\!z_4(\cC)}\,,\EE
where $z_i(\cC)$ is the $i$-th marked point of the marked curve $\cC$.
For $\ell\!\in\!\Z^+$ and a quadruple 
\hbox{$\fq\!\equiv\!(i,j,k,m)$} of distinct elements of~$[\ell]$, let 
\BE{CRftdfn_e}\CR_{\fq}\!=\!\CR_{ijkm}: \ov\cM_{\ell}\xlra{\ff_{ijkm}}\ov\cM_4\xlra{~\CR~}\C\P^1\EE
be the composition of~$\CR$ with the forgetful morphism~$\ff_{ijkm}$
sending the marked points $z_i,z_j,z_k,z_m$ to the marked points indexed by $1,2,3,4$, respectively.
We note~that
\BE{CRprp_e}\begin{split} 
&\CR_{kmij}=\CR_{ijkm}, \quad \CR_{jikm},\CR_{ijmk}=1/\CR_{ijkm}, \quad
\CR_{mjki},\CR_{ikjm}=1-\CR_{ijkm},\\
&\hspace{.8in}\CR_{kjim},\CR_{imkj}=-\frac{\CR_{ijkm}}{1\!-\!\CR_{ijkm}}, \quad
\CR_{ijkn}=\CR_{ijkm}\CR_{ijmn}
\end{split}\EE
for all $i,j,k,m,n\!\in\![\ell]$ distinct.
Denote the set of quadruples~$\fq$ as above by~$\cQ_{\ell}$.\\

\noindent
We call a function $f\!:W\!\lra\!\C\P^1$ from a subset $W$ of $(\C\P^1)^N$ \sf{rational}
if there are polynomials
$$g,h\!:(\C^2)^N\lra \C$$
with {\it integer} coefficients
that are homogeneous of the same degree in the two components of each $\C^2$-factor
(the degrees could be different for different $\C^2$-factors) and 
$$f\big(\!\big([w_i]\big)_{i\in[N]}\big)=
\big[g\big(\!(w_i)_{i\in[N]}\big),h\big(\!(w_i)_{i\in[N]}\big)\big]
\qquad\forall~\big([w_i]\big)_{i\in[N]}\!\in\!W.$$
If $W$ is an open subset (resp.~a submanifold), then a rational function 
from $W$ is holomorphic (resp.~smooth).

\subsection{The combinatorial structure of $\ov\cM_{\ell}$}
\label{cMellstrC_subs0}

\noindent
For $\ell\!\in\!\Z^+$ with $\ell\!\ge\!3$,
the topological type of an element of $\ov\cM_{\ell}$ is described
by a \sf{trivalent $\ell$-marked tree}.
A \sf{graph} consists of a finite set~$\Ver$ of \sf{vertices} and 
a collection $\Edg$ of \sf{edges}, i.e.~of two-element subsets of~$\Ver$.
A graph $(\Ver,\Edg)$ is a \sf{tree} if it is connected and contains no~\sf{loops}, 
i.e.~for every pair of distinct vertices $v,v'\!\in\!\Ver$ 
there exist a unique \sf{path of vertices from~$v$ to~$v'$},
$$v_0\!=\!v',~v_1,~\ldots,~v_{n-1},~v_n\!=\!v'$$
with $v_1,\ldots,v_{n-1}\!\in\!\Ver$ and 
$\{v_0,v_1\},\{v_1,v_2\},\ldots,\{v_{n-1},v_n\}\!\in\!\Edg$.
An \sf{$\ell$-marked tree} is a tuple \hbox{$\Ga\!\equiv\!(\Ver,\Edg,\mu)$} so that
$(\Ver,\Edg)$ is a tree and $\mu\!:[\ell]\!\lra\!\Ver$ is a map.
For such a tuple~$\Ga$ and $v\!\in\!\Ver$, the \sf{valence} of~$v$ in~$\Ga$
is the number
$$\val_{\Ga}(v)=\big|\mu^{-1}(v)\!\sqcup\!\{e\!\in\!\Edg\!:e\!\ni\!v\}\big|.$$
An $\ell$-marked tree $\Ga\!\equiv\!(\Ver,\Edg,\mu)$ is \sf{trivalent}
if $\val_{\Ga}(v)\!\ge\!3$ for every $v\!\in\!\Ver$.
The center diagram of Figure~\ref{markedtree_fig} shows a trivalent 11-marked tree,
representing each vertex by a~dot and each edge by the line segment 
between the corresponding vertices.
The elements of the set~$[11]$ are linked by short line segments to their images under
the map~$\mu$.\\

\noindent
For an $\ell$-marked tree $\Ga\!\equiv\!(\Ver,\Edg,\mu)$, let
$$\wt\Edg\equiv\big\{vv'\!\equiv\!(v,v')\!:\{v,v'\}\!\in\!\Edg\big\}$$
be the set of \sf{oriented edges}.
For each $v\!\in\!\Ver$, let
$$\nE_{\Ga}(v)\equiv\big\{\!\{v,v'\}\!\in\!\Edg\} \qquad\hbox{and}\qquad
\wt\nE_{\Ga}(v)\equiv\big\{vv'\!\in\!\wt\Edg\}$$
be the set of edges containing~$v$ and 
the set of oriented edges leaving~$v$, respectively. 
The removal of an edge $e\!\equiv\!\{v,v'\}$ from~$\Ga$ breaks this tree into
two subtrees.
For an oriented edge $\wt{e}\!\equiv\!vv'$, we denote by 
$(\Ver_{\wt{e}},\Edg_{\wt{e}})$ and $(\Ver_{\wt{e}}^c,\Edg_{\wt{e}}^c)$
the two resulting subtrees so that $\Ver_{\wt{e}}\!\ni\!v$ and $\Ver_{\wt{e}}^c\!\ni\!v'$;
see the center and right diagrams in Figure~\ref{markedtree_fig} for an example.\\

\noindent
A stable $\ell$-marked rational curve~$\cC$ determines
a trivalent $\ell$-marked tree $\Ga\!\equiv\!(\Ver,\Edg,\mu)$
so~that
\begin{enumerate}[label=$\bullet$,leftmargin=*]

\item $\Ver$ is the set of the irreducible components of~$\cC$,

\item $\Edg$ is the set of pairs $\{v,v'\}\!\subset\!\Ver$
so that the irreducible components~$\C\P^1_v$ and~$\C\P^1_{v'}$
corresponding to~$v$ and~$v'$ are joined by a node, and

\item
$\mu\!:[\ell]\!\lra\!\Ver$ is the map so that the $i$-th marked point~$z_i(\cC)$ of~$\cC$
lies on~$\C\P^1_{\mu(i)}$.

\end{enumerate}
We call the trivalent $\ell$-marked tree $\Ga$ obtained in this way
the \sf{dual graph} of~$\cC$.
As an example, the dual graph of the stable 11-marked rational curve in the left diagram 
of Figure~\ref{markedtree_fig} is represented by the center diagram.\\

\begin{figure}
\begin{pspicture}(-1,-3.3)(10,3.6)
\psset{unit=.4cm}
\pscircle(3,2){1.5}\rput(3.1,2){\smsize{$\C\P^1_{v'}$}}
\pscircle(1,4){1.33}\rput(1.1,4){\smsize{$\C\P^1_a$}}\pscircle*(1.94,3.06){.2}
\pscircle*(-.33,4){.2}\pscircle*(1,5.33){.2}\rput(-.9,4){\smsize{$4$}}\rput(1,6){\smsize{$5$}}
\pscircle(5,4){1.33}\rput(5.1,4){\smsize{$\C\P^1_b$}}\pscircle*(4.06,3.06){.2}
\pscircle*(5.94,3.06){.2}\rput(6.5,2.9){\smsize{$6$}}
\pscircle(3.2,5.8){1.22}\rput(3.2,5.8){\smsize{$\C\P^1_c$}}\pscircle*(4.06,4.94){.2}
\pscircle*(2.34,6.66){.2}\pscircle*(4.06,6.66){.2}
\rput(2.3,7.3){\smsize{$7$}}\rput(4.1,7.3){\smsize{$8$}}
\pscircle(6.8,5.8){1.22}\rput(6.8,5.8){\smsize{$\C\P^1_d$}}\pscircle*(5.94,4.94){.2}
\pscircle*(5.94,6.66){.2}\pscircle*(7.66,6.66){.2}\pscircle*(7.66,4.94){.2}
\rput(5.9,7.3){\smsize{$9$}}\rput(7.77,7.3){\smsize{$10$}}\rput(8.1,4.4){\smsize{$11$}}
\pscircle(3,-1){1.5}\rput(3.1,-1){\smsize{$\C\P^1_v$}}\pscircle*(3,.5){.2}
\pscircle*(1.5,-1){.2}\pscircle*(3,-2.5){.2}\pscircle*(4.5,-1){.2}
\rput(1,-1){\smsize{$1$}}\rput(3,-3.2){\smsize{$2$}}\rput(5,-1){\smsize{$3$}}
\rput(4,-5){$\cC$}
\pscircle*(15,2){.2}\rput(15.2,2.8){\smsize{$v'$}}
\psline[linewidth=.04](12,4)(15,2)\pscircle*(12,4){.2}\rput(12,4.7){\smsize{$a$}}
\psline[linewidth=.02](12,4)(10.5,5.5)\rput(10.5,6){\smsize{$4$}}
\psline[linewidth=.02](12,4)(13.5,5.5)\rput(13.5,6){\smsize{$5$}}
\psline[linewidth=.04](18,4)(15,2)\pscircle*(18,4){.2}\rput(18,4.7){\smsize{$b$}}
\psline[linewidth=.02](18,4)(19.5,2.5)\rput(20,2.5){\smsize{$6$}}
\psline[linewidth=.04](18,4)(15,6)\pscircle*(15,6){.2}\rput(15,6.7){\smsize{$c$}}
\psline[linewidth=.02](15,6)(13.5,7.5)\rput(13.5,8){\smsize{$7$}}
\psline[linewidth=.02](15,6)(16.5,7.5)\rput(16.5,8){\smsize{$8$}}
\psline[linewidth=.04](18,4)(21,6)\pscircle*(21,6){.2}\rput(21.5,5.9){\smsize{$d$}}
\psline[linewidth=.02](21,6)(19.5,7.5)\rput(19.5,8){\smsize{$9$}}
\psline[linewidth=.02](21,6)(21,8)\rput(21,8.5){\smsize{$10$}}
\psline[linewidth=.02](21,6)(22.5,7.5)\rput(22.5,8){\smsize{$11$}}
\psline[linewidth=.04,arrowsize=.5]{<-}(15,1)(15,-1)
\psline[linewidth=.04](15,2)(15,-1)\rput(16.8,.5){\smsize{$\wt{e}\!=\!vv'$}}
\pscircle*(15,-1){.2}\rput(14.5,-.7){\smsize{$v$}}
\psline[linewidth=.02](15,-1)(13.5,-2.5)\rput(13.5,-3){\smsize{$1$}}
\psline[linewidth=.02](15,-1)(15,-3)\rput(15,-3.5){\smsize{$2$}}
\psline[linewidth=.02](15,-1)(16.5,-2.5)\rput(16.5,-3){\smsize{$3$}}
\rput(16,-5){$\Ga\!=\!(\Ver,\Edg,\mu)$}
\pscircle*(28,2){.2}\rput(28.2,2.8){\smsize{$v'$}}
\psline[linewidth=.04](25,4)(28,2)\pscircle*(25,4){.2}\rput(25,4.7){\smsize{$a$}}
\psline[linewidth=.02](25,4)(23.5,5.5)\rput(23.5,6){\smsize{$4$}}
\psline[linewidth=.02](25,4)(26.5,5.5)\rput(26.5,6){\smsize{$5$}}
\psline[linewidth=.04](31,4)(28,2)\pscircle*(31,4){.2}\rput(31,4.7){\smsize{$b$}}
\psline[linewidth=.02](31,4)(32.5,2.5)\rput(33,2.5){\smsize{$6$}}
\psline[linewidth=.04](31,4)(28,6)\pscircle*(28,6){.2}\rput(28,6.7){\smsize{$c$}}
\psline[linewidth=.02](28,6)(26.5,7.5)\rput(26.5,8){\smsize{$7$}}
\psline[linewidth=.02](28,6)(29.5,7.5)\rput(29.5,8){\smsize{$8$}}
\psline[linewidth=.04](31,4)(34,6)\pscircle*(34,6){.2}\rput(34.5,5.9){\smsize{$d$}}
\psline[linewidth=.02](34,6)(32.5,7.5)\rput(32.5,8){\smsize{$9$}}
\psline[linewidth=.02](34,6)(34,8)\rput(34,8.5){\smsize{$10$}}
\psline[linewidth=.02](34,6)(35.5,7.5)\rput(35.5,8){\smsize{$11$}}
\pscircle*(33,-1){.2}\rput(32.5,-.7){\smsize{$v$}}
\psline[linewidth=.02](33,-1)(31.5,-2.5)\rput(31.5,-3){\smsize{$1$}}
\psline[linewidth=.02](33,-1)(33,-3)\rput(33,-3.5){\smsize{$2$}}
\psline[linewidth=.02](33,-1)(34.5,-2.5)\rput(34.5,-3){\smsize{$3$}}
\rput(28,1){$(\Ver_{\wt{e}}^c,\Edg_{\wt{e}}^c)$}
\rput(33,-5){$(\Ver_{\wt{e}},\Edg_{\wt{e}})$}
\rput(20,-7){\smsize{$[\Ga]_{\eta;v}\!=\!\{1,2,3,6\},~~[\Ga]_{\eta;v'}\!=\!\{2,4,6\},~~
[\Ga]_{\eta;a}\!=\!\{2,4,5\}$}}
\rput(20,-8.2){\smsize{$[\Ga]_{\eta;b}\!=\!\{4,6,8,9\},~~
[\Ga]_{\eta;c}\!=\!\{6,7,8\},~~[\Ga]_{\eta;d}\!=\!\{8,9,10,11\}$}}
\end{pspicture}
\caption{A stable 11-marked rational curve~$\cC$, its dual 11-marked tree~$\Ga$,
the subtrees $(\Ver_{\wt{e}},\Edg_{\wt{e}})$ and $(\Ver_{\wt{e}}^c,\Edg_{\wt{e}}^c)$
of~$\Ga$ determined by the oriented edge~$\wt{e}$ of~$\Ga$,
and the subsets~\eref{ellGavdfn_e} of~[11]
determined by a systematic marking map~$\eta$.}
\label{markedtree_fig}
\end{figure}

\noindent
For a trivalent $\ell$-marked tree $\Ga\!\!\equiv\!(\Ver,\Edg,\mu)$, 
denote~by \hbox{$\cM_{\Ga}\!\subset\!\ov\cM_{\ell}$} 
the stratum of the elements~$[\cC]$ with dual graph~$\Ga$.
Let $\cT(\Ga)$ be the collection of (trivalent) $\ell$-marked  trees 
\hbox{$\Ga'\!\equiv\!(\Ver',\Edg',\mu')$} so that there exists a surjection 
$\ka\!:\Ver\!\lra\!\Ver'$ such~that
\BE{Gakacond_e}\begin{split}
&\quad\Edg'=\big\{\!\{\ka(v),\ka(v')\}\!:\!\{v,v'\}\!\in\!\Edg,\,\ka(v)\!\neq\!\ka(v')\!\big\},
\quad \mu'\!=\!\ka\!\circ\!\mu\!:[\ell]\lra\Ver',\\
&\{v_1,v_1'\},\{v_2,v_2'\}\in\Edg,~\ka(v_1)\!=\!\ka(v_2)\!\neq\!\ka(v_1')=\!\ka(v_2')
\quad\Lra\quad v_1\!=\!v_2,~v_1'\!=\!v_2';
\end{split}\EE
the last condition implies that the fibers of~$\ka$ are subtrees of~$(\Ver,\Edg)$.
By the first and last conditions in~\eref{Gakacond_e}, the~map
$$\big\{\!\{v,v'\}\!\in\!\Edg\!:\ka(v)\!\neq\!\ka(v')\!\big\}\lra\Edg',\quad
\{v,v'\}\lra\big\{\!\ka(v),\ka(v')\!\big\},$$
is a bijection.
The trees $\Ga'\!\in\!\cT(\Ga)$ are obtained by collapsing some of the edges of~$\Ga$
and combining the endpoints of each such edge into a single vertex.  
The elements of~$\cM_{\Ga'}$ with $\Ga'\!\in\!\cT(\Ga)$ are obtained from
the elements of~$\cM_{\Ga}$ by smoothing some of the nodes of the latter.
In particular,
\BE{cWGadfn_e}\cW_{\Ga}\equiv \bigcup_{\Ga'\in\cT(\Ga)}\!\!\!\!\!\!\cM_{\Ga'}\subset\ov\cM_{\ell}\EE 
is an open neighborhood of~$\cM_{\Ga}$ in~$\ov\cM_{\ell}$.
For $\vr\!\subset\![\ell]$, 
\BE{cMGaDvr_e}\cW_{\Ga}\!\cap\!D_{\ell;\vr}\neq\eset \qquad\Llra\qquad
\mu^{-1}\big(\Ver_{\wt{e}_{\vr}}\big)=\vr
~~\hbox{for some}~\wt{e}_{\vr}\!\in\!\wt\Edg.\EE
In such a case, $\cM_{\Ga}\!\subset\!D_{\ell;\vr}$.
We denote by $e_{\vr}\!\in\!\Edg$ the (unordered) edge determined by~$\wt{e}_{\vr}$.\\

\noindent
With $\Ga$ as above and $v_+\!\in\!\Ver$, define
\BE{wtGav_e}\Ga v_+=(\Ver,\Edg,\mu_{v_+}), \quad
\mu_{v_+}\!\!:[\ell\!+\!1]\lra\Ver, ~~
\mu_{v_+}(i)=\begin{cases}\mu(i),&\hbox{if}~i\!\in\![\ell];\\
v_+,&\hbox{if}~i\!=\!\ell\!+\!1.
\end{cases}\EE
Thus, $\Ga v_+$ is a trivalent $(\ell\!+\!1)$-marked tree.
The corresponding stratum \hbox{$\cM_{\Ga v_+}\!\subset\!\ov\cM_{\ell+1}$}
consists of the elements $\wt\cC\!\in\!\ov\cM_{\ell+1}$ with the marked point~$z_{\ell+1}(\wt\cC)$
on the irreducible component~$\C\P^1_{v_+}$ of an $\ell$-marked curve
$\cC\!\in\!\cM_{\Ga}$.\\

\noindent
For $e_+\!\equiv\!\{v_+,v_+'\}\in\!\Edg$, define
\BE{wtGae_e}\begin{split}
&\Ver_{e_+}=\Ver\!\sqcup\!\{v_{\bu}\}, \quad
\Edg_{e_+}=\big(\Edg\!-\!\{e_+\!\}\!\big)\!\sqcup\!
\big\{\!\{v_+,v_{\bu}\},\{v_{\bu},v_+'\}\!\big\},\\
&\hspace{.5in}\mu_{e_+}\!:[\ell\!+\!1]\lra\Ver_{\bu}, \quad
\mu_{e_+}(i)=\begin{cases}\mu(i),&\hbox{if}~i\!\in\![\ell];\\
v_{\bu},&\hbox{if}~i\!=\!\ell\!+\!1.
\end{cases}\end{split}\EE
Thus, $\Ga e_+\!\equiv\!(\Ver_{e_+},\Edg_{e_+},\mu_{e_+})$ 
is a  trivalent $(\ell\!+\!1)$-marked tree.
If $e_+\!=\!e_{\vr}$ with $\vr\!\in\!\cA_{\ell}$,
the corresponding stratum 
$$\cM_{\Ga e_+}=\ff_{\ell+1}^{~-1}(\cM_{\Ga})
\!\cap\!D_{\ell+1;\vr}\!\cap\!D_{\ell+1;\vr_{\ell}^c} 
\subset \ov\cM_{\ell+1}$$
is the section of~$\ff_{\ell+1}$ over $\cM_{\Ga}\!\subset\!\ov\cM_{\ell}$
with the $(\ell\!+\!1)$-th marked point ``at" the node~$e_+$ 
of the elements of~$\cM_{\Ga}$.\\

\noindent
For any trivalent $\ell$-marked tree $\Ga\!\!\equiv\!(\Ver,\Edg,\mu)$,  
$$\ff_{\ell+1}^{\,-1}(\cM_{\Ga})=
\bigsqcup_{v_+\in\Ver}\!\!\!\!\!\cM_{\Ga v_+}\sqcup
\bigsqcup_{e_+\in\Edg}\!\!\!\!\!\cM_{\Ga e_+}\,.$$

\subsection{The complex structure of $\ov\cM_{\ell}$}
\label{cMellstrC_subs}

\noindent
We call $\eta\!:\wt\Edg\!\lra\![\ell]$ a \sf{marking map} 
for an $\ell$-marked tree $\Ga\!\equiv\!(\Ver,\Edg,\mu)$
if \hbox{$\mu(\eta(\wt{e})\!)\!\in\!\Ver_{\wt{e}}^c$} for every $\wt{e}\!\in\!\wt\Edg$.
For each $v\!\in\!\Ver$, the restriction of such a map~$\eta$ to~$\wt\nE_{\Ga}(v)$
is injective into~$[\ell]\!-\!\mu^{-1}(v)$.
Thus, the cardinality of the subset
\BE{ellGavdfn_e} 
[\Ga]_{\eta;v}\equiv\mu^{-1}(v)\!\sqcup\!\eta\big(\wt\nE_{\Ga}(v)\!\big)\subset [\ell]\EE
is $\val_{\Ga}(v)$.
We call a marking map~$\eta$ as above \sf{systematic} if 
$$\eta(vv')\in  [\Ga]_{\eta;v'} \quad\forall~vv'\!\in\!\wt\Edg.$$
For example, the marking map given~by
$$\eta\!:\wt\Edg\lra[\ell], \qquad \eta(\wt{e})=\min\mu^{-1}(\Ver_{\wt{e}}^c),$$
is systematic.
Figure~\ref{markedtree_fig} provides another example of a systematic marking map~$\eta$ 
by listing the subsets $[\Ga]_{\eta;v}\!\subset\![11]$, which encode
the map~$\eta$.\\

\noindent
Let $\Ga\!\equiv\!(\Ver,\Edg,\mu)$ be a trivalent $\ell$-marked tree.
If $v\!\in\!\Ver$, we call two distinct elements \hbox{$i,j\!\in[\ell]$}
\sf{$(\Ga,v)$-independent} if either 
$$\{i,j\}\!\cap\!\mu^{-1}(v)\neq\eset \quad\hbox{or}\quad
\exists~\wt{e}_i,\wt{e}_j\in\wt\nE_{\Ga}(v)~\hbox{s.t.}~\wt{e}_i\!\neq\!\wt{e}_j,\,
i\!\in\!\mu^{-1}\big(\Ver_{\wt{e}_i}^c\big),\,j\!\in\!\mu^{-1}\big(\Ver_{\wt{e}_j}^c\big).$$
If $\eta\!:\wt\Edg\!\lra\![\ell]$ is a marking map for~$\Ga$,
then $[\Ga]_{\eta;v}\!\subset\![\ell]$ is a maximal subset of (pairwise)
$(\Ga,v)$-independent elements.
If $\Ga'\!\in\!\cT(\Ga)$, $\ka$ is as in~\eref{Gakacond_e},
and $i,j\!\in\![\ell]$ are $(\Ga,v)$-independent,
then $i,j$ are also $(\Ga',\ka(v)\!)$-independent.
For $i,j\!\in\![\ell]$ distinct, the collection of vertices $v\!\in\!\Ver$
such that $i,j$ are $(\Ga,v)$-independent forms the unique path in~$\Ga$ 
from~$\mu(i)$ to~$\mu(j)$.
For \hbox{$i,j,k\!\in\![\ell]$} distinct, there is a unique vertex $v_{ijk}\!\in\!\Ver$
so that $i,j,k$ are pairwise $(\Ga,v_{ijk})$-independent.\\

\noindent
For each $v\!\in\!\Ver$, choose a maximal subset $[\Ga]_v\!\subset\![\ell]$
of $(\Ga,v)$-independent elements and an ordering \hbox{$i_1(v),\ldots,i_{\val_{\Ga}(v)}(v)$}
of its elements and let
\BE{cQetavGadfn_e}
\cQ_{\Ga}(v)=\big\{\!(i_1(v),i_2(v),i_2(v),i_r(v)\!)\!:r\!=\!4,\ldots,\val_{\Ga}(v)\!\big\}
\subset \cQ_{\ell}\,.\EE
For each $e\!\equiv\!\{v,v'\}\!\in\!\Edg$, choose 
an oriented edge representative $\wt{e}\!\equiv\!vv'$ and
\BE{fqedfn_e}\fq_e\equiv \big(i_e,j_e,k_e,m_e\big)\in \cQ_{\ell}\EE
so that $i_e,k_e\!\in\!\mu^{-1}(\Ver_{vv'})$ are $(\Ga,v)$-independent
and $j_e,m_e\!\in\!\mu^{-1}(\Ver_{v'v})$ are $(\Ga,v')$-independent.
Let
\BE{cQetaGadfn_e}
\cQ_{\Ga}=\bigcup_{v\in\Ver}\!\!\!\cQ_{\Ga}(v)\cup\big\{\fq_e\!:e\!\in\!\Edg\big\}
\subset\cQ_{\ell}.\EE
We call such a subset $\cQ_{\Ga}\!\subset\!\cQ_{\ell}$ a \sf{$\Ga$-basis}
for~$\cQ_{\ell}$; see Corollary~\ref{cMspan_crl}.
We call such a subset~$\cQ_{\Ga}$ for~$\cQ_{\ell}$ \sf{compatible} 
with a systematic marking map~$\eta$ for~$\Ga$ if 
\BE{MarkSyst_e}\begin{aligned}
\hspace{1in} [\Ga]_v&=[\Ga]_{\eta;v}&~~
&\forall\,v\!\in\!\Ver \qquad\hbox{and}\\
i_e=\eta(v'v),\,j_e=\eta(vv'),&\,\,k_e\!\in\![\Ga]_{\eta;v},\,
m_e\!\in\![\Ga]_{\eta;v'}
&~~
&\forall\,e\!\equiv\!\{v,v'\}\!\in\!\Edg.
\end{aligned}\EE
This in particular implies that each of $[\Ga]_v,[\Ga]_{v'}$ contains both $\eta(vv'),\eta(v'v)$
whenever \hbox{$\{v,v'\}\!\in\!\Edg$}.\\

\noindent
For $\fq\!\equiv\!(i,j,k,m)\!\in\!\cQ_{\ell}$ and $v\!\in\!\Ver$,
the cardinality of the subset
$$[\Ga]_{\fq;v}\equiv
\big(\{i,j,k,m\}\!\cap\!\mu^{-1}(v)\!\big)\!\cup\!
\big\{\wt{e}\!\in\!\wt\nE_v(\Ga)\!:\{i,j,k,m\}\!\cap\!\Ver_{\wt{e}}^c\!\neq\!\eset\big\}
\subset \mu^{-1}(v)\!\sqcup\!\wt\nE_{\Ga}(v)$$
is the cardinality of a maximal subset of $\{i,j,k,m\}$ of $(\Ga,v)$-independent elements.
For every $\fq\!\in\!\cQ_{\ell}$, there exists $v\!\in\!\Ver$ so that 
$|[\Ga]_{\fq;v}|\!\ge\!3$.
If there exists $v\!\in\!\Ver$ such that $|[\Ga]_{\fq;v}|\!=\!4$,
then $|[\Ga]_{\fq;v'}|\!\le\!2$ for all $v'\!\in\!\Ver\!-\!\{v\}$.
Otherwise, there exist exactly two distinct elements $v,v'\!\in\!\Ver$ 
such that \hbox{$|[\Ga]_{\fq;v}|,|[\Ga]_{\fq;v'}|\!=\!3$}.
For example,
\BE{qtypes_e}
\big|[\Ga]_{\fq;v}\big|=4~~\forall~\fq\!\in\!\cQ_{\Ga}(v),\,v\!\in\!\Ver, 
\qquad
\big|[\Ga]_{\fq_e;v}\big|,\big|[\Ga]_{\fq_e;v'}\big|=3~~
\forall~e\!\equiv\!\{v,v'\}\!\in\!\Edg.\EE
Thus, all elements $\fq_e\!\in\!\cQ_{\ell}$ with $e\!\in\!\Edg$ are pairwise distinct,
all subsets $\cQ_{\Ga}(v)\!\subset\!\cQ_{\ell}$ with $v\!\in\!\Ver$ are pairwise disjoint,
and none of the latter contains any of the former.
It follows~that
\BE{cMdimcount_e} 
\big|\cQ_{\Ga}\big|
=\sum_{v\in\Ver}\!\!\!\big(\val_{\Ga}(v)\!-\!3)\!+\!|\Edg|
=\ell\!-\!3=\dim_{\C}\ov\cM_{\ell}\,.\EE
Lemma~\ref{cMspan_lmm} below is essentially the proof of \cite[Theorem~D.4.2(iii)]{MS},
establishing that the subspace
\hbox{$\ov{M}_{\ell}\!\subset\!(\C\P^1)^{\cQ_{\ell}}$}
in the proof of Corollary~\ref{cMspan_crl}
is a complex submanifold. 

\begin{lmm}\label{cMspan_lmm} Suppose $\ell\!\in\!\Z^+$ with $\ell\!\ge\!3$,
\hbox{$\Ga\!\equiv\!(\Ver,\Edg,\mu)$} is a trivalent $\ell$-marked tree,
\hbox{$\eta\!:\wt\Edg\!\lra\![\ell]$} is a systematic marking map for~$\Ga$,
and $\cQ_{\Ga}\!\subset\!\cQ_{\ell}$ as in~\eref{cQetaGadfn_e} 
is a $\Ga$-basis compatible with~$\eta$.
Let $\cW_{\Ga}\!\subset\!\ov\cM_{\ell}$ be as in~\eref{cWGadfn_e}.
For every $\fq\!\in\!\cQ_{\ell}$, the cross ratio~\eref{CRftdfn_e} is a rational function
of the cross ratios~$\CR_{\fq'}$ with $\fq'\!\in\!\cQ_{\Ga}$
on the image~$W_{\Ga}$ of~$\cW_{\Ga}$ under $(\CR_{\fq'})_{\fq'\in\cQ_{\Ga}}$.
\end{lmm}

\begin{proof} Let $[\Ga]_v\!\subset\![\ell]$ and
$\cQ_{\Ga}(v)\!\subset\!\cQ_{\ell}$ with $v\!\in\!\Ver$ and
$\fq_e\!\in\!\cQ_{\ell}$ with $e\!\in\!\Edg$ be as in~\eref{cQetavGadfn_e} 
and~\eref{fqedfn_e}, subject to the conditions~\eref{MarkSyst_e}.
For a subtree $(\Ver',\Edg')$ of $(\Ver,\Edg)$, let
\BE{cMspan_e1a}[\Ga]_{\Ver'}=\bigcup_{v\in\Ver'}\!\!\![\Ga]_v
\subset[\ell].\EE 
In particular, \hbox{$[\Ga]_{\Ver}\!=\![\ell]$}.
By~\eref{qtypes_e},
\BE{cMspan_e1}
\CR_{\fq}(\cW_{\Ga})\subset\C^*\!-\!\{1\}~~\forall~\fq\!\in\!\cQ_{\Ga}(v),\,v\!\in\!\Ver,
\quad
\CR_{\fq_e}(\cW_{\Ga})\subset\C~~\forall~e\!\in\!\Edg.\EE
We show that the subcollection $\cQ_{\ell}^{\star}\!\subset\!\cQ_{\ell}$
for which the claim of the lemma holds is all of~$\cQ_{\ell}$.
By the first line in~\eref{CRprp_e}, 
$\cQ_{\ell}^{\star}\!\subset\![\ell]^4$ is preserved by the permutations of 
the components of its elements.
By the last relation in~\eref{CRprp_e} and the first equation in~\eref{cMspan_e1},
\BE{cMspan_e2}
\cQ_{\ell}^{\star}\supset\cQ_{\ell}\!\cap\![\Ga]_v^{\,4}
\quad\hbox{and}\quad
\CR_{\fq}(\cW_{\Ga})\in\C^*\!-\!\{1\}~~
\forall~\fq\!\in\!\cQ_{\ell}\!\cap\![\Ga]_v^{\,4}\EE
for every $v\!\in\!\Ver$.\\

\noindent
Suppose $(\Ver',\Edg')$ is proper subtree of $(\Ver,\Edg)$
so that 
$\cQ_{\ell}^{\star}$ contains \hbox{$\cQ_{\ell}\!\cap\![\Ga]_{\Ver'}^{\,4}$}.
Choose \hbox{$v\!\in\!\Ver\!-\!\Ver'$} so that $e\!\equiv\!\{v,v'\}\!\in\!\Edg$
for some $v'\!\in\!\Ver'$.
We can assume that \hbox{$i_e,k_e\!\in\!\mu^{-1}(\Ver_{vv'})$} and 
$k_e,m_e\!\in\!\mu^{-1}(\Ver_{v'v})$.
For \hbox{$k\!\in\![\Ga]_v$} and $m\!\in\![\Ga]_{\Ver'}$ distinct from~$i_e,j_e$,
the last equation in~\eref{CRprp_e} gives
\BE{cMspan_e5}
\CR_{i_ej_ekm}
=\CR_{i_ej_em_em}\CR_{i_ej_ekm_e}\\
=\CR_{i_ej_em_em}\CR_{i_ej_ekk_e}
\CR_{\fq_e},\EE
with $\CR_{i_ej_em_em}\!\equiv\!1$ if $m_e\!=\!m$
and $\CR_{i_ej_ekk_e}\!\equiv\!1$ if $k\!=\!k_e$.
Since $[\Ga]_v$ and $[\Ga]_{\Ver'}$ contain both $i_e,j_e$,
$\cQ_{\ell}^{\star}$ contains $(i_e,j_e,k,k_e)$ if $k\!\neq\!k_e$
(by~\eref{cMspan_e2}) and $(i_e,j_e,m_e,m)$ if $m_e\!\neq\!m$  
(by the assumption on~$\Ver'$ above).
The middle factor on the right side of~\eref{cMspan_e5} lies in~$\C^*\!-\!1$
if $k\!\neq\!k_e$ (when evaluated on~$\cW_{\Ga}$).
Since \hbox{$j_e,m_e\!\in\![\Ga]_{v'}$} are distinct,
\hbox{$i_e\!\in\!\mu^{-1}(\Ver_{vv'})$}, and 
\hbox{$j_e,m_e,m\!\in\!\mu^{-1}(\Ver_{v'v})$},
the first factor on the right side of~\eref{cMspan_e5} lies in~$\C$.
Since $\cQ_{\ell}^{\star}$ also contains~$q_e$, 
\eref{cMspan_e5} implies that $\cQ_{\ell}^{\star}$ contains 
$(i_e,j_e,k,m)$ as~well.\\

\noindent
For $k\!\in\![\Ga]_v$ and $j,m\!\in\![\Ga]_{\Ver'}$ 
distinct from each other and from~$i_e,j_e$,
we similarly obtain
\BE{cMspan_e7}
\CR_{i_ekjm}
=\CR_{i_ekjj_e}\CR_{i_ekj_em}.\EE
If $i,k\!\in\![\Ga]_v$ are distinct from each other and from~$i_e,j_e$,
and $j,m\!\in\![\Ga]_{\Ver'}$ are distinct from each other and from~$i_e$, 
we obtain
\BE{cMspan_e8}
\CR_{jmik}=\CR_{jmii_e}\CR_{jmi_ek}.\EE
Suppose $j,m,m'\!\in\![\Ga]_{\Ver'}$ are distinct from each other and
from~$i_e$ and there exists no $\wt{e}\!\in\!\wt\Edg'$ so that 
$i_e,j\!\in\!\mu^{-1}(\Ver_{\wt{e}})$ and 
$m,m'\!\in\!\mu^{-1}(\Ver_{\wt{e}}^c)$;
the last condition is satisfied at least after a permutation of~$j,m,m'$.
For $k\!\in\![\Ga]_v$ distinct from $i_e,j_e$, 
the last equation in~\eref{CRprp_e} then gives
\BE{cMspan_e9}
\CR_{jmm'k}=\CR_{jmm'i_e}\CR_{jmi_ek}.\EE
By assumption, $(i_e,j,m,m')\!\in\![\Ga]_{\Ver'}$.
All six factors on the right-hand sides of \eref{cMspan_e7}-\eref{cMspan_e9}
lie in~$\C$.\\

\noindent
By \eref{cMspan_e2}, the inductive assumption, and 
 \eref{cMspan_e7}-\eref{cMspan_e9}, $\cQ_{\ell}^{\star}$ contains every element~$\fq$
of \hbox{$\cQ_{\ell}\!\cap\![\Ga]_{\Ver'\cup\{v\}}^{\,4}$}  
that has 0,2,3, or~4 components in~$\mu^{-1}(\Ver_{v'v})$.
Analogously to~\eref{cMspan_e7} and~\eref{cMspan_e9}, without the concern for the edge condition, 
$\cQ_{\ell}^{\star}$ also contains every element~$\fq$
of \hbox{$\cQ_{\ell}\!\cap\![\Ga]_{\Ver'\cup\{v\}}^{\,4}$} 
that has~3 components in~$\mu^{-1}(\Ver_{vv'})$.
Thus, $\cQ_{\ell}^{\star}$ contains all of 
\hbox{$\cQ_{\ell}\!\cap\![\Ga]_{\Ver'\cup\{v\}}^{\,4}$}.
This implies that \hbox{$\cQ_{\ell}^{\star}\!=\!\cQ_{\ell}$}.
\end{proof}

\begin{crl}\label{cMspan_crl} Suppose $\ell\!\in\!\Z^+$ with $\ell\!\ge\!3$,
\hbox{$\Ga\!\equiv\!(\Ver,\Edg,\mu)$} is a trivalent $\ell$-marked tree,
and $\cQ_{\Ga}\!\subset\!\cQ_{\ell}$ as in~\eref{cQetaGadfn_e} 
is a $\Ga$-basis for~$\cQ_{\ell}$.
\begin{enumerate}[label=(\arabic*),leftmargin=*]

\item\label{cMfqgenC_it} The cross ratios~$\CR_{\fq'}$ with $\fq'\!\in\!\cQ_{\Ga}$
are coordinates on a neighborhood~$\cW_{\Ga}'$ of~$\cM_{\Ga}$ in~$\cW_{\Ga}$.
For every $\fq\!\in\!\cQ_{\ell}$, the cross ratio~\eref{CRftdfn_e} is a rational function
of these cross ratios 
on the image~$W_{\Ga}'$ of~$\cW_{\Ga}'$ under $(\CR_{\fq'})_{\fq'\in\cQ_{\Ga}}$.

\item\label{cMfqvrC_it} For every $\vr\!\subset\![\ell]$ satisfying
either side of~\eref{cMGaDvr_e},
\BE{cMGaDvr_e2} D_{\ell;\vr}\!\cap\!\cW_{\Ga}'
=\big\{[\cC]\!\in\!\cW_{\Ga}'\!:\CR_{\fq_{e}}([\cC])\!=\!0\big\}.\EE

\end{enumerate}
\end{crl}

\begin{proof} We denote by $\ov{M}_{\ell}\!\subset\!(\C\P^1)^{\cQ_{\ell}}$ the subspace of tuples
$(c_{\fq})_{\fq\in\cQ_{\ell}}$ satisfying the relations in~\eref{CRprp_e}
with~$\CR$ replaced by~$c$.
By \cite[Theorem~D.5.2]{MS}, the~map 
\BE{MSemd_e1} (\CR_{\fq})_{\fq\in\cQ_{\ell}}\!:
\ov\cM_{\ell}\lra \ov{M}_{\ell}\EE
is a homeomorphism with respect to Gromov's topology on its domain.
Let 
$$M_{\ell},M_{\Ga}\subset W_{\Ga}\subset \ov{M}_{\ell}\subset\!(\C\P^1)^{\cQ_{\ell}}$$ 
be the images of $\cM_{\ell},\cM_{\Ga},\cW_{\Ga}\!\subset\!\ov\cM_{\ell}$,
respectively, under~\eref{MSemd_e1}.\\

\noindent
If $\Ga$ is a one-vertex graph, i.e.~$\cM_{\Ga}\!=\!\cM_{\ell}$, 
the~map
$$(\CR_{\fq})_{\fq\in\cQ_{\Ga}}\!:
\cM_{\ell}\lra \big(\C^*\!-\!\{1\}\!\big)^{\cQ_{\Ga}}\subset(\C\P^1)^{\cQ_{\Ga}}$$
is surjective.
Thus, so is the projection
$$M_{\ell}\lra  \big(\C^*\!-\!\{1\}\!\big)^{\cQ_{\Ga}},
\quad
\big(c_{\fq}\big)_{\fq\in\cQ_{\ell}}\lra \big(c_{\fq}\big)_{\fq\in\cQ_{\Ga}}\,.$$
Along with the previous paragraph and Lemma~\ref{cMspan_lmm}, 
this implies that the components~$c_{\fq}$ with $\fq\!\in\!\cQ_{\Ga}$
are holomorphic coordinates on the neighborhood~$W_{\Ga}$ of~$M_{\Ga}$ in~$\ov{M}_{\ell}$
if $\cQ_{\Ga}$ is  compatible with
a systematic marking map \hbox{$\eta\!:\wt\Edg\!\lra\![\ell]$} for~$\Ga$. 
Thus, $\ov{M}_{\ell}$ is a complex submanifold of~$(\C\P^1)^{\cQ_{\ell}}$.
Since the map~\eref{MSemd_e1} is holomorphic with respect to the standard complex structure on its domain,
it follows that it is a biholomorphism.
This in turn means that the cross ratios~$\CR_{\fq}$ with $\fq\!\in\!\cQ_{\Ga}$
are coordinates on the neighborhood~$\cW_{\Ga}$ of~$\cM_{\Ga}$ in~$\ov\cM_{\ell}$.\\

\noindent
By the above and the proof of Lemma~\ref{cMspan_lmm}, 
the collection $(c_{\fq})_{\fq\in\cQ_{\Ga}(v),v\in\Ver}$ 
provides coordinates on the submanifold $M_{\Ga}\!\subset\!\ov{M}_{\ell}$
for any choices of maximal subsets $[\Ga]_v\!\subset\![\ell]$
of $(\Ga,v)$-independent elements.
For any $\fq_e\!\in\!\cQ_{\ell}$ as in~\eref{fqedfn_e}, 
\hbox{$\CR_{\fq_e}(\cM_{\Ga})\!=\!\{0\}$}.
For any collection of such elements, the induced differential 
$$\nd(\CR_{\fq_e})_{e\in\Edg}\!:\cN_{\ov\cM_{\ell}}\cM_{\Ga}\lra\C^{\Edg}$$
on the normal bundle of~$\cM_{\Ga}$ in~$\ov\cM_{\ell}$ is an isomorphism;
this follows from Lemma~\ref{cMspan_lmm} and \eref{cMspan_e5}-\eref{cMspan_e8}.
This implies the first claim in~\ref{cMfqgenC_it}.
The second claim in~\ref{cMfqgenC_it} follows from the first and Lemma~\ref{cMspan_lmm}.
The last claim of the corollary is immediate from the definition of~$D_{\ell;\vr}$.
\end{proof}

\subsection{Coordinate charts for~$X_{\vr^*}$}
\label{CCharts_subs}

\noindent
Let $\ell\!\in\!\Z^+$ with $\ell\!\ge\!3$.
For $\vr\!\subset\![\ell]$ and $(i,j,k,m)\!\in\![\ell\!+\!1]^4$, define
$$\vr\!\cap\!(i,j,k,m)=\vr\!\cap\!\{i,j,k,m\}\subset[\ell].$$ 
For $\fq\!\in\!\cQ_{\ell+1}$ and $\vr\!\in\!\cA_{\ell}$, 
the map~$\CR_{\fq}$ from~$\ov\cM_{\ell+1}$ 
is constant on every fiber of the forgetful morphism 
$$\ff_{\ell+1}\!: D_{\ell+1;\vr}\lra D_{\ell;\vr}$$
if and only if 
either $\fq\!\in\!\cQ_{\ell}$ or $|\vr\!\cap\!\fq|\!\ge\!2$.
For each $\vr^*\!\in\!\{0\}\!\sqcup\!\cA_{\ell}$,
this map thus descends to a continuous~map
\BE{CRdescC_e}\begin{aligned}
\CR_{\vr^*;\fq}\!: X_{\vr^*}&\lra\C\P^1 &\quad&\hbox{if}~~ \fq\!\in\!\cQ_{\ell} 
\qquad\hbox{and}\\
\CR_{\vr^*;\fq}\!: X_{\vr^*}-
\bigcup_{\begin{subarray}{c}\vr\in\cA_{\ell}(\vr^*)\\ 
|\vr\cap\fq|\le1\end{subarray}}\!\!\!\!\!\!Y_{\vr^*;\vr}^0
&\lra\C\P^1
&\quad&\hbox{if}~~ \fq\!\in\!\cQ_{\ell+1}\!-\!\cQ_{\ell}\,.
\end{aligned}\EE

\vspace{.15in}

\noindent
For $\vr^*\!\in\!\{0\}\!\sqcup\!\cA_{\ell}$ and
 a trivalent $\ell$-marked tree \hbox{$\Ga\!\equiv\!(\Ver,\Edg,\mu)$}, let
\begin{alignat}{1}
\notag
\cA_{\Ga}(\vr^*)&=\big\{\mu^{-1}(\Ver_{\wt{e}})\!:\wt{e}\!\in\!\wt\Edg,\,
\mu^{-1}(\Ver_{\wt{e}})\!\in\!\cA_{\ell}(\vr^*)\!\big\},\\
\notag
\cA_{\Ga}^{\star}(\vr^*)&=\big\{\vr\!\in\!\cA_{\Ga}(\vr^*)\!:
\vr\!\not\supset\!\vr'~\forall\,\vr'\!\in\!\cA_{\Ga}(\vr^*)\!-\!\{\vr\}\!\big\},\\
\label{wtcQGavrdfn_e}
\wt\cQ_{\Ga}(\vr^*)&=\cQ_{\ell}\!\cup\!\big\{\fq\!\in\!\cQ_{\ell+1}\!-\!\cQ_{\ell}\!:
|\vr\!\cap\!\fq|\!\ge\!2~\forall\,\vr\!\in\!\cA_{\Ga}(\vr^*)\!\big\}\subset \cQ_{\ell+1}.
\end{alignat}
The condition $\vr\!\in\!\cA_{\Ga}(\vr^*)$ in the definition of~$\wt\cQ_{\Ga}(\vr^*)$
can be equivalently replaced by $\vr\!\in\!\cA_{\Ga}^{\star}(\vr^*)$.
If \hbox{$\Ga'\!\in\cT(\Ga)$}, then $\cA_{\Ga'}(\vr^*)\!\subset\!\cA_{\Ga}(\vr^*)$
and thus $\wt\cQ_{\Ga'}(\vr^*)\!\supset\!\wt\cQ_{\Ga}(\vr^*)$.\\

\noindent
By~\eref{ndCprp_e0}, 
$\vr\!\supsetneq\!\vr'^c_{\ell}$ for any $\vr,\vr'\!\in\!\cA_{\Ga}^{\star}(\vr^*)$ distinct and
\BE{CblowupInp_e} 
\nV_{\Ga}(\vr^*)\equiv
\bigcap_{\vr\in\cA_{\Ga}(\vr^*)}\hspace{-.18in}\Ver_{\wt{e}_{\vr}}
=\bigcap_{\vr\in\cA_{\Ga}^{\star}(\vr^*)}\hspace{-.18in}\Ver_{\wt{e}_{\vr}}\neq\eset
\subset\Ver\,.\EE
The above intersection in particular contains every $v_+\!\in\!\Ver$
such that $\wt{e}_{\vr}\!=\!v_+v_0$ for some \hbox{$\vr\!\in\!\cA_{\Ga}^{\star}(\vr^*)$}
and $v_0\!\in\!\Ver$.
Furthermore, $\nV_{\Ga}(\vr^*)\!\subset\!\Ver$ is a subtree of~$\Ga$ and
\BE{CblowupInp_e2}\ka^{-1}(v')\!\cap\!\nV_{\Ga}(\vr^*)\neq\eset
\qquad\forall~v'\!\in\!\nV_{\Ga'}(\vr^*),\,\Ga'\!\in\!\cT(\Ga),\EE
with $\ka\!:\Ver\!\lra\!\Ver'$ as in~\eref{Gakacond_e}.
If $i,j,k\!\in\![\ell]$ are distinct, then 
\hbox{$(i,j,k,\ell\!+\!1)\!\in\!\wt\cQ_{\Ga}(\vr^*)$} 
if and only~if the unique vertex $v_{ijk}\!\in\!\Ver$
such that $i,j,k$ are pairwise $(\Ga,v_{ijk})$-independent
lies in~$\nV_{\Ga}(\vr^*)$.\\

\noindent
For $v\!\in\!\Ver$, let
\BE{wtcWdf_e}\begin{split}
\cA_{\Ga}(v)&=\big\{\mu^{-1}(\Ver_{\wt{e}})\!:\wt{e}\!\in\!\wt\Edg,\,
v\!\in\!\Ver_{\wt{e}}\!\big\},\\
\wt\cW_{\vr^*;\Ga,v}&=q_{\vr^*}\big(\ff_{\ell+1}^{~-1}(\cW_{\Ga})\!\big)
\!-\!\bigcup_{\begin{subarray}{c}\vr\in\cA_{\Ga}(v)\\
\vr\not\in\cA_{\Ga}(\vr^*)\end{subarray}}\hspace{-.18in}Y_{\vr^*;\vr}^0
\subset X_{\vr^*}.
\end{split}\EE 
Since 
$$q_{\vr^*}^{\,-1}\big(\wt\cW_{\vr^*;\Ga,v}\big)=\ff_{\ell+1}^{~-1}(\cW_{\Ga})
\!-\!\bigcup_{\begin{subarray}{c}\vr\in\cA_{\Ga}(v)\\
\vr\not\in\cA_{\Ga}(\vr^*)\end{subarray}}
\hspace{-.18in}D_{\ell+1;\vr}
\!-\!\bigcup_{\begin{subarray}{c}\vr\in\cA_{\Ga}(\vr^*)\\
\vr_{\ell}^c\in\cA_{\Ga}(v)\end{subarray}}
\hspace{-.18in}D_{\ell+1;\vr},$$
by~\eref{ndCprp_e},
the subspace $\wt\cW_{\vr^*;\Ga,v}\!\subset\!X_{\vr^*}$ is open.
If $v\!\in\!\nV_{\Ga}(\vr^*)$, then 
$$\cA_{\Ga}(\vr^*)=\cA_{\ell}(\vr^*)\!\cap\!\cA_{\Ga}(v)$$
and the last union above is empty.\\

\noindent
Let $v_+\!\in\!\nV_{\Ga}(\vr^*)$.
Choose \hbox{$[\Ga]_v\!\subset\![\ell]$} and \hbox{$\cQ_{\Ga}(v)\!\subset\!\cQ_{\ell}$} 
for each $v\!\in\!\Ver$ as in~\eref{cQetavGadfn_e} and $\fq_e\!\in\!\cQ_{\ell}$ 
for each $e\!\in\!\Edg$ as in~\eref{fqedfn_e}.
With $\cQ_{\Ga}\!\subset\!\cQ_{\ell}$ as in~\eref{cQetaGadfn_e},
let \hbox{$\cW_{\Ga}'\!\subset\!\cW_{\Ga}$} be as in Corollary~\ref{cMspan_crl} and
\BE{wtcWdf2_e}\wt\cW_{\vr^*;\Ga,v_+}'=q_{\vr^*}\big(\ff_{\ell+1}^{~-1}(\cW_{\Ga}')\!\big)
\!\cap\!\wt\cW_{\vr^*;\Ga,v_+}
=q_{\vr^*}\big(\ff_{\ell+1}^{~-1}(\cW_{\Ga}')\!\big)
\!-\!\bigcup_{\begin{subarray}{c}\vr\in\cA_{\Ga}(v)\\
\vr\not\in\cA_{\Ga}(\vr^*)\end{subarray}}\hspace{-.18in}Y_{\vr^*;\vr}^0
\subset X_{\vr^*};\EE
the subspace~$\wt\cW_{\vr^*;\Ga,v_+}'\!\subset\!X_{\vr^*}$ is again open.
Choose $i_{v_+},j_{v_+},k_{v_+}\!\in\![\ell]$ $(\Ga,v_+)$-independent and define
\BE{wtfqvrdfn_e}
\wt\fq_{v_+}=\big(i_{v_+},j_{v_+},k_{v_+},\ell\!+\!1\big)\in\wt\cQ_{\Ga}(\vr^*),
\quad
\wt\cQ_{\Ga,v_+}=\cQ_{\Ga}\!\cup\!\big\{\wt\fq_{v_+}\big\}
\subset \wt\cQ_{\Ga}(\vr^*).\EE

\begin{prp}\label{Cblowup_prp}
Suppose $\ell\!\in\!\Z^+$ with $\ell\!\ge\!3$, 
$\vr^*\!\in\!\{0\}\!\sqcup\!\cA_{\ell}$, 
\hbox{$\Ga\!\equiv\!(\Ver,\Edg,\mu)$} is a trivalent $\ell$-marked tree,
and $v_+\!\in\!\nV_{\Ga}(\vr^*)$.
With the notation and assumptions as in~\eref{wtcWdf2_e} and~\eref{wtfqvrdfn_e},
the~map 
\BE{cMspan2_e}
\wt\CR_{\vr^*;\Ga,v_+}\!\equiv\!
(\CR_{\vr^*;\fq})_{\fq\in\wt\cQ_{\Ga,v_+}}\!:
\wt\cW_{\vr^*;\Ga,v_+}'\lra (\C\P^1)^{\wt\cQ_{\Ga,v_+}}\EE
is a well-defined homeomorphism onto an open subset 
$\wt{W}_{\vr^*;\Ga,v_+}'\!\subset\!(\C\P^1)^{\wt\cQ_{\Ga,v_+}}$.
For every $\fq\!\in\!\wt\cQ_{\Ga}(\vr^*)$, $\CR_{\vr^*;\fq}$ is a rational function
on~$\wt{W}_{\vr^*;\Ga,v_+}'$
of the cross ratios~$\CR_{\vr^*;\fq'}$ with $\fq'\!\in\!\wt\cQ_{\Ga,v_+}$.
\end{prp}

\begin{proof}
By~\eref{cMGaDvr_e}, $q_{\vr^*}(\ff_{\ell+1}^{~-1}(\cW_{\Ga})\!)\!\cap\!Y_{\vr^*;\vr}^0\!=\!\eset$
if $\vr\!\in\!\cA_{\ell}(\vr^*)\!-\!\cA_{\Ga}(\vr^*)$.
Thus, the map~\eref{cMspan2_e} is well-defined (and continuous).\\

\noindent
Suppose $\wt\cC_1,\wt\cC_2\!\in\!\wt\cW_{\vr^*;\Ga,v_+}'$ are such that
\hbox{$[\wt\cC_1]_{\vr^*},[\wt\cC_2]_{\vr^*}\!\in\!X_{\vr^*}$} are distinct, but  
$$\CR_{\fq}(\wt\cC_1)=\CR_{\fq}(\wt\cC_2) \qquad\forall\,
\fq\!\in\!\cQ_{\ell}\subset\wt\cQ_{\Ga,v_+}.$$
Since the map~\eref{MSemd_e1} is injective, it follows that 
$$\cC\equiv \ff_{\ell+1}(\wt\cC_1)=\ff_{\ell+1}(\wt\cC_2)\in\cM_{\Ga'}$$ 
for some $\Ga'\!\in\!\cT(\Ga)$.
Below we show~that 
\BE{Cblowup_e2}\CR_{\wt\fq_{v_+}}(\wt\cC_1)\!\neq\!\CR_{\wt\fq_{v_+}}(\wt\cC_2),\EE
thus establishing the injectivity of~\eref{cMspan2_e}.\\

\noindent
Let \hbox{$\Ga'\!=\!(\Ver',\Edg',\mu')$} and $\ka\!:\Ver\!\lra\!\Ver'$ 
be as in~\eref{Gakacond_e}.
Let $r\!=\!1,2$.  By~\eref{ndCprp_e}, the set 
\BE{Cblowup_e8}\cA_{\wt\cC_r}(\vr^*)\equiv\big\{\vr\!\in\!\cA_{\ell}(\vr^*)\!:
\wt\cC_r\!\in\!D_{\ell+1;\vr}\big\}\subset\cA_{\Ga'}(\vr^*)
\subset\cA_{\Ga}(\vr^*)\EE
is ordered by the inclusion of subsets of~$[\ell]$.
If $\cA_{\wt\cC_r}(\vr^*)\!\neq\!\eset$, then
$$\vr_r\!\equiv\!\min\cA_{\wt\cC_r}(\vr^*)\in \cA_{\Ga'}^{\star}(\vr^*).$$
Let $\wt{e}_{\vr_r}\!\equiv\!v_r^+v_r^0$ and
$z_{\ell+1}'(\wt\cC_r)\!\in\!\C\P^1_{v_r^+}$ be the node of~$\wt\cC_r$
determined by~$\vr_r'$.
By~\eref{wtcWdf2_e}, $v_r^+\!=\!\ka(v_+)$.
Since the set  $q_{\vr^*}(\ff_{\ell+1}^{~-1}(\cC)\!\cap\!D_{\ell+1;\vr_r}^-)$ consists of 
a single point, 
\BE{Cblowup_e4}z_{\ell+1}'(\wt\cC_1)\neq z_{\ell+1}'(\wt\cC_2)\in 
\C\P^1_{\ka(v_+)}\EE
if $\cA_{\wt\cC_1}(\vr^*),\cA_{\wt\cC_2}(\vr^*)\!\neq\!\eset$.\\

\noindent
Suppose $r\!=\!1,2$ and $\cA_{\wt\cC_r}(\vr)\!=\!\eset$.
Let $\C\P^1_v$ be the irreducible component of~$\wt\cC_r$ carrying
the marked point~$z_{\ell+1}(\wt\cC_r)$ of~$\wt\cC_r$.
There are two possibilities, which are
illustrated by the left and middle diagrams of Figure~\ref{Cblowup_fig}:
\begin{enumerate}[label=($\wt\cC\arabic*$),leftmargin=*]

\item\label{wtCmain_it} $z_{\ell+1}(\wt\cC_r)$ is a smooth unmarked 
point of the irreducible component~$\C\P^1_{\ka(v_+)}$ of~$\cC$.
We then set \hbox{$z_{\ell+1}'(\wt\cC_r)\!\!=\!z_{\ell+1}(\wt\cC_r)\!\in\!\C\P^1_{\ka(v_+)}$}.

\item\label{wtCmarked_it} $z_{\ell+1}(\wt\cC_r)\!=\!z_m(\cC)$ for some 
$m\!\in\!\mu'^{-1}(\ka(v_+)\!)$,
i.e.~$\C\P^1_v$ also carries exactly one other marked point~$z_m(\wt\cC_r)$ 
of~$\wt\cC_r$ and contains exactly one node, which it shares with 
the irreducible component~$\C\P^1_{\ka(v_+)}$ of~$\cC$.
We then set $z_{\ell+1}'(\wt\cC_r)\!=\!z_m(\cC)\!\in\!\C\P^1_{\ka(v_+)}$.

\setcounter{temp}{\value{enumi}}

\end{enumerate}
Since $[\wt\cC_1]_{\vr^*}\!\neq\![\wt\cC_2]_{\vr^*}$, \eref{Cblowup_e4} holds
in all cases.\\

\begin{figure}
\begin{pspicture}(-2.5,-1.2)(10,2.5)
\psset{unit=.4cm}
\pscircle(3,2){1.5}\rput(2.9,2){\smsize{$\ka(v_+)$}}
\pscircle*(1.94,3.06){.2}\pscircle*(4.06,3.06){.2}
\pscircle(3,-1){1.5}
\pscircle*(1.5,-1){.2}\pscircle*(3,-2.5){.2}\pscircle*(4.5,-1){.2}
\pscircle*(4.5,2){.2}\rput(5.6,2){\smsize{$\ell\!+\!1$}}
\pscircle(13,2){1.5}\rput(13.1,2){\smsize{$\ka(v_+)$}}
\pscircle*(11.94,3.06){.2}
\pscircle(15,4){1.33}\rput(15.1,4){\smsize{$v'$}}\pscircle*(14.06,3.06){.2}
\pscircle*(15.94,3.06){.2}\rput(16.7,2.9){\smsize{$m$}}
\pscircle*(15.94,4.94){.2}\rput(17.1,4.94){\smsize{$\ell\!+\!1$}}
\pscircle(13,-1){1.5}
\pscircle*(11.5,-1){.2}\pscircle*(13,-2.5){.2}\pscircle*(14.5,-1){.2}
\pscircle(23,2){1.5}\rput(23.1,2){\smsize{$v$}}
\pscircle*(21.94,6.06){.2}\pscircle*(24.06,6.06){.2}
\pscircle(23,-1){1.5}\rput(23.1,-1){\smsize{$v'$}}\pscircle*(23,.5){.2}
\pscircle*(21.5,-1){.2}\pscircle*(23,-2.5){.2}\pscircle*(24.5,-1){.2}
\pscircle*(23,3.5){.2}
\pscircle(23,5){1.5}\rput(23.1,5){\smsize{$v''$}}
\pscircle*(24.5,2){.2}\rput(25.6,2){\smsize{$\ell\!+\!1$}}
\end{pspicture}
\caption{The three possibilities for $z_{\ell+1}(\wt\cC_r)\!\in\!\C\P^1_v$
in the proof of Proposition~\ref{Cblowup_prp} and in Remark~\ref{Cblowup_rmk}}
\label{Cblowup_fig}
\end{figure}

\noindent
If $\mu'(i_{v_+})\!=\!\ka(v_+)$, let $z_i'(\cC)\!\in\!\C\P^1_{\ka(v_+)}$ be 
the $i_{v_+}$-th marked point~$z_{\wt{i}_{v_+}}(\cC)$ of~$\cC$.
Otherwise, let $z_i'(\cC)\!\in\!\C\P^1_{\ka(v_+)}$ be the node of this irreducible component
separating it from the component $\C\P^1_{\mu'(i_{v_+})}$ 
carrying~$z_{i_{v_+}}(\cC)$.
Define $z_j'(\cC),z_k'(\cC)\!\in\!\C\P^1_{\ka(v_+)}$ similarly.
Since $i,j,k$ are $(\Ga',v_+)$-independent,
these three points of~$\C\P^1_{\ka(v_+)}$ are distinct.\\

\noindent
In summary, each of the two collections
$$z_i'(\cC),z_j'(\cC),z_k'(\cC)\in\C\P^1_{\ka(v_+)}
\qquad\hbox{and}\qquad z_{\ell+1}'(\wt\cC_1),z_{\ell+1}'(\wt\cC_2)\in\C\P^1_{\ka(v_+)}$$
consists of distinct points (the two collections may overlap). 
Thus, the cross ratios
$$\CR_{\wt\fq_{v_+}}(\wt\cC_r)\equiv
\frac{z_i'(\cC)\!-\!z_k'(\cC)}{z_i'(\cC)\!-\!z_{\ell+1}'(\wt\cC_r)}\!:
\!\frac{z_j'(\cC)\!-\!z_k'(\cC)}{z_j'(\cC)\!-\!z_{\ell+1}'(\wt\cC_r)}
\qquad\hbox{with}~~r=1,2$$
are distinct.
This establishes~\eref{Cblowup_e2}.\\

\noindent
We next show that the map~\eref{cMspan2_e} is open.
By Corollary~\ref{cMspan_crl}\ref{cMfqgenC_it},
it is sufficient to establish this  under 
the assumption that the subcollection $\cQ_{\Ga}\!\subset\!\cQ_{\ell}$ 
is  compatible with
a systematic marking map \hbox{$\eta\!:\wt\Edg\!\lra\![\ell]$} for~$\Ga$. 
Since $(\CR_{\fq})_{\fq\in\wt\cQ_{\Ga}}$ is then a chart on the dense open subset 
$\cW_{\Ga}\!\subset\!\ov\cM_{\ell}$,
\BE{cMspan2_e2} 
\big\{\!(\CR_{\fq})_{\fq\in\wt\cQ_{\Ga}}\big\}\!\big(\cW_{\Ga}\big)
\!\cap\!
\big\{\!(\CR_{\fq})_{\fq\in\wt\cQ_{\Ga}}\big\}\!\big(\ov\cM_{\ell}\!-\!\cW_{\Ga}\big)
=\eset.\EE
Since $\CR_{\wt\fq_{v_+}}(\ff_{\ell+1}^{\,-1}(\cC)\!)\!=\!\C\P^1$
for every $\cC\!\in\!\cM_{\ell}$,
$$\big\{\!(\CR_{\fq})_{\fq\in\wt\cQ_{\Ga,v_+}}\!\big\}
\big(\ff_{\ell+1}^{\,-1}(W)\!\big)=W\!\times\!\C\P^1$$
for every open subset $W\!\subset\!\ov\cM_{\ell}$.\\

\noindent
Let $[\wt\cC]_{\vr^*}\!\in\!\wt\cW_{\vr^*;\Ga,v_+}'$ and
$\wt{C}$ be its image under~\eref{cMspan2_e}.
We show that the image of any neighborhood $\wt{W}_{\vr^*}\!\subset\!\wt\cW_{\vr^*;\Ga,v_+}'$ 
of~$[\wt\cC]_{\vr^*}$
under~\eref{cMspan2_e} contains a neighborhood of~$\wt{C}$.
If not, there exists a sequence $\wt\cC_r\!\in\!\ov\cM_{\ell+1}$ such~that 
$$\wt\cC_r\not\in q_{\vr^*}^{\,-1}\big(\wt\cW_{\vr^*;\Ga,v_+}'\big)
\quad\hbox{and}\quad
\lim_{r\lra\i}
\big\{\!(\CR_{\fq})_{\fq\in\wt\cQ_{\Ga,v_+}}\big\}(\wt\cC_r)=\wt{C}.$$
After passing to a subsequence, we can assume that this sequence converges
to some $\wt\cC'\!\in\!\ov\cM_{\ell+1}$ with 
\BE{cMspan2_e4}[\wt\cC']_{\vr^*}\neq[\wt\cC]_{\vr^*} \quad\hbox{and}\quad
\big\{\!(\CR_{\fq})_{\fq\in\wt\cQ_{\Ga,v_+}}\big\}(\wt\cC')=\wt{C}\,.\EE
By~\eref{cMspan2_e2} and the last property in~\eref{cMspan2_e4},
$$\ff_{\ell+1}(\wt\cC')=\ff_{\ell+1}(\wt\cC)\in\ov\cM_{\ell}\,.$$
By the proof the injectivity of~\eref{cMspan2_e} above
with $\wt\cC_1\!=\!\wt\cC$ and~$\wt\cC_2\!=\!\wt\cC'$,
$$\CR_{\wt\fq_{v_+}}(\wt\cC)\neq \CR_{\wt\fq_{v_+}}(\wt\cC').$$
However, this contradicts the last property in~\eref{cMspan2_e4}.\\

\noindent
Let $\wt\cQ_{\Ga}^{\star}(\vr^*)\!\subset\!\wt\cQ_{\Ga}(\vr^*)$ 
be the subcollection of the tuples~$\fq$ for which 
the rational function claim holds.
As in the proof of Lemma~\ref{cMspan_lmm},
$\wt\cQ_{\Ga}^{\star}(\vr^*)\!\subset\![\ell\!+\!1]^4$ is preserved by the permutations of 
the components of its elements and
\BE{cMspan2_e1}\wt\cQ_{\Ga}^{\star}(\vr^*)\supset 
\wt\cQ_{\Ga}(\vr^*)\!\cap\!\big([\Ga]_{v_+}\!\cup\!\{\ell\!+\!1\}\big)^{\!4}.\EE
By Corollary~\ref{cMspan_crl}\ref{cMfqgenC_it}, 
$\wt\cQ_{\Ga}^{\star}(\vr^*)\!\supset\!\cQ_{\ell}$.\\

\noindent
If $(i,j,k,\ell\!+\!1)\!\in\!\wt\cQ_{\Ga}(\vr^*)$ and $m\!\in\![\ell]$
is $(\Ga,v_+)$-independent of $i,j$, then
$(i,j,m,\ell\!+\!1)\!\in\!\wt\cQ_{\Ga}(\vr^*)$ as~well.
By the last equation in~\eref{CRprp_e}, 
\BE{cMspan2_e3} 
\CR_{\vr^*;ijk,\ell+1}=\CR_{\vr^*;ijkm}\CR_{\vr^*;ijm,\ell+1}
\qquad\hbox{on}~~\wt\cW_{\vr^*;\Ga,v_+}',\EE
if the product on the right side of~\eref{cMspan2_e3} makes sense, i.e.~if
\BE{cMspan2_e5}\begin{split}
\CR_{\vr^*;ijkm}^{~-1}(\i)\!\cap\!
\CR_{\vr^*;ijm,\ell+1}^{~-1}(0)\!\cap\!\wt\cW_{\vr^*;\Ga,v_+}'&=\eset,\\
\CR_{\vr^*;ijkm}^{~-1}(0)\!\cap\!
\CR_{\vr^*;ijm,\ell+1}^{~-1}(\i)\!\cap\!\wt\cW_{\vr^*;\Ga,v_+}'&=\eset.
\end{split}\EE
If the first intersection above is not empty, 
then there exists $\vr\!\in\!\cA_{\Ga}(v_+)$ so that~either
\BE{cMspan2_e6}j,k\in\vr,~~i,m\in\vr_{\ell}^c  \qquad\hbox{or}\qquad
i,m\in\vr,~~j,k\in\vr_{\ell}^c,~~
Y_{\vr^*;\vr}^0\!\cap\!\wt\cW_{\vr^*;\Ga,v_+}'\neq\eset.\EE
If $i,m\!\in\!\vr_{\ell}^c$, then $i,m$ would not be $(\Ga,v_+)$-independent.
If $j,k\!\in\!\vr_{\ell}^c$, then $\vr\!\not\in\!\cA_{\Ga}(\vr^*)$ by the definition 
of~$\wt\cQ_{\Ga}(\vr^*)$.
However, this contradicts the last condition in~\eref{cMspan2_e6}.
Thus, the first set in~\eref{cMspan2_e5} is empty.
By the same reasoning with~$i$ and~$j$ interchanged,
the second set in~\eref{cMspan2_e5} is empty as~well.
It follows that the product on the right side of~\eref{cMspan2_e3} makes sense
and \hbox{$\wt\cQ_{\Ga}^{\star}(\vr^*)\!=\!\wt\cQ_{\Ga}(\vr^*)$}.
\end{proof}

\begin{rmk}\label{Cblowup_rmk}
By the reasoning in the proof of Proposition~\ref{Cblowup_prp}, the~map
$$(\CR_{\vr^*;\fq})_{\fq\in\wt\cQ_{\Ga}(\vr^*)}\!:
q_{\vr^*}\big(\ff_{\ell+1}^{~-1}(\cW_{\Ga})\!\big)\lra (\C\P^1)^{\wt\cQ_{\Ga}(\vr^*)}$$
is also well-defined, continuous, and injective.
In addition to~\ref{wtCmain_it} and~\ref{wtCmarked_it} 
with $v_+\!\in\!\nV_{\Ga}(\vr^*)$ dependent on $r\!=\!1,2$,
a third situation, illustrated by the right diagram of Figure~\ref{Cblowup_fig},
can now arise:
\begin{enumerate}[label=($\wt\cC\arabic*$),leftmargin=*]

\setcounter{enumi}{\value{temp}}

\item\label{wtCnode_it} $z_{\ell+1}(\wt\cC_r)$ is ``at" a node $\{v',v''\}$ of~$\cC$
with \hbox{$v',v''\!\in\!\nV_{\Ga'}(\vr^*)$},
i.e.~$\C\P^1_v$ carries no other marked points of~$\wt\cC_r$ 
and contains exactly two nodes, which it shares with irreducible components 
$\C\P^1_{v'},\C\P^1_{v''}$ of~$\cC$.

\setcounter{temp}{\value{enumi}}

\end{enumerate}
This situation does not arise in the proof of Proposition~\ref{Cblowup_prp}
because $v',v''$ can be chosen so that \hbox{$\vr\!\equiv\!\mu^{-1}(\Ver_{v'v''})$} 
lies in~$\cA_{\Ga}(v_+)$ and thus \hbox{$[\wt\cC_r]_{\vr^*}\!\in\!Y_{\vr^*;\vr}^0$}.
On the other hand, let
$$\wt\cQ_{\ell}(\vr^*)=\cQ_{\ell}\!\cup\!\big\{\fq\!\in\!\cQ_{\ell+1}\!-\!\cQ_{\ell}\!:
|\vr\!\cap\!\fq|\!\ge\!2~\forall\,\vr\!\in\!\cA_{\ell}(\vr^*)\!\big\}\subset \cQ_{\ell+1}.$$
This is the maximal subcollection of~$\cQ_{\ell+1}$ so that the~map
$$(\CR_{\vr^*;\fq})_{\fq\in\wt\cQ_{\ell}(\vr^*)}\!:
X_{\vr^*}\lra (\C\P^1)^{\wt\cQ_{\ell}(\vr^*)}$$
is well-defined (and thus continuous).
However, this map is usually not injective,
as the intersection in~\eref{CblowupInp_e} with~$\cA_{\Ga}(\vr^*)$ 
replaced by $\cA_{\ell}(\vr^*)$ could be empty.
Thus, the cross ratios do not induce an embedding of the quotient~$X_{\vr^*}$
into a product of copies of~$\C\P^1$ for arbitrary values of~$\vr^*$.
\end{rmk}

\vspace{.1in}

\noindent
With the assumptions as in Proposition~\ref{Cblowup_prp},
let $\Ga v_+$ be as in~\eref{wtGav_e}.
By Corollary~\ref{cMspan_crl}\ref{cMfqgenC_it} with~$\ell$ replaced by~$\ell\!+\!1$, 
the collection~$(\CR_{\fq})_{\fq\in\wt\cQ_{\Ga,v_+}}$ of Proposition~\ref{Cblowup_prp}
is a coordinate chart 
on a neighborhood 
$$\cW_{\Ga v_+}'\subset \cW_{\Ga v_+}\!\cap \ff_{\ell+1}^{~-1}(\cW_{\Ga}')
\subset \ov\cM_{\ell+1}-
\bigcup_{\vr\in\cA_{\Ga}(v_+)}\hspace{-.2in}D_{\ell+1;\vr}$$ 
of~$\cM_{\Ga v_+}$ in~$\ov\cM_{\ell+1}$.
If $v_+\!\in\!\nV_{\Ga}(\vr^*)$, the restriction of the quotient~map
$$q_{\vr^*}\!:\cW_{\Ga v_+}'\lra 
\cW_{\vr^*;\Ga v_+}'\!\equiv\!q_{\vr^*}(\cW_{\Ga v_+}')\subset
\wt\cW_{\vr^*;\Ga,v_+}'
\!\equiv\!q_{\vr^*}\big(\ff_{\ell+1}^{~-1}(\cW_{\Ga}')\!\big)
\!-\!\bigcup_{\begin{subarray}{c}\vr\in\cA_{\Ga}(v_+)\\
\vr\not\in\cA_{\Ga}(\vr^*)\end{subarray}}\hspace{-.2in}Y_{\vr^*;\vr}^0
\subset X_{\vr^*}$$
is a homeomorphism onto an open subset.
Proposition~\ref{Cblowup_prp} extends the induced chart $(\CR_{\vr^*;\fq})_{\fq\in\wt\cQ_{\Ga,v_+}}$
on~$\cW_{\vr^*;\Ga v_+}'$ to~$\wt\cW_{\vr^*;\Ga,v_+}'$.
The latter open subset of~$X_{\vr^*}$ includes the points~$[\wt\cC]_{\vr^*}$ 
consisting of an $\ell$-marked curve $\cC\!\!\in\!\cM_{\Ga'}\!\cap\!\cW_{\Ga}'$ 
for some $\Ga'\!\in\!\cT(\Ga)$ with the marked point
$z_{\ell+1}(\wt\cC)\!\in\!\C\P^1_{\ka(v_+)}$, with~$\ka$ as in~\eref{Gakacond_e}.
This marked point
may coincide with the marked point~$z_m(\cC)$ for some $m\!\in\!\mu^{-1}(\ka^{-1}(v_+)\!)$,
as in~\ref{wtCmarked_it} in the proof of Proposition~\ref{Cblowup_prp}, 
or with a node~$\wt{e}_{\vr}\!\in\!\wt\nE_{\Ga'}(\ka(v_+)\!)$ 
corresponding to some $\vr\!\in\!\cA_{\Ga'}(\vr^*)$,
similarly to~\ref{wtCnode_it} in Remark~\ref{Cblowup_rmk}.

\subsection{Proof of Theorem~\ref{Cblowup_thm}}
\label{CblowupPf_subs}

\noindent
\ref{Cspaces_it} 
Since the quotient map~$q_{\vr^*}$ restricts to a homeomorphism
\BE{CblowupPf_e0}
q_{\vr^*}\!:\ov\cM_{\ell+1}-
\bigcup_{\vr\in\cA_{\ell}(\vr^*)}\!\!\!\!\!\!\!D_{\ell+1;\vr}\lra 
X_{\vr^*}-\bigcup_{\vr\in\cA_{\ell}(\vr^*)}\!\!\!\!\!\!Y_{\vr^*;\vr}^0\,,\EE
the right-hand side above inherits a complex structure from~$\ov\cM_{\ell+1}$.
We combine it with the coordinate charts of Proposition~\ref{Cblowup_prp} 
to obtain a complex structure on~$X_{\vr^*}$ with the stated properties.\\

\noindent
Let \hbox{$\Ga\!\equiv\!(\Ver,\Edg,\mu)$} be a trivalent $\ell$-marked tree.
For $\vr\!\in\!\cA_{\Ga}(\vr^*)\!\subset\!\cA_{\ell}(\vr^*)$, define
$$Y_{\vr^*;\Ga,\vr}^0=q_{\vr^*}\big(\ff_{\ell+1}^{~-1}(\cM_{\Ga})
\!\cap\!D_{\ell+1;\vr}\big)
\subset Y_{\vr^*;\vr}^0\,.$$
If $\vr'\!\in\!\cA_{\Ga}(\vr^*)$ and either $\vr'\!\subset\!\vr$ or $\vr'\!\supset\!\vr$,
then $Y_{\vr^*;\Ga,\vr'}^0\!=\!Y_{\vr^*;\Ga,\vr}^0$.
Thus,
\BE{CblowupPf_e2}
\bigcup_{\vr\in\cA_{\Ga}(\vr^*)}\hspace{-.18in}Y_{\vr^*;\Ga,\vr}^0
=\bigcup_{\vr\in\cA_{\Ga}^{\star}(\vr^*)}\hspace{-.18in}Y_{\vr^*;\Ga,\vr}^0
=\bigcup_{v_+\in\nV_{\Ga}(\vr^*)}
\bigcup_{\begin{subarray}{c}\vr\in\cA_{\Ga}^{\star}(\vr^*)\\
v_+\in e_{\vr}\end{subarray}}
\hspace{-.18in}Y_{\vr^*;\Ga,\vr}^0.\EE
The subspaces \hbox{$Y_{\vr^*;\Ga',\vr}^0\!\subset\!Y_{\vr^*;\vr}^0$} 
partition~$Y_{\vr^*;\vr}^0$ as $\Ga'$ runs over the trivalent $\ell$-marked trees 
with $\vr\!\in\!\cA_{\Ga'}(\vr^*)$.\\

\noindent
With $[\Ga]_v\!\subset\![\ell]$, $\cQ_{\Ga}(v)\!\subset\!\cQ_{\ell}$, 
$\fq_e\!\in\!\cQ_{\ell}$, and $\cQ_{\Ga}\!\subset\!\cQ_{\ell}$
as above~\eref{wtcWdf2_e} and  $v_+\!\in\!\nV_{\Ga}(\vr^*)$,
let 
$$\wt\cW_{\vr^*;\Ga,v_+}'\subset X_{\vr^*}, \qquad
\wt\fq_{v_+}\in\wt\cQ_{\Ga}(\vr^*),  \quad\hbox{and}\quad
\wt\cQ_{\Ga,v_+}\!\subset\!\wt\cQ_{\Ga}(\vr^*)$$ 
be as in~\eref{wtcWdf2_e} and~\eref{wtfqvrdfn_e}.
By~\eref{ndCprp_e0}, \hbox{$Y_{\vr^*;\Ga,\vr}^0\!\subset\!\wt\cW_{\vr^*;\Ga,v_+}'$}
for every $\vr\!\in\!\cA_{\Ga}(\vr^*)$ such that $v_+\!\in\!e_{\vr}$.
Along with~\eref{CblowupPf_e2}, this implies that 
the open subsets $\wt\cW_{\vr^*;\Ga',v_+}'\!\subset\!X_{\vr^*}$  cover 
the union on the right-hand side of~\eref{CblowupPf_e0}
as~$\Ga'$ runs over the trivalent $\ell$-marked trees 
and $v_+$ runs over the elements of~$\nV_{\Ga'}(\vr^*)$.\\

\noindent
If $\Ga'\!\equiv\!(\Ver',\Edg',\mu')\!\in\!\cT(\Ga)$, 
$\wt\cQ_{\Ga'}(\vr^*)\!\supset\!\wt\cQ_{\Ga}(\vr^*)$. 
If in addition $v_+'\!\in\!\Ver'$, the overlap~map
$$\wt\CR_{\vr^*;\Ga,v_+}\!\circ\!\wt\CR_{\vr^*;\Ga',v_+'}^{~-1}\!:
\wt\CR_{\vr^*;\Ga',v_+'}
\big(\wt\cW_{\vr^*;\Ga,v_+}'\!\cap\!\wt\cW_{\vr^*;\Ga',v_+'}'\big)
\lra 
\wt\CR_{\vr^*;\Ga,v_+}
\big(\wt\cW_{\vr^*;\Ga,v_+}'\!\cap\!\wt\cW_{\vr^*;\Ga',v_+'}'\big)$$
is thus holomorphic by the rational function claim of Proposition~\ref{Cblowup_prp}
with~$(\Ga,v_+)$ replaced by~$(\Ga',v_+')$.
If $\Ga''\!\equiv\!(\Ver'',\Edg',\mu'')$  is any trivalent $\ell$-marked tree
and $v_+''\!\in\!\Ver''$,
$\wt\cW_{\vr^*;\Ga,v_+}'\!\cap\!\wt\cW_{\vr^*;\Ga'',v_+''}'$
is contained in the union of the open subsets 
$\wt\cW_{\vr^*;\Ga',v_+'}'\!\subset\!X_{\vr^*}$ with 
\hbox{$\Ga'\!\in\!\cT(\Ga),\cT(\Ga'')$} as above and $v_+'\!\in\!\Ver'$.
Thus, the coordinate charts of Proposition~\ref{Cblowup_prp}
determine a complex structure on~$X_{\vr^*}$.\\

\noindent
We next show that each subspace $Y_{\vr^*;\vr}^+,Y_{\vr^*;\vr}^0\!\subset\!X_{\vr^*}$
is a slice in some coordinate chart of Proposition~\ref{Cblowup_prp}.
Let $v_+\!\in\!\nV_{\Ga}(\vr^*)$ and $\vr\!\subset\![\ell]$ with 
$|\vr|,|\vr_{\ell}^c|\!\ge\!2$.
By~\eref{Yvrdfn_e}, 
$$\big(Y_{\vr^*;\vr}^+\!\cap\!\wt\cW_{\vr^*;\Ga,v_+}'\big)
\!\cup\! \big(Y_{\vr^*;\vr}^0\!\cap\!\wt\cW_{\vr^*;\Ga,v_+}'\big)
= q_{\vr^*}\big(f_{\ell+1}^{~-1}(D_{\ell;\vr})\!\big)\!\cap\!\wt\cW_{\vr^*;\Ga,v_+}'.$$
By~\eref{wtcWdf2_e} and~\eref{cMGaDvr_e}, 
$\vr\!\in\!\cA_{\Ga}(v_+)$ or $\vr_{\ell}^c\!\in\!\cA_{\Ga}(v_+)$
if the right-hand side above is non-empty.
In such a case,
\eref{cMGaDvr_e2} gives
\BE{Csubspaces_e8c}\begin{split}
\big(Y_{\vr^*;\vr}^+\!\cap\!\wt\cW_{\vr^*;\Ga,v_+}'\big)
\!\cup\! \big(Y_{\vr^*;\vr}^0\!\cap\!\wt\cW_{\vr^*;\Ga,v_+}'\big)
&\equiv q_{\vr^*}\big(f_{\ell+1}^{~-1}(D_{\ell;\vr})\!\big)\!\cap\!\wt\cW_{\vr^*;\Ga,v_+}'\\
&=\big\{x\!\in\!\wt\cW_{\vr^*;\Ga,v_+}'\!\!:
\CR_{\vr^*;\fq_{e_{\vr}}}\!(x)\!=\!0\big\}.
\end{split}\EE
By~\eref{wtcWdf2_e},
\BE{Csubspaces_e8d}\begin{aligned}
Y_{\vr^*;\vr}^0\!\cap\!\wt\cW_{\vr^*;\Ga,v_+}'&=\eset
&\quad&\hbox{if}~~
\vr\!\in\!\cA_{\Ga}(v_+)\!-\!\cA_{\ell}(\vr^*),\\
Y_{\vr^*;\vr}^+\!\cap\!\wt\cW_{\vr^*;\Ga,v_+}'
=Y_{\vr^*;\vr^c_{\ell}}^0\!\cap\!\wt\cW_{\vr^*;\Ga,v_+}'
&=\eset &\quad&\hbox{if}~~\vr\!\in\!\cA_{\ell},~
\vr^c_{\ell}\!\in\!\cA_{\Ga}(v_+).
\end{aligned}\EE

\vspace{.15in}

\noindent
By~\eref{Csubspaces_e8c}, \eref{Csubspaces_e8d}, and~\eref{Yvrmp_e},  
\BE{Csubspaces_e5b}
Y_{\vr^*;\vr}^+\!\cap\!\wt\cW_{\vr^*;\Ga,v_+}'
=\begin{cases}\{x\!\in\!\wt\cW_{\vr^*;\Ga,v_+}'\!\!:
\CR_{\vr^*;\fq_{e_{\vr}}}(x)\!=\!0\},
&\hbox{if}~\vr\!\in\!\cA_{\ell}\!\cap\!\cA_{\Ga}(v_+);\\
\eset,&\hbox{if}~\vr\!\in\!\cA_{\ell}\!-\!\cA_{\Ga}(v_+).
\end{cases}\EE
By~\eref{Csubspaces_e8c} and~\eref{Csubspaces_e8d},
\BE{Csubspaces_e5a}
Y_{\vr^*;\vr}^0\!\cap\!\wt\cW_{\vr^*;\Ga,v_+}'
=\begin{cases}
\{x\!\in\!\wt\cW_{\vr^*;\Ga,v_+}'\!\!:
\CR_{\vr^*;\fq_{e_{\vr}}}(x)\!=\!0\},
&\hbox{if}~\vr\!\in\!\cA_{\ell},\,\vr_{\ell}^c\!\in\!\cA_{\Ga}(v_+);\\
\eset,
&\hbox{if}~\vr\!\in\!\cA_{\Ga}(v_+)\!-\!\cA_{\ell}(\vr^*),\\
\eset,
&\hbox{if}~\vr\!\in\!\cA_{\ell},~\vr,\vr^c_{\ell}\!\not\in\!\cA_{\Ga}(v_+).
\end{cases}\EE
Thus, the subspaces $Y_{\vr^*;\vr}^+\!\subset\!X_{\vr^*}$ with $\vr\!\in\!\cA_{\ell}$ and
$Y_{\vr^*;\vr}^0\!\subset\!X_{\vr^*}$ with $\vr\!\in\!\cA_{\ell}\!-\!\cA_{\ell}(\vr^*)$
are complex submanifolds.
Since $v_+\!\in\!\nV_{\Ga}(\vr^*)$, the conditions in the first case of~\eref{Csubspaces_e5a}
imply that $\vr\!\in\!\cA_{\ell}\!-\!\cA_{\ell}(\vr^*)$.
The remaining case of the intersection in~\eref{Csubspaces_e5a}, 
with $\vr\!\in\!\cA_{\Ga}(\vr^*)$
and thus $Y_{\vr^*;\vr}^0\!\subset\!Y_{\vr^*;\vr}^+$, is 
treated in the next paragraph.\\

\noindent
Suppose either $\vr\!\in\!\cA_{\Ga}(\vr^*)$ or $\vr\!=\![\ell]\!-\{j\}$
for some $j\!\in\![\ell]$.
We can then choose~$\wt\fq_{v_+}$ in~\eref{wtfqvrdfn_e} 
so that $j_{v_+}\!\in\!\vr_{\ell}^c$.
Since $v_+\!\in\!\nV_{\Ga}(\vr^*)$, then $i_{v_+},k_{v_+}\!\in\!\vr$.
Thus, $\CR_{\wt\fq_{v_+}}$ vanishes on~$D_{\ell+1;\vr}$ and
\BE{Csubspaces_e8a}Y_{\vr^*;\vr}^0\!\cap\!\wt\cW_{\vr^*;\Ga,v_+}'
\subset\big\{x\!\in\!\wt\cW_{\vr^*;\Ga,v_+}'\!\!:
\CR_{\vr^*;\wt\fq_{v_+}}\!(x)\!=\!0\big\}\,.\EE
If $[\wt\cC]_{\vr^*}\!\in\!\wt\cW_{\vr^*;\Ga,v_+}'$ with
$$  \CR_{\wt\fq_{v_+}}(\wt\cC)=0, \quad
\ff_{\ell+1}(\wt\cC)\in\cM_{\Ga'},
\quad \Ga'\!\equiv\!(\Ver',\Edg',\mu')\in\cT(\Ga),$$
and $\ka\!:\Ver\!\lra\!\Ver'$ as in~\eref{Gakacond_e}, then
either $\wt\cC\!\in\!D_{\ell+1;\vr}$ or there exists 
$$\vr'\in\cA_{\Ga'}\big(\ka(v_+)\!\big)\subset\cA_{\Ga}(v_+) \qquad\hbox{s.t.}\quad
i_{v_+},k_{v_+}\in\vr',~~j_{v_+}\not\in\vr',~~
\wt\cC\in D_{\ell+1;\vr'}\,.$$
In the latter case, $\vr'\!\in\!\cA_{\ell}(\vr^*)$ by~\eref{wtcWdf2_e}
and there exists 
$$\wt\cC'\in D_{\ell+1;\vr'}\!\cap\!D_{\ell+1;[\ell]-\{j_{v_+}\}}
\quad\hbox{s.t.}~~\ff_{\ell+1}(\wt\cC')=\ff_{\ell+1}(\wt\cC),$$
i.e.~$\wt\cC'\!\sim_{\vr'}\!\wt\cC$ and so $[\wt\cC']_{\vr^*}\!=\![\wt\cC]_{\vr^*}$.
Thus,
\BE{Csubspaces_e3b}
Y_{\vr^*;[\ell]-\{j_{v_+}\}}^0\!\cap\!\wt\cW_{\vr^*;\Ga,v_+}'
=\big\{x\!\in\!\wt\cW_{\vr^*;\Ga,v_+}'\!\!:
\CR_{\vr^*;\wt\fq_{v_+}}\!(x)\!=\!0\big\}.\EE
Along with~\eref{Yvrmp_e} and~\eref{Csubspaces_e5b}, this gives
\BE{Csubspaces_e3}
Y_{\vr^*;\vr}^0\!\cap\!\wt\cW_{\vr^*;\Ga,v_+}'
=\big\{x\!\in\!\wt\cW_{\vr^*;\Ga,v_+}'\!\!:
\CR_{\vr^*;\fq_{e_{\vr}}}(x)\!=\!0,\,\CR_{\vr^*;\wt\fq_{v_+}}\!(x)\!=\!0\big\}\EE
for $\vr\!\in\!\cA_{\Ga}(\vr^*)$.
Thus, the subspaces $Y_{\vr^*;\vr}^0\!\subset\!X_{\vr^*}$ 
with $\vr\!\in\!\wt\cA_{\ell}(\vr^*)$
are complex submanifolds as~well.\\

\noindent
\ref{Cblowup_it} The holomorphic map~$\pi_{\vr^*}$ in~\eref{bldowndfn_e} restricts to 
a bijection from~$X_{\vr^*}\!-\!Y_{\vr^*;\vr^*}^0$ to \hbox{$X_{\vr^*-1}\!-\!Y_{\vr^*-1;\vr^*}^0$}
and sends $Y_{\vr^*;\vr^*}^0\!\subset\!X_{\vr^*}$ to~$Y_{\vr^*-1;\vr^*}^0\!\subset\!X_{\vr^*-1}$.
Below we cover~$Y_{\vr^*-1;\vr^*}^0$ and~$Y_{\vr^*;\vr^*}^0$ by charts as
in Lemma~\ref{StanBl_lmm}.\\

\noindent
Suppose \hbox{$\Ga\!\equiv\!(\Ver,\Edg,\mu)$} is a trivalent $\ell$-marked tree
with \hbox{$\vr^*\!\in\!\cA_{\Ga}(\vr^*\!-\!1)$} 
and $\wt{e}_{\vr^*}\!=\!v_+v_+'$.
In particular, 
$$v_+\in\nV_{\Ga}(\vr^*\!-\!1), \qquad v_+'\not\in\nV_{\Ga}(\vr^*\!-\!1),
\quad\hbox{and}\quad
e_+\!\equiv\!e_{\vr^*}\!\equiv\!\{v_+,v_+'\}\subset\nV_{\Ga}(\vr^*).$$
Let \hbox{$\wt\cQ_{\Ga}(\vr^*\!-\!1),\wt\cQ_{\Ga}(\vr^*)\!\subset\!\cQ_{\ell+1}$} 
be as in~\eref{wtcQGavrdfn_e}
and \hbox{$\cA_{\Ga}(v_+),\cA_{\Ga}(v_+')\!\subset\![\ell]$} as in~\eref{wtcWdf_e}.
Define
$$\cA_{\Ga}(e_+)=\cA_{\Ga}(v_+)\!\cup\!\cA_{\Ga}(v_+')-\big\{\vr^*,(\vr^*)_{\ell}^c\big\}.$$
Take 
\hbox{$[\Ga]_v\!\subset\![\ell]$}, \hbox{$\cQ_{\Ga}(v)\!\subset\!\cQ_{\ell}$}, 
$\fq_e\!\in\!\cQ_{\ell}$, and $\cQ_{\Ga}\!\subset\!\cQ_{\ell}$ be as above~\eref{wtcWdf2_e}
with \hbox{$i_{e_+},k_{e_+}\!\in\!\vr^*$}.
Let
\begin{alignat*}{2}
\wt\fq_{v_+},\wt\fq_{e_+}&=
\big(i_{e_+},j_{e_+},k_{e_+},\ell\!+\!1\big)\in\wt\cQ_{\Ga}(\vr^*\!-\!1), 
&\qquad
\wt\fq_{v_+'}&=\big(i_{e_+},j_{e_+},m_{e_+},\ell\!+\!1\big)\in\wt\cQ_{\Ga}(\vr^*),\\
\wt\fq_{e_+}'&=\big(i_{e_+},j_{e_+},\ell\!+\!1,m_{e_+}\big)\in\wt\cQ_{\Ga}(\vr^*),
&\qquad 
\wt\cQ_{\Ga;\vr^*}&=\cQ_{\Ga}\!-\!\big\{\fq_{e_+}\big\}\subset\cQ_{\ell}\,.
\end{alignat*}
The coordinates in our charts will be indexed by the collections
\BE{Cblowup_e23}
\wt\cQ_{\Ga,v_+}=\wt\cQ_{\Ga;\vr^*}\!\cup\!\big\{\fq_{e_+},\wt\fq_{v_+}\!\big\},~~
\wt\cQ_{\Ga,v_+'}=\wt\cQ_{\Ga;\vr^*}\!\cup\!\big\{\fq_{e_+},\wt\fq_{v_+'}\!\big\},~~
\wt\cQ_{\Ga,e_+}=\wt\cQ_{\Ga;\vr^*}\!\cup\!\big\{\wt\fq_{v_+},\wt\fq_{e_+}'\!\big\}.\EE

\vspace{.15in}

\noindent
Let $\cW_{\Ga}'\!\subset\!\ov\cM_{\ell}$ and 
\hbox{$\wt\cW_{\vr^*;\Ga,v_+'}'\!\!\subset\!X_{\vr^*}$} be as in~\eref{wtcWdf2_e}.
Define
\begin{equation*}\begin{split}
\wt\cW_{\vr^*-1;\Ga,v_+}''&=q_{\vr^*-1}\big(\ff_{\ell+1}^{~-1}(\cW_{\Ga}')\!\big)
\!-\!\bigcup_{\begin{subarray}{c}\vr\in\cA_{\Ga}(v_+)\\ \vr\not\supset\vr^*\end{subarray}}
\hspace{-.18in}Y_{\vr^*-1;\vr}^0
\!-\!\CR_{\vr^*-1;\wt\fq_{e_+}}^{-1}\!(\i)\subset \wt\cW_{\vr^*-1;\Ga,v_+}'
\subset X_{\vr^*-1};\\
\wt\cW_{\vr^*;\Ga,e_+}''&=q_{\vr^*}\big(\ff_{\ell+1}^{~-1}(\cW_{\Ga}')\!\big)
\!-\!\bigcup_{\vr\in\cA_{\Ga}(e_+)}
\hspace{-.18in}Y_{\vr^*;\vr}^0
\!-\!\CR_{\vr^*;\wt\fq_{e_+}}^{-1}\!(\i)
\!-\!\CR_{\vr^*;\wt\fq_{e_+}'}^{-1}\!(\i)\subset X_{\vr^*}.
\end{split}\end{equation*}
Since $i_{e_+},j_{e_+},k_{e_+}$ are $(\Ga,v_+)$-independent
and $v_+\!\in\!\nV_{\Ga}(\vr^*\!-\!1)$,
$\wt\fq_{e_+}\!\in\!\wt\cQ_{\Ga}(\vr^*\!-\!1)$ and so 
$\CR_{\vr^*-1;\wt\fq_{e_+}}$ is well-defined on~$\wt\cW_{\vr^*-1;\Ga,v_+}'$.
By the definitions of the open subsets $\wt\cW_{\vr^*-1;\Ga,v_+}''$,
$\wt\cW_{\vr^*;\Ga,v_+'}'$, and~$\wt\cW_{\vr^*;\Ga,e_+}''$,
\begin{equation*}\begin{split}
\CR_{\vr^*-1;\fq_{e_+}}\big(\wt\cW_{\vr^*-1;\Ga,v_+}''\big),
\CR_{\vr^*;\fq_{e_+}}\big(\wt\cW_{\vr^*;\Ga,v_+'}'\big),
\CR_{\vr^*;\wt\fq_{e_+}'}\big(\wt\cW_{\vr^*;\Ga,e_+}''\big)
&\subset\C,\\
\CR_{\vr^*-1;\wt\fq_{v_+}}\big(\wt\cW_{\vr^*-1;\Ga,v_+}''\big),
\CR_{\vr^*;\wt\fq_{v_+'}}\big(\wt\cW_{\vr^*;\Ga,v_+'}'\big),
\CR_{\vr^*;\wt\fq_{v_+}}\big(\wt\cW_{\vr^*;\Ga,e_+}''\big)&\subset\C.
\end{split}\end{equation*}
By the definition of the equivalence relation~$\sim_{\vr^*}$ in Section~\ref{Cintro_subs},
\BE{Cblowup_e25}
\CR_{\vr^*;\fq}\!=\!\CR_{\vr^*-1;\fq}\!\circ\!\pi_{\vr^*}\!:
\pi_{\vr^*}^{\,-1}\big(\wt\cW_{\vr^*-1;\Ga,v_+}'')
\!=\!\wt\cW_{\vr^*;\Ga,v_+'}'\!\cup\!\wt\cW_{\vr^*;\Ga,e_+}''\lra\C\P^1
\quad\forall\,\fq\!\in\!\wt\cQ_{\Ga,v_+}\,.\EE
Combining this equation with~\eref{CRprp_e}, we obtain
\BE{Cblowup_e23b}\begin{split}
\CR_{\vr^*-1;\wt\fq_{v_+}}\!\!\circ\!\pi_{\vr^*}
&=\!\CR_{\vr^*;\wt\fq_{v_+}}
\!=\!\CR_{\vr^*;\wt\fq_{v_+'}}\!\cdot\CR_{\vr^*;\fq_{e_+}}\!\!:
\wt\cW_{\vr^*;\Ga,v_+'}'\lra\C, \\
\CR_{\vr^*-1;\fq_{e_+}}\!\!\circ\!\pi_{\vr^*}
&=\!\CR_{\vr^*;\fq_{e_+}}
\!=\!\CR_{\vr^*;\wt\fq_{e_+}'}\!\cdot\CR_{\vr^*;\wt\fq_{v_+}}\!\!:
\wt\cW_{\vr^*;\Ga,e_+}''\lra\C.
\end{split}\EE

\vspace{.15in}

\noindent
Since the quotient map~$q_{\vr^*}$ restricts to a homeomorphism
$$q_{\vr^*}\!:
\ff_{\ell+1}^{~-1}(\cW_{\Ga}')
\!-\!\bigcup_{\vr\in\cA_{\Ga}(e_+)}\hspace{-.22in}D_{\ell+1;\vr}
-\!\CR_{\wt\fq_{e_+}}^{-1}\!(\i)\!-\!\CR_{\wt\fq_{e_+}'}^{-1}\!(\i)
\lra \wt\cW_{\vr^*;\Ga,e_+}''\,,$$
Corollary~\ref{cMspan_crl}\ref{cMfqgenC_it} and~\ref{Cspaces_it} above imply that 
the holomorphic coordinates 
$\{\CR_{\fq}\}_{\fq\in\wt\cQ_{\Ga,e_+}}$ on~$\cW_{\Ga e_+}'$,
with $\Ga e_+$ as in~\eref{wtGae_e}, 
descend to the holomorphic coordinates
$\{\CR_{\vr^*;\fq}\}_{\fq\in\wt\cQ_{\Ga,e_+}}$ on~$\wt\cW_{\vr^*;\Ga,e_+}''$.
By Proposition~\ref{Cblowup_prp}  and~\ref{Cspaces_it} above again,
the collections $\{\CR_{\vr^*-1;\fq}\}_{\fq\in\wt\cQ_{\Ga,v_+}}$ and
$\{\CR_{\vr^*;\fq}\}_{\fq\in\wt\cQ_{\Ga,v_+'}}$ are
holomorphic coordinates on~$\wt\cW_{\vr^*-1;\Ga,v_+}''$ and
$\wt\cW_{\vr^*;\Ga,v_+'}'$, respectively.
Since the open subsets \hbox{$\wt\cW_{\vr^*-1;\Ga,v_+}'\!\subset\!X_{\vr^*-1}$}
cover the complex submanifold~$Y_{\vr^*-1;\vr^*}^0$
as~$\Ga$ runs over the trivalent $\ell$-marked trees with 
$\vr^*\!\in\!\cA_{\Ga}(\vr^*\!-\!1)$,
\eref{Csubspaces_e3} with~$(\vr^*,\vr)$ replaced~by $(\vr^*\!-\!1,\vr^*)$ 
implies that so do the open subsets~$\wt\cW_{\vr^*-1;\Ga,v_+}''$ and 
\BE{Csubspaces_e3d}
Y_{\vr^*-1;\vr^*}^0\!\cap\!\wt\cW_{\vr^*-1;\Ga,v_+}''
=\big\{x\!\in\!\wt\cW_{\vr^*-1;\Ga,v_+}'\!\!:
\CR_{\vr^*-1;\fq_{e_+}}(x)=\!0,\,\CR_{\vr^*-1;\wt\fq_{v_+}}\!(x)=\!0\big\}.\EE
By Lemma~\ref{StanBl_lmm} with $(\bF,\fc)\!=\!(\C,2)$ and \eref{Cblowup_e23}-\eref{Csubspaces_e3d},
the map~$\pi_{\vr^*}$ in~\eref{bldowndfn_e}
is the holomorphic blowup of~$X_{\vr^*-1}$ along~$Y_{\vr^*-1;\vr^*}^0$
with the exceptional divisor~$Y_{\vr^*;\vr^*}^0$.\\

\noindent
\ref{Csubspaces_it} Let $\vr\!\in\!\cA_{\ell}(\vr^*)$.
If 
$$Y_{\vr^*-1;\vr^*}^0\!\cap\!Y_{\vr^*-1;\vr}^0\!\cap\!\wt\cW_{\vr^*-1;\Ga,v_+}'\neq\eset,$$
then $\vr\!\supsetneq\!\vr^*$ and $\vr\!\in\!\cA_{\Ga}(v_+),\cA_{\Ga}(v_+')$.
We then continue with the setup in~\ref{Cspaces_it},
choosing~$\wt\fq_{v_+}$ in~\eref{wtfqvrdfn_e} so that 
\hbox{$j_{v_+}\!\!\in\!\vr_{\ell}^c\!\subset\!(\vr^*)_{\ell}^c$}.
By~\eref{Csubspaces_e3} with $(\vr^*,\vr)$ replaced by~$(\vr^*\!-\!1,\vr^*)$ 
and~$(\vr^*\!-\!1,\vr)$,
\begin{equation*}\begin{split}
Y_{\vr^*-1;\vr^*}^0\!\cap\!Y_{\vr^*-1;\vr}^0\!\cap\!\wt\cW_{\vr^*-1;\Ga,v_+}''
=\big\{&x\!\in\!\wt\cW_{\vr^*-1;\Ga,v_+}''\!\!:\\
&\CR_{\vr^*-1;\fq_{e_+}}(x),\CR_{\vr^*-1;\fq_{e_{\vr}}}(x),
\CR_{\vr^*-1;\wt\fq_{v_+}}\!(x)\!=\!0\big\}.
\end{split}\end{equation*}
Along with~\eref{Csubspaces_e3} with $(\vr^*,\vr)$ replaced by~$(\vr^*\!-\!1,\vr^*)$ 
and~$(\vr^*\!-\!1,\vr)$ again, this implies that 
\hbox{$Y_{\vr^*-1;\vr^*}^0,Y_{\vr^*-1;\vr}^0\!\subset\!X_{\vr^*-1}$} intersect cleanly
for all $\vr\!\in\!\cA_{\ell}(\vr^*)$.\\

\noindent
Let $\vr\!\in\!\cA_{\Ga}(v_+),\cA_{\Ga}(v_+')$ as before.
By~\eref{Csubspaces_e3} as stated and with~$\vr$ replaced by~$\vr^*$,
\begin{equation*}\begin{split}
Y_{\vr^*;\vr^*}^0\!\cap\!Y_{\vr^*;\vr}^0\!\cap\!
\pi_{\vr^*}^{\,-1}\big(\wt\cW_{\vr^*-1;\Ga,v_+'}'\big)
=\big\{&\wt{x}\!\in\!\wt\cW_{\vr^*;\Ga,v_+'}''\!\!:\\
&\CR_{\vr^*;\fq_{e_+}}(\wt{x}),\CR_{\vr^*;\fq_{e_{\vr}}}(\wt{x}),
\CR_{\vr^*;\wt\fq_{v_+'}}\!(\wt{x})\!=\!0\big\}.
\end{split}\end{equation*}
Along with~\eref{Cblowup_e25}, the first line in~\eref{Cblowup_e23b},
and~\eref{Csubspaces_e3} as stated and with~$\vr$ replaced by~$\vr^*$ again, 
this implies that 
\hbox{$Y_{\vr^*;\vr^*}^0,Y_{\vr^*;\vr}^0\!\subset\!X_{\vr^*}$} intersect cleanly
for all $\vr\!\in\!\cA_{\ell}(\vr^*)$,
$$T_{\wt{x}}Y_{\vr^*;\vr}^0\cap\ker\nd_{\wt{x}}\big\{\pi_{\vr^*}|_{Y_{\vr^*;\vr^*}^0}\big\}=\{0\}
\quad\forall\,\wt{x}\!\in\!Y_{\vr^*;\vr^*}^0\!\cap\!Y_{\vr^*;\vr}^0,$$
and the restriction~$\pi_{\vr^*}|_{Y_{\vr^*;\vr}^0}$ is a bijection
onto~$Y_{\vr^*-1;\vr}^0$.
This establishes the $\pi_{\vr^*}$-equivalent claim.\\

\noindent
For any $\vr\!\in\!\wt\cA_{\ell}$, the restrictions
$$\pi_{\vr^*}\!:Y_{\vr^*;\vr}^+\!-\!Y_{\vr^*;\vr^*}^0\lra 
Y_{\vr^*-1;\vr}^+\!-\!Y_{\vr^*-1;\vr^*}^0 
\quad\hbox{and}\quad
\pi_{\vr^*}\!:Y_{\vr^*;\vr}^0\!-\!Y_{\vr^*;\vr^*}^0\lra 
Y_{\vr^*-1;\vr}^0\!-\!Y_{\vr^*-1;\vr^*}^0$$
are holomorphic bijections.
The submanifolds $Y_{\vr^*;\vr}^+\!\subset\!X_{\vr^*}$ with $\vr\!\in\!\wt\cA_{\ell}$
and $Y_{\vr^*;\vr}^0\!\subset\!X_{\vr^*}$ with $\vr\!\not\in\!\cA_{\ell}(\vr^*\!-\!1)$
are complex divisors.
By definition, 
\begin{equation*}\begin{split}
Y_{\vr^*;\vr}^{\bu}\!\cap\!Y_{\vr^*;\vr^*}^0
&=q_{\vr^*}\big(Y_{\vr_{\max};\vr}^{\bu}\!\cap\!Y_{\vr_{\max};\vr^*}^0 \big)
\subset X_{\vr^*} \qquad\hbox{and}\\
Y_{\vr^*-1;\vr}^{\bu}\!\cap\!Y_{\vr^*-1;\vr^*}^0
&=q_{\vr^*-1}\big(Y_{\vr_{\max};\vr}^{\bu}\!\cap\!Y_{\vr_{\max};\vr^*}^0 \big)
\subset X_{\vr^*-1}.
\end{split}\end{equation*}
Since the intersections of distinct divisors $D_{\ell+1;\wt\vr}\!\subset\!\ov\cM_{\ell+1}$
with $\wt\vr\!\in\!\cA_{\ell+1}$ are transverse, 
$$Y_{\vr_{\max};\vr}^{\bu}\!\cap\!Y_{\vr_{\max};\vr^*}^0\subset Y_{\vr_{\max};\vr}^{\bu}$$
is a divisor whenever $(\vr,\bu)\!\neq\!(\vr^*,0)$.
Thus,
\begin{alignat*}{2}
\big\{\pi_{\vr_*}\big\}_*\big([Y_{\vr^*;\vr}^+]\big)&=[Y_{\vr^*-1;\vr}^+]\in 
H_*\big(Y_{\vr^*-1;\vr}^+;\Z\big)
&\quad&\forall~\vr\!\in\!\wt\cA_{\ell} \qquad\hbox{and}\\
\big\{\pi_{\vr_*}\big\}_*\big([Y_{\vr^*;\vr}^0]\big)&=[Y_{\vr^*-1;\vr}^0]
\in H_*\big(Y_{\vr^*-1;\vr}^0;\Z\big)
&\quad&\forall~\vr\!\in\!\wt\cA_{\ell}\!-\!\cA_{\ell}(\vr^*\!-\!1)\,.
\end{alignat*}
This establishes the $(\pi_{\vr^*},\Z)$-related claim.

\section{The real quotient sequence}
\label{Rblowup_sec}

\noindent
For $\ell\!\in\!\Z^+$ with $\ell\!\ge\!2$, 
we denote by $\ov\cM_{\ell}^{\pm}$ the Deligne-Mumford moduli space
of stable rational $[\ell^{\pm}]$-marked curves and 
by~$\cQ_{\ell}^{\pm}$ the set of quadruples~$\fq$
of distinct elements of~$[\ell^{\pm}]$.
For $\fq\!\in\!\cQ_{\ell}^{\pm}$, we define the cross ratios
$$\CR_{\fq}\!: \ov\cM_{\ell}^{\pm}\lra\C\P^1 $$
as in~\eref{CRftdfn_e}.
The relations~\eref{CRprp_e} hold for all distinct elements $i,j,k,m,n$ of~$[\ell^{\pm}]$.
For $\fq\!\in\!\cQ_{\ell}^{\pm}$, $i\!\in\![\ell^{\pm}]$, and $I\!\subset\![\ell^{\pm}]$,
we define $i\!\in\!\fq$ and $|\fq\!\cap\!I|$ to mean that $i$ is one of 
the four components of~$\fq$ and the number of components of~$\fq$ contained in the set~$I$,
respectively.
For $i\!\in\!\Z^+$, let $\ov{i^{\pm}}\!=\!i^{\mp}$.
We denote the induced involution on~$\cQ_{\ell}^{\pm}$ in the same way.
For $\vr\!\in\!\cA_{\ell}^{\pm}$, we define the divisor
$$D_{\ell;\vr}^{\pm}\subset \ov\cM_{\ell}^{\pm}$$
analogously to $D_{\ell;\vr}\!\subset\!\ov\cM_{\ell}$ in Section~\ref{Cintro_subs}.

\subsection{The combinatorial structure of $\R\ov\cM_{0,\ell}$}
\label{RcMellstrC_subs0}

\noindent
The topological type of an element of $\R\ov\cM_{0,\ell}$ is described
by a \sf{trivalent real $\ell$-marked tree}.
An \sf{$[\ell^{\pm}]$-marked tree} is a tuple \hbox{$\Ga\!\equiv\!(\Ver,\Edg,\mu)$} so that
$(\Ver,\Edg)$ is a tree and \hbox{$\mu\!:[\ell^{\pm}]\!\lra\!\Ver$} is a map.
For such a tuple~$\Ga$ and $v\!\in\!\Ver$, the \sf{valence} of~$v$ in~$\Ga$
is the number
$$\val_{\Ga}(v)=\big|\mu^{-1}(v)\!\sqcup\!\{e\!\in\!\Edg\!:e\!\ni\!v\}\big|.$$
An $[\ell^{\pm}]$-marked tree $\Ga\!\equiv\!(\Ver,\Edg,\mu)$ is \sf{trivalent}
if $\val_{\Ga}(v)\!\ge\!3$ for every $v\!\in\!\Ver$.\\

\noindent
A \sf{real} $[\ell^{\pm}]$-marked tree is a tuple $(\Ver,\Edg,\mu,\phi)$
so that $(\Ver,\Edg,\mu)$ is an $[\ell^{\pm}]$-marked tree and
$\phi$ is an involution on the set~$\Ver$ so~that
$$\big\{\phi(v),\phi(v')\!\big\}\in\Edg~~\forall\,\{v,v'\}\!\in\!\!\Edg
\quad\hbox{and}\quad
\phi\big(\mu(i)\!\big)=\mu(\ov{i})~~\forall\,i\!\in\![\ell^{\pm}].$$
A stable real rational curve~$\cC$ with $\ell$~conjugate pairs of marked points 
determines a trivalent real $[\ell^{\pm}]$-marked tree $\R\Ga\!\equiv\!(\Ver,\Edg,\mu,\phi)$
so~that $\Ga\!\equiv\!(\Ver,\Edg,\mu)$ is the dual graph of~$\cC$ without its involution.
The stability conditions on~$\cC$ determine the involution~$\phi$ on~$\Ver$.
We call the trivalent real $[\ell^{\pm}]$-marked tree $\R\Ga$ obtained in this way
the \sf{dual graph} of~$\cC$.
For a trivalent real $[\ell^{\pm}]$-marked tree~$\R\Ga$, 
denote~by \hbox{$\cM_{\R\Ga}\!\subset\!\R\ov\cM_{0,\ell}$} 
the stratum of the elements~$[\cC]$ with dual graph~$\R\Ga$.\\

\noindent
For a real $[\ell^{\pm}]$-marked tree $\R\Ga\!\equiv\!(\Ver,\Edg,\mu,\phi)$, 
we define the collection~$\wt\Edg$ of oriented edges and 
the~subsets 
$$\nE_{\R\Ga}(v)\subset\Edg, \quad
\wt\nE_{\R\Ga}(v)\subset\wt\Edg, \quad
\Ver_{\wt{e}},\Ver_{\wt{e}}^c\subset\Ver, \quad\hbox{and}\quad
\Edg_{\wt{e}},\Edg_{\wt{e}}^c\subset\Edg,$$
for each $v\!\in\!\Ver$ and $\wt{e}\!\in\!\wt\Edg$
as in Section~\ref{cMellstrC_subs0}.
Let
\begin{gather*}
\Ver^{\R}=\big\{v\!\in\!\Ver\!:\phi(v)\!=\!v\big\}, ~~
\Ver^{\C}=\big\{v\!\in\!\Ver\!:\phi(v)\!\neq\!v\big\},~~
\Edg^H=\big\{\!\{v,v'\}\!\in\!\Edg\!:v,v'\!\in\!\Ver^{\R}\big\}, \\
\Edg^E=\big\{\!\{v,v'\}\!\in\!\Edg\!:\phi(v)\!=\!v'\big\},\quad
\Edg^{\C}=\big\{\!\{v,v'\}\!\in\!\Edg\!:\{\phi(v),\phi(v')\}\!\neq\!\{v,v'\}\!\big\}.
\end{gather*}

\vspace{.15in}

\noindent
Let $\cT(\R\Ga)$ be the collection of (trivalent) real $[\ell^{\pm}]$-marked  trees 
\hbox{$\R\Ga'\!\equiv\!(\Ver',\Edg',\mu',\phi')$} so that 
$$\Ga'\in\cT(\Ga) \qquad\hbox{and}\qquad
\ka\!\circ\!\phi=\phi'\!\circ\!\ka,$$
with $\ka\!:\Ver\!\lra\!\Ver'$ as in~\eref{Gakacond_e}.
Similarly to~\eref{cWGadfn_e},
\BE{cWGadfnR_e}
\cW_{\Ga}\equiv \bigcup_{\Ga'\in\cT(\Ga)}\!\!\!\!\!\cM_{\Ga'}\subset\ov\cM_{\ell}^{\pm}
\quad\hbox{and}\quad
\cW_{\R\Ga}\equiv \bigcup_{\R\Ga'\in\cT(\R\Ga)}\!\!\!\!\!\!\!\!\!\cM_{\R\Ga'}
=\cW_{\Ga}\!\cap\!\R\ov\cM_{\ell}
\subset\ov\cM_{\ell}^{\pm}\EE 
are open neighborhoods of~$\cM_{\Ga}$ in~$\ov\cM_{\ell}$ and~$\cM_{\R\Ga}$ in~$\R\ov\cM_{0,\ell}$
satisfying~\eref{cMGaDvr_e}.
Given a $\Ga$-basis \hbox{$\cQ_{\Ga}\!\subset\!\cQ_{\ell}^{\pm}$} as in~\eref{cQetaGadfn_e},
let \hbox{$\cW_{\Ga}'\!\subset\!\cW_{\Ga}$} be as in Corollary~\ref{cMspan_crl}
and define
\BE{cWGadfnR_e2}\cW_{\R\Ga}'=\cW_{\Ga}'\!\cap\!\R\ov\cM_{0,\ell}
\subset\cW_{\R\Ga}.\EE

\vspace{.15in}

\noindent
With $\R\Ga$ as above and $e_+\!\equiv\!\{v_+,v_+'\}\in\!\Edg^H$, define
\BE{RCwtGae_e}\begin{split}
&\Ver_{e_+}=\Ver\!\sqcup\!\{v_{\bu}\}, ~~
\phi_{e_+}\!:\Ver_{e_+}\lra\Ver_{e_+}, \quad
\phi_{e_+}(v)=\begin{cases}\phi(v),&\hbox{if}~v\!\in\!\Ver;\\
v_{\bu},&\hbox{if}~v\!=\!v_{\bu};
\end{cases}\\
&~\Edg_{e_+}=\big(\Edg\!-\!\{e_+\}\!\big)\!\sqcup\!
\big\{\!\{v_+,v_{\bu}\},\{v_{\bu},v_+'\}\!\big\},\\
&\hspace{.5in}\mu_{e_+}\!:\big[(\ell\!+\!1)^{\pm}\big]\lra\Ver_{\bu}, \quad
\mu_{e_+}(i)=\begin{cases}\mu(i),&\hbox{if}~i\!\in\![\ell^{\pm}];\\
v_{\bu},&\hbox{if}~i\!=\!(\ell\!+\!1)^{\pm}.
\end{cases}
\end{split}\EE
Thus, $\R\Ga e_+\!\equiv\!(\Ver_{e_+},\Edg_{e_+},\mu_{e_+},\phi_{e_+})$ 
is a  trivalent real $[(\ell\!+\!1)^+]$-marked tree.
If $e_+\!=\!e_{\vr}$ with \hbox{$\vr\!\in\!\cA_{\ell}^H$},
$$\ff_{(\ell+1)^{\pm}}^{\R}\!:
\cM_{\R\Ga e_+}^{\dag}\!\equiv\!\big\{\ff_{(\ell+1)^{\pm}}^{\R}\big\}^{\!-1}(\cM_{\R\Ga})\!\cap\!
\wt{D}_{\vr}''\!\cap\!\wt{D}_{\vr_{\ell^{\pm}}^c}'' \lra \cM_{\R\Ga}
\subset \R\ov\cM_{0,\ell}$$
is an $S^1$-fiber bundle containing the stratum $\cM_{\R\Ga e_+}\!\subset\!\R\ov\cM_{0,\ell}$.
The former
 consists of the $[(\ell\!+\!1)^{\pm}]$-marked curves obtained from an element of~$\cM_{\R\Ga}$
by replacing its node~$e_+$ with $\C\P^1$ or a wedge of two copies of~$\C\P^1$,
as in the middle diagram of Figure~\ref{MRblowup_fig}, and placing
the last conjugate pair of marked points onto this replacement.
The complement of~$\cM_{\R\Ga e_+}$ in~$\cM_{\R\Ga e_+}^{\dag}$ 
is a section of~$\ff_{(\ell+1)^{\pm}}^{\R}$;
it is represented by the middle diagram in Figure~\ref{MRblowup_fig}.

\subsection{The smooth structure of $\R\ov\cM_{0,\ell}$}
\label{RcMellstrC_subs}

\noindent
The real moduli space~$\R\ov\cM_{0,\ell}$ is the fixed locus of
the anti-holomorphic involution
\BE{Psielldfn_e}\begin{split}
\Psi_{\ell}\!:\ov\cM_{\ell}^{\pm}&\lra\ov\cM_{\ell}^{\pm}, \\
\Psi_{\ell}\big([\Si,\fJ,(z_1^+,z_1^-),\ldots,(z_{\ell}^+,z_{\ell}^-)\big]\big)
&=\big[\Si,-\fJ,(z_1^-,z_1^+),\ldots,(z_{\ell}^-,z_{\ell}^+)\big],
\end{split}\EE
where $\Si$ is a nodal surface and $\fJ$ is a complex structure on~$\Si$.
Along with \cite[Theorems~D.4.2/5.2]{MS} and \cite[Lemma~3]{Meyer81},
this implies that $\R\ov\cM_{0,\ell}$ inherits a smooth structure from
the complex structure of~$\ov\cM_{\ell}^{\pm}$.
The forgetful morphism~$\ff_{(\ell+1)^{\pm}}^{\R}$ in~\eref{ffRdfn_e}
is the restriction of the forgetful morphism
$$\ff_{\ell+1}^{\pm}\!:\ov\cM_{\ell+1}^{\pm}\lra\ov\cM_{\ell}^{\pm}$$
dropping the marked points $z_{\ell+1}^+,z_{\ell+1}^-$.\\

\noindent
The involution~preserves the divisors \hbox{$D_{\ell;\vr}^{\pm}\!\subset\!\ov\cM_{\ell}^{\pm}$}
with \hbox{$\vr\!\in\!\cA_{\ell}^E\!\cup\!\cA_{\ell}^H$} and interchanges
the divisors $D_{\ell;\vr}^{\pm},D_{\ell;\ov\vr}^{\pm}$
with \hbox{$\vr\!\in\!\cA_{\ell}^{\pm}\!-\!\cA_{\ell}^D$}.
Furthermore,
\BE{CRfqPsi_e}\CR_{\fq}\big(\Psi_{\ell}(\cC)\!\big)=\ov{\CR_{\ov\fq}(\cC)}
\qquad\forall~\cC\!\in\!\ov\cM_{\ell}^{\pm},\,\fq\!\in\!\cQ_{\ell}^{\pm}\,.\EE 
Along with \cite[Theorem~D.4.2(ii)]{MS}, 
the above fixed locus property and~\eref{CRfqPsi_e} imply~that 
\BE{RDMstr_e}\R\ov\cM_{0,\ell}=\big\{\cC\!\in\!\ov\cM_{\ell}^{\pm}\!:
\ov{\CR_{\fq}(\cC)}\!=\!\CR_{\ov\fq}(\cC)
~\forall\,\fq\!\in\!\cQ_{\ell}^{\pm}\big\}.\EE
With the notation as in~\eref{WFdfn_e}, we note the following.

\begin{lmm}\label{cMspanR_lmm} 
Suppose $\ell\!\in\!\Z^+$ with $\ell\!\ge\!2$,
\hbox{$\R\Ga\!\equiv\!(\Ver,\Edg,\mu,\phi)$} is a trivalent real $[\ell^{\pm}]$-marked tree,
and $\cQ_{\Ga}\!\subset\!\cQ_{\ell}^{\pm}$ is a $\Ga$-basis as in~\eref{cQetaGadfn_e}.
Let \hbox{$\cW_{\R\Ga}'\!\subset\!\cW_{\Ga}'$} be as in~\eref{cWGadfnR_e2}.
There exists a smooth involution~$F_{\R\Ga}$ on an open subset $W_{\Ga}'\!\subset\!(\C\P^1)^{\cQ_{\Ga}}$
so that the restriction
\BE{cMspanR_e}
(\CR_{\fq})_{\fq\in\cQ_{\Ga}}\!:
\cW_{\R\Ga}'\lra(W_{\Ga}')_{F_{\R\Ga}}\subset W_{\Ga}'\subset
(\C\P^1)^{\cQ_{\Ga}}\EE
is a diffeomorphism and
for every $\fq\!\in\!\wt\cQ_{\ell}$, $\CR_{\fq}$ is a rational function on~$W_{\Ga}'$
of the cross ratios~$\CR_{\fq'}$ with $\fq'\!\in\!\cQ_{\Ga}$.
\end{lmm}

\begin{proof}
Replacing $\cW_{\Ga}'$ by $\cW_{\Ga}'\!\cap\!\Psi_{\ell}(\cW_{\Ga}')$,
we can assume that~$\cW_{\Ga}'$ is $\Psi_{\ell}$-invariant.
Let 
$$W_{\Ga}'=\big\{\!(\CR_{\fq})_{\fq\in\cQ_{\Ga}}\big\}\big(\cW_{\Ga}'\big)
\subset (\C\P^1)^{\cQ_{\Ga}} \quad\hbox{and}\quad
\ov{\cQ_{\Ga}}=\big\{\ov\fq\!:\fq\!\in\!\cQ_{\Ga}\big\}
\subset\cQ_{\ell}^{\pm}.$$
By Corollary~\ref{cMspan_crl}\ref{cMfqgenC_it}, the subspace 
$W_{\Ga}'\!\subset\!(\C\P^1)^{\cQ_{\Ga}}$ is open and the~map
\BE{cMspanR_e5}(\CR_{\fq})_{\fq\in\cQ_{\Ga}}\!:\cW_{\Ga}'\lra W_{\Ga}'\EE
is biholomorphic.
Furthermore, for each $\fq\!\in\!\ov{\cQ_{\Ga}}$ there exists 
a rational function~$F_{\fq}$ on~$W_{\Ga}'$ so~that 
\BE{cMspanR_e6}
\CR_{\fq}\!=\!F_{\fq}\!\circ\!(\CR_{\fq'})_{\fq'\in\cQ_{\Ga}}\!:\cW_{\Ga}'\lra\C\P^1\,.\EE
Define
$$F_{\R\Ga}\!\equiv\!(F_{\R\Ga;\fq})_{\fq\in\cQ_{\Ga}}:
W_{\Ga}'\lra (\C\P^1)^{\cQ_{\Ga}}, \qquad
F_{\R\Ga;\fq}(w)=\ov{F_{\ov\fq}(w)}.$$
By~\eref{cMspanR_e6} and~\eref{CRfqPsi_e},
$$F_{\R\Ga}\!\circ\!(\CR_{\fq})_{\fq\in\cQ_{\Ga}}
=(\CR_{\fq})_{\fq\in\cQ_{\Ga}}\!\circ\!\Psi_{\ell}\!:
\cW_{\Ga}'\lra W_{\Ga}'\,.$$
Thus, $F_{\R\Ga}$ is a smooth involution on~$W_{\Ga}'$.\\

\noindent
Since the subcollection $\cQ_{\Ga}\!\cup\!\ov{\cQ_{\Ga}}\!\subset\!\cQ_{\ell}^{\pm}$
is preserved the conjugation~$\ov{\,\cdot\,}$ and the map
$$(\CR_{\fq})_{\fq\in\cQ_{\Ga}\cup\ov{\cQ_{\Ga}}}\!:
\cW_{\Ga}'\lra (\C\P^1)^{\cQ_{\Ga}\cup\ov{\cQ_{\Ga}}}$$
is injective, \eref{CRfqPsi_e} implies that
$$\cW_{\Ga}\!\cap\!\R\ov\cM_{0,\ell}
=\big\{\cC\!\in\!\cW\!:\ov{\CR_{\fq}(\cC)}\!=\!\CR_{\ov\fq}(\cC)
~\forall\,\fq\!\in\!\cQ_{\Ga}\big\}.$$
Combining this with~\eref{cMspanR_e6} and the map~\eref{cMspanR_e5} being a diffeomorphism,
we obtain the claim.
\end{proof}

\subsection{Local embeddings for $X_{\vr^*}$}
\label{RlocalEmbedd_subs}

\noindent
Let $\ell\!\in\!\Z^+$ with $\ell\!\ge\!2$.
We continue with the notation and setup of Section~\ref{Rthm_subs}.
For \hbox{$\fq\!\in\!\cQ_{\ell+1}^{\pm}\!-\!\cQ_{\ell}^{\pm}$}, let
$$\cA_{\ell}^{\R}(\fq)=\begin{cases}
\{\vr\!\in\!\cA_{\ell}^{\R}\!:|\vr\!\cap\!\fq|\!\le\!1\!\},
&\hbox{if}~(\ell\!+\!1)^+\!\in\!\fq,~(\ell\!+\!1)^-\!\not\in\!\fq;\\
\{\vr\!\in\!\cA_{\ell}^{\R}\!:|\vr\!\cap\!\ov\fq|\!\le\!1\!\},
&\hbox{if}~(\ell\!+\!1)^-\!\in\!\fq,~(\ell\!+\!1)^+\!\not\in\!\fq;\\
\{\vr\!\in\!\cA_{\ell}^{\R}\!:|\vr\!\cap\!\ov\vr\!\cap\fq|\!\le\!1\},
&\hbox{if}~(\ell\!+\!1)^+,(\ell\!+\!1)^-\!\in\!\fq.
\end{cases}$$
For $\fq\!\in\!\cQ_{\ell+1}^{\pm}$ and $\vr\!\in\!\cA_{\ell}^{\R}$,
the map~$\CR_{\fq}$  from~$\R\ov\cM_{0,\ell+1}$ is constant on every fiber~of the forgetful morphism
$$\ff_{(\ell+1)^{\pm}}^{\R}\!: 
\wt{D}_{\vr}''\!\equiv\!\wt{D}_{\vr}^0\!\cup\!\wt{D}_{\vr}^-\lra 
D_{\ell;\vr}\subset\R\ov\cM_{0,\ell}$$
if and only if either $\vr\!\in\!\cQ_{\ell}^{\pm}$ or $\vr\!\not\in\!\cA_{\ell}^{\R}(\fq)$.
For each $\vr^*\!\in\!\{0\}\!\sqcup\!\cA_{\ell}^{\R}$,
this map thus descends to a continuous~map
\BE{CRdescR_e}\begin{aligned}
\CR_{\vr^*;\fq}\!: X_{\vr^*}&\lra\C\P^1 &\quad&\hbox{if}~~ \fq\!\in\!\cQ_{\ell}^{\pm}
\qquad\hbox{and}\\
\CR_{\vr^*;\fq}\!: X_{\vr^*}-
\bigcup_{\begin{subarray}{c}\vr\in\cA_{\ell}^{\R}(\vr^*)\\ 
\vr\in\cA_{\ell}^{\R}(\fq)\end{subarray}}\!\!\!\!\!\!Y_{\vr^*;\vr}^0&\lra\C\P^1
&\quad&\hbox{if}~~ \fq\!\in\!\cQ_{\ell+1}^{\pm}\!-\!\cQ_{\ell}^{\pm}\,.
\end{aligned}\EE
For any $\vr\!\subset\![\ell^{\pm}]$, define
\begin{gather*}
\wt{D}_{\vr}'=D_{\ell+1;\vr\cup\{(\ell+1)^+\}}\!\cup\!
D_{\ell+1;\vr\cup\{(\ell+1)^+,(\ell+1)^-\}}\subset\R\ov\cM_{0,\ell+1}, \\
Y_{\vr^*;\vr}'=q_{\vr^*}\big(\wt{D}_{\vr}'\big)\subset X_{\vr^*},\qquad
Y_{\vr^*;\vr}''=q_{\vr^*}\big(\wt{D}_{\vr}''\big)
=Y_{\vr^*;\vr}^0\!\cup\!Y_{\vr^*;\vr}^-\subset X_{\vr^*}.
\end{gather*}
If $\vr\!\in\!\cA_{\ell}^{\R}(\vr^*)$,  $Y_{\vr^*;\vr}''\!=\!Y_{\vr^*;\vr}^0$.\\

\noindent
For $\vr^*\!\in\!\{0\}\!\sqcup\!\cA_{\ell}^{\R}$ 
and a trivalent real $[\ell^{\pm}]$-marked tree $\R\Ga\!\equiv\!(\Ver,\Edg,\mu,\phi)$, let 
\begin{alignat}{1}
\notag
\cA_{\R\Ga}(\vr^*)&=\big\{\mu^{-1}(\Ver_{\wt{e}})\!:\wt{e}\!\in\!\wt\Edg,\,
\mu^{-1}(\Ver_{\wt{e}})\!\in\!\cA_{\ell}^{\R}(\vr^*)\!\big\},\\
\notag
\cA_{\R\Ga}^{\star}(\vr^*)&=\big\{\vr\!\in\!\cA_{\R\Ga}(\vr^*)\!:
\vr\!\not\supset\!\vr'~\forall\,\vr'\!\in\!\cA_{\R\Ga}(\vr^*)\!-\!\{\vr\}\!\big\},\\
\label{wtcQGavrdfnR_e}
\wt\cQ_{\R\Ga}(\vr^*)&=\cQ_{\ell}^{\pm}\!\cup\!\big\{
\fq\!\in\!\cQ_{\ell+1}^{\pm}\!:\big|\fq\!\cap\![\ell^{\pm}]\big|\!=\!3,\,
\cA_{\R\Ga}(\vr^*)\!\cap\!\cA_{\ell}^{\R}(\fq)\!=\!\eset\big\}\subset \cQ_{\ell+1}^{\pm},\\
\notag
\nV_{\R\Ga}(\vr^*)&=
\bigcap_{\vr\in\cA_{\R\Ga}(\vr^*)}\hspace{-.2in}\Ver_{\wt{e}_{\vr}}
=\bigcap_{\vr\in\cA_{\R\Ga}^{\star}(\vr^*)}\hspace{-.2in}\Ver_{\wt{e}_{\vr}}\subset\Ver.
\end{alignat}
The subcollection $\wt\cQ_{\R\Ga}(\vr^*)\!\subset\!\cQ_{\ell+1}^{\pm}$ 
is preserved by the conjugation~$\ov{\,\cdot\,}$.
As with~\eref{CblowupInp_e}, \hbox{$\nV_{\R\Ga}(\vr^*)\!\neq\!\eset$}.\\

\noindent
For $v\!\in\!\Ver$, define
\BE{wtcWdf_e2}\begin{split}
\cA_{\R\Ga}(v)&=\big\{\mu^{-1}(\Ver_{\wt{e}})\!:\wt{e}\!\in\!\wt\Edg,\,
v\!\in\!\Ver_{\wt{e}}\!\big\},\\
\wt\cW_{\vr^*;\R\Ga,v}&=q_{\vr^*}\big(\{\ff_{(\ell+1)^{\pm}}^{\R}\}^{-1}(\cW_{\R\Ga})\!\big)
\!-\!\bigcup_{\begin{subarray}{c}\vr\in\cA_{\R\Ga}(v)\\
\vr\not\in\cA_{\R\Ga}(\vr^*)\end{subarray}}
\hspace{-.2in}Y_{\vr^*;\vr}''\subset X_{\vr^*}.
\end{split}\EE
Since 
\begin{equation*}\begin{split}
q_{\vr^*}^{\,-1}\big(\wt\cW_{\vr^*;\R\Ga,v}\big)=\{\ff_{(\ell+1)^{\pm}}^{\R}\}^{-1}(\cW_{\R\Ga})
&-\!\bigcup_{\begin{subarray}{c}\vr\in\cA_{\R\Ga}(v)\\
\vr\not\in\cA_{\R\Ga}(\vr^*)\end{subarray}}
\hspace{-.2in}\wt{D}_{\vr}''
\!-\!\bigcup_{\begin{subarray}{c}\vr\in\cA_{\R\Ga}(\vr^*)\\
\vr_{\ell^{\pm}}^c\in\cA_{\R\Ga}(v)\end{subarray}}\hspace{-.24in}\wt{D}_{\vr}''
\end{split}\end{equation*}
by~\eref{ndRprp_e},
the subspace $\wt\cW_{\vr^*;\R\Ga,v}\!\subset\!X_{\vr^*}$ is open;
the last union above is empty if $v\!\in\!\nV_{\R\Ga}(\vr^*)$.\\

\noindent
Let $v_+\!\in\!\nV_{\R\Ga}(\vr^*)$.
Choose \hbox{$[\Ga]_v\!\subset\![\ell^{\pm}]$} and 
\hbox{$\cQ_{\Ga}(v)\!\subset\!\cQ_{\ell}^{\pm}$} 
for each $v\!\in\!\Ver$ as in~\eref{cQetavGadfn_e} and $\fq_e\!\in\!\cQ_{\ell}^{\pm}$ 
for each $e\!\in\!\Edg$ as in~\eref{fqedfn_e}.
With $\cQ_{\Ga}\!\subset\!\cQ_{\ell}^{\pm}$ as in~\eref{cQetaGadfn_e},
let \hbox{$\cW_{\R\Ga}'\!\subset\!\cW_{\R\Ga}$} be as in~\eref{cWGadfnR_e2} and
\BE{CRwtcWdf2_e}\wt\cW_{\vr^*;\R\Ga,v_+}'
=q_{\vr^*}\big(\{\ff_{(\ell+1)^{\pm}}^{\R}\}^{\!-1}(\cW_{\R\Ga}')\!\big)
\!\cap\!\wt\cW_{\vr^*;\R\Ga,v_+}
=q_{\vr^*}\big(\{\ff_{(\ell+1)^{\pm}}^{\R}\}^{\!-1}(\cW_{\R\Ga}')\!\big)
\!-\!\bigcup_{\begin{subarray}{c}\vr\in\cA_{\R\Ga}(v)\\
\vr\not\in\cA_{\R\Ga}(\vr^*)\end{subarray}}
\hspace{-.2in}Y_{\vr^*;\vr}'';\EE
the subspace~$\wt\cW_{\vr^*;\R\Ga,v_+}'\!\subset\!X_{\vr^*}$ is again open.
Choose $i_{v_+},j_{v_+},k_{v_+}\!\in\![\ell^{\pm}]$ 
$(\R\Ga,v_+)$-independent and define
\BE{CRwtfqvrdfn_e}
\wt\fq_{v_+}=\big(i_{v_+},j_{v_+},k_{v_+},(\ell\!+\!1)^+\big)\in\wt\cQ_{\R\Ga}(\vr^*),\quad
\wt\cQ_{\Ga,v_+}=\cQ_{\Ga}\!\cup\!\big\{\wt\fq_{v_+},\ov{\wt\fq_{v_+}}\big\}
\subset \wt\cQ_{\R\Ga}(\vr^*).\EE

\begin{prp}\label{Rblowup_prp} 
Suppose $\ell\!\in\!\Z^+$ with $\ell\!\ge\!2$,
$\vr^*\!\in\!\{0\}\!\sqcup\!\cA_{\ell}^{\R}$, 
\hbox{$\R\Ga\!\equiv\!(\Ver,\Edg,\mu,\phi)$} is a trivalent real $[\ell^{\pm}]$-marked tree,
and $v_+\!\in\!\nV_{\R\Ga}(\vr^*)$.
Let~$W_{\Ga}'$ and $F_{\R\Ga}$ be as in Lemma~\ref{cMspanR_lmm}.
With the notation and assumptions as in~\eref{CRwtcWdf2_e} and~\eref{CRwtfqvrdfn_e},
the~map
\BE{cMRspan2_e}
\wt\CR_{\vr^*;\R\Ga,v_+}\!\equiv\!
(\CR_{\vr^*;\fq})_{\fq\in\wt\cQ_{\Ga,v_+}}\!:
\wt\cW_{\vr^*;\R\Ga,v_+}'\lra(\C\P^1)^{\wt\cQ_{\Ga,v_+}}
\!\equiv\!(\C\P^1)^{\cQ_{\Ga}}\!\times\!(\C\P^1)^{\{\wt\fq_{v_+},\ov{\wt\fq_{v_+}}\}}\EE
is a well-defined homeomorphism onto its image~$\wt{W}_{\vr^*;\R\Ga,v_+}'$.
Furthermore,
there exists an open subset
\hbox{$\wt{W}_{\vr^*;\Ga,v_+}'\!\subset\!W_{\Ga}'\!\times\!
(\C\P^1)^{\{\wt\fq_{v_+},\ov{\wt\fq_{v_+}}\}}$} such~that 
\BE{MRblowup_e0b}
\wt{W}_{\vr^*;\R\Ga,v_+}'=
\big\{\!(c_{\fq})_{\fq\in\wt\cQ_{\R\Ga,v_+}}\!\!\in\!(W_{\Ga}')_{F_{\R\Ga}}
\!\times\!(\C\P^1)^{\{\wt\fq_{v_+},\ov{\wt\fq_{v_+}}\}}\!:
\ov{c_{\wt\fq_{v_+}}}\!=\!c_{\ov{\wt\fq_{v_+}}}\big\}\EE 
and for every $\fq\!\in\!\wt\cQ_{\R\Ga}(\vr^*)$, 
$\CR_{\vr^*;\fq}$ is a rational function on~$\wt{W}_{\vr^*;\Ga,v_+}'$
of the cross ratios~$\CR_{\vr^*;\fq'}$ with $\fq'\!\in\!\wt\cQ_{\R\Ga,v_+}$.
\end{prp}

\begin{proof}
By~\eref{cMGaDvr_e},
$\wt\cW_{\vr^*;\R\Ga,v_+}'\!\subset\!\{\ff_{(\ell+1)^{\pm}}^{\R}\}^{-1}(\cW_{\R\Ga})$
is disjoint from $Y_{\vr^*;\vr}^0\!=\!Y_{\vr^*;\vr}''$
if \hbox{$\vr\!\in\!\cA_{\ell}^{\R}(\vr^*)\!-\!\cA_{\R\Ga}(\vr^*)$}.
Thus, the map~\eref{cMRspan2_e} is well-defined (and continuous).\\

\noindent
The proof of the injectivity of~\eref{cMRspan2_e} proceeds as 
in the case of Proposition~\ref{Cblowup_prp}.
We can start with \hbox{$\wt\cC_1,\wt\cC_2\!\in\!\R\ov\cM_{0,\ell+1}$} such~that 
$$[\wt\cC_1]_{\vr^*},[\wt\cC_2]_{\vr^*}\in\wt\cW_{\vr^*;\R\Ga,v_+}', ~~ 
[\wt\cC_1]_{\vr^*}\neq[\wt\cC_2]_{\vr^*}, ~~
\cC\equiv\ff_{(\ell+1)^{\pm}}^{\R}(\wt\cC_1)
=\ff_{(\ell+1)^{\pm}}^{\R}(\wt\cC_2)
\in\cM_{\R\Ga'}\subset\R\ov\cM_{0,\ell}$$
for some $\R\Ga'\!\equiv\!(\Ver',\Edg',\mu',\phi')\!\in\!\cT(\R\Ga)$.
Let $\ka\!:\Ver\!\lra\!\Ver'$ be as in~\eref{Gakacond_e}.
For $r\!=\!1,2$, the sets
$$\cA_{\wt\cC_r}^{\R}(\vr^*)\equiv\big\{\vr\!\in\!\cA_{\ell}^{\R}(\vr^*)\!:
\wt\cC_r\!\in\!\wt{D}_{\vr}''\big\}\subset
\cA_{\R\Ga'}(\vr^*)\subset\cA_{\R\Ga}(\vr^*) $$
are ordered by the inclusion of subsets of~$[\ell^{\pm}]$.
If $\cA_{\wt\cC_r}^+(\vr^*)\!\neq\!\eset$, we define
\hbox{$v_r^+\!\in\!\nV_{\R\Ga'}(\vr^*)$} and \hbox{$z_{\ell+1}'(\wt\cC_r)\!\in\!\C\P^1_{v_r^+}$} 
as above~\eref{Cblowup_e4} with $\ell\!+\!1$ replaced by~$(\ell\!+\!1)^+$.
As before, $v_r^+\!=\!\ka(v_+)$.\\

\noindent
Suppose $r\!=\!1,2$ and $\cA_{\wt\cC_r}^{\R}(\vr^*)\!=\!\eset$.
Let $\C\P^1_v$ be the irreducible component of~$\wt\cC_r$ carrying
the marked point~$z_{\ell+1}^+(\wt\cC_r)$ of~$\wt\cC_r$.
In addition to the situations~\ref{wtCmain_it} and~\ref{wtCmarked_it} 
in the proof of Proposition~\ref{Cblowup_prp}
with $z_{\ell+1}(\wt\cC_r)$ replaced by~$z_{\ell+1}^+(\wt\cC_r)$, 
another possibility, 
illustrated by the left diagram in Figure~\ref{MRblowup_fig},
can now arise if $v_+\!\in\!\Ver^{\R}$:
\begin{enumerate}[label=($\wt\cC\arabic*$),leftmargin=*]

\setcounter{enumi}{\value{temp}}

\item $\C\P^1_v$ also carries the marked point $z_{\ell+1}^-(\wt\cC_r)$ of~$\wt\cC_r$ only 
and contains exactly one (real) node, which it shares with the irreducible component 
$\C\P^1_{\ka(v_+)}$ of~$\cC$.
We then set \hbox{$z_{\ell+1}'(\wt\cC_r)\!\in\!\C\P^1_{\ka(v_+)}$} to be 
the (real) node shared with~$\C\P^1_v$.

\setcounter{temp}{\value{enumi}}

\end{enumerate}
In all cases, \eref{Cblowup_e4} still holds.
The reasoning in the proof of Proposition~\ref{Cblowup_prp} with $\ell\!+\!1$ replaced
by $(\ell\!+\!1)^+$ shows~that 
$$\CR_{\wt\fq_{v_+}}(\wt\cC_1)\neq\CR_{\wt\fq_{v_+}}(\wt\cC_2).$$
A situation as in~\ref{wtCnode_it} in Remark~\ref{Cblowup_rmk}, 
represented by the right diagram in Figure~\ref{Cblowup_fig} and
the middle and right diagrams in Figure~\ref{MRblowup_fig}, 
cannot arise for the same reason as in Remark~\ref{Cblowup_rmk}.
This establishes the injectivity of~\eref{cMRspan2_e}.\\

\begin{figure}
\begin{pspicture}(-1,-2.6)(10,2.1)
\psset{unit=.4cm}
\pscircle(8,-1){1.5}\rput(11.1,-1){\smsize{$v$}}
\rput(8.1,-1){\smsize{$\ka(v_+)$}}
\pscircle*(6.94,3.06){.2}\pscircle*(9.06,3.06){.2}
\pscircle*(6.94,-5.06){.2}\pscircle*(9.06,-5.06){.2}
\pscircle*(6.94,.06){.2}\pscircle*(6.94,-2.06){.2}
\pscircle(8,-4){1.5}\pscircle*(8,-2.5){.2}
\pscircle*(8,.5){.2}\pscircle(8,2){1.5}
\pscircle*(9.5,-1){.2}\pscircle(11,-1){1.5}
\pscircle*(11,.5){.2}\pscircle*(11,-2.5){.2}
\rput(11.5,1.2){\smsize{$(\ell\!+\!1)^+$}}
\rput(11.5,-3.2){\smsize{$(\ell\!+\!1)^-$}}
\pscircle(18,-1){1.5}\rput(21.1,-1){\smsize{$v$}}
\pscircle*(16.94,3.06){.2}\pscircle*(19.06,3.06){.2}
\pscircle*(16.94,-5.06){.2}\pscircle*(19.06,-5.06){.2}
\pscircle(18,-4){1.5}\rput(18.1,-4){\smsize{$v'$}}\pscircle*(18,-2.5){.2}
\pscircle*(18,.5){.2}
\pscircle(18,2){1.5}\rput(18.1,2){\smsize{$v''$}}
\pscircle*(19.5,-1){.2}\pscircle(21,-1){1.5}
\pscircle*(21,.5){.2}\pscircle*(21,-2.5){.2}
\rput(21.5,1.2){\smsize{$(\ell\!+\!1)^+$}}
\rput(21.5,-3.2){\smsize{$(\ell\!+\!1)^-$}}
\pscircle(28,3.5){1.5}\rput(28.1,3.5){\smsize{$v''$}}
\pscircle(28,.5){1.5}\pscircle*(28,2){.2}
\pscircle(28,-2.5){1.5}\pscircle*(28,-1){.2}\rput(28,-2.5){\smsize{$v$}}
\pscircle(28,-5.5){1.5}\rput(28.1,-5.5){\smsize{$v'$}}\pscircle*(28,-4){.2}
\pscircle*(26.94,4.56){.2}\pscircle*(29.06,4.56){.2}
\pscircle*(26.94,-6.56){.2}\pscircle*(29.06,-6.56){.2}
\pscircle*(29.5,.5){.2}\rput(31.2,.5){\smsize{$(\ell\!+\!1)^-$}}
\pscircle*(29.5,-2.5){.2}\rput(31.2,-2.5){\smsize{$(\ell\!+\!1)^+$}}
\end{pspicture}
\caption{The three possibilities for $z_{\ell+1}^{\pm}(\wt\cC_r)\!\in\!\C\P^1_v$
in the proof of Proposition~\ref{Rblowup_prp} and in Remark~\ref{Rblowup_rmk} 
beyond those of Figure~\ref{Cblowup_fig}.}
\label{MRblowup_fig}
\end{figure}

\noindent
With $\cW_{\Ga}'\!\subset\!\cW_{\Ga}$ as in~\eref{cWGadfnR_e2}, let 
$$\wt\cW_{\Ga,v_+}'\subset
\big\{\ff_{\ell+1}^{\pm}\big\}^{\,-1}(\cW_{\Ga}')
-\!\bigcup_{\begin{subarray}{c}\vr\in\cA_{\R\Ga}(v)\\
\vr\not\in\cA_{\R\Ga}(\vr^*)\end{subarray}}
\hspace{-.22in}
\big(D_{\ell+1;\vr}^{\pm}\!\cup\!D_{\ell+1;\vr\cup\{(\ell+1)^-\}}^{\pm}
\!\cup\!D_{\ell+1;\ov\vr}^{\pm}\!\cup\!D_{\ell+1;\ov\vr\cup\{(\ell+1)^+\}}^{\pm}\big)
\subset\ov\cM_{\ell+1}^{\pm}$$ 
be an open subset invariant under the involution~$\Psi_{\ell+1}$ as in~\eref{Psielldfn_e}
such~that 
$$q_{\vr^*}^{\,-1}\big(\wt\cW_{\vr^*;\R\Ga,v_+}'\big)=
\R\ov\cM_{0,\ell+1}\!\cap\!\wt\cW_{\Ga,v_+}'
\subset\ov\cM_{\ell+1}^{\pm}\,.$$
By the proof of Proposition~\ref{Cblowup_prp}, the~subset
\BE{MRblowup_e11}\wt{W}_{\vr^*;\Ga,v_+}'\equiv
(\CR_{\fq})_{\fq\in\wt\cQ_{\R\Ga,v_+}}
\big(\wt\cW_{\Ga,v_+}'\big) \subset(\C\P^1)^{\wt\cQ_{\R\Ga,v_+}}\EE
is open.
By Lemma~\ref{cMspanR_lmm} and~\eref{RDMstr_e},
\BE{MRblowup_e15}\begin{split}
\wt{W}_{\vr^*;\R\Ga,v_+}'&\equiv\!
\wt\CR_{\vr^*;\R\Ga,v_+}\!\big(\wt\cW_{\vr^*;\R\Ga,v_+}'\big)\\
&\subset
\big\{\!(c_{\fq})_{\fq\in\wt\cQ_{\R\Ga,v_+}}\!\!\in\!(W_{\Ga}')_{F_{\R\Ga}}
\!\times\!(\C\P^1)^{\{\wt\fq_{v_+},\ov{\wt\fq_{v_+}}\}}\!:
\ov{c_{\wt\fq_{v_+}}}\!=\!c_{\ov{\wt\fq_{v_+}}}\big\}.
\end{split}\EE

\vspace{.15in}

\noindent
Suppose $\wt\cC\!\in\!\wt\cW_{\Ga,v_+}'$ with 
$(\CR_{\fq}(\wt\cC))_{\fq\in\wt\cQ_{\R\Ga,v_+}}$ lying in 
the right-hand side of~\eref{MRblowup_e15}.
By Lemma~\ref{cMspanR_lmm}, 
$$\cC\equiv\ff_{\ell+1}^{\pm}(\wt\cC)\in\cM_{\R\Ga'}\subset\R\ov\cM_{0,\ell}$$
for some $\R\Ga'$ and $\ka$ as above.
Below we construct 
\BE{MRblowup_e8}
\wt\cC'\in \big\{\ff_{\ell+1}^{\R}\big\}^{\!-1}(\cC)\!\cap\!
q_{\vr^*}^{\,-1}\big(\wt\cW_{\vr^*;\R\Ga,v_+}'\big)
\subset\R\ov\cM_{0,\ell+1}
\quad\hbox{s.t.}\quad
\ov{\CR_{\wt\fq_{v_+}}(\wt\cC')}
=\CR_{\ov{\wt\fq_{v_+}}}(\wt\cC)
=\ov{\CR_{\wt\fq_{v_+}}(\wt\cC)}\EE
by changing the position of the marked point $z_{\ell+1}^+(\cC)$.
The above implies~that 
$$\CR_{\fq}(\wt\cC')\!=\!\CR_{\fq}(\wt\cC)
\qquad\forall\,\fq\!\in\!\big\{\wt\fq_{v_+},\ov{\wt\fq_{v_+}}\big\}$$
and confirms~\eref{MRblowup_e0b}.\\

\noindent
Let $\si$ be the involution on~$\cC$.
Thus, $\si$ interchanges the marked points and nodes of~$\cC$.
If \hbox{$z_{\ell+1}^+(\wt\cC')\!\!=\!\si(z_{\ell+1}^-(\wt\cC))$}, while all other marked
points of~$\wt\cC$ and~$\wt\cC'$ are the same in a suitable sense,
$\wt\cC'$ satisfies the first equality in~\eref{MRblowup_e8}.
We also note that the cross ratios $\CR_{\wt\fq_{v_+}}$ and $\CR_{\ov{\wt\fq_{v_+}}}$
are constant on each fiber~of
$$D_{\ell+1;\vr}^{\pm}\!\cup\!D_{\ell+1;\vr\cup\{(\ell+1)^-\}}^{\pm}\lra D_{\vr}^{\pm}
\quad\hbox{and}\quad
D_{\ell+1;\ov\vr}^{\pm}\!\cup\!D_{\ell+1;\ov\vr\cup\{(\ell+1)^+\}}^{\pm}
\lra D_{\ell;\ov\vr}^{\pm}\,,$$
respectively, for every $\vr\!\in\!\cA_{\R\Ga}(v)$.
These constants are conjugate for the fibers over $D_{\ell;\vr}\!\subset\!\ov\cM_{0,\ell}$.\\

\noindent
Let $\C\P^1_v$ be the irreducible component of~$\wt\cC$ carrying
the marked point~$z_{\ell+1}^-(\wt\cC)$ of~$\wt\cC$.
There are four possibilities now:
\begin{enumerate}[label=($\wt\cC\arabic*\R$),leftmargin=*]

\item\label{wtCmainR_it} $z_{\ell+1}^-(\wt\cC)$ is a smooth unmarked point of~$\cC$.
We then set
$z_{\ell+1}^+(\wt\cC')\!\!=\!\si(z_{\ell+1}^-(\wt\cC))\!\in\!\C\P^1_{\phi(v)}$.

\item\label{wtCmarkedR_it} $z_{\ell+1}^-(\wt\cC)\!=\!z_m(\cC)$ for some $m\!\in\![\ell^{\pm}]$.
We then set $z_{\ell+1}^+(\wt\cC')\!\!=\!z_{\ov{m}}(\wt\cC)$.

\item $z_{\ell+1}^-(\wt\cC)$ is ``at" a node $e\!\in\!\Edg^{\C}$ of~$\cC$
(with $\C\P^1_v$ possibly carrying the marked point $z_{\ell+1}^+(\wt\cC)$ as~well).
We then take $z_{\ell+1}^+(\wt\cC')$ to be the node $\phi(e)$ of~$\cC$,
so that~$\wt\cC'$ is obtained from~$\cC$ by replacing the node~$e$ 
by a~$\C\P^1$ carrying the marked point $z_{\ell+1}^-(\wt\cC')\!=\!z_{\ell+1}^-(\wt\cC)$
and the node~$\phi(e)$ 
by a~$\C\P^1$ carrying the marked point~$z_{\ell+1}^+(\wt\cC')$.

\item $z_{\ell+1}^-(\wt\cC)$ is ``at" a node $e\!\in\!\Edg^H\!\cup\!\Edg^E$ of~$\cC$.
We then take~$\wt\cC'$ to be the $[(\ell^{\pm}\!+\!1)^{\pm}]$-marked curve
obtained from~$\cC$ by replacing the node~$e$ by a~$\C\P^1$ 
carrying 
the marked points~$z_{\ell+1}^+(\wt\cC')$ and~$z_{\ell+1}^-(\wt\cC')$
at conjugate positions, along with two real nodes if  $e\!\in\!\Edg^H$
and a conjugate pair of nodes if $e\!\in\!\Edg^E$.

\end{enumerate} 
The first two situations above are analogous to \ref{wtCmain_it} and~\ref{wtCmarked_it}
in the proof of Proposition~\ref{Cblowup_prp},
illustrated by the left and middle diagrams in Figure~\ref{Cblowup_fig}.
The last two situations are analogous to~\ref{wtCnode_it}
in Remark~\ref{Cblowup_rmk}, illustrated by the right diagram in Figure~\ref{Cblowup_fig}
and the middle and right diagrams in Figure~\ref{MRblowup_fig}.
This yields an $[(\ell\!+\!1)^{\pm}]$-marked curve~$\wt\cC'$
satisfying~\eref{MRblowup_e8}.\\

\noindent
By~\eref{MRblowup_e0b}, \eref{MRblowup_e11}, and 
the open map part of the proof of Proposition~\ref{Cblowup_prp},
the map~\eref{cMRspan2_e} is open onto its image as well.
The rational function claim is obtained as in the proof of Proposition~\ref{Cblowup_prp}
with $(v_+,\ell\!+\!1)$ replaced by $(v_+,(\ell\!+\!1)^+)$
and $(\phi(v_+),(\ell\!+\!1)^-)$.
\end{proof}

\begin{rmk}\label{Rblowup_rmk}
In contrast to the situation in Remark~\ref{Cblowup_rmk}, 
the well-defined continuous map 
$$(\CR_{\vr^*;\fq})_{\fq\in\wt\cQ_{\Ga}(\vr^*)}\!:
q_{\vr^*}\big(\!\{\ff_{\ell+1}\}^{\,-1}(\cW_{\Ga})\!\big)
\lra(\C\P^1)^{\wt\cQ_{\Ga}(\vr^*)}$$
is not injective in the setting of Proposition~\ref{Rblowup_rmk}
if $\mu^{-1}(\Ver_{\wt{e}})\!\in\!\cA_{\ell}^H\!\cup\!\cA_{\ell}^E\!-\!\cA_{\ell}(\vr^*)$
for some \hbox{$\wt{e}\!\in\!\wt\Edg$}.
\end{rmk}

\subsection{Proof of Theorem~\ref{Rblowup_thm}\eref{Rspaces_it}}
\label{RblowupPf_subs1}

\noindent
Since the quotient map $q_{\vr^*}$ restricts to a homeomorphism
\BE{RblowupPf_e0}
q_{\vr^*}\!:\R\ov\cM_{0,\ell+1}-\bigcup_{\vr\in\cA_{\ell}^{\R}(\vr^*)}\!\!\!\!\!\!\!\!
\wt{D}_{\vr}''
\lra X_{\vr^*}-\bigcup_{\vr\in\cA_{\ell}^{\R}(\vr^*)}\!\!\!\!\!\!Y_{\vr^*;\vr}^0\,,\EE
the right-hand side above inherits a smooth structure from~$\R\ov\cM_{0,\ell+1}$.
We combine it with the embeddings of Proposition~\ref{Rblowup_prp}
to obtain a smooth structure on~$X_{\vr^*}$ with the stated properties.\\

\noindent
Let \hbox{$\R\Ga\!\equiv\!(\Ver,\Edg,\mu,\phi)$} be a trivalent real $[\ell^{\pm}]$-marked tree.
With $[\Ga]_v\!\subset\![\ell^{\pm}]$, \hbox{$\cQ_{\Ga}(v)\!\subset\!\cQ_{\ell}^{\pm}$}, 
$\fq_e\!\in\!\cQ_{\ell}^{\pm}$, and $\cQ_{\Ga}\!\subset\!\cQ_{\ell}^{\pm}$
as above~\eref{CRwtcWdf2_e} and $v_+\!\in\!\nV_{\R\Ga}(\vr^*)$, let 
\begin{gather*}
\cW_{\vr^*;\R\Ga,v_+}'\subset X_{\vr^*}, \qquad
\wt\fq_{v_+},\wt\fq_{v_-}\in\wt\cQ_{\R\Ga}(\vr^*),\qquad
\wt\cQ_{\R\Ga,v_+}\!\subset\!\wt\cQ_{\R\Ga}(\vr^*),\\
\hbox{and}\qquad
\wt{W}_{\vr^*;\R\Ga,v_+}'\subset
\wt{W}_{\vr^*;\Ga,v_+}'\subset (\C\P^1)^{\wt\cQ_{\R\Ga,v_+}}
\end{gather*}
be as in Proposition~\ref{Rblowup_prp}.
By the first statement of Proposition~\ref{Rblowup_prp}, \eref{MRblowup_e0b},
and \cite[Lemma~3]{Meyer81}, 
the map~\eref{cMRspan2_e} determines a smooth structure on 
the open subset \hbox{$\cW_{\vr^*;\R\Ga,v_+}'\!\subset\!X_{\vr^*}$}.
By the same reasoning as in the complex case at the beginning of Section~\ref{CblowupPf_subs},
the open subsets $\cW_{\vr^*;\R\Ga',v_+}'\!\subset\!X_{\vr^*}$  cover 
the union on the right-hand side of~\eref{RblowupPf_e0}
as~$\R\Ga'$ runs over the trivalent real $[\ell^{\pm}]$-marked trees 
and $v_+$ runs over the elements of~$\nV_{\R\Ga'}(\vr^*)$.\\

\noindent
If $\R\Ga'\!\equiv\!(\Ver',\Edg',\mu',\phi')\!\in\!\cT(\R\Ga)$, 
$\wt\cQ_{\R\Ga'}(\vr^*)\!\supset\!\wt\cQ_{\R\Ga}(\vr^*)$. 
Suppose in addition $v_+'\!\in\!\Ver'$.
By Proposition~\ref{Rblowup_prp} with~$(\R\Ga,v_+)$ replaced by~$(\R\Ga',v_+')$, 
there exists a $\wt\cQ_{\R\Ga,v_+}$-tuple 
$$F\!\equiv\!(F_{\fq})_{\fq\in\wt\cQ_{\R\Ga,v_+}}\!:
\wt{W}_{\vr^*;\Ga',v_+'}'\lra (\C\P^1)^{\wt\cQ_{\R\Ga,v_+}}$$
of rational functions such that 
\begin{equation*}\begin{split}
\wt\CR_{\vr^*;\R\Ga,v_+}\!=\!F\!\circ\!\wt\CR_{\vr^*;\R\Ga',v_+'}\!:
\,&\wt\cW_{\vr^*;\R\Ga,v_+}'\!\cap\!\wt\cW_{\vr^*;\R\Ga',v_+'}'\\
&\hspace{.5in}\lra \wt\CR_{\vr^*;\R\Ga,v_+}
\big(\wt\cW_{\vr^*;\R\Ga,v_+}'\!\cap\!\wt\cW_{\vr^*;\R\Ga',v_+'}'\big)
\subset \wt{W}_{\vr^*;\R\Ga',v_+'}'\,.
\end{split}\end{equation*}
Since $F$ is a bijection on a neighborhood of 
$$\wt\CR_{\vr^*;\R\Ga',v_+'}
\big(\wt\cW_{\vr^*;\R\Ga,v_+}'\!\cap\!\wt\cW_{\vr^*;\R\Ga',v_+'}'\big)
\subset  \wt{W}_{\vr^*;\Ga',v_+'}'\subset (\C\P^1)^{\wt\cQ_{\R\Ga',v_+'}},$$
it follows that the smooth structures on $\wt\cW_{\vr^*;\R\Ga,v_+}'$ and
$\wt\cW_{\vr^*;\R\Ga',v_+'}'$ induced by the corresponding maps~\eref{cMRspan2_e}
agree on the overlap.
Similarly to~\ref{Cspaces_it} in Section~\ref{CblowupPf_subs},
we conclude that such charts determine a smooth structure on~$X_{\vr^*}$.\\

\noindent
With $F_{\R\Ga}$ and $\wt{W}_{\vr^*;\Ga,v_+}'\!\subset\!(\C\P^1)^{\wt\cQ_{\Ga,v_+}}$
as in Lemma~\ref{cMspanR_lmm} and Proposition~\ref{Rblowup_prp}, respectively,
the~map
$$\wt{F}_{\R\Ga}\!:\wt{W}_{\vr^*;\Ga,v_+}'\lra\wt{W}_{\vr^*;\Ga,v_+}'
\subset W_{\Ga}'\!\times\!(\C\P^1)^{\{\wt\fq_{v_+},\ov{\wt\fq_{v_+}}\}},
\quad
\wt{F}_{\R\Ga}(c,c_+,c_-)=\big(F_{\R\Ga}(c),\ov{c_-},\ov{c_+}\big),$$
is a smooth involution cutting out~\eref{MRblowup_e0b}.
We next show that each subspace $Y_{\vr^*;\vr}^{\bu}\!\subset\!X_{\vr^*}$
with $\vr\!\in\!\cA_{\ell}^{\R}$ and $\bu\!\in\!\{+,0,-\}$
corresponds to the fixed locus of the restriction of~$\wt{F}_{\R\Ga}$ 
to a coordinate slice in~$\wt{W}_{\vr^*;\Ga,v_+}'$ via~\eref{cMRspan2_e}.\\

\noindent
Let $\vr\!\in\!\cA_{\ell}^{\R}$ and $v_+\!\in\!\nV_{\Ga}(\vr^*)$.
By the relevant definitions in Sections~\ref{Rthm_subs} 
and~\ref{RlocalEmbedd_subs},
$$\big(Y_{\vr^*;\vr}'\!\cap\!\wt\cW_{\vr^*;\R\Ga,v_+}'\big)
\!\cup\!\big(Y_{\vr^*;\vr}''\!\cap\!\wt\cW_{\vr^*;\R\Ga,v_+}'\big)
= q_{\vr^*}\big(\{f_{(\ell+1)^{\pm}}^{\R}\}^{\!-1}
(D_{\ell;\vr})\!\big)\!\cap\!\wt\cW_{\vr^*;\R\Ga,v_+}'.$$
By~\eref{CRwtcWdf2_e} and~\eref{cMGaDvr_e}, 
$\vr\!\in\!\cA_{\R\Ga}(v_+)$ or $\vr_{\ell^{\pm}}^c\!\in\!\cA_{\R\Ga}(v_+)$
if the right-hand side above is non-empty.
In such a case,
\eref{cMGaDvr_e2} gives
\BE{Rsubspaces_e8c}\begin{split}
\big(Y_{\vr^*;\vr}'\!\cap\!\wt\cW_{\vr^*;\R\Ga,v_+}'\big)
\!\cup\!\big(Y_{\vr^*;\vr}''\!\cap\!\wt\cW_{\vr^*;\R\Ga,v_+}'\big)
&\equiv q_{\vr^*}\big(\{f_{(\ell+1)^{\pm}}^{\R}\}^{\!-1}
(D_{\ell;\vr})\!\big)\!\cap\!\wt\cW_{\vr^*;\R\Ga,v_+}'\\
&=\big\{x\!\in\!\wt\cW_{\vr^*;\R\Ga,v_+}'\!\!:\CR_{\vr^*;\fq_{e_{\vr}}}\!(x)\!=\!0\big\}.
\end{split}\EE
By~\eref{CRwtcWdf2_e},
\BE{Rsubspaces_e8d}\begin{aligned}
Y_{\vr^*;\vr}''\!\cap\!\wt\cW_{\vr^*;\R\Ga,v_+}'&=\eset
&\quad&\hbox{if}~~
\vr\!\in\!\cA_{\R\Ga}(v_+)\!-\!\cA_{\ell}^{\R}(\vr^*),\\
Y_{\vr^*;\vr}'\!\cap\!\wt\cW_{\vr^*;\R\Ga,v_+}'
\equiv Y_{\vr^*;\vr^c_{\ell^{\pm}}}''\!\cap\!\wt\cW_{\vr^*;\R\Ga,v_+}'
&=\eset &\quad&\hbox{if}~~\vr\!\in\!\cA_{\ell}^{\R},~
\vr^c_{\ell^{\pm}}\!\in\!\cA_{\R\Ga}(v_+).
\end{aligned}\EE

\vspace{.15in}

\noindent
By~\eref{Rsubspaces_e8c}, \eref{Rsubspaces_e8d}, and~\eref{RYvrmp_e},  
\BE{Rsubspaces_e5b}
Y_{\vr^*;\vr}'\!\cap\!\wt\cW_{\vr^*;\R\Ga,v_+}'
=\begin{cases}\{x\!\in\!\wt\cW_{\vr^*;\R\Ga,v_+}'\!\!:
\CR_{\vr^*;\fq_{e_{\vr}}}(x)\!=\!0\},
&\hbox{if}~\vr\!\in\!\cA_{\ell}^{\R}\!\cap\!\cA_{\R\Ga}(v_+);\\
\eset,&\hbox{if}~\vr\!\in\!\cA_{\ell}^{\R}\!-\!\cA_{\R\Ga}(v_+).
\end{cases}\EE
By~\eref{Rsubspaces_e8c} and~\eref{Rsubspaces_e8d},
\BE{Rsubspaces_e5a}
Y_{\vr^*;\vr}''\!\cap\!\wt\cW_{\vr^*;\R\Ga,v_+}'
=\begin{cases}
\{x\!\in\!\wt\cW_{\vr^*;\R\Ga,v_+}'\!\!:\CR_{\vr^*;\fq_{e_{\vr}}}(x)\!=\!0\},
&\hbox{if}~\vr\!\in\!\cA_{\ell}^{\R},\,\vr_{\ell^{\pm}}^c\!\in\!\cA_{\R\Ga}(v_+);\\
\eset,
&\hbox{if}~\vr\!\in\!\cA_{\R\Ga}(v_+)\!-\!\cA_{\ell}^{\R}(\vr^*),\\
\eset,
&\hbox{if}~\vr\!\in\!\cA_{\ell}^{\R},~\vr,\vr^c_{\ell^{\pm}}\!\not\in\!\cA_{\R\Ga}(v_+).
\end{cases}\EE
Since $v_+\!\in\!\nV_{\R\Ga}(\vr^*)$, the conditions in the first case of~\eref{Rsubspaces_e5a}
imply that $\vr\!\in\!\cA_{\ell}^{\R}\!-\!\cA_{\ell}^{\R}(\vr^*)$.
The remaining case of the intersection in~\eref{Rsubspaces_e5a} is 
thus $\vr\!\in\!\cA_{\R\Ga}(\vr^*)$
and so $Y_{\vr^*;\vr}''\!\subset\!Y_{\vr^*;\vr}^+$.\\

\noindent
Since the involution~$\Psi_{\ell}$ in~\eref{Psielldfn_e} preserves 
the submanifold 
$D_{\ell;\vr}^{\pm}\!\cap\!D_{\ell;\ov\vr}^{\pm}\!\subset\!\ov\cM_{\ell}^{\pm}$,
the analogue of~\eref{cMGaDvr_e2} for $\ov\cM_{\ell}^{\pm}$ implies that
the nonempty subspaces in~\eref{Rsubspaces_e5b} and~\eref{Rsubspaces_e5a} correspond
to the fixed locus of
the restriction of the involution~$\wt{F}_{\R\Ga}$ to the coordinate slice
\BE{Rsubspaces_e15}\big\{\!(c_{\fq})_{\fq\in\wt\cQ_{\R\Ga,v_+}}\!\!\in\!\wt{W}_{\vr^*;\Ga,v_+}'\!:
c_{\fq_{e_{\vr}}},c_{\fq_{e_{\ov\vr}}}\!=\!0\big\}
\subset \wt{W}_{\vr^*;\Ga,v_+}'\EE
via~\eref{cMRspan2_e}
(if $\vr\!\in\!\cA_{\ell}^H\!\cup\!\cA_{\ell}^E$, $e_{\vr}\!=\!e_{\ov\vr}$).
Thus,
$$Y_{\vr^*;\vr}'\!\cap\!\wt\cW_{\vr^*;\R\Ga,v_+}',
Y_{\vr^*;\vr}''\!\cap\!\wt\cW_{\vr^*;\R\Ga,v_+}'\subset X_{\vr^*}$$
are smooth submanifolds for all $\vr\!\in\!\cA_{\ell}^{\R}$ and
$\vr\!\in\!\cA_{\ell}^{\R}\!-\!\cA_{\R\Ga}(\vr^*)$, respectively.
Each of these submanifolds either is empty or corresponds to
the fixed locus of the involution~$\wt{F}_{\R\Ga}$ to the coordinate 
slice~\eref{Rsubspaces_e15} via~\eref{cMRspan2_e}.\\

\noindent
If $\vr\!\in\!\cA_{\ell}^H$, $Y_{\vr^*;\vr}^+\!=\!Y_{\vr^*;\vr}'$ and
$Y_{\vr^*;\vr}^0,Y_{\vr^*;\vr}^-\!=\!Y_{\vr^*;\vr}''$.
If $\vr\!\in\!\cA_{\ell}^E$, then
$$Y_{\vr^*;\vr}'=Y_{\vr^*;\vr}^+\!\cup\!Y_{\vr^*;\vr}^0, \quad
Y_{\vr^*;\vr}''=Y_{\vr^*;\vr}^-\!\cup\!Y_{\vr^*;\vr}^0, \quad
Y_{\vr^*;\vr}'\!\cap\!Y_{\vr^*;\vr}''=Y_{\vr^*;\vr}^0.$$
Along with~\eref{Rsubspaces_e8d} and~\eref{RYvrmp_e}, this implies that 
\begin{gather*}
Y_{\vr^*;\vr}^0\!\cap\!\wt\cW_{\vr^*;\R\Ga,v_+}'=\eset 
~~\hbox{if}~\vr\!\in\!\cA_{\ell}^E\!-\!\cA_{\R\Ga}(\vr^*),\\
Y_{\vr^*;\vr}^+\!\cap\!\wt\cW_{\vr^*;\R\Ga,v_+}'
=Y_{\vr^*;\vr}'\!\cap\!\wt\cW_{\vr^*;\R\Ga,v_+}',~
Y_{\vr^*;\vr}^-\!\cap\!\wt\cW_{\vr^*;\R\Ga,v_+}'
=Y_{\vr^*;\vr}''\!\cap\!\wt\cW_{\vr^*;\R\Ga,v_+}' 
~~\hbox{if}~\vr\!\in\!\cA_{\ell}^E.
\end{gather*}
If $\vr\!\in\!\cA^D_{\ell;1}$, then 
$\vr\!\subset\!\ov\vr^c_{\ell^{\pm}}\!\in\!\cA^D_{\ell;2}$ and
$$Y_{\vr^*;\vr}^+=Y_{\vr^*;\vr}', \quad
Y_{\vr^*;\vr}^-,
Y_{\vr^*;\ov\vr^c_{\ell^{\pm}}}^0,Y_{\vr^*;\ov\vr^c_{\ell^{\pm}}}^-
=Y_{\vr^*;\ov\vr^c_{\ell^{\pm}}}'', \quad
Y_{\vr^*;\vr}^0,Y_{\vr^*;\ov\vr^c_{\ell^{\pm}}}^+
=Y_{\vr^*;\vr}''\!\cap\!Y_{\vr^*;\ov\vr^c_{\ell^{\pm}}}'\,.$$
If $\vr\!\in\!\cA^D_{\ell;3}$, then $\ov\vr\!\in\!\cA^D_{\ell;3}$ as well 
and
$$Y_{\vr^*;\vr}^0,Y_{\vr^*;\vr}^-=Y_{\vr^*;\vr}'', \quad
Y_{\vr^*;\vr}^+=Y_{\vr^*;\vr}'\!\cap\!Y_{\vr^*;\ov\vr}'\,.$$
Combining these statements with the conclusion of the previous paragraph,
we verify~that
\begin{gather*}
Y_{\vr^*;\vr}^+\subset X_{\vr^*} 
~~\hbox{with}~\vr\!\in\!\cA_{\ell}^{\R}\!-\!
\big\{\vr\!\in\!\cA_{\ell;2}^D\!:\ov\vr_{\ell^{\pm}}^c\!\in\!\cA_{\ell}^{\R}(\vr^*)\!\big\},
\quad
Y_{\vr^*;\vr}^0\subset X_{\vr^*} 
~~\hbox{with}~\vr\!\in\!\cA_{\ell}^{\R}\!-\!\cA_{\ell}^{\R}(\vr^*),\\
Y_{\vr^*;\vr}^-\subset X_{\vr^*} 
~~\hbox{with}~\vr\!\in\!\cA_{\ell}^{\R}\!-\!\cA_{\ell}^{\R}(\vr^*)
\!-\!
\big\{\vr\!\in\!\cA_{\ell;1}^D\!:\ov\vr_{\ell^{\pm}}^c\!\in\!\cA_{\ell}^{\R}(\vr^*)\!\big\}
\end{gather*}
are smooth submanifolds.\\

\noindent
If either $\vr\!\in\!\cA_{\R\Ga}(\vr^*)$ or $\vr\!=\![\ell^{\pm}]\!-\{j\}$
for some $j\!\in\![\ell^{\pm}]$,
we can choose~$\wt\fq_{v_+}$ in~\eref{CRwtfqvrdfn_e} 
so that $j_{v_+}\!\in\!\vr_{\ell^{\pm}}^c$.
Combining the reasoning in~\ref{Cspaces_it} in Section~\ref{CblowupPf_subs} with $\ell\!+\!1$ replaced
by~$(\ell\!+\!1)^+$
and~\ref{wtCmainR_it} and~\ref{wtCmarkedR_it} in the proof 
of Proposition~\ref{Rblowup_prp} with the roles 
of~$(\ell\!+\!1)^+$ and~$(\ell\!+\!1)^-$ reversed,
we then obtain
\begin{alignat}{1}\label{Rsubspaces_e3b}
&Y_{\vr^*;[\ell]-\{j_{v_+}\}}^0\!\cap\!\wt\cW_{\vr^*;\R\Ga,v_+}'
=\big\{x\!\in\!\wt\cW_{\vr^*;\R\Ga,v_+}'\!\!:\CR_{\vr^*;\wt\fq_{v_+}}\!(x)\!=\!0\big\},\\
\label{Rsubspaces_e3}
&Y_{\vr^*;\vr}^0\!\cap\!\wt\cW_{\vr^*;\R\Ga,v_+}'
=\big\{x\!\in\!\wt\cW_{\vr^*;\R\Ga,v_+}'\!\!:
\CR_{\vr^*;\fq_{e_{\vr}}}(x)\!=\!0,\,\CR_{\vr^*;\wt\fq_{v_+}}\!(x)\!=\!0\big\}
~~\forall\,\vr\!\in\!\cA_{\R\Ga}(\vr^*).
\end{alignat}
Via~\eref{cMRspan2_e},
the subspaces in~\eref{Rsubspaces_e3b} and~\eref{Rsubspaces_e3} correspond
to the fixed loci of
the restrictions of the involution~$\wt{F}_{\R\Ga}$ to the coordinate slices
\begin{equation*}\begin{split}
\big\{\!(c_{\fq})_{\fq\in\wt\cQ_{\R\Ga,v_+}}\!\!\in\!\wt{W}_{\vr^*;\Ga,v_+}'\!:
c_{\wt\fq_{v_+}},c_{\ov{\wt\fq_{v_+}}}\!=\!0\big\}
&\subset \wt{W}_{\vr^*;\Ga,v_+}',\\
\big\{\!(c_{\fq})_{\fq\in\wt\cQ_{\R\Ga,v_+}}\!\!\in\!\wt{W}_{\vr^*;\Ga,v_+}'\!:
c_{\fq_{e_{\vr}}},c_{\fq_{e_{\ov\vr}}},c_{\wt\fq_{v_+}},c_{\ov{\wt\fq_{v_+}}}\!=\!0\big\}
&\subset \wt{W}_{\vr^*;\Ga,v_+}'
\end{split}\end{equation*}
respectively.
Thus, the subspaces $Y_{\vr^*;\vr}^0\!=Y_{\vr^*;\vr}^-\!\!\subset\!X_{\vr^*}$ 
with $\vr\!\in\!\wt\cA_{\ell}(\vr^*)$
are smooth submanifolds as~well.

\subsection{Proof of Theorem~\ref{Rblowup_thm}\eref{Rblowup_it}}
\label{RblowupPf_sub2}

\noindent
The smooth map~$\pi_{\vr^*}$ in~\eref{Rbldowndfn_e} restricts to 
a bijection from~$X_{\vr^*}\!-\!Y_{\vr^*;\vr^*}''$ to 
\hbox{$X_{\vr^*-1}\!-\!Y_{\vr^*-1;\vr^*}^0$}
and sends $Y_{\vr^*;\vr^*}''\!\subset\!X_{\vr^*}$ 
to~$Y_{\vr^*-1;\vr^*}^0\!\subset\!X_{\vr^*-1}$.
Below we cover~$Y_{\vr^*-1;\vr^*}^0$ and~$Y_{\vr^*;\vr^*}''$ by charts as
in Lemmas~\ref{StanBl_lmm}-\ref{Rblowup_lmm},
depending on whether $\vr^*$~is an element of~$\cA_{\ell}^D$, $\cA_{\ell}^E$,
or~$\cA_{\ell}^H$.
Let
$$\ff_{\vr^*-1}\!:X_{\vr^*-1}\lra\R\ov\cM_{0,\ell} \quad\hbox{and}\quad
\ff_{\vr^*}\!=\!\ff_{\vr^*-1}\!\circ\!\pi_{\vr^*}\!:X_{\vr^*}\lra\R\ov\cM_{0,\ell}$$
be the maps  induced by $\ff_{(\ell+1)^{\pm}}^{\R}$;
these maps are smooth by Lemma~\ref{cMspanR_lmm} and Section~\ref{RblowupPf_subs1}.\\

\noindent
Suppose \hbox{$\R\Ga\!\equiv\!(\Ver,\Edg,\mu,\phi)$} is a trivalent $\ell$-marked tree
with \hbox{$\vr^*\!\in\!\cA_{\R\Ga}(\vr^*\!-\!1)$} 
and $\wt{e}_{\vr^*}\!=\!v_+v_+'$.
In particular, 
$$v_+\in\nV_{\R\Ga}(\vr^*\!-\!1), \qquad v_+'\not\in\nV_{\R\Ga}(\vr^*\!-\!1),
\quad\hbox{and}\quad
e_+\!\equiv\!e_{\vr^*}\!\equiv\!\{v_+,v_+'\}\subset\nV_{\R\Ga}(\vr^*).$$
Let  $\wt\cQ_{\R\Ga}(\vr^*\!-\!1),\wt\cQ_{\R\Ga}(\vr^*)\!\subset\!\cQ_{\ell+1}^{\pm}$ 
be as in~\eref{wtcQGavrdfnR_e} and
\hbox{$\cA_{\R\Ga}(v_+),\cA_{\R\Ga}(v_+')\!\subset\![\ell^{\pm}]$} 
as in~\eref{wtcWdf_e2}.
Define
$$\cA_{\R\Ga}(e_+)=\cA_{\R\Ga}(v_+)\!\cup\!\cA_{\R\Ga}(v_+')
-\big\{\vr^*,(\vr^*)_{\ell^{\pm}}^c\big\}.$$
Take
\hbox{$[\Ga]_v\!\subset\![\ell^{\pm}]$}, \hbox{$\cQ_{\Ga}(v)\!\subset\!\cQ_{\ell}^{\pm}$}, 
$\fq_e\!\in\!\cQ_{\ell}^{\pm}$, and $\cQ_{\Ga}\!\subset\!\cQ_{\ell}^{\pm}$ 
as above~\eref{CRwtcWdf2_e} with \hbox{$i_{e_+},k_{e_+}\!\in\!\vr^*$}.
Let 
\begin{alignat*}{2}
\wt\fq_{v_+},\wt\fq_{e_+}
&=\big(i_{e_+},j_{e_+},k_{e_+},(\ell\!+\!1)^+\big)\in\wt\cQ_{\R\Ga}(\vr^*\!-\!1),
&~~
\wt\fq_{v_+'}&=\big(i_{e_+},j_{e_+},m_{e_+},(\ell\!+\!1)^+\big)\in\wt\cQ_{\R\Ga}(\vr^*),\\
\wt\fq_{e_+}'&=\big(i_{e_+},j_{e_+},(\ell\!+\!1)^+,m_{e_+}\big)\in\wt\cQ_{\R\Ga}(\vr^*),
&~~
\wt\cQ_{\Ga;\vr^*}&=\cQ_{\Ga}\!-\!\big\{\fq_{e_+},\fq_{\phi(e_+)}\big\}\subset\cQ_{\ell}^{\pm}.
\end{alignat*}

\vspace{.15in}

\noindent
Let $\cW_{\R\Ga}'\!\subset\!\R\ov\cM_{0,\ell}$ be as in~\eref{cWGadfnR_e2} and
\hbox{$\wt\cW_{\vr^*;\R\Ga,v_+'}'\!\!\subset\!X_{\vr^*}$} be as in~\eref{CRwtcWdf2_e}.
Similarly to~\ref{Cblowup_it} in Section~\ref{CblowupPf_subs}, define
\begin{equation*}\begin{split}
\wt\cW_{\vr^*-1;\R\Ga,v_+}''&=\ff_{\vr^*-1}^{\,-1}(\cW_{\R\Ga}')
\!-\!\bigcup_{\begin{subarray}{c}\vr\in\cA_{\R\Ga}(v_+)\\ \vr\not\supset\vr^*\end{subarray}}
\hspace{-.22in}Y_{\vr^*-1;\vr}''
\!-\!\CR_{\vr^*-1;\wt\fq_{e_+}}^{-1}\!(\i)\subset X_{\vr^*-1};\\
\wt\cW_{\vr^*;\R\Ga,e_+}''&
=\ff_{\vr^*}^{\,-1}(\cW_{\R\Ga}')
\!-\!\bigcup_{\vr\in\cA_{\R\Ga}(e_+)}
\hspace{-.22in}Y_{\vr^*;\vr}''
\!-\!\CR_{\vr^*;\wt\fq_{e_+}}^{-1}\!(\i)
\!-\!\CR_{\vr^*;\wt\fq_{e_+}'}^{-1}\!(\i)\subset X_{\vr^*}.
\end{split}\end{equation*}
As before, 
\begin{equation*}\begin{split}
\CR_{\vr^*-1;\fq_{e_+}}\big(\wt\cW_{\vr^*-1;\R\Ga,v_+}''\big),
\CR_{\vr^*;\fq_{e_+}}\big(\wt\cW_{\vr^*;\R\Ga,v_+'}'\big),
\CR_{\vr^*;\wt\fq_{e_+}'}\big(\wt\cW_{\vr^*;\R\Ga,e_+}''\big)
&\subset\C,\\
\CR_{\vr^*-1;\wt\fq_{v_+}}\big(\wt\cW_{\vr^*-1;\R\Ga,v_+}''\big),
\CR_{\vr^*;\wt\fq_{v_+'}}\big(\wt\cW_{\vr^*;\R\Ga,v_+'}'\big),
\CR_{\vr^*;\wt\fq_{v_+}}\big(\wt\cW_{\vr^*;\R\Ga,e_+}''\big)&\subset\C.
\end{split}\end{equation*}

\vspace{.15in}

\noindent
As in the complex case in Section~\ref{CblowupPf_subs}, 
the quotient map~$q_{\vr^*}$ restricts to a diffeomorphism
\BE{Rblowup_e24}
\wt\cW_{\R\Ga,e_+}''\!\equiv\!
q_{\vr^*}\!:\{\ff_{(\ell+1)^{\pm}}^{\R}\}^{\!-1}(\cW_{\R\Ga}')
\!-\!\bigcup_{\vr\in\cA_{\R\Ga}(e_+)}
\hspace{-.22in}\wt{D}_{\vr}''
\!-\!\CR_{\wt\fq_{e_+}}^{-1}\!(\i)
\!-\!\CR_{\wt\fq_{e_+}'}^{-1}\!(\i)\lra\wt\cW_{\vr^*;\R\Ga,e_+}'' .\EE
By the definition of the equivalence relation~$\sim_{\vr^*}$ in Section~\ref{Rintro_subs},
\BE{Rblowup_e25}\begin{split}
&\hspace{.8in}
\pi_{\vr^*}^{\,-1}\big(\wt\cW_{\vr^*-1;\R\Ga,v_+}'')
\!=\!\wt\cW_{\vr^*;\R\Ga,v_+'}'\!\cup\!\wt\cW_{\vr^*;\R\Ga,e_+}'' \qquad\hbox{and}\\
&\CR_{\vr^*;\fq}\!=\!\CR_{\vr^*-1;\fq}\!\circ\!\pi_{\vr^*}\!:
\pi_{\vr^*}^{\,-1}\big(\wt\cW_{\vr^*-1;\R\Ga,v_+}'')\lra\C\P^1
~~\forall\,\fq\!\in\!\wt\cQ_{\R\Ga}(\vr^*\!-\!1)\,.
\end{split}\EE
Combining the last equation with~\eref{CRprp_e}, we obtain
\BE{Rblowup_e23b}\begin{split}
\CR_{\vr^*-1;\wt\fq_{v_+}}\!\!\circ\!\pi_{\vr^*}
&=\!\CR_{\vr^*;\wt\fq_{v_+}}
\!=\!\CR_{\vr^*;\wt\fq_{v_+'}}\!\cdot\CR_{\vr^*;\fq_{e_+}}\!\!:
\wt\cW_{\vr^*;\R\Ga,v_+'}'\lra\C, \\
\CR_{\vr^*-1;\fq_{e_+}}\!\!\circ\!\pi_{\vr^*}
&=\!\CR_{\vr^*;\fq_{e_+}}
\!=\!\CR_{\vr^*;\wt\fq_{e_+}'}\!\cdot\CR_{\vr^*;\wt\fq_{v_+}}\!\!:
\wt\cW_{\vr^*;\R\Ga,e_+}''\lra\C.
\end{split}\EE

\vspace{.15in}

\noindent
{\bf{\emph{Complex blowup.}}} Suppose $\vr^*\!\in\!\cA_{\ell}^D$ and 
thus $e_+\!\neq\!\phi(e_+)$.
By Lemma~\ref{cMspanR_lmm} and~\eref{cMGaDvr_e2}, there exists a collection 
$$\big\{\vph_{\al}\!:
U_{\al}\lra\R^{\cQ_{\Ga;\vr^*}}\big\}_{\al\in\cI_{\Ga;\vr^*}}$$ 
of smooth maps from open subsets of~$\cW_{\R\Ga}'$ covering
$D_{\ell;\vr^*}\!\cap\!\cW_{\R\Ga}'$ so~that 
$$\big(\CR_{\fq_{e_+}},\vph_{\al}\big)\!:
U_{\al}\lra \C\!\times\!\R^{\cQ_{\Ga;\vr^*}}$$
is a chart for~$D_{\ell;\vr^*}$ in~$\R\ov\cM_{0,\ell}$ for each $\al\!\in\!\cI_{\Ga;\vr^*}$.
By Proposition~\ref{Rblowup_prp} and~\eref{Rsubspaces_e3}
with~$(\vr^*,\vr)$ replaced~by $(\vr^*\!-\!1,\vr^*)$, the collection
$$\big\{\!\big(\CR_{\vr^*-1;\fq_{e_+}},\CR_{\vr^*-1;\wt\fq_{v_+}},
\vph_{\al}\!\circ\!\ff_{\vr^*-1}\big)\!:
\wt{U}_{\al}\!\equiv\!\wt\cW_{\vr^*-1;\R\Ga,v_+}''\big|_{U_{\al}}
\lra \C\!\times\!\C\!\times\!\R^{\cQ_{\Ga;\vr^*}}\big\}_{\al\in\cI_{\Ga;\vr^*}}$$
is then a $\C$-atlas for $Y_{\vr^*-1;\vr^*}^0\!\cap\!\wt\cW_{\vr^*-1;\R\Ga,v_+}'$
in~$\wt\cW_{\vr^*-1;\R\Ga,v_+}'$.
Similarly,
$$\wt\phi_{\al;1}\!\equiv\!\big(\CR_{\vr^*;\fq_{e_+}},\CR_{\vr^*;\wt\fq_{v_+'}},
\vph_{\al}\!\circ\!\ff_{\vr^*}\big)\!:
\wt{U}_{\al;1}\!\equiv\!\wt\cW_{\vr^*;\R\Ga,v_+'}'\big|_{U_{\al}}
\lra \C\!\times\!\C\!\times\!\R^{\cQ_{\Ga;\vr^*}}$$
is a chart on~$X_{\vr^*}$.\\

\noindent
Since the restriction~\eref{Rblowup_e24} is a diffeomorphism,
Lemma~\ref{cMspanR_lmm} and Section~\ref{RblowupPf_subs1} imply that 
$$\wt\phi_{\al;2}\!\equiv\!\big(\CR_{\vr^*;\wt\fq_{e_+}'},\CR_{\vr^*;\wt\fq_{v_+}},
\vph_{\al}\!\circ\!\ff_{\vr^*}\big)\!:
\wt{U}_{\al;2}\!\equiv\!\wt\cW_{\vr^*;\R\Ga,e_+}''\big|_{U_{\al}}
\lra \C\!\times\!\C\!\times\!\R^{\cQ_{\Ga;\vr^*}}$$
is a chart on~$X_{\vr^*}$.
Since the open subsets \hbox{$\wt\cW_{\vr^*-1;\R\Ga,v_+}'\!\subset\!X_{\vr^*-1}$}
cover the submanifold~$Y_{\vr^*-1;\vr^*}^0$
as~$\Ga$ runs over the trivalent real $[\ell^{\pm}]$-marked trees with 
$\vr^*\!\in\!\cA_{\R\Ga}(\vr^*\!-\!1)$,
so do the open subsets~$\wt{U}_{\al}$ with $\al\!\in\!\cI_{\Ga;\vr^*}$.
By Lemma~\ref{StanBl_lmm} with $(\bF,\fc)\!=\!(\C,2)$,
the previous paragraph, \eref{Rblowup_e25}, and~\eref{Rblowup_e23b},
the map~$\pi_{\vr^*}$ in~\eref{Rbldowndfn_e}
is a complex blowup of~$X_{\vr^*-1}$ along~$Y_{\vr^*-1;\vr^*}^0$
with the exceptional divisor $Y_{\vr^*;\vr^*}''\!=\!Y_{\vr^*;\vr^*}^0$.\\

\noindent
{\bf{\emph{Augmented blowup.}}} Suppose $\vr^*\!\in\!\cA_{\ell}^E$ 
and thus $e_+\!=\!\phi(e_+)$.
In this case, we take $\fq_{e_+}\!\in\!\cQ_{\ell}^{\pm}$ so that 
$\ov{i_{e_+}}\!=\!j_{e_+}$ and $\ov{k_{e_+}}\!=\!m_{e_+}$.
In light of~\eref{CRprp_e}, this implies that 
$$\CR_{\fq_{e_+}}\!=\!\CR_{\ov{\fq_{e_+}}}\!:\ov\cM_{\ell}^{\pm}\lra\C\P^1\,.$$
By Lemma~\ref{cMspanR_lmm} and~\eref{cMGaDvr_e2}, there then exists a collection 
$$\big\{\vph_{\al}\!:
U_{\al}\lra\R^{\cQ_{\Ga;\vr^*}}\big\}_{\al\in\cI_{\Ga;\vr^*}}$$ 
of smooth maps from open subsets of~$\cW_{\R\Ga}'$ covering
$D_{\ell;\vr^*}\!\cap\!\cW_{\R\Ga}'$ so~that 
$$\big(\CR_{\fq_{e_+}},\vph_{\al}\big)\!:
U_{\al}\lra \R\!\times\!\R^{\cQ_{\Ga;\vr^*}}$$
is a chart for~$D_{\ell;\vr^*}$ in~$\R\ov\cM_{0,\ell}$ for each $\al\!\in\!\cI_{\Ga;\vr^*}$.
By Proposition~\ref{Rblowup_prp} and~\eref{Rsubspaces_e3}
with~$(\vr^*,\vr)$ replaced~by $(\vr^*\!-\!1,\vr^*)$, the collection
$$\big\{\!\big(\CR_{\vr^*-1;\fq_{e_+}},\CR_{\vr^*-1;\wt\fq_{v_+}},
\vph_{\al}\!\circ\!\ff_{\vr^*-1}\big)\!:
\wt{U}_{\al}\!\equiv\!\wt\cW_{\vr^*-1;\R\Ga,v_+}''\big|_{U_{\al}}
\lra \R\!\times\!\C\!\times\!\R^{\cQ_{\Ga;\vr^*}}\big\}_{\al\in\cI_{\Ga;\vr^*}}$$
is then a 1-augmented atlas for $Y_{\vr^*-1;\vr^*}^0\!\cap\!\wt\cW_{\vr^*-1;\R\Ga,v_+}'$
in~$\wt\cW_{\vr^*-1;\R\Ga,v_+}'$.
Similarly,
$$\wt\phi_{\al;1}^1\!\equiv\!\big(\CR_{\vr^*;\fq_{e_+}},\CR_{\vr^*;\wt\fq_{v_+'}},
\vph_{\al}\!\circ\!\ff_{\vr^*}\big)\!:
\wt{U}_{\al;1}^1\!\equiv\!\wt\cW_{\vr^*;\R\Ga,v_+'}'\big|_{U_{\al}}
\lra \R\!\times\!\C\!\times\!\R^{\cQ_{\Ga;\vr^*}}$$
is a chart on~$X_{\vr^*}$.\\

\noindent
With the assumptions as before, let 
\begin{alignat*}{2}
\wt{q}_{e_+}^{\pm}&=\big(i_{e_+},j_{e_+},(\ell\!+\!1)^+,(\ell\!+\!1)^-\big)
\in\cQ_{\ell+1}^{\pm}, &\quad
\wt\cW_{\vr^*;\R\Ga,e_+}^{\pm}&=
\wt\cW_{\vr^*;\R\Ga,e_+}''\!-\!\CR_{\vr^*;\wt{q}_{e_+}^{\pm}}^{\,-1}(\i),\\
\wt{q}_{e_+}^{\mp}&=\big(i_{e_+},j_{e_+},(\ell\!+\!1)^-,(\ell\!+\!1)^+\big)
\in\cQ_{\ell+1}^{\pm},
&\quad
\wt\cW_{\vr^*;\R\Ga,e_+}^{\mp}&=
\wt\cW_{\vr^*;\R\Ga,e_+}''\!-\!\CR_{\vr^*;\wt{q}_{e_+}^{\mp}}^{\,-1}(\i).
\end{alignat*}
By~\eref{CRprp_e} and~\eref{Rblowup_e25}, $\wt\cW_{\vr^*;\R\Ga,e_+}''\!=\!
\wt\cW_{\vr^*;\R\Ga,e_+}^{\pm}\!\cup\!\wt\cW_{\vr^*;\R\Ga,e_+}^{\mp}$
and
\BE{Rblowup_e33a}
\CR_{\vr^*-1;\fq_{v_+}}\!\!\circ\!\pi_{\vr^*}
\!=\!\CR_{\vr^*;\fq_{v_+}}
\!=\!\CR_{\ov{\wt\fq_{e_+}'}}\!\cdot\!\CR_{\vr^*;\wt\fq_{e_+}^{\mp}}\!:
\wt\cW_{\vr^*;\R\Ga,e_+'}^{\mp}\lra\C.\EE
By \eref{RDMstr_e} and 
the assumptions $\ov{i_{e_+}}\!=\!j_{e_+}$ and $\ov{k_{e_+}}\!=\!m_{e_+}$, 
\begin{equation*}\begin{split}
\big|\CR_{\wt\fq_{v_+}}\big|^2&=\CR_{i_{e_+}j_{e_+}k_{e_+}(\ell+1)^+}\!
\cdot\!\CR_{i_{e_+}j_{e_+}(\ell+1)^-m_{e_+}}\!:
\R\ov\cM_{0,\ell+1}\lra \R\P^1\subset\C\P^1,\\
\big|\CR_{\ov{\wt\fq_{e_+}'}}\big|^2&=\CR_{i_{e_+}j_{e_+}k_{e_+}(\ell+1)^-}\!
\cdot\!\CR_{i_{e_+}j_{e_+}(\ell+1)^+m_{e_+}}\!:
\R\ov\cM_{0,\ell+1}\lra \R\P^1\subset\C\P^1.
\end{split}\end{equation*}
Along with~\eref{CRprp_e} and~\eref{Rblowup_e25}, these identities give
\BE{Rblowup_e33b}\begin{split}
\CR_{\vr^*-1;\fq_{e_+}}\!\!\circ\!\pi_{\vr^*}
&=\!\big|\CR_{\vr^*;\wt\fq_{v_+}}\big|^2\CR_{\vr^*;\wt\fq_{e_+}^{\pm}}\!:
\wt\cW_{\vr^*;\R\Ga,e_+'}^{\pm}\lra\C, \\
\CR_{\vr^*-1;\fq_{e_+}}\!\!\circ\!\pi_{\vr^*}
&=\!\big|\CR_{\vr^*;\ov{\wt\fq_{e_+}'}}\big|^2\CR_{\vr^*;\wt\fq_{e_+}^{\mp}}\!:
\wt\cW_{\vr^*;\R\Ga,e_+'}^{\mp}\lra\C.
\end{split}\EE

\vspace{.15in}

\noindent
Since $\ov{i_{e_+}}\!=\!j_{e_+}$, \eref{CRprp_e} also implies that 
$$\CR_{\wt\fq_{e_+}^{\pm}}\!=\!\CR_{\ov{\wt\fq_{e_+}^{\pm}}}\!:\ov\cM_{\ell}^{\pm}\lra\C\P^1
\quad\hbox{and}\quad
\CR_{\wt\fq_{e_+}^{\mp}}\!=\!\CR_{\ov{\wt\fq_{e_+}^{\mp}}}\!:\ov\cM_{\ell}^{\pm}\lra\C\P^1\,.$$
Since the restriction~\eref{Rblowup_e24} is a diffeomorphism,
Lemma~\ref{cMspanR_lmm} and Section~\ref{RblowupPf_subs1} thus imply that 
\begin{equation*}\begin{split}
\wt\phi_{\al;0}^2\!\equiv\!\big(\CR_{\vr^*;\wt\fq_{e_+}^{\pm}},\CR_{\vr^*;\wt\fq_{v_+}},
\vph_{\al}\!\circ\!\ff_{\vr^*}\big)\!:
\wt{U}_{\al;0}^2\!\equiv\!\wt\cW_{\vr^*;\R\Ga,e_+}^{\pm}\big|_{U_{\al}}
\lra \R\!\times\!\C\!\times\!\R^{\cQ_{\Ga;\vr^*}},\\
\wt\phi_{\al;1}^2\!\equiv\!\big(\CR_{\vr^*;\wt\fq_{e_+}^{\mp}},
\CR_{\vr^*;\ov{\wt\fq_{e_+}'}},
\vph_{\al}\!\circ\!\ff_{\vr^*}\big)\!:
\wt{U}_{\al;1}^2\!\equiv\!\wt\cW_{\vr^*;\R\Ga,e_+}^{\mp}\big|_{U_{\al}}
\lra \R\!\times\!\C\!\times\!\R^{\cQ_{\Ga;\vr^*}}
\end{split}\end{equation*}
are charts on~$X_{\vr^*}$.
Since the open subsets \hbox{$\wt\cW_{\vr^*-1;\R\Ga,v_+}'\!\subset\!X_{\vr^*-1}$}
cover the submanifold~$Y_{\vr^*-1;\vr^*}^0$
as~$\Ga$ runs over the trivalent real $[\ell^{\pm}]$-marked trees with 
$\vr^*\!\in\!\cA_{\R\Ga}(\vr^*\!-\!1)$,
so do the open subsets~$\wt{U}_{\al}$ with $\al\!\in\!\cI_{\Ga;\vr^*}$.
By Lemma~\ref{AugBl_lmm} with $(\fc,\fc_1)\!=\!(3,1)$, \eref{Rblowup_e25},
\eref{Rblowup_e33a}, and \eref{Rblowup_e33b},
the map~$\pi_{\vr^*}$ in~\eref{Rbldowndfn_e}
is a 1-augmented blowup of~$X_{\vr^*-1}$ along~$Y_{\vr^*-1;\vr^*}^0$
with the exceptional divisor $Y_{\vr^*;\vr^*}''\!=\!Y_{\vr^*;\vr^*}^0\!\cup\!Y_{\vr^*;\vr^*}^-$.\\

\noindent
{\bf{\emph{Real blowup.}}} Suppose $\vr^*\!\in\!\cA_{\ell}^H$ 
and thus $e_+\!=\!\phi(e_+)$ again.
With \hbox{$\cW_{\Ga}'\!\subset\!\ov\cM_{\ell}^{\pm}$} 
and \hbox{$W_{\Ga}'\!\subset\!(\C\P^1)^{\cQ_{\Ga}}$} as in Lemma~\ref{cMspanR_lmm},
let $F$ be the $\C$-valued rational function on~$W_{\Ga}'$ so~that 
$$\CR_{\ov{\fq_{e_+}}}\!=\!F\!\circ\!(\CR_{\fq})_{\fq\in\cQ_{\Ga}}\!:\cW_{\Ga}'\lra\C\P^1\,.$$
Since $\CR_{\ov{\fq_{e_+}}}$ can be used in place of $\CR_{\fq_{e_+}}$ in~$\cQ_{\Ga}$
and both vanish on \hbox{$D_{\ell;\vr^*}^{\pm}\!\subset\!\ov\cM_{\ell}^{\pm}$},
there exists a $\C^*$-valued rational function~$g_{\Ga;\vr^*}$ on~$W_{\Ga}'$ with 
\BE{Rblowup_e41}F\big(\!(c_{\fq})_{\fq\in\cQ_{\Ga}}\big)
=c_{\fq_{e_+}}g_{\Ga;\vr^*}\big(\!(c_{\fq})_{\fq\in\cQ_{\Ga}}\big)
~~\forall~(c_{\fq})_{\fq\in\cQ_{\Ga}}\!\in\!W_{\Ga}' .\EE
By Lemma~\ref{cMspanR_lmm} and~\eref{cMGaDvr_e2}, there  exists a collection 
$$\big\{\!\big(\vph_{\al}\!:U_{\al}\!\lra\!\R^{\cQ_{\Ga;\vr^*}},
\io_{\al}\!:\R^{\cQ_{\Ga;\vr^*}}\!\lra\!\C^{\cQ_{\Ga;\vr^*}}\big)\!\big\}_{\al\in\cI_{\Ga;\vr^*}}$$ 
of pairs consisting of a smooth map from an open subset $U_{\al}\!\subset\!\cW_{\R\Ga}'$
and a smooth map from~$\R^{\cQ_{\Ga;\vr^*}}$ so that 
the sets~$U_{\al}$ cover $D_{\ell;\vr^*}\!\cap\!\cW_{\R\Ga}'$,
$$(\CR_{\fq})_{\fq\in\cQ_{\Ga;\vr^*}}
\!=\!\io_{\al}\!\circ\vph_{\al}\!:U_{\al}\!\lra\!\C^{\cQ_{\Ga;\vr^*}}\,,$$
the~map
$$F_{\al}\!:W_{\al}\!\equiv\!
\big\{\id_{\C}\!\times\!\io_{\al}\big\}^{-1}(W_{\Ga}')\lra\C\!\times\!\R^{\cQ_{\Ga;\vr^*}}, 
\quad
F_{\al}(z,r)=\big(\ov{zg_{\Ga;\vr^*}(z,\io_{\al}(r)\!)},r\big),$$
is an involution, and the~map
$$\big(\CR_{\fq_{e_+}},\vph_{\al}\big)\!:
U_{\al}\lra \big\{\!(z,r)\!\in\!W_{\al}\!:F_{\al}(z,r)\!=\!(z,r)\big\}
\subset \C\!\times\!\R^{\cQ_{\Ga;\vr^*}}$$ 
is an open embedding which identifies $U_{\al}\!\cap\!D_{\ell;\vr^*}$
with the locus $z\!=\!0$ in the image.
Define
$$g_{\al}\!:W_{\al}\lra\C^*,\quad g_{\al}(z,r)=g_{\Ga;\vr^*}\big(z,\io_{\al}(r)\!\big)
~~\forall~(z,r)\in W_{\al}\!\subset\!\C\!\times\!\R^{\cQ_{\Ga;\vr^*}}\,.$$
By Proposition~\ref{Rblowup_prp} and~\eref{Rsubspaces_e3}
with~$(\vr^*,\vr)$ replaced~by $(\vr^*\!-\!1,\vr^*)$, 
the collection
\begin{equation*}\begin{split}
\big\{\!\big((\CR_{\vr^*-1;\wt\fq_{v_+}},\CR_{\vr^*-1;\fq_{e_+}},
\vph_{\al}\!\circ\!\ff_{\vr^*-1})\!:
\wt{U}_{\al}\!\equiv\!\wt\cW_{\vr^*-1;\R\Ga,v_+}''\big|_{U_{\al}}
\lra \C^2\!\times\!\R^{\cQ_{\Ga;\vr^*}},&\\
g_{\al}\!:W_{\al}\!\lra\!\C^*&\big)\!\big\}_{\al\in\cI_{\Ga;\vr^*}}
\end{split}\end{equation*}
is a quasi-atlas for $Y_{\vr^*-1;\vr^*}^0\!\cap\!\wt\cW_{\vr^*-1;\R\Ga,v_+}'$
in~$\wt\cW_{\vr^*-1;\R\Ga,v_+}'$.
Similarly, the~pair
$$\big(\wt\phi_{\al;2}\!\equiv\!(\CR_{\vr^*;\wt\fq_{v_+'}},\CR_{\vr^*;\fq_{e_+}},
\vph_{\al}\!\circ\!\ff_{\vr^*})\!:
\wt{U}_{\al;2}\!\equiv\!\wt\cW_{\vr^*;\R\Ga,v_+'}'\big|_{U_{\al}}
\lra\C^2\!\times\!\R^{\cQ_{\Ga;\vr^*}},
g_{\al}\!:W_{\al}\!\lra\!\C^*\big)$$
is a type~1 pre-chart on~$X_{\vr^*}$.\\

\noindent
Let $\R\Ga_+$ be as in~\eref{RCwtGae_e}  and 
$$\wt\cQ_{\Ga,e_+}=\cQ_{\Ga;\vr^*}\!\cup\!
\big\{\wt\fq_{v_+},\wt\fq_{e_+}',\wt\fq_{e_+}^{\mp}\!\big\}
\subset \cQ_{\ell}^{\pm}\,.$$
By the definition of $\wt\cW_{\R\Ga,e_+}''$ in~\eref{Rblowup_e24},
$\CR_{\wt\fq_{e_+}^{\mp}}\!(\wt\cW_{\R\Ga,e_+}'')\!\subset\!\C^*$.
Similarly to the proof of Lemma~\ref{cMspanR_lmm},
there exists an open subset \hbox{$\wt\cW_{\Ga,e_+}''\!\subset\!\ov\cM_{0,\ell+1}^{\pm}$}
so~that 
\BE{Rblowup_e42}
\wt\cW_{\R\Ga,e_+}''=\wt\cW_{\Ga,e_+}''\!\cap\!\R\ov\cM_{0,\ell+1}\,,
\quad \ff_{\ell+1}^{\pm}\big(\wt\cW_{\Ga,e_+}''\big)\subset\cW_{\Ga}', 
\quad \CR_{\wt\fq_{e_+}^{\mp}}\!\big(\wt\cW_{\Ga,e_+}''\big)\subset\C^*,\EE
the holomorphic map 
$$\wt\CR_{\Ga,e_+}\!\equiv\!
(\CR_{\fq})_{\fq\in\wt\cQ_{\Ga,e_+}}\!:\wt\cW_{\Ga,e_+}''
\lra(\C\P^1)^{\wt\cQ_{\Ga,e_+}}$$
is an open holomorphic embedding, and 
for every $\fq\!\in\!\cQ_{\ell+1}^{\pm}$, $\CR_{\fq}$ is a rational function~$F_{\fq}$
of the cross ratios~$\CR_{\fq'}$ with $\fq'\!\in\!\wt\cQ_{\Ga,e_+}$
on the image~$\wt{W}_{\Ga,e_+}''$ of~$\wt\cW_{\Ga,e_+}''$
under~$\wt\CR_{\Ga,e_+}$.\\

\noindent
By the last condition in~\eref{Rblowup_e42}, 
there exists a $\C^*$-valued rational function~$\wt{g}_{\Ga;\vr^*;0}$
on~$\wt{W}_{\Ga,e_+}''$  so that 
\BE{Rblowup_e44}
F_{\ov{\wt\fq_{e_+}^{\mp}}}\big(\!(c_{\fq})_{\wt\fq\in\cQ_{\Ga,e_+}}\big)
=c_{\wt\fq_{e_+}^{\mp}}^{\,-1}\wt{g}_{\Ga;\vr^*;0}
\big(\!(c_{\fq})_{\fq\in\wt\cQ_{\Ga,e_+}}\big)\quad
~~\forall~(c_{\fq})_{\fq\in\wt\cQ_{\Ga,e_+}}\!\in\!\wt{W}_{\Ga,e_+}''.\EE
Since $\CR_{\ov{\wt\fq_{e_+}^{\mp}}}\!=\!1/\CR_{\wt\fq_{e_+}^{\mp}}$ on
$D_{\ell+1;\vr}^{\pm}\!\cap\!D_{\ell+1;\vr_{\ell^{\pm}}^c}^{\pm}\!\subset\!\ov\cM_{\ell+1}^{\pm}$
by~\eref{CRprp_e} and $\wt\CR_{\Ga,e_+}$ identifies
$$\wt\cW_{\Ga,e_+}''\!\cap\!D_{\ell+1;\vr}^{\pm}\!\cap\!
D_{\ell+1;\vr_{\ell^{\pm}}^c}^{\pm}\subset\ov\cM_{\ell+1}^{\pm}$$
with the locus $c_{\wt\fq_{v_+}},c_{\wt\fq_{e_+}'}\!=\!0$ in~$\wt{W}_{\Ga,e_+}''$,
\BE{Rblowup_e49}
\wt{g}_{\Ga;\vr^*;0}\big(\!(c_{\fq})_{\fq\in\wt\cQ_{\Ga,e_+}}\big)=1
\quad\forall~(c_{\fq})_{\fq\in\wt\cQ_{\Ga,e_+}}\!\in\!\wt{W}_{\Ga,e_+}''
~\hbox{s.t.}~c_{\wt\fq_{v_+}},c_{\wt\fq_{e_+}'}\!=\!0.\EE

\vspace{.15in}

\noindent
By the same reasoning as for~\eref{Rblowup_e41},
there exist $\C^*$-valued rational functions~$\wt{g}_{\Ga;\vr^*;1}$ and~$\wt{g}_{\Ga;\vr^*;2}$
on~$\wt{W}_{\Ga,e_+}''$  so that 
\BE{Rblowup_e45}\begin{split}
F_{\ov{\wt\fq_{v_+}}}\big(\!(c_{\fq})_{\wt\fq\in\wt\cQ_{\Ga,e_+}}\big)
&=c_{\wt\fq_{v_+}}\wt{g}_{\Ga;\vr^*;1}\big(\!(c_{\fq})_{\fq\in\wt\cQ_{\Ga,e_+}}\big),\\
F_{\ov{\wt\fq_{e_+}'}}\big(\!(c_{\fq})_{\wt\fq\in\wt\cQ_{\Ga,e_+}}\big)
&=c_{\wt\fq_{e_+}'}\wt{g}_{\Ga;\vr^*;2}\big(\!(c_{\fq})_{\fq\in\wt\cQ_{\Ga,e_+}}\big)
\end{split}
\quad
~~\forall~(c_{\fq})_{\fq\in\wt\cQ_{\Ga,e_+}}\!\in\!\wt{W}_{\Ga,e_+}''.\EE
Since $\CR_{\fq_{e_+}}\!=\!\CR_{\wt\fq_{v_+}}\CR_{\wt\fq_{e_+}'}$ by~\eref{CRprp_e},
\eref{Rblowup_e41} and \eref{Rblowup_e45} imply that 
\begin{gather}\label{Rblowup_e47}
\wt{g}_{\Ga;\vr^*;1}\big(\!(c_{\fq})_{\fq\in\wt\cQ_{\Ga,e_+}}\big)
\wt{g}_{\Ga;\vr^*;2}\big(\!(c_{\fq})_{\fq\in\wt\cQ_{\Ga,e_+}}\big)
=g_{\Ga;\vr^*}\big(\!(c_{\fq}')_{\fq\in\cQ_{\Ga}}\big),\\
\notag
\hbox{with}\quad
c_{\fq}'=\begin{cases}
c_{\wt\fq_{v_+}}c_{\wt\fq_{e_+}'},&\hbox{if}~\fq\!=\!\fq_{e_+};\\
c_{\fq},&\hbox{if}~\fq\!\in\!\cQ_{\Ga;\vr^*}\,.
\end{cases}
\end{gather}
On
$D_{\ell+1;\vr}^{\pm}\!\cap\!D_{\ell+1;\vr_{\ell^{\pm}}^c}^{\pm}\!\subset\!\ov\cM_{\ell+1}^{\pm}$,
$$\CR_{\ov{\wt\fq_{v_+}}}=\CR_{i_{e_+}j_{e_+}k_{e_+}(\ell+1)^-}
=\CR_{\wt\fq_{v_+}}\CR_{\wt\fq_{e_+}^{\mp}}$$
by~\eref{CRprp_e}.
Thus, 
\BE{Rblowup_e49b}
\wt{g}_{\Ga;\vr^*;1}\big(\!(c_{\fq})_{\fq\in\wt\cQ_{\Ga,e_+}}\big)=1
\quad\forall~(c_{\fq})_{\fq\in\wt\cQ_{\Ga,e_+}}\!\in\!\wt{W}_{\Ga,e_+}''
~\hbox{s.t.}~c_{\wt\fq_{v_+}},c_{\wt\fq_{e_+}'}\!=\!0.\EE

\vspace{.15in}

\noindent
For each $\al\!\in\!\cI_{\Ga;\vr^*}$ and $j\!=\!0,1,2$, define
\begin{gather*}
\wt{g}_{\al;j}\!:
\wt{W}_{\al}\!\equiv\!\big\{\id_{\C^*}\!\times\!\id_{\C^2}\!\times\!\io_{\al}\big\}^{-1}
\big(\wt{W}_{\Ga,e_+}''\big)\lra\C^*,  \\
\wt{g}_{\al;j}(z,r)=\wt{g}_{\Ga;\vr^*;j}\big(z,\io_{\al}(r)\!\big)
~~\forall~(z,r)\in W_{\al}\!\subset\!\C\!\times\!\R^{\cQ_{\Ga;\vr^*}}\,.
\end{gather*}
Let $\wt{g}_{\al}\!=\!(\wt{g}_{\al;j})_{j\in\dbsqbr{2}}$.
By~\eref{Rblowup_e47}, \eref{Rblowup_e49b}, and~\eref{Rblowup_e49}, 
\begin{alignat}{2}
\label{Rblowup_e51b}
\wt{g}_{\al;1}\big(\!(z_j)_{j},r\big)
\wt{g}_{\al;2}\big(\!(z_j)_{j},r\big)
&=g_{\al}\big(z_1z_2,r\big)
&\quad&\forall~\big(\!(z_j)_{j},r\big)\!\in\!
\wt{W}_{\al}\!\subset\!\C^*\!\times\!\C^2\!\times\!\R^{\cQ_{\Ga;\vr^*}},\\
\label{Rblowup_e51a}
\wt{g}_{\al;0}\big(\!(z_j)_{j},r\big),\wt{g}_{\al;1}\big(\!(z_j)_{j},r\big)&=1
&\quad&\forall~\big(\!(z_j)_{j},r\big)\!\in\!
\wt{W}_{\al}\!\subset\!\C^*\!\times\!0^2\!\times\!\R^{\cQ_{\Ga;\vr^*}}.
\end{alignat}
Since the restriction~\eref{Rblowup_e24} is a diffeomorphism,
Lemma~\ref{cMspanR_lmm}, Section~\ref{RblowupPf_subs1}, and~\eref{Rblowup_e51a} 
imply that  the~pair
\begin{equation*}\begin{split}
\big(\wt\phi_{\al;1}\!\equiv\!\big(\CR_{\vr^*;\wt\fq_{e_+}^{\mp}},
\CR_{\vr^*;\wt\fq_{v_+}},\CR_{\vr^*;\wt\fq_{e_+}'},
\vph_{\al}\!\circ\!\ff_{\vr^*}\big)\!:
\wt{U}_{\al;1}\!\equiv\!\wt\cW_{\vr^*;\R\Ga,e_+}''\big|_{U_{\al}}
\lra \C^*\!\times\!\C^2\!\times\!\R^{\cQ_{\Ga;\vr^*}},&\\
\wt{g}_{\al}\!:\wt{W}_{\al}\lra(\C^*)^3&\big)
\end{split}\end{equation*}
is a type $(1,2)$ pre-chart on~$X_{\vr^*}$.
Furthermore,
$$\pi_{\vr^*}^{\,-1}\big(\wt{U}_{\al}\big)\!-\!\wt{U}_{\al;2}
=\big\{\wt{x}\!\in\!\wt{U}_{\al;1}\!:
\CR_{\vr^*;\wt\fq_{v_+}}(\wt{x}),\CR_{\vr^*;\wt\fq_{e_+}'}(\wt{x})\!=\!0\big\}.$$
Since the open subsets \hbox{$\wt\cW_{\vr^*-1;\R\Ga,v_+}'\!\subset\!X_{\vr^*-1}$}
cover the submanifold~$Y_{\vr^*-1;\vr^*}^0$
as~$\Ga$ runs over the trivalent real $[\ell^{\pm}]$-marked trees with 
$\vr^*\!\in\!\cA_{\R\Ga}(\vr^*\!-\!1)$,
so do the open subsets~$\wt{U}_{\al}$ with $\al\!\in\!\cI_{\Ga;\vr^*}$.
By Lemma~\ref{Rblowup_lmm}, \eref{Rblowup_e33b}, and~\eref{Rblowup_e51b},
the map~$\pi_{\vr^*}$ in~\eref{Rbldowndfn_e}
is thus real blowup of~$X_{\vr^*-1}$ along~$Y_{\vr^*-1;\vr^*}^0$
with the exceptional divisor $Y_{\vr^*;\vr^*}''\!=\!Y_{\vr^*;\vr^*}^0$.

\subsection{Proof of Theorem~\ref{Rblowup_thm}\eref{Rsubspaces_it}}
\label{RblowupPf_sub3}

\noindent
We continue with the notation in Sections~\ref{RblowupPf_subs1} and~\ref{RblowupPf_sub2}.
Suppose $\vr\!\in\!\cA_{\ell}^{\R}(\vr^*)$.
By the same reasoning as in~\ref{Csubspaces_it} in Section~\ref{CblowupPf_subs}
with $\Ga$, \eref{wtfqvrdfn_e}, \eref{Csubspaces_e3}, \eref{Cblowup_e25},
and~\eref{Cblowup_e23b},
replaced by $\R\Ga$, \eref{CRwtfqvrdfn_e}, \eref{Rsubspaces_e3}, \eref{Rblowup_e25},
and~\eref{Rblowup_e23b}, respectively, 
the submanifolds \hbox{$Y_{\vr^*;\vr}^0\!\subset\!X_{\vr^*}$} 
and  \hbox{$Y_{\vr^*-1;\vr}^0\!\subset\!X_{\vr^*-1}$} 
are $\pi_{\vr^*}$-equivalent.
Furthermore, the restriction 
$$\ff_0\!:Y_{0;\vr}^0\lra D_{\ell;\vr}\subset \R\ov\cM_{0,\ell}$$
is a diffeomorphism for every $\vr\!\in\!\cA_{\ell}^{\R}$.
Since there are identifications
\BE{Rsubspaces_e55}D_{\ell;\vr}\approx\begin{cases}\ov\cM_{\ell+1},&
\hbox{if}~\vr\!\in\!\cA_{\ell}^E;\\
\R\ov\cM_{|\ov\vr-\vr^c_{\ell^{\pm}}|/2+1}\!\times\!\ov\cM_{|\vr|+1},&
\hbox{if}~\vr\!\in\!\cA_{\ell;1}^D;\\
\R\ov\cM_{|\vr^c_{\ell^{\pm}}-\ov\vr|/2+1}\!\times\!\ov\cM_{|\vr|+1},&
\hbox{if}~\vr\!\in\!\cA_{\ell;3}^D.
\end{cases}\EE
Thus, the submanifolds $Y_{\vr^*;\vr}^0\!=\!Y_{\vr^*;\vr}^-\!\subset\!X_{\vr^*}$
with $\vr\!\in\!\cA_{\ell}^{\R}(\vr^*)\!\cap\!(\cA_{\ell}^E\!\cup\!\cA_{\ell}^D)$
are orientable by the upward induction on~$\vr^*$.
Along with the $\pi_{\vr^*}$-equivalence conclusion above,
$\pi_{\vr^*}$-related claim of Theorem~\ref{Rblowup_thm}\eref{Rsubspaces_it}
whenever $\vr\!\in\!\cA_{\ell}^{\R}(\vr^*)$ and $\bu\!=\!0,-$.\\

\noindent
For any $\vr\!\in\!\wt\cA_{\ell}^{\R}$ and $\bu\!=\!+,0,-$, the restriction
$$\pi_{\vr^*}\!:Y_{\vr^*;\vr}^{\bu}\!-\!
Y_{\vr^*;\vr^*}^0\!\cup\!Y_{\vr^*;\vr^*}^-\!\lra 
Y_{\vr^*-1;\vr}^{\bu}\!-\!Y_{\vr^*-1;\vr^*}^0$$
is a diffeomorphism.
The submanifolds 
$Y_{\vr^*;\vr}^+,Y_{\vr^*;\vr}^0,Y_{\vr^*;\vr}^-\subset X_{\vr^*}$
with 
\BE{Rsubspaces_e57}\vr\in \wt\cA_{\ell}^{\R}\!-\!\cA_{\ell;2}^D,\quad
\vr\!\in\!\wt\cA_{\ell}^{\R}\!-\!\cA_{\ell}(\vr^*\!-\!1),
\quad\hbox{and}\quad
\vr\in\wt\cA_{\ell}^{\R}\!-\!\cA_{\ell}(\vr^*\!-\!1)\!-\!\cA_{\ell;1}^D\,,\EE
respectively, have the same dimensions as the corresponding submanifolds
\hbox{$Y_{\vr_{\max};\vr}^{\bu}\!\subset\!\R\ov\cM_{0,\ell+1}$}.
By definition, 
\BE{Rsubspaces_e53}\begin{split}
Y_{\vr^*;\vr}^{\bu}\!\cap\!(Y_{\vr^*;\vr^*}^0\!\cup\!Y_{\vr^*;\vr^*}^-\big)
&=q_{\vr^*}\big(Y_{\vr_{\max};\vr}^{\bu}\!\cap\!Y_{\vr_{\max};\vr^*}^0 \big)
\!\cup\!q_{\vr^*}\big(Y_{\vr_{\max};\vr}^{\bu}\!\cap\!Y_{\vr_{\max};\vr^*}^-\big)
\subset X_{\vr^*} \quad\hbox{and}\\
Y_{\vr^*-1;\vr}^{\bu}\!\cap\!Y_{\vr^*-1;\vr^*}^0
&=q_{\vr^*-1}\big(Y_{\vr_{\max};\vr}^{\bu}\!\cap\!Y_{\vr_{\max};\vr^*}^0 \big)
\subset X_{\vr^*-1}.
\end{split}\EE
Since the intersections of distinct submanifolds 
$D_{\ell+1;\wt\vr}\!\subset\!\R\ov\cM_{0,\ell+1}$
with $\wt\vr\!\in\!\cA_{\ell+1}^{\R}$ are transverse, 
$$Y_{\vr_{\max};\vr}^{\bu}\!\cap\!Y_{\vr_{\max};\vr^*}^0
\subset Y_{\vr_{\max};\vr}^{\bu}$$
is a proper submanifold if $(\vr,\bu)\!\neq\!(\vr^*,0)$, 
$(\vr,\bu)\!\neq\!(\ov{\vr^*}_{\ell^{\pm}}^c,+)$ if $\vr\!\in\!\cA_{\ell;1}^D$,
and $(\vr,\bu)\!\neq\!(\ov{\vr^*}_{\ell^{\pm}}^c,-)$ if $\vr\!\in\!\cA_{\ell;2}^D$.
Similarly, 
$$Y_{\vr_{\max};\vr}^{\bu}\!\cap\!Y_{\vr_{\max};\vr^*}^-
\subset Y_{\vr_{\max};\vr}^{\bu}$$
is a proper submanifold if $(\vr,\bu)\!\neq\!(\vr^*,-)$,
$(\vr,\bu)\!\neq\!(\ov{\vr^*}_{\ell^{\pm}}^c,0),(\ov{\vr^*}_{\ell^{\pm}}^c,0)$ 
if $\vr\!\in\!\cA_{\ell;1}^D$,
and \hbox{$(\vr,\bu)\!\neq\!(\ov{\vr^*}_{\ell^{\pm}}^c,-)$} if $\vr\!\in\!\cA_{\ell;2}^D$.\\

\noindent
The conclusion of the previous paragraph completes the proof of
the $\pi_{\vr^*}$-related claim for $\vr\!\in\!\cA_{\ell}^H$.
Since 
$$Y_{\vr^*;\vr}^+=Y_{\vr^*;\ov\vr_{\ell^{\pm}}^c}^0
~~\hbox{if}~\vr\!\in\!\cA_{\ell;2}^D
\quad\hbox{and}\quad
Y_{\vr^*;\vr}^-=Y_{\vr^*;\ov\vr_{\ell^{\pm}}^c}^0
~~\hbox{if}~\vr\!\in\!\cA_{\ell;1}^D,$$
the above conclusion would also complete the proof of 
this claim for $\vr\!\in\!\wt\cA_{\ell}^{\R}\!-\!\cA_{\ell}^H$
if each submanifold \hbox{$Y_{\vr^*;\vr}^{\bu},Y_{\vr^*-1;\vr}^{\bu}\!\subset\!X_{\vr^*}$} 
with~$\vr$ satisfying the respective condition in~\eref{Rsubspaces_e57} is orientable.
By the reasoning as in~\eref{Rsubspaces_e55},
this is the case for~$Y_{\vr_{\max};\vr}^{\bu}$.
By~\eref{Rsubspaces_e5b} and~\eref{Rsubspaces_e5a}, 
the two paragraphs afterward, and~\eref{Rsubspaces_e3},
all with~$\vr^*$ replaced by~$\vr^*\!-\!1$,
the intersection on the second line in~\eref{Rsubspaces_e53}
is a submanifold of~$Y_{\vr^*-1;\vr}^{\bu}$
of codimension at least~2 if $\vr\!\in\!\cA_{\ell}^{\R}$
and the respective condition in~\eref{Rsubspaces_e57} holds.
This submanifold thus does not affect the orientability of~$Y_{\vr^*-1;\vr}^{\bu}$.
If $\vr\!=\![\ell^{\pm}]\!-\!\{i\}$ for some $i\!\in\![\ell^{\pm}]$,
\hbox{$Y_{\vr^*;\vr}^+,Y_{\vr^*-1;\vr}^{\bu}\!=\!\eset$} and
the restriction 
$$\pi_{\vr^*}\!:Y_{\vr^*;\vr}^0\!=\!Y_{\vr^*;\vr}^-\lra 
Y_{\vr^*-1;\vr}^0\!=\!Y_{\vr^*-1;\vr}^-$$
is a diffeomorphism by~\eref{Rsubspaces_e3b} and~\eref{Rsubspaces_e3},
both as stated and with~$\vr^*$ replaced by~$\vr^*\!-\!1$.
Thus, the orientability of~$Y_{\vr^*;\vr}^{\bu}$ implies the orientability
of~$Y_{\vr^*-1;\vr}^{\bu}$ whenever $(\vr,\bu)$ satisfies~\eref{Rsubspaces_e57}.
By the downward induction on~$\vr^*$, 
we then obtain the desired orientability claim.\\

\vspace{.2in}

\noindent
{\it Department of Mathematics, Harvard University, Cambridge, MA 02138\\
xujiachen@g.harvard.edu}\\

\noindent
{\it Department of Mathematics, Stony Brook University, Stony Brook, NY 11794\\
azinger@math.stonybrook.edu}


\begin{thebibliography}{99}

\bibitem{AkKi85} S.~Akbulut and H.~King,
{\it Submanifolds and homology of nonsingular real algebraic varieties},
Amer.~J.~Math.~107 (1985), no.~1, 45–-83

\bibitem{ArKa} G.~Arone and M.~Kankaanrinta, 
{\it On the functoriality of the blow-up construction},
Bull.~Belg.~Math.~Soc.~Simon Stevin 17 (2010), no.~5, 821–-832

\bibitem{Cey}  O.~Ceyhan,  
{\it On moduli of pointed real curves of genus zero}, 
Proceedings of G\"okova Geometry-Topology Conference 2006, 1--38

\bibitem{RDMhomol} X.~Chen, P.~Georgieva, and A.~Zinger,
{\it The cohomology ring of the Deligne-Mumford space
of real rational curves with conjugate marked points}, math/2305.08798

\bibitem{EHKR} P.~Etingof, A.~Henriques, J.~Kamnitzer, and E.~Rains,
{\it The cohomology ring of the real locus of the moduli space of stable curves 
of genus~0 with marked points}, Ann.~of Math.~171 (2010), no.~2, 731-–777

\bibitem{Penka2} P.~Georgieva,
{\it Open Gromov-Witten invariants in the presence of an anti-symplectic involution},
 Adv.~Math.~301 (2016), 116-–160

\bibitem{RealEnum} P.~Georgieva and A.~Zinger,
{\it Enumeration of real curves in $\C\P^{2n-1}$ 
and a WDVV relation for real Gromov-Witten invariants}, 
Duke Math.~J.~166 (2017), no.~17, 3291-–3347

\bibitem{GH}  P.~Griffiths and J.~Harris,
{\it Principles of Algebraic Geometry}, Wiley, 1994

\bibitem{Keel} S.~Keel, 
{\it Intersection theory of moduli spaces of stable $n$-pointed curves of genus zero},
Trans.~AMS 330 (1992), no.~2, 545--574 

\bibitem{MS} D.~McDuff and D.~Salamon, 
{\it $J$-holomorphic Curves and Symplectic Topology},
Colloquium Publications~52, AMS, 2012

\bibitem{Meyer81} K.~Meyer, 
{\it Hamiltonian systems with a discrete symmetry},
J.~Diff.~Eqs 41 (1981), no.~2, 228–-238

\bibitem{blowups} A.~Zinger,
{\it Smooth blowups: global vs.~local perspectives}, math/2312.16112

\end{thebibliography}
\end{document}